\documentclass[12pt]{amsart}
\usepackage{amsmath}
\usepackage{amsxtra}
\usepackage{amscd}
\usepackage{amsthm}
\usepackage{amsfonts}
\usepackage{amssymb}
\usepackage{eucal}
\usepackage{epsfig}
\usepackage{graphics}
\usepackage{graphicx}
\usepackage{accents}
\usepackage{dynkin-diagrams}
\usepackage{comment}
\usepackage{enumitem}
\usepackage{mathtools}
\usepackage{tikz-cd}
\usepackage{mathtools}
\usepackage{pgf}
\usetikzlibrary{arrows, matrix}
\usepackage{here} 
\usepackage{simplewick}
\usepackage{ulem}
\numberwithin{equation}{section}
\newtheorem{thm}{Theorem}[section]
\newtheorem{prop}[thm]{Proposition}
\newtheorem{lem}[thm]{Lemma}
\newtheorem{rem}[thm]{Remark}

\newcommand{\bh}{{\bf h}_N({\bf u}_0,\ldots, {\bf u}_l; \vartheta)}
\theoremstyle{definition}

\newcommand{\ket}[1]{{| #1 \rangle}}

\newcommand{\C}{{\mathbb C}}
\newcommand{\Z}{{\mathbb Z}}

\newcommand{\R}{{\mathbb R}}

\newcommand{\E}{{\mathcal E}}
\newcommand{\F}{\mathcal F}

\newcommand{\K}{\mathcal{K}}

\newcommand{\gl}{\mathfrak{gl}}

\newcommand{\bLa}{{\boldsymbol \Lambda}}

\newcommand{\sfq}{{\sf q}}

\newcommand{\ssq}{\textsf{q}}

\newcommand{\ssu}{\textsf{u}}

\newcommand{\End}{\mathop{\rm End}}

\newcommand{\res}{{\rm res}}

\newcommand{\Sym}{\mathrm{Sym}}

\newcommand{\wt}{{\rm wt}\,}

\newcommand{\on}{\operatorname}
\newcommand{\mc}{\mathcal}
\newcommand{\al}{\alpha}

\newcommand{\bLambda}{{\boldsymbol \varLambda}}

\newcommand{\g}{\mathfrak{g}}

\newcommand{\bT}{\boldsymbol T}

\newcommand{\uq}{U_\ssq(\ddot{D_l})}
\newcommand{\FF}{\mathbb{F}}
\newcommand{\omb}{{\bar{\omega}}}
\newcommand{\alb}{\bar{\alpha}}
\newcommand{\jm}{\textcolor{red}}

\newcommand{\com}[2]{\left\lVert\begin{array}{c}#1\\
#2\end{array}\right\lVert}
\newcommand{\comb}[2]{\left(\begin{array}{c}#1\\
#2\end{array}\right)}
\newcommand{\ol}[1]{\overline{#1}}
\newcommand{\ul}[1]{\underline{#1}}

\begin{document}
\begin{title}
{Integrals of motion for
the deformed ${\mathcal W}$-algebra ${\mathcal W}_{q,t}(D_l)$}
\end{title}
\author{Michio Jimbo and Takeo Kojima}

\address{MJ: 
Professor Emeritus,
Rikkyo University, Toshima-ku, Tokyo 171-8501, Japan}
\email{jimbomm@rikkyo.ac.jp}
\address{TK:
Yamagata University, Jonan 4-3-16, Yonezawa 992-8510, Japan}
\email{kojima@yz.yamagata-u.ac.jp}

\begin{abstract} 
The deformed $W$-algebra $\mathcal{W}_{q,t}(D_l)$ of type $D_l$ possesses an 
infinite family of commuting elements that can be viewed as a two-parameter 
deformation of the classical Hamiltonians for the KdV hierarchy of type $D_l$; 
these conserved charges are termed local integrals of motion (IMs). Concurrently, 
the quantum toroidal algebra $U_{\mathsf{q}}(\ddot D_l)$ of type $D_l$ admits four 
distinct representations via free field realizations, each corresponding to one 
of the four end nodes of the affine Dynkin diagram. Utilizing these representations, 
we associate an infinite set of commuting operators with each node, which we define 
as non-local integrals of motion. Both local and non-local integrals of motion act 
simultaneously on a shared Fock space consisting of $l+1$ free fields. We provide a 
direct, computational 
proof demonstrating that these local and non-local integrals of 
motion mutually commute.
\end{abstract}

\hfill \today

\vspace{5mm}

\dedicatory{To the memory of Kenji Iohara}



\maketitle

\section{Introduction}

One of the remarkable features of $W$-algebras is the presence of large 
commutative subalgebras that quantize the Hamiltonians of integrable 
hierarchies of the Korteweg-de Vries (KdV) type. Triggered by the seminal 
papers 
\cite{Bazhanov-Lukyanov-Zamolodchikov1}--\cite{ Bazhanov-Lukyanov-Zamolodchikov3}, 
considerable progress has been achieved in analyzing their spectra and 
exploring their applications to quantum integrable models. These Hamiltonians 
are defined via spatial integrals of local field densities and are conventionally 
termed local integrals of motion (IMs). To extend this framework, Bazhanov, 
Lukyanov, and Zamolodchikov introduced non-local IMs, which take the form of 
multiple integrals of fields. Specifically, they constructed spectral-parameter-dependent 
transfer matrices $\mathbf{T}(\lambda)$ whose asymptotic expansions around 
$\lambda \to \infty$ and $\lambda \to 0$ systematically generate the local 
and non-local IMs, respectively. This formulation renders the 
mutual commutativity of local IMs with non-local ones immediately manifest.

Deformed $W$-algebras $\mathcal{W}_{q,t}(\mathfrak{g})$ provide a two-parameter 
deformation of the ordinary $W$-algebras $\mathcal{W}(\mathfrak{g})$. They 
were initially introduced and studied in \cite{Shiraishi-Kubo-Awata-Odake, 
Awata-Kubo-Odake-Shiraishi, Feigin-Frenkel} for the $\mathfrak{sl}_n$ case, and subsequently 
generalized to arbitrary simple Lie algebras $\mathfrak{g}$ in 
\cite{Frenkel-Reshetikhin}. For type $A$ algebras, the deformed analogs 
of both local and non-local IMs were successfully constructed in 
\cite{Feigin-Kojima-Shiraishi-Watanabe, Kojima-Shiraishi}. These formulations 
critically rely on the free field realization. Upon $q$-deformation, 
both types of IMs become non-local. Nevertheless, for the sake 
of book-keeping, we shall retain the designations 
``local IM'' and ``non-local IM'' to refer to these deformed counterparts.

Quantum toroidal algebras offer an
algebraic framework for a uniform 
description of these deformed structures. First, the $\mathfrak{gl}_1$ quantum 
toroidal algebra $\mathcal{E}_1$ intrinsically encodes information about 
$\mathcal{W}_{q,t}(\mathfrak{sl}_n)$ for all $n$. Specifically, the 
generating currents $e(z)$ and $f(z)$ of $\mathcal{E}_1$, when acting on an 
$n$-fold tensor product of Fock representations, yield the fundamental current of 
$\mathcal{W}_{q,t}(\mathfrak{sl}_n)$ enhanced by a free boson 
\cite{Feigin-Hashizume-Hoshino-Shiraishi-Yanagida}. 

In the deformed setting, the transfer matrix construction of 
\cite{Bazhanov-Lukyanov-Zamolodchikov1} was carried out
 in \cite{Feigin-Jimbo-Mukhin3}. Consider a Fock representation 
of the quantum toroidal $\mathfrak{gl}_n$ algebra $\mathcal{E}_n$. Each 
of its weight spaces can be identified as vector-space 
with an $n$-fold tensor 
product of Fock representations of $\mathcal{E}_1$. On this space, two distinct 
families of commuting operators act simultaneously: (i) the expansion 
coefficients of the trace of the universal $R$-matrix of $\mathcal{E}_1$, and 
(ii) the corresponding coefficients for $\mathcal{E}_n$. Provided that the 
algebraic parameters of $\mathcal{E}_1$ and $\mathcal{E}_n$ are appropriately matched, 
the operators from (i) and (ii) coincide exactly with the local and non-local 
IMs of \cite{Feigin-Kojima-Shiraishi-Watanabe, Kojima-Shiraishi}, respectively. 

The mutual commutativity of these two sets of operators is a highly non-trivial 
property. It can be interpreted as a $(\mathfrak{gl}_1, \mathfrak{gl}_n)$-duality, 
closely analogous to the dualities established in the affine setting 
\cite{Mukhin-Tarasov-Varchenko1, Mukhin-Tarasov-Varchenko2}. 
A generalization to full $(\mathfrak{gl}_m, \mathfrak{gl}_n)$-duality within the 
toroidal framework was subsequently established in \cite{Feigin-Jimbo-Mukhin2}.

For types other than $A$, the situation is less clear. 
The construction of $\mathcal{W}_{q,t}(\mathfrak{g})$ and its associated local 
IMs was presented in \cite{Feigin-Jimbo-Mukhin-Vilkoviskiy} by employing 
a certain comodule algebra $\mathcal{K}_1$ instead of $\mathcal{E}_1$. 
We will briefly review this formulation in the next section. 

For a simply laced Lie algebra $\mathfrak{g}$, the corresponding non-local 
IMs are expected to arise from the transfer matrix of the quantum toroidal 
algebra $U_{\sfq}(\ddot\g)$ 
of type $\mathfrak{g}$ 
\footnote{For types other than $A$, quantum toroidal algebras depend on a 
single parameter $\mathsf{q}$.}. There are, however, substantial technical 
difficulties. At the time of writing, we are unaware of any literature 
establishing the existence of the universal $R$-matrix for types other than $A$ 
(specifically regarding the non-degeneracy of the Hopf pairing in the quantum 
double construction). Even assuming that the universal $R$-matrix exists, explicit 
methods to compute its trace remain unavailable, except when the underlying 
representation admits a free field realization.

In the present paper, we present partial results regarding the duality of IMs 
for the deformed $W$-algebra $\mathcal{W}_{q,t}(D_l)$ of type $D_l$. 
The quantum toroidal algebra $U_{\mathsf{q}}(\ddot D_l)$ possesses four distinct 
representations that admit a free field realization, corresponding to the four 
end nodes of the affine $D_l$ Dynkin diagram. We derive explicit formulas 
for the associated non-local IMs and establish their mutual commutativity 
via direct computation. Furthermore, we provide a direct proof that these 
non-local IMs commute with the local IMs constructed in 
\cite{Feigin-Jimbo-Mukhin-Vilkoviskiy}. Characterizing the non-local IMs 
associated with the remaining Dynkin nodes, as well as establishing the duality 
in its full generality, are left as subjects for future investigation.

The remainder of this paper is structured as follows. 

In Section~\ref{sec:gl1 K1}, we collect basic facts concerning the 
quantum toroidal $\mathfrak{gl}_1$ algebra $\mathcal{E}_1$ and its comodule 
algebra $\mathcal{K}_1$. Section~\ref{section:deformed-W} reviews the 
construction of the deformed $W$-algebra $\mathcal{W}_{q,t}(D_l)$ utilizing 
the quantum toroidal algebras $\mathcal{E}_1$ and its comodule algebra
$\mathcal{K}_1$. In 
Section~\ref{section:main-result}, we introduce two infinite families of 
commuting operators, denoted by $\{{\bf I}_M\}_{M=1}^\infty$ and 
$\{{\bf G}_N(\vartheta)\}_{N=1}^\infty$,
where $\vartheta$ ranges over a certain four-dimensional vector 
space of theta functions, see Lemma \ref{lem:theta1} below.
Sections~\ref{section:commutativity} 
and~\ref{section:duality} are devoted to establishing the commutativity relations 
$[{{\bf G}_M(\vartheta_1), {\bf G}_N(\vartheta_2)}]=0$ and 
$[{{\bf I}_M, {\bf G}_N(\vartheta)}]=0$, respectively, through direct, 
explicit evaluations.

The main text is supplemented by five appendices. Appendix~\ref{appendix:bosonization-quantum toroidal D} 
presents the free field realizations of the relevant modules for the quantum 
toroidal algebra $U_{\mathsf{q}}(\ddot{D}_l)$. Appendix~\ref{appendix:normal-ordering} 
tabulates the normal ordering prescriptions for the fundamental $W$-currents 
and screening currents. Finally, Appendices~\ref{appendix:IM-GN}, 
\ref{appendix:I1-GN}, and~\ref{appendix:I2-GN} illustrate the specific pole 
configurations of the integrands appearing in the product 
${\bf I}_M \cdot {\bf G}_N(\vartheta)$.

\bigskip

\noindent\textit{Notation.}\quad
Throughout this paper, we fix the parameters $q, t \in \mathbb{C}$ such that 
$|q| < 1$ and $|t| < 1$. We parameterize these variables as
\begin{align*}
t = q^\beta \,, \quad q = e^{-2\pi\sqrt{-1}/\tau} \,,
\end{align*}
where $\operatorname{Im}\tau > 0$. 
We adopt the standard notation for the multiple $q$-Pochhammer symbols, defined by
\begin{align*}
(z_1, \ldots, z_r ; q)_\infty = \prod_{s=1}^r \prod_{j=0}^\infty (1 - z_s q^j) \,.
\end{align*}
Furthermore, we introduce the theta functions
\begin{align*}
\theta(u) = q^{u(u-1)/2} \Theta_q(q^u) \,, \quad 
\Theta_q(x) = (x, qx^{-1}, q ; q)_\infty \,,
\end{align*}
which obey the quasi-periodicity relations
\begin{align*}
\theta(u+1) = -\theta(u) = \theta(-u) \,, \quad 
\theta(u+\tau) = -e^{-2\pi\sqrt{-1}u - \pi\sqrt{-1}\tau} \theta(u) \,.
\end{align*}

\section{Quantum  algebras.}\label{sec:gl1 K1}

The quantum toroidal $\gl_1$ algebra $\E_1$ \cite{Miki, Ding-Iohara, Schiffmann-Vasserot, Feigin-Feigin-Jimbo-Miwa-Mukhin,
Feigin-Jimbo-Mukhin-Vilkoviskiy}
is a two-parameter deformation of the ${\mathcal W}_{1+\infty}$-algebra and has various names:
the Ding-Iohara-Miki algebra or the elliptic Hall algebra.
In this section we collect basic facts concerning 
$\E_1$ and
its comodule algebra $\K_1$.

\subsection{Quantum toroidal $\gl_1$ algebra ${\E_1}$.}\label{sec:gl1}

We set
\begin{align}
q_1=q^{1-\beta}\,,\quad
q_2=q^\beta\,,\quad q_3=q^{-1},\notag
\end{align}
so that $q_1q_2q_3=1$. 
Let
\begin{align*}
& g(z,w)=(z-q_1w)(z-q_2w)(z-q_3w)\,,\\
&\kappa_r=(1-q_1^r)(1-q_2^r)(1-q_3^r), \qquad  r\in \Z\,.
\end{align*} 

The quantum toroidal $\gl_1$ algebra $\E_1$ is a unital associative algebra
generated by elements $e_j,f_j, h_r$,
where $j\in\Z$, $r\in\Z\backslash\{0\}$,
and invertible central elements $C$, $\psi_0$. 
In terms of the generating series
\begin{align*}
&e(z) =\sum_{j\in \Z} e_{j}z^{-j}, \quad 
f(z) =\sum_{j\in\Z} f_{j}z^{-j}, \quad 
\psi^{\pm}(z) = \psi_0^{\pm1}
\exp\bigl(\sum_{r=1}^\infty \kappa_r h_{\pm r}z^{\mp r}\bigr)\,,
\end{align*}
the defining relations read as follows:
\begin{align*}
&[h_{r},h_{s}]=
\delta_{r+s,0}\,\frac{1}{r} \frac{C^{r}-C^{-r}}{\kappa_r}\,,
\\ 
&g(C^{(1\pm1)/2}z,w)\psi^\pm(z)e(w)+g(w,C^{(1\pm1)/2}z)e(w)\psi^\pm(z)=0
\,,\\ 
&g(w,C^{(1\mp1)/2}z) \psi^\pm(z)f(w)+g(C^{(1\mp1)/2}z,w)f(w)\psi^\pm(z)=0\,,
\\
&[e(z),f(w)]=\frac{1}{\kappa_1}
\bigl(\delta(C w/z)\psi^+(w)
-\delta(Cz/w) \psi^-(z)\bigr),\\
&g(z,w)e(z)e(w)+g(w,z)e(w)e(z)=0, \\
&g(w,z)f(z)f(w)+g(z,w)f(w)f(z)=0,\\
&\mathop{\on{Sym}}_{z_1,z_2,z_3}\, z_1 z_2^2 [e(z_1),[e(z_2),e(z_3)]]=0,\\
&\mathop{\on{Sym}}_{z_1,z_2,z_3}\, z_1 z_2^2 [f(z_1),[f(z_2),f(z_3)]]=0\,,
\end{align*}
where $\Sym$ stands for symmetrization in the variables indicated.

We use the coproduct as in \cite{Feigin-Jimbo-Mukhin3}:
\begin{align*}
&\Delta e(z)=e(C_2^{-1}z) \otimes \psi^+(C_2^{-1}z)+1 \otimes e(z),
\\
&\Delta f(z)=f(z) \otimes 1+\psi^-(C_1^{-1} z)\otimes f(C_1^{-1}z),\\
&\Delta \psi^+(z)=\psi^+(z) \otimes  \psi^+(C_1 z),\\
&\Delta \psi^+(z)=\psi^-(C_2 z)\otimes \psi^-(z),
\end{align*}
where $C_1=C\otimes 1$ and $C_2=1 \otimes C$.

Algebra $\E_1$ has a $\Z$ grading, denoted deg, given by
\begin{align*}
\deg e_j=j\,,\quad \deg f_j=j\,,\quad \deg h_r=r\,,\quad \deg C=\deg\psi_0=0\,. 
\end{align*}

Quite generally, we say that a graded vector space $W=\oplus_{d\in\Z}W_d$ is admissible if
\begin{enumerate}
\item $\dim W_d<\infty$ for all $d$,
\item there exists a $d_0$ such that $W_d=0$ for $d>d_0$.
\end{enumerate}
A graded current on $W$ is a formal series $a(z)=\sum_{j\in\Z}a_j z^{-j}$,
$a_j\in\End(W)$, satisfying $a_j W_d\subset W_{d+j}$ for all $d,j$. 
We say that an $\E_1$ module is admissible 
if it is admissible as graded vector space and $e(z)$, $f(z)$, $\psi^\pm(z)$ 
act on it as graded
currents.

On admissible modules, the current
\begin{align*}
\tilde{f}(z)=-{\psi^-(z)}^{-1}f(Cz)\,
\end{align*}
has a well-defined action. We have
\begin{align*}
\Delta\tilde{f}(z)=\tilde{f}(C_2z)\otimes{\psi^-(z)}^{-1}+1\otimes\tilde{f}(z)\,.
\end{align*}
\medskip

We say that an $\E_1$ module $W$ has level $\ell\in\C^\times$ if $C$ acts as the scalar $\ell$.

We shall use the 
Fock representation $\F_2(\ssu)$ 
of color $2$ and parameter $\ssu\in\C^\times$
defined as follows. 
Consider a 
Heisenberg algebra with generators $\{h_r\}_{r \in \Z\backslash\{0\}}$ 
satisfying
\begin{align*}
&[h_{r},h_{s}]=\delta_{r+s,0}\,\frac{1}{r} \frac{q_2^{r/2}-q_2^{-r/2}}{\kappa_r}\,.
\end{align*}
Let $\F_2(\ssu)=\C[h_{-r}\mid r>0]
|\ssu, q_2^{{1}/{2}}\rangle$ be its Fock representation
with vacuum vector
$|\ssu, q_2^{{1}/{2}} \rangle$  
such that $h_r |\ssu, q_2^{{1}/{2}}\rangle=0$, $r>0$.  
We extend
it to an admissible $\E_1$ module of level $q_2^{1/2}$
by letting the generators act by vertex operators:
\begin{align*}
&e(z)\mapsto b^{-1}\, V_2(z;\ssu)\,,
\quad
\tilde{f}(z)\mapsto b^{-1}\, :V_2(z;\ssu)^{-1}:\,,
\quad \psi_0\mapsto 1\,,
\end{align*}
where $b=q_2^{1/2}(1-q_1)(1-q_3)$ and 
\begin{align}
&V_2(z;\ssu)=\ssu
\exp\Bigl(\sum_{r>0}(1-q_1^r)(1-q_3^r)h_{-r}z^r\Bigr)
\exp\Bigl(\sum_{r>0}(1-q_1^{r})(1-q_3^{r})q_2^{r/2} h_{r}z^{-r}\Bigr)\,.
\label{Fock VO}
\end{align}

\subsection{Quantum algebra ${\mathcal K}_1$}

We recall algebra $\K_1$ introduced in
\cite{Feigin-Jimbo-Mukhin-Vilkoviskiy}.
The quantum algebra $\K_1$ is a unital associative algebra 
generated by elements
$E_j$, $j\in\Z$,  
$H_r$, $r\in\Z\backslash\{0\}$,
and an invertible central element $C$.
In terms of the generating currents
\begin{align*}
E(z)=\sum_{j \in {\mathbf Z}}E_j z^{-j},\quad
K^{\pm}(z)=\exp\Bigl(\sum_{r>0}\kappa_r H_{\pm r} z^{\mp r}\Bigr),
\end{align*}
the defining relations are as follows.
\begin{align*}
&[H_r, H_s]=-\delta_{r+s,0}\frac{1}{r} \frac{1+C^{2r}}{\kappa_r}\,,
\qquad \text{for $r>0$},
\\
&g(z,w)K^\pm(z)E(w)+g(w,z)E(w)K^\pm(z)=0\,,
\\
&g(z,w)E(z)E(w)+g(w,z)E(w)E(z)&
\\
=&\frac{1}{\kappa_1}
\left(g(z,w)\delta(C^2z/w)K^-(z)K^+(C^2z)+g(w,z)\delta(C^2 w/z)K^-(w)K^+(C^2 w)\right)\,,
\\
&\underset{z_1, z_2, z_3}{\rm Sym} z_2 z_3^{-1}
[E(z_1), [E(z_2), E(z_3)]]
=\frac{1}{\kappa_1} \underset{z_1, z_2, z_3}{\rm Sym}
X(z_1, z_2, z_3)\delta(C^2 z_1/z_3)K^-(z_1)E(z_2)K^+(z_3)\,.
\end{align*}
Here
\begin{align}
X(z_1,z_2,z_3)=
\frac{(z_1+z_2)(z_3^2-z_1z_2)}{z_1 z_2 z_3}G(z_2/z_3)
+\frac{(z_2+z_3)(z_1^2-z_2z_3)}{z_1z_2z_3}G(z_1/z_2)
+\frac{(z_3+z_1)(z_2^2-z_3z_1)}{z_1z_2z_3},\notag
\end{align}
and $G(w/z)$ stands for the power series expansion of 
$-g(w,z)/g(z,w)$
in $w/z$.

Since the right hand sides contain infinite sums, a justification is needed. 
We consider these relations on an admissible graded vector space $W$ such that
$E(z)$, $K^\pm(z)$ act as graded currents. Then the above relations are well defined
in $\End(W)$. 
In this case we say that $W$ is an admissible $\K_1$ module.

We say that $W$ has level $\ell\in\C^\times$ if $C$ acts as the scalar $\ell$.

Consider a Heisenberg algebra with generators $\{H_r\}_{r \in \Z\backslash\{0\}}$ 
satisfying
\begin{align*}
[H_r, H_s]=-\delta_{r+s,0}\frac{1}{r} \frac{1+q_2^{-r}}{\kappa_r}\,.
\end{align*}
Let  $\F_2^D=\C[H_{-r}|r>0]\ket{q_2^{-{1}/{2}}}$ be its Fock representation 
{with vacuum vector $|q_2^{-{1}/{2}}\rangle$}.
The following assignment makes $\F_2^D$  an admissible $\K_1$ module: 
\begin{align*}
E(z) \mapsto 0\,,\quad K^\pm(z) \mapsto K^\pm(z)\,,\quad C\mapsto q_2^{-1/2}.
\end{align*}
We call $\F_2^D$ the boundary Fock module of type $D$ and color $2$.

The following is a special case of Proposition 4.1 in 
\cite{Feigin-Jimbo-Mukhin-Vilkoviskiy}, see also the proof of 
Corollary 3.7 in \cite{Feigin-Jimbo-Mukhin1}:
\begin{prop}\label{prop:coproK}
Let $M$ be an admissible $\E_1$ module of level $\ell$, such that
$q_2^{-1}\ell^2\neq \pm1,q_2^{\pm1}$. 
Then the following assignment makes $M\otimes \F_2^D$ an admissible 
$\K_1$ module of level $\ell q_2^{-1/2}$:
\begin{align*}
& E(z)\mapsto
e(C_2^{-1}z)\otimes K^+(z)+\tilde{f}(C_2 z) \otimes K^-(z)\,,
\\
&K^+(z)\mapsto
\psi^+(C_1^{-1}C_2^{-1}z)\otimes K^+(z)\,,
\\
&K^-(z)\mapsto 
\psi^-(C_2 z)^{-1} \otimes K^-(z)\,,
\\
&C\mapsto C \otimes C\,.
\end{align*}
\end{prop}
In particular, by choosing $M$ to be the $\E_1$ module
${\mathcal F}_{2}(\ssu_1) \otimes \cdots \otimes {\mathcal F}_{2}(\ssu_l)$
with $l\ge3$, 
we obtain an admissible ${\mathcal K}_1$ module of level $q_2^{(l-1)/2}$,
\begin{align}
\mathcal{F}^D(\ssu_1,\ldots,\ssu_l)
=
{\mathcal F}_{2}(\ssu_1) \otimes \cdots \otimes {\mathcal F}_{2}(\ssu_l)
\otimes {\mathcal F}_2^D\,.
\label{l+1 Fock}
\end{align}
The space \eqref{l+1 Fock} 
is the Fock space of the Heisenberg algebra generated by 
$l+1$ species of free bosons
\begin{align}
1\otimes\cdots\otimes
\overset{\underset{\smile}{i}}{h_r}\otimes\cdots\otimes 1\,, 
\quad
1\otimes\cdots\otimes 1\otimes H_r\,, 
\label{osc}
\end{align}
where $i=1,\ldots,l$ and $r\in \Z\setminus \{0\}$. 
We shall refer to them as ``oscillators''. 
We fix the vacuum vector $\ket{\ssu_1,\ldots,\ssu_l}
{=|\ssu_1, q_2^{{1}/{2}}\rangle \otimes 
|\ssu_2, q_2^{{1}/{2}}\rangle
\otimes \cdots
\otimes |\ssu_l, q_2^{{1}/{2}}\rangle \otimes
|q_2^{-{1}/{2}}\rangle}$
of \eqref{l+1 Fock} annihilated by oscillators with $r>0$.

In the following sections we discuss 
integrals of motion which act on this space. 

\section{Deformed $W$-algebra ${\mathcal W}_{q, t}(D_l)$} 
\label{section:deformed-W}

In this section,
we review the construction in \cite{Feigin-Jimbo-Mukhin-Vilkoviskiy}
of the deformed $W$-algebra 
${\mathcal W}_{q, t}(D_l)$
via the quantum algebras ${\mathcal E}_1$ and ${\mathcal K}_1$.
Compared with \cite{Frenkel-Reshetikhin} 
we have made some notational changes, see Remark \ref{rem:FR} below. 
Throughout this section we assume $l\ge 4$.

We use the following convention for the root system.

Let $\{\epsilon_j\}_{1\le j\le l}$ be an orthonormal basis of $\R^l$
relative to the standard inner product $(\epsilon_i,\epsilon_j)=\delta_{i,j}$. 
We set
\begin{align*}
&\alb_0=-\epsilon_1-\epsilon_2\,,\quad
\alb_i=\epsilon_i-\epsilon_{i+1}\quad 1\le i\le l-1\,,
\quad
\alb_l=\epsilon_{l-1}+\epsilon_l\,,
\\
&\omb_0=0\,,\quad 
\omb_i=\epsilon_1+\cdots+\epsilon_i\quad 1\le i\le l-2\,,
\\
&\omb_{l-1}=\frac{1}{2}(\epsilon_1+\cdots+\epsilon_{l-1}-\epsilon_l)\,,
\quad
\omb_{l}=\frac{1}{2}(\epsilon_1+\cdots+\epsilon_{l-1}+\epsilon_l)\,,
\\
& \bar{P}=\oplus_{i=1}^l\Z\omb_i\,, 
\quad
 \bar{Q}=\oplus_{i=1}^l\Z\alb_i\,.
\end{align*}
The Cartan matrix $(C_{i,j})_{i,j=0}^l$
of type $D^{(1)}_l$ is given by 
\begin{align*}
 C_{i,j}=(\bar{\alpha}_i,\bar{\alpha}_j)\,. 
\end{align*}
We have $\sum_{i=0}^la_i\alb_i=0$, where
\begin{align}
a_i=\begin{cases}
     1&\text{for $i=0,1,l-1,l$\,,}\\
     2&\text{for $2\le i\le l-2$}\,,\\
    \end{cases} 
\label{def:Kac}
\end{align}
are called Kac label.

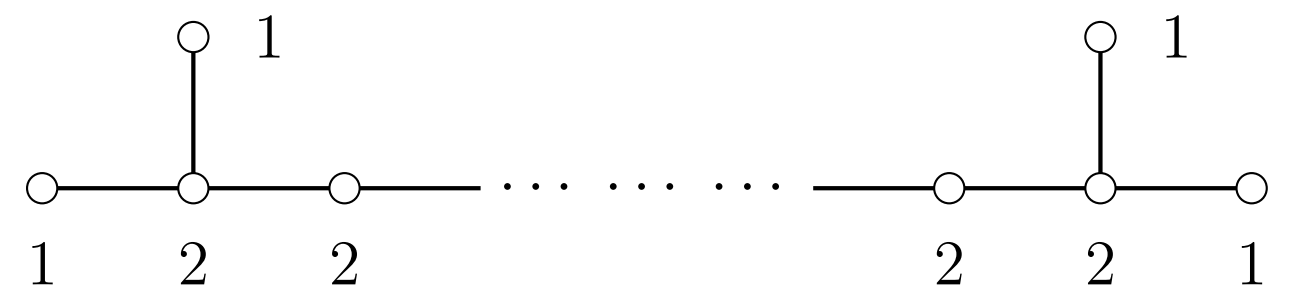
\begin{figure}[ht]
\begin{align}
\begin{tikzpicture}{scale=1}
\draw(-40mm,0mm)circle(1mm);
\draw(-30mm,0mm)circle(1mm);
\draw(-20mm,0mm)circle(1mm);
\draw(-30mm,10mm)circle(1mm);
\draw[thick](-39mm,0mm)--(-31mm,0mm);
\draw[thick](-29mm,0mm)--(-21mm,0mm);
\draw[thick](-19mm,0mm)--(-11mm,0mm);
\draw[thick](-30mm,1mm)--(-30mm,9mm);
\draw[thick](39mm,0mm)--(31mm,0mm);
\draw[thick](29mm,0mm)--(21mm,0mm);
\draw[thick](19mm,0mm)--(11mm,0mm);
\draw[thick](30mm,1mm)--(30mm,9mm);
\draw(40mm,0mm)circle(1mm);
\draw(30mm,0mm)circle(1mm);
\draw(20mm,0mm)circle(1mm);
\draw(30mm,10mm)circle(1mm);
\draw(-7mm,0mm)node{$\cdots$};
\draw(0mm,0mm)node{$\cdots$};
\draw(7mm,0mm)node{$\cdots$};
\draw(-40mm,-5mm)node{$1$};
\draw(-30mm,-5mm)node{$2$};
\draw(-20mm,-5mm)node{$2$};
\draw(-25mm,10mm)node{$1$};
\draw(40mm,-5mm)node{$1$};
\draw(30mm,-5mm)node{$2$};
\draw(20mm,-5mm)node{$2$};
\draw(35mm,10mm)node{$1$};
\end{tikzpicture}
\notag
\end{align}
\caption{Dynkin diagram and Kac label for $D_l^{(1)}$, $l\ge4$}
\label{fig:Dynkin}
\end{figure}

\subsection{Basic $W$-currents}\label{sec:basicW}

For an element $x\in \mathcal{K}_1$, we shall write its action 
on \eqref{l+1 Fock} as $\pi (x)$.

The current $T_1(z)=\pi(E(z))$ acts as a sum of $2l$ vertex operators in 
$l+1$ oscillators \eqref{osc}
\begin{align}
T_1(z)=\sum_{i=1}^l \Lambda_i(z)+\sum_{i=1}^l \Lambda_{\bar{i}}(z)\,,
\end{align}
where for $i=1,2,\ldots,l$
\begin{align*}
&\Lambda_i(z)=1 \otimes \cdots \otimes 1 \otimes 
\overset{\text{$i$-th}}{\overbrace{
V_2(c_i z; \ssu_i)}}
\otimes 
\psi^+(q_2^{-1/2}c_{i+1}z)\otimes \cdots \otimes \psi^+(q_2^{-1/2}c_l z)\otimes K^+(z)\,,
\\
&\Lambda_{\bar{i}}(z)=1 \otimes \cdots \otimes 1 \otimes 
:\overset{\text{$i$-th}}{\overbrace{
V_2^{-1}(c_i^{-1} z; \ssu_i)}}:
\otimes 
\psi^-(c_{i+1}^{-1} z)^{-1} \otimes \cdots \otimes \psi^-(c_l^{-1} z)^{-1} \otimes K^-(z)\,,
\end{align*}
and we set 
\begin{align*}
c_i=q_2^{(-l+i+1)/2}\,. 
\end{align*}

In what follows, we order the indices as 
\begin{align*}
1\prec\cdots\prec l\prec \bar{l}\prec\cdots\prec \bar{1} \,
\end{align*}
and set $I=\{1,\ldots,l,\bar{l},\ldots,\bar{1}\}$.

Using the dressing operator 
\begin{align*}
&K_p^+(z)=\prod_{s=0}^\infty K^+(p^{-s}z)\,,\quad
p=q^{(l-2)\beta}\,,
\end{align*}
we define the basic $W$-current ${\mathbf T}_1(z)$ of 
${\mathcal W}_{q,t}(D_l)$ as follows.
\begin{align*}
{\mathbf T}_1(z)=\pi(E(z) K_p^+(z)^{-1})
=\sum_{i=1}^l \bLambda_i(z)+\sum_{i=1}^l \bLambda_{\bar{i}}(z)\,,
\end{align*}
where
\begin{align*}
\bLambda_i(z)=\Lambda_i(z) \pi(K_p^+(z))^{-1}\,,
\quad
\bLambda_{\bar{i}}(z)=\Lambda_{\bar{i}}(z) \pi(K_p^+(z))^{-1}\,.
\end{align*}

Quite generally, when the product of vertex operators $V(z),W(z)$ 
has the form
\begin{align*}
&V(z)W(w)=\varphi_{V,W}(w/z):V(z)W(w):\,,\\
&\varphi_{V,W}(z)
=\prod_{i=1}^M \frac{(\beta_i z; q_i)_\infty}{(\alpha_i z; p_i)_\infty}
=\exp\left(\sum_{i=1}^M \sum_{n>0}\frac{1}{n}
\left(\frac{\alpha_i^n}{1-p_i^n}-\frac{\beta_i^n}{1-q_i^n}\right)z^n \right)\,,
\end{align*}
we write this as 
\begin{align*}
{\mathcal C}(V(z), W(w))=\sum_{i=1}^M
\left(\frac{\alpha_i}{1-p_i}-\frac{\beta_i}{1-q_i}\right)\,.
\end{align*}
In this shorthand notation we have 
\begin{align*}
&{\mathcal C}(\Lambda_i(z), \Lambda_j(w))=
\begin{cases}
-\kappa_1 &i,j\in I\,,
i \prec  j, i\neq \bar{j},\\ 
0&i,j\in I\,, 
i \succ j, i \neq \bar{j},
\end{cases}
\\
&{\mathcal C}(\Lambda_i(z), \Lambda_i(w))
=-(1-q_1)(1-q_3)\,,\quad 
i\in I\,,
\\
&{\mathcal C}(\Lambda_i(z), \Lambda_{\bar{i}}(w))
=(1-q_1)(1-q_3)(q_2^{l-i}-1+q_2)\,,
\quad i=1,\ldots,l\,,
\\
&{\mathcal C}(\Lambda_{\bar{i}}(z),\Lambda_i(w))
=(1-q_1)(1-q_3)q^{-l+i+1}\,,
\quad i=1,\ldots,l\,,
\end{align*}
and
\begin{align*}
{\mathcal C}(\bLambda_i(z), \bLambda_j(w)))
={\mathcal C}(\Lambda_i(z), \Lambda_j(w))+\frac{\kappa_1}{1-p}\,, 
\quad i,j\in I\,.
\end{align*}
Note the symmetry
\begin{align*}
{\mathcal C}(\bLambda_i(z), \bLambda_j(w))
=
\left.{\mathcal C}(\bLambda_{\bar{i}}(z), \bLambda_{\bar{j}}(w))\right|_{q \to q^{-1}},
\quad i,j\in I\,.
\end{align*}
More explicit formulas can be found in
Appendix \ref{appendix:normal-ordering}.

\subsection{Root current}

With each neighboring pair of Fock spaces appearing in \eqref{l+1 Fock}, 
we associate currents ${A}_i(z)$ given by
\begin{align*}
&{ A}_0(z)=:\frac{\bLambda_{\bar{2}}(p^{-\frac{1}{2}}z)}{\bLambda_{1}(p^{\frac{1}{2}}z)}:
=
:\frac{\bLambda_{\bar{1}}(p^{\frac{1}{2}}z)}{\bLambda_{2}(p^{-\frac{1}{2}}z)}:\,,
\\
&
{A}_i(z)=:\frac{\bLambda_i(c_i^{-1}z)}{\bLambda_{i+1}(c_i^{-1}z)}:=:\frac{\bLambda_{\overline{i+1}}(c_i z)}{\bLambda_{\overline{i}}(c_i z)}:
\quad \text{for $1\le i \le l-1$}\,,
\\
&{A}_l(z)=:\frac{\bLambda_{l-1}(z)}{\bLambda_{\overline{l}}(z)}:
=
:\frac{\bLambda_{l}(z)}{\bLambda_{\overline{l-1}}(z)}:\,.
\end{align*}
Their contractions read as 
\begin{align*}
{\mathcal C}({A}_i(z), {A}_j(w))=
- (q_1^{1/2}-q_1^{-1/2})  (q_3^{1/2}-q_3^{-1/2}) 
C_{i,j}(q_1,q_2,q_3)\,,
\end{align*}
where
$C(q_1,q_2,q_3)=(C_{i,j}(q_1,q_2,q_3))_{i,j=0}^l$  
is the deformed Cartan matrix 
\begin{align*}
C_{i,j}(q_1,q_2,q_3)=
\begin{cases}
q_2^{1/2}+q_2^{-1/2}&\text{if $C_{i,j}=2$}\,,\\ 
-1&\text{if $C_{i,j}=-1$}\,,\\
0&\text{if $C_{i,j}=0$}\,.
\end{cases}
\end{align*}

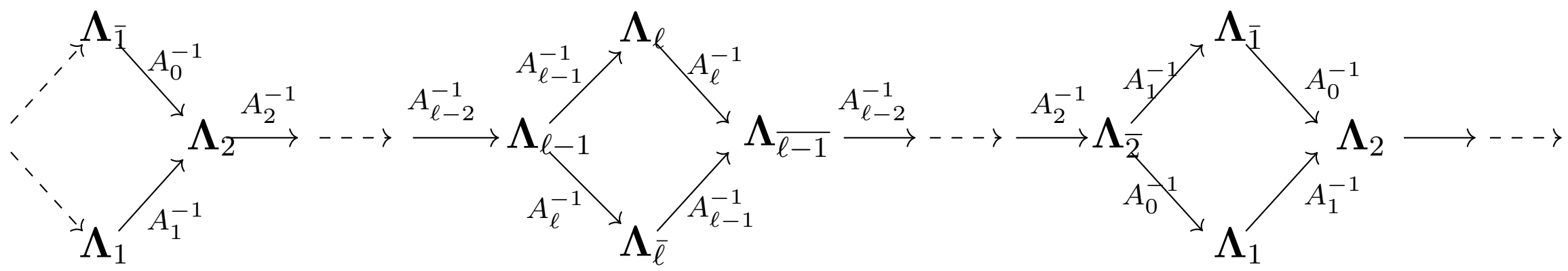
\begin{figure}[ht]
\begin{align*}
\begin{tikzpicture}[scale=0.7]
\draw[->,dashed](-6.8,0.2)--(-5.8,1.3);
\draw[->,dashed](-6.8,-0.2)--(-5.8,-1.3);
\draw(-5.5,1.5)node{$\bLa_{\bar 1}$};
\draw(-5.5,-1.5)node{$\bLa_1$};
\draw[->](-5.3,1.3)--(-4.4,0.3);
\draw(-4.5,0.6)node[above]{$\scriptstyle A_0^{-1}$};
\draw(-4.5,-0.7)node[below]{$\scriptstyle  A_1^{-1}$};
\draw[->](-5.3,-1.3)--(-4.4,-0.3);
\draw(-4,0)node{$\bLa_{2}$};
\draw[->](-3.8,0)--(-2.8,0);
\draw(-3.2,0)node[above]{$\scriptstyle  A_2^{-1}$};
\draw[->,dashed](-2.5,0)--(-1.5,0);
\draw[->](-1.2,0)--(-0,0);
\draw(-0.8,0)node[above]{$\scriptstyle  A_{\ell-2}^{-1}$};
\draw(0.7,0)node{$\bLa_{\ell-1}$};
\draw[->](0.7,0.2)--(1.7,1.2);
\draw[->](0.7,-0.2)--(1.7,-1.2);
\draw(1.4,1)node[left]{$\scriptstyle  A_{\ell-1}^{-1}$};
\draw(1.4,-1)node[left]{$\scriptstyle  A_{\ell}^{-1}$};
\draw(2,1.5)node{$\bLa_{\ell}$};
\draw(2,-1.5)node{$\bLa_{\bar\ell}$};
\draw[->](2.2,1.3)--(3.2,0.2);
\draw[->](2.2,-1.3)--(3.2,-0.2);
\draw(2.4,1)node[right]{$\scriptstyle  A_{\ell}^{-1}$};
\draw(2.4,-1)node[right]{$\scriptstyle  A_{\ell-1}^{-1}$};
\draw(4,0)node{$\bLa_{\overline{\ell-1}}$};
\draw[->](4.8,0)--(5.8,0);
\draw(5.2,0)node[above]{$\scriptstyle  A_{\ell-2}^{-1}$};
\draw[->,dashed](6,0)--(7,0);
\draw[->](7.2,0)--(8.2,0);
\draw(7.8,0)node[above]{$\scriptstyle  A_{2}^{-1}$};
\draw(8.6,0)node{$\bLa_{\overline{2}}$};
\draw[->](8.8,0.2)--(9.8,1.3);
\draw[->](8.8,-0.2)--(9.8,-1.3);
\draw(10.3,1.5)node{$\bLa_{\bar 1}$};
\draw(10.3,-1.5)node{$\bLa_{1}$};
\draw[->](10.4,1.3)--(11.4,0.2);
\draw[->](10.4,-1.3)--(11.4,-0.2);
\draw(9.7,0.8)node[left]{$\scriptstyle  A_{1}^{-1}$};
\draw(9.7,-0.8)node[left]{$\scriptstyle  A_{0}^{-1}$};
\draw(11,0.8)node[right]{$\scriptstyle  A_{0}^{-1}$};
\draw(11,-0.8)node[right]{$\scriptstyle  A_{1}^{-1}$};
\draw(12,0)node{$\bLa_{2}$};
\draw[->](12.6,0)--(13.6,0);
\draw[->,dashed](13.8,0)--(14.8,0);
\end{tikzpicture}
\end{align*}
\caption{Extended vector representation of $D_l$}
\label{fig:vector}
\end{figure}

The basic $W$-current ${\bf T}_1(z)$ is pictured as 
the $qq$ character of the vector representation of $D_l$
 \cite{Frenkel-Reshetikhin} 
extended periodically. 
See Figure \ref{fig:vector}.

\subsection{Screening current}

Introduce ``zero mode'' operators for $1\le i\le l$
\begin{align*}
&e^{Q_{\bar{\alpha}_i}}:
\mathcal{F}^D(\ssu_1,\ldots,\ssu_i,\ldots,\ssu_l)
\longrightarrow
\mathcal{F}^D(\ssu_1,\ldots,q^{1-\beta}\ssu_i,\ldots,\ssu_l)\,,
\\
&P_{\bar{\omega}_i}:\mathcal{F}^D(\ssu_1,\ldots,\ssu_l)
\longrightarrow
\mathcal{F}^D(\ssu_1,\ldots,\ssu_l)\,,
\end{align*}
which commute with the oscillators and 
\begin{align*}
&e^{Q_{\bar{\alpha}_i}}\ket{\ssu_1,\ldots,\ssu_i,\ldots,\ssu_l}
 =\ket{\ssu_1,\ldots,q^{1-\beta}\ssu_i,\ldots,\ssu_l} \,,\\
&
P_{\bar{\omega}_i}\ket{\ssu_1,\ldots,\ssu_l}
=\frac{\log \ssu_i}{\log q} 
\ket{\ssu_1,\ldots,\ssu_l}\,.
\end{align*}
We set
$e^{Q_{\bar{\alpha}_0}}=\prod_{i=1}^l(e^{Q_{\bar{\alpha}_i}})^{-a_i}$ and
$P_{\bar{\alpha}_i}=\sum_{i=1}^lC_{i,j}P_{\bar{\omega}_j}$ for $0\le i\le l$.
We have 
\begin{align*}
[P_{\bar{\alpha}_i},e^{Q_{\bar{\alpha}_j}}]=
(1-\beta)C_{i,j} e^{Q_{\bar{\alpha}_j}}\,,
\quad 0\le i,j\le l\,.
\end{align*}

For $0\le i\le l$, define the screening current $S_i(z)$ by 
\begin{align*}
&S_i(z)=
e^{Q_{\alpha_i}}z^{P_{\alpha_i}+1-\beta}S_i^{osc}(z)\,,\\
&{A}_i(z)=q^{1-\beta}:\frac{S_i(q^{1/2}z)}{S_i(q^{-1/2}z)}:\,.
\end{align*}
Here the oscillator part $S_i^{osc}(z)$ 
is a vertex operator consisting of oscillators.
It is uniquely fixed by the second requirement.

We have for $0\le i,j\le l$
\begin{align}\label{cont-screening}
S_i(z)S_j(w)=:S_i(z)S_j(w):\times
\begin{cases}
z^{2(1-\beta)}{\displaystyle \left(1-\frac{w}{z}\right)
\frac{\left(q^\beta w/z;q\right)_\infty}
{\left(q^{1-\beta} w/z; q\right)_\infty}}
&\text{if $C_{i,j}=2$}\,,
\\ 
z^{-(1-\beta)}{\displaystyle 
\frac{\left(q^{1-\beta/2}w/z; q\right)_\infty}
{\left(q^{\beta/2}w/z; q\right)_\infty}}&
\text{if $C_{i,j}=-1$}\,,
\\
1&\text{if $C_{i,j}=0$}\,.
\end{cases}
\end{align}
The following commutation relations hold as meromorphic functions:
\begin{align}
\theta\left(u-v+\beta\right)\theta(u-v)S_i(z)S_i(w)
&=\theta(v-u+\beta)\theta(v-u)S_i(w)S_i(z)\,,
\label{SS1}\\
\theta\left(u-v+\beta/2\right)^{-1}S_i(z)S_j(w)
&=\theta\left(v-u+\beta/2\right)^{-1}S_j(w)S_i(z)
\qquad \text{if $C_{i,j}=-1$},
\label{SS2}\\
S_i(z)S_j(w)&=S_j(w)S_i(z)
\qquad \text{if $C_{i,j}=0$},
\label{SS3}
\end{align}
where $0\le i,j\le l$ and $z=q^u, w=q^v$.

Contraction rules between $\bLambda_i(w)$ and $S_j(z)$ are summarized in 
TABLE \ref{table:pole1} and \ref{table:pole2} in Appendix \ref{appendix:normal-ordering}.
\medskip

\begin{rem}\label{rem:FR}
{\rm
Our notation for $\mathcal{W}_{q,t}(D_l)$
is related to the one for $\mathcal{W}^{FR}_{q^{FR},t^{FR}}(D_l)$
in \cite{Frenkel-Reshetikhin}
as follows:
$q^{FR}=q_3^{1/2}$, $t^{FR}=q_1^{-1/2}$, $\beta^{FR}=1-\beta$.
Up to  overall constants and shifts of $z$, the currents
$A_i(z), S_i(z)$ 
are $A_i^{FR}(z),S_i^{+,FR}(z)$ 
in \cite{Frenkel-Reshetikhin} for $i=1,\ldots,l$. 

In comparison with \cite{Frenkel-Reshetikhin}, 
we have in addition $A_0(z),S_0(z)$ which involve an
 extra free boson $\pi(K^{\pm}(z))$ 
commuting with $A_i(z), S_i(z)$ with $i=1,\ldots,l$. 
They play a key role in constructing 
local IMs as explained in the next section. 
}
\end{rem}

\section{Main results}
\label{section:main-result}

In this section we present two infinite sets of operators  
$\{{\bf I}_M\}_{M=1}^{\infty}$ and $\{{\bf G}_N(\vartheta)\}_{N=1}^{\infty}$ 
which act on the space $\mathcal{F}^D(\ssu_1,\ldots,\ssu_l)$ in 
\eqref{l+1 Fock}. 
We then state main theorems.

\subsection{Local integrals of motion}

Following \cite{Feigin-Jimbo-Mukhin-Vilkoviskiy},
we give an infinite set of operators which we call local integrals of motion.
For each $M \geq 1$, we define ${\bf I}_M$ as follows.
\begin{align}
{\bf I}_M=\prod_{i=1}^M \int_{|w_i|=R}
\underline{dw_i}
\overset{\curvearrowright}
{\prod_{1\leq i \leq M}}{\bf T}_1(w_i) \cdot
\omega_M(w_1, w_2, \ldots, w_M).
\label{def:local}
\end{align}
Here and after, we use the abbreviation
\begin{align*}
\underline{dw}=\frac{dw}{2\pi\sqrt{-1}w} 
\end{align*}
and symbols for ordered products
\begin{align*}
\overset{\curvearrowright}
{\prod_{1\leq i \leq M}}A_i=A_1A_2\cdots A_M\,,
\quad
\overset{\curvearrowleft}
{\prod_{1\leq i \leq M}}A_i=A_M\cdots A_2A_1\,.
\end{align*}
The kernel function $\omega_M$ is defined as
\begin{align}
\omega_M(w_1, w_2,\ldots, w_M)=\prod_{1\leq i<j \leq M}
\omega(w_j/w_i)\,,
\quad 
\omega(w)=\frac{\Theta_p(w)\Theta_p(q_2^{-1}w)}
{\Theta_p(q_1 w)\Theta_p(q_3 w)}\,. 
\label{def:omega}
\end{align}
We set $\omega_1(w)=1$.

The integral is taken on a common circle $|w_i|=R>0$.
The definition is independent of the choice of $R$.

\subsection{Non-local integrals of motion}

We give another set of operators acting on \eqref{l+1 Fock}
which we call non-local integrals of motion. 

Let $\nu_i\in\C$, $i=1,\ldots,l$, and $\nu_0=0$. Set
\begin{align}
\mu_i=\sum_{j=1}^l C_{i,j}\nu_j, \quad {i=0,1,\ldots,l}. \label{munu}
\end{align}
Denote by $X_l=X_l(\nu_1,\ldots,\nu_l)$ the space of entire functions
 $\vartheta(v_0,v_1,\ldots, v_l)$ that satisfies
\begin{align*}
&\vartheta(v_0,\ldots, v_i+1,\ldots, v_l)
=\vartheta(v_0,\ldots, v_i, \ldots, v_l)\,,\quad
0\leq i \leq l\,,
\\
&\vartheta(v_0,\ldots, v_i+\tau,\ldots, v_l)
=e^{-2\pi\sqrt{-1}(\sum_{j=0}^l C_{i,j}v_j-\mu_{i}+\tau)}
\vartheta(v_0,\ldots, v_i, \ldots, v_l)\,,\quad 
0\leq i \leq l\,.
\end{align*}
\begin{lem}\label{lem:theta1}
(i) If ${\vartheta} \in X_l$ then 
\begin{align*}
{\vartheta}(u_0+a_0 u,\ldots, u_i+a_i u,\ldots, u_l +a_l u)
={\vartheta}(u_0,\ldots, u_i,\ldots,u_l)\,,\quad u \in \C\,.
\end{align*}
(ii) We have $\dim X_l=|\bar{P}/\bar{Q}|=4$. A basis of $X_l$ is given by
$\{\vartheta_{\bar\omega_i}\mid i=0,1,l-1,l\}$, 
where
\begin{align*}
\vartheta_{\bar\omega}(v_0,\ldots,v_l) 
=\sum_{\alpha\in \bar{Q}+\bar{\omega}
}e^{-\frac{\pi\sqrt{-1}}{\tau}(\alpha-v,\alpha-v)}\,,
\quad v=\sum_{j=0}^l\bar{\alpha}_j(v_j-\nu_j)\,.
\end{align*}
\end{lem}
\begin{proof}
(i) Let 
$\varphi(u)={\vartheta}(u_0+a_0 u,\ldots, u_i+a_i u,\ldots, u_l +a_l u)$.
Then $\varphi(u)$ is entire and 
satisfies $\varphi(u+1)=\varphi(u)$, $\varphi(u+\tau)=\varphi(u)$.
Hence $\varphi(u)$ is a constant: $\varphi(u)=\varphi(0)$.

(ii) By (i), $\vartheta\in X_l$ can be written as
$$
\vartheta(v_0,\ldots,v_i,\ldots, v_l)=\bar{\vartheta}(v_1-a_1v_0-\nu_1,\ldots, v_i-a_i v_0-\nu_i,\ldots, v_l-a_l v_0-\nu_l),
$$
where $\bar{\vartheta}(v_1,\ldots,v_l)$ is entire, 
$1$-periodic in $v_1,\ldots,v_l$, and satisfies
\begin{align*}
&\bar{\vartheta}(v_1,\ldots, v_i+\tau,\ldots, v_l)
=e^{-2\pi\sqrt{-1}(\sum_{j=1}^l C_{i,j}v_j+\tau)}
\bar{\vartheta}(v_1,\ldots, v_i, \ldots, v_l)\,,\quad 
1\leq i \leq l\,.
\end{align*}
Since $(C_{i,j})_{i,j=1}^l$ is positive definite, 
the assertion about the dimension and the basis is a well-known classical fact.
\end{proof}

Choose
\begin{align}
\nu_j=\frac{\log \ssu_j}{\log q}\,,\quad j=1,\ldots,l\,. 
\label{choice-nu}
\end{align}
For each $N \geq 1$ and $\vartheta \in X_l(\nu_1,\ldots,\nu_l)$, following \cite{Jimbo-Kojima}, we define
${\bf G}_N(\vartheta)$ as follows.
\begin{align}
{\bf G}_N(\vartheta)&
=\prod_{i=0}^l \prod_{a=1}^{N a_i} \int_{|z_{i,a}|=R}
\underline{dz_{i,a}}
~\overset{\curvearrowright}{\prod_{0\leq i \leq l}}~
\overset{\curvearrowright}{\prod_{1\leq a \leq N a_i}}
S_i(z_{i,a})
\cdot 
{\bf h}_N({\bf u}_0, {\bf u}_1,\ldots, {\bf u}_l; \vartheta).
\label{def:non-local}
\end{align}
Here $a_i$ are given by \eqref{def:Kac}, $z_{i,a}=q^{u_{i,a}}$, and 
\begin{align*}
{\bf u}_i=(u_{i,1}, u_{i,2},\ldots, u_{i, Na_i})\,,
\quad
\overline{u}_i=\sum_{a=1}^{N a_i}u_{i,a}\,.
\end{align*}
The kernel function $\mathbf{h}_N$ is given by
\begin{align*}
{\bf h}_N({\bf u}_0, {\bf u}_1,\ldots, {\bf u}_l; \vartheta)=
\frac{\displaystyle
\prod_{i=0}^l
\prod_{1\leq a<b \leq N a_i}
\theta\left(u_{i,a}-u_{i,b}+\beta \right)\theta\left( u_{i,a}-u_{i,b} \right)
}{
\prod\limits_{\substack{\scriptscriptstyle 0\leq i < j \leq l \\{C{i j}=-1}}}
\prod_{a=1}^{N a_i}
\prod_{b=1}^{N a_j}
\theta\bigl(u_{i,a}-u_{j,b}+\beta/2\bigr)}
~\vartheta\left(
\overline{u}_0,\ldots,\overline{u}_i,\ldots,\overline{u}_l\right)\,.
\end{align*}
The condition \eqref{choice-nu} ensures that the integrand is 
a single-valued function in each $z_{i,a}$. 
Using the contraction rules \eqref{cont-screening} for $S_i(z)$, one
can check that the only poles of the integrand are 
$q^{\beta/2+n}z_{i,a}/z_{j,b}=1$ for some $a,b$, $i\neq j$, 
and $n\ge0$. 
Provided $|q|<1,|q^\beta|<1$, the integral
taken on a common circle $|z_{i,a}|=R$ is well-defined and
is independent of the choice of $R>0$.

The product of screening operators in \eqref{def:non-local}
is so chosen that the shift
operators $e^{Q_{\bar{\alpha}_i}}$ cancel out. Hence $\mathbf{G}_N(\vartheta)$
preserves the space $\F^D(\ssu_1,\ldots,\ssu_l)$. 

\subsection{Main theorems}

The following result has been known.
\begin{thm}\cite{Feigin-Jimbo-Mukhin-Vilkoviskiy}
The local integrals of motion $\{{\bf I}_M\}$ mutually commute,
\begin{align}
[{\mathbf I}_{M_1}, {\mathbf I}_{M_2}]=0\quad
{\it for}\quad M_1, M_2 \geq 1\,.
\notag
\end{align}
\end{thm}

We state the main results of this paper. 
\begin{thm}\label{thm:commutativity}
The non-local integrals of motion $\{{\bf G}_N (\vartheta)\}$ 
mutually commute, 
\begin{align*}
[{\mathbf G}_{N_1}(\vartheta_1), {\mathbf G}_{N_2}(\vartheta_2)]=0
\quad
\text{for $N_1,N_2 \geq 1$ and
$\vartheta_1, \vartheta_2 \in X_l(\nu_1,\ldots,\nu_l)$}\,.
\end{align*}
\end{thm}

\vspace{5mm}

\begin{thm}\label{thm:duality}
The local integrals of motion $\{{\bf I}_M\}$ 
and the non-local integrals of motion $\{{\bf G}_N (\vartheta)\}$ 
mutually commute,
\begin{align*}
[{\mathbf I}_M, {\mathbf G}_N(\vartheta)]=0\quad
{\it for}\quad M,N \geq 1 \quad {\it and}\quad \vartheta \in 
X_l(\nu_1,\ldots,\nu_l)\,.
\end{align*}
\end{thm}

The rest of the paper is devoted to proofs of Theorem
\ref{thm:commutativity} and Theorem \ref{thm:duality}.

\section{Proof of Theorem \ref{thm:commutativity}}
\label{section:commutativity}

In this section we show $[{\bf G}_M(\vartheta_1), {\bf G}_N(\vartheta_2)]=0$ by a direct calculation.

\subsection{An identity of theta functions}

We prove the commutativity assuming $0<\beta<1$. The general 
case follows by analytic continuation in $\beta$.
Under this assumption, the product $S_i(z)S_j(w)$ converges absolutely 
in the domain $|w/z|\le 1$. Therefore, in the product of operators
\begin{align*}
{\bf G}_M(\vartheta_1){\bf G}_N(\vartheta_2)
&=\prod_{i=0}^l \prod_{a=1}^{M a_i} \int_{|z_{i,a}|=R}
\underline{dz_{i,a}}
\prod_{j=0}^l \prod_{b=1}^{N a_j} \int_{{|w_{j,b}|}=R}
\underline{dw_{j,b}}
\\
&\times
~\overset{\curvearrowright}{\prod_{0\leq i \leq l}}~
\overset{\curvearrowright}{\prod_{1\leq a \leq M a_i}}
S_i(z_{i,a})
~\overset{\curvearrowright}{\prod_{0\leq j \leq l}}~
\overset{\curvearrowright}{\prod_{1\leq b \leq N a_i}}
S_j(w_{j,b})\\
&\times
{\bf h}_M({\bf u}_0, {\bf u}_1,\ldots, {\bf u}_l; \vartheta_1)
{\bf h}_N({\bf v}_0, {\bf v}_1,\ldots, {\bf v}_l; \vartheta_2)\,,
\end{align*}
one can choose the common radius $R$ for the circles of integration. 

Using the commutation relations \eqref{SS1}-\eqref{SS3} one can further 
symmetrize the integrand. 

Hence, Theorem \ref{thm:commutativity} is reduced to
the following identity of theta functions:

\begin{prop}\label{prop:theta1}
Let $\nu_1, \ldots, \nu_l\in\C$ and $l \geq 3$.
The following identity holds for all $M,N\ge 1$
and ${\vartheta}_1, {\vartheta}_2 \in X_l(\nu_1,\ldots,\nu_l)$:
\begin{align}
&
\sum_{I_0, J_0}\cdots \sum_{I_l, J_l} 
{\vartheta}_1(\overline{u}_{0,I_0},\ldots, \overline{u}_{l, I_l})
{\vartheta}_2(\overline{u}_{0,J_0},\ldots, \overline{u}_{l, J_l})
\notag\\
&\times
\frac{
\prod\limits_{\substack{\scriptscriptstyle 0\leq i<j \leq l\\ \scriptscriptstyle C_{ij}=-1}}
\Bigl(
\prod\limits_{\substack{\scriptscriptstyle 
a\in I_i \\ \scriptscriptstyle b \in J_j}}
\theta\bigl(u_{i,a}-u_{j,b}+\beta/2\bigr)
\prod\limits_{\substack{\scriptscriptstyle a\in I_j \\ \scriptscriptstyle  b \in J_i}}
\theta\bigl(u_{j,a}-u_{i,b}+\beta/2\bigr)
\Bigr)}
{
\prod\limits_{i=0}^l
\prod\limits_{\substack{\scriptscriptstyle a \in I_i \\ \scriptscriptstyle  
b \in J_i}}
\theta(u_{i,a}-u_{i,b}+\beta)
\theta(u_{i,a}-u_{i,b})
}
\notag\\
=&
\sum_{I_0, J_0}\cdots \sum_{I_l, J_l} 
{\vartheta}_1(\overline{u}_{0,I_0},\ldots, \overline{u}_{l, I_l})
{\vartheta}_2(\overline{u}_{0,J_0},\ldots, \overline{u}_{l, J_l})
\notag\\
&\times
\frac{
\prod\limits_{\substack{\scriptscriptstyle  0\leq i<j \leq l \\ 
\scriptscriptstyle  C_{ij}=-1}}
\Bigl(
\prod\limits_{\substack{\scriptscriptstyle  a\in J_i \\
\scriptscriptstyle  b \in I_j}}
\theta\bigl(u_{i,a}-u_{j,b}+\beta/2\bigr)
\prod\limits_{\substack{\scriptscriptstyle  a\in J_j \\ \scriptscriptstyle  b \in I_i}}
\theta\bigl(u_{j,a}-u_{i,b}+\beta/2\bigr)
\Bigr)}{
\prod\limits_{i=0}^l
\prod\limits_{\substack{\scriptscriptstyle  
a \in J_i
\\ \scriptscriptstyle  
b \in I_i}}
\theta(u_{i,a}-u_{i,b}+\beta)
\theta(u_{i,a}-u_{i,b})}.
\label{eqn:theta1}
\end{align}
Here the sum is taken over all partitions
$I_i \cup J_i=\{1,2,\ldots, (M+N)a_i\}$
satisfying
$|I_i|=M a_i$, $|J_i|=N a_i$, $0\leq i \leq l$.
We set $\overline{u}_{i,I}=\sum_{a \in I} u_{i,a}$.
\end{prop}

\subsection{Proof of Proposition \ref{prop:theta1}}

The proof of Proposition \ref{prop:theta1} 
is an adaptation of the one 
in \cite{Kojima-Shiraishi, Feigin-Jimbo-Mukhin1}. 
We use induction on $l$. The base of induction is $l=3$, 
in which case the identity (\ref{eqn:theta1}) 
reduces 
to that of $A_3^{(1)}$ proved in \cite{Kojima-Shiraishi, Feigin-Jimbo-Mukhin1}.

We begin with a simple fact which 
appears as Lemma A.3 in \cite{Feigin-Jimbo-Mukhin1}. 
For completeness we include its proof.
\begin{lem}\label{lem:theta2}
Let $f(u_0,\ldots,u_L)$ be an entire function 
with the quasi-periodicity properties
\begin{align}
&f(u_0,\ldots, u_i+1,\ldots, u_L)
=f(u_0,\ldots, u_i,\ldots, u_L)\,,\quad 0\leq i \leq L\,,
\notag\\
&f(u_0,\ldots, u_i+\tau,\ldots,u_L)=
e^{-2\pi\sqrt{-1}m \sum_{i=0}^L u_i+k}
f(u_0,\ldots, u_i,\ldots, u_L)\,,\quad 0\leq i \leq L,
\notag
\end{align}
with some $m\in\Z_{\ge0}$ and $k \in {\bf C}$.
Then $f$ is a function of $\overline{u}=\sum_{i=0}^L u_i$.
\end{lem}
\begin{proof}
Fix $0\leq i \neq j \leq L$ and 
set $\varphi(u)=f(u_0,\ldots, u_i+u, \ldots, u_j-u,\ldots,u_L)$.
Then $\varphi(u)$ is an entire function 
and satisfies $\varphi(u+1)=\varphi(u)$, $\varphi(u+\tau)=\varphi(u)$.
Hence $\varphi(u)$ is a constant $\varphi(u)=\varphi(0)$.
We obtain $f(\overline{u},0,\ldots,0)=f(u_0,u_1,\ldots, u_L)$.
\end{proof}

Denote each side of the identity (\ref{eqn:theta1}) by 
${\rm LHS}^{(l)}_{M,N}({\bf u}_0, {\bf u}_1,\ldots, {\bf u}_l |\vartheta_1, \vartheta_2)$, 
${\rm RHS}^{(l)}_{M,N}({\bf u}_0, {\bf u}_1,\ldots, {\bf u}_l |\vartheta_1, \vartheta_2)$, 
respectively.
Let 
\begin{align}
{\Phi}_{M,N}^{(l)}({\bf u}_0, {\bf u}_1,\ldots, {\bf u}_l|\vartheta_1, \vartheta_2)=
{\rm LHS}_{M,N}^{(l)}({\bf u}_0, {\bf u}_1, \ldots, {\bf u}_l |\vartheta_1, \vartheta_2)-
{\rm RHS}_{M,N}^{(l)}({\bf u}_0, {\bf u}_1,\ldots, {\bf u}_l |\vartheta_1, \vartheta_2)\,.
\notag
\end{align}
We set ${\Phi}_{M,N}^{(l)}({\bf u}_0, {\bf u}_1, \ldots, {\bf u}_l |\vartheta_1, \vartheta_2)=0$ if $M=0$ or $N=0$.

We retain the notation $\mu_i$ in \eqref{munu}. 
\begin{lem}\label{lem:theta5}
The function 
${\Phi}_{M,N}^{(l)}({\bf u}_0,  \ldots, {\bf u}_l|\vartheta_1, \vartheta_2)$ 
satisfies
the quasi-periodicity 
in each  $u_{i,a}$ with $0\leq i \leq l$, $1 \leq a \leq (M+N)a_i$:
\begin{align}
&
{\Phi}_{M,N}^{(l)}(\cdots, u_{i,a}+1,\ldots|\vartheta_1, \vartheta_2)=
{\Phi}_{M,N}^{(l)}(\ldots, u_{i,a} \ldots |\vartheta_1, \vartheta_2),\notag
\\
&
{\Phi}_{M,N}^{(l)}(\cdots, u_{i,a}+\tau,\ldots|\vartheta_1, \vartheta_2)=
e^{-2\pi\sqrt{-1}\left(\sum_{j=0}^l {C_{i,j}}\overline{u}_j-{\mu_i}
+\tau\right)}
{\Phi}_{M,N}^{(l)}(\ldots, u_{i,a} \ldots |\vartheta_1, \vartheta_2).\notag
\end{align}
Modulo ${\bf Z} \oplus {\bf Z}\tau$, the only poles of 
$\Phi_{M,N}^{(l)}({\bf u}_0,\ldots, {\bf u}_l | \vartheta_1, \vartheta_2)$ 
are $u_{i,a}=u_{i,b}+\beta$ with some $a,b,i$.
\end{lem}
\begin{proof}
The quasi-periodicity is satisfied by each summand separately.
Since $\Phi_{M,N}^{(l)}$ is symmetric in each group of variables $\{u_{i,a}\}_{a=1}^{(M+N)a_i}$,
it cannot have a simple pole at $u_{i,a}=u_{i,b}$.
\end{proof}

\begin{lem}\label{lem:recursion1}
Let $\mathrm{Res}$ stand for the residue on the intersection of the lines
\begin{align*}
&u_{i,2}=u_{i,1}+\beta\quad \text{for $i=0,1,l-1,l$}\,,\\ 
&u_{i,3}=u_{i,1}+\beta\,,\quad u_{i,4}=u_{i,2}+\beta\,,
\quad \text{for $2\le i\le l-2$}\,.
\end{align*}
Then the following recursion relation holds for $M,N \geq 1$.
\begin{align}
&
\mathop{\mathrm{Res}}
{\Phi}_{M,N}^{(l)}({\bf u}_0, {\bf u}_1, \ldots, {\bf u}_l|\vartheta_1, \vartheta_2)=
{c^{2(l-1)}}
\Phi_{M-1,N-1}^{(l)}({\bf u}_0', {\bf u}_1',\ldots, {\bf u}_l'|{\vartheta_1'}, \vartheta_2')
\notag\\
\times&\frac{
\prod\limits_{\substack{\scriptscriptstyle 
0\leq i<j \leq l \\ \scriptscriptstyle 
{C_{i,j}=-1}}}
\Bigl(
\prod\limits_{\substack{\scriptscriptstyle a=1 \\
\scriptscriptstyle a \neq a_{i}+1, 2 a_i}}^{(M+N)a_i}
\prod\limits_{b=1}^{a_j}
\theta\bigl(u_{i,a}-u_{j,b}-\beta/2\bigr)
\prod\limits_{a=1}^{a_i}
\prod\limits_{\substack{\scriptscriptstyle b=1
\\ \scriptscriptstyle b \neq a_{j}+1, 2a_j}}^{(M+N)a_j}
\theta\bigl(u_{i,a}-u_{j,b}+\beta/2\bigr)
\Bigr)
}{
\prod\limits_{i=0}^l 
{\prod\limits_{a=1}^{a_i}}
{\prod\limits_{\substack{\scriptscriptstyle b=1\\ \scriptscriptstyle
b\neq a,a_i+1,2a_i}}^{(M+N)a_i}}
\theta(u_{i,a}-u_{i,b}+\beta)\theta(u_{i,a}-u_{i,b})}.
\notag
\end{align}
Here we set
\begin{align*}
&{{\bf u}_i=(u_{i, 1}, u_{i,2},\ldots, u_{i, (M+N)a_i}),}\quad
{\bf u}_i'=(u_{i, 2a_i+1}, u_{i,2a_i+2},\ldots, u_{i, (M+N)a_i}),~~~0\leq i \leq l\,,
\\
&
\vartheta_j'(v_0,\ldots, v_i,\ldots, v_l)=
\vartheta_j\Bigl(
v_0+u_{0,1}, \ldots, v_i+\sum_{a=1}^{a_i} u_{i,a},\ldots, v_l+u_{l,1}
\Bigr){\in X_l(\nu_1',\ldots, \nu_l')}\,,\quad j=1,2,
\end{align*}
with ${\displaystyle \nu_i'=\nu_i-\sum_{a=1}^{a_i} u_{i,a},~ 0\leq i \leq l}$,
and $c=\frac{1}{\theta(\beta)}\underset{u=0}{\res}\frac{du}{\theta(u)}$.
\end{lem}

\begin{proof}
In the sum ${\rm LHS}_{M,N}^{(l)}$, the iterated residue is non-zero only for terms such that $1 \in I_i$, $2 \in J_i$
for $i=0,1,l-1,l$, and $1,2 \in I_i$, $3,4 \in J_i$ for $2\leq i \leq l-2$.
Write $I_i=\{1\} \cup \overline{I}_i$, $J_i=\{2\} \cup \overline{J}_i$ for $i=0,1,l-1,l$, and
$I_i=\{1, 2\} \cup \overline{I}_i$, $J_i=\{3,4 \} \cup \overline{J}_i$ for $2\leq i \leq l-2$.
Similarly, for ${\rm RHS}_{M,N}^{(l)}$, the iterated residue is non-zero only for terms such that 
$2 \in I_i$, $1 \in J_i$
for $i=0,1,l-1,l$, and $3,4 \in I_i$, $1,2 \in J_i$ for $2\leq i \leq l-2$.
Write $I_i=\{2\} \cup \overline{I}_i$, $J_i=\{1\} \cup \overline{J}_i$ for $i=0,1,l-1,l$, and
$I_i=\{3, 4\} \cup \overline{I}_i$, $J_i=\{1,2 \} \cup \overline{J}_i$ for $2\leq i \leq l-2$.

In either case the functions $\vartheta_j$ specialize to 
$\vartheta'_j$ due to Lemma \ref{lem:theta1}.

The rest is a straightforward calculation. 
\end{proof}

\begin{lem}\label{lem:recursion2}
Let $\vartheta_1, \vartheta_2 \in X_l$.
We have the following recursion relations under specialization.

(1) For $l=4$ and $0\leq i \neq j \leq 4$, {$i,j\neq 2$}, 
\begin{align}
&
\Phi_{M,N}^{(4)}({\bf u}_0, {\bf u}_1, {\bf u}_2, {\bf u}_3, {\bf u}_4| \vartheta_1, \vartheta_2)
\Bigl|_{u_{2,a}=u_{i,a}\pm \frac{\beta}{2},\,
u_{2,M+N+a}=u_{j,a}\pm \frac{\beta}{2},\ 1\leq a \leq M+N}
\notag\\
=&
\Phi_{M,N}^{(3)}({\bf u}_0, {\bf u}_1, {\bf u}_3, {\bf u}_4|\vartheta_1'', \vartheta_2''')
\Bigl|_{{u_{i,a} \to u_{i,a}\pm \frac{\beta}{2},\,
u_{j,a} \to u_{j,a}\pm \frac{\beta}{2},\ 1\leq a \leq M+N}}
\notag
\end{align}
where
\begin{align}
\vartheta_1''(v_0,v_1,v_3,v_4)=\vartheta_1(v_0, v_1, v_i+v_j\pm M\beta, v_3, v_4),\notag\\
\vartheta_2'''(v_0,v_1,v_3,v_4)=\vartheta_2(v_0, v_1, v_i+v_j\pm N\beta, v_3, v_4).\notag
\end{align}
Here we used the Cartan matrix $\left(C_{i,j}\right)_{i,j=0}^l$ for $l=3$ given as
$\left(
\begin{array}{cccc}
2&0&-1&-1\\
0&2&-1&-1\\
-1&-1&2&0\\
-1&-1&0&2
\end{array}
\right)$ 
for $(i,j)=(0,1), (3,4)$,
$\left(
\begin{array}{cccc}
2&-1&-1&0\\
-1&2&0&-1\\
-1&0&2&-1\\
0&-1&-1&2
\end{array}
\right)$ 
for $(i,j)=(0,4), (1,3)$,
and
$\left(
\begin{array}{cccc}
2&-1&0&-1\\
-1&2&-1&0\\
0&-1&2&-1\\
-1&0&-1&2
\end{array}
\right)$ 
for $(i,j)=(0,3), (1,4)$.

\vspace{2mm}

(2) For $l \geq 5$
\begin{align}
&
\Phi_{M,N}^{(l)}({\bf u}_0, {\bf u}_1, {\bf u}_2, \ldots, {\bf u}_l| \vartheta_1, \vartheta_2)
\Bigl|_{u_{2,a}=u_{0,a}\pm \frac{\beta}{2},\,
u_{2,M+N+a}=u_{1,a}\pm \frac{\beta}{2},\ 
1\leq a \leq M+N}
\notag\\
=&
\Phi_{M,N}^{(l-1)}({\bf u}_0, {\bf u}_1, {\bf u}_3, \ldots, {\bf u}_l|\vartheta_1'', \vartheta_2''')
\Bigl|_{
{u_{0,a} \to u_{0,a}\pm \frac{\beta}{2},\,
u_{1,a} \to u_{1,a}\pm \frac{\beta}{2},\ 1\leq a \leq M+N}}
\notag
\end{align}
where
\begin{align}
\vartheta_1''(v_0,v_1,v_3,\ldots,v_l)=\vartheta_1(v_0, v_1, v_0+v_1\pm M\beta, v_3, \ldots, v_l),\notag\\
\vartheta_2'''(v_0,v_1,v_3,\ldots, v_l)=\vartheta_2(v_0, v_1, v_0+v_1\pm N\beta, v_3, \ldots, v_l).\notag
\end{align}

\begin{align}
&
\Phi_{M,N}^{(l)}({\bf u}_0, \ldots, {\bf u}_{l-3}, {\bf u}_{l-2}, {\bf u}_{l-1}, {\bf u}_{l}| \vartheta_1, \vartheta_2)
\Bigl|_{u_{l-2,a}=u_{l-1,a}\pm \frac{\beta}{2},\,
u_{l-2,M+N+a}=u_{l,a}\pm \frac{\beta}{2},\,
1\leq a \leq M+N}
\notag\\
=&
\Phi_{M,N}^{(l-1)}({\bf u}_0, \ldots, {\bf u}_{l-3}, {\bf u}_{l-1}, {\bf u}_l|\vartheta_1'', \vartheta_2''')
\Bigl|_{
{u_{0,a} \to u_{0,a}\pm \frac{\beta}{2},\,
u_{1,a} \to u_{1,a}\pm \frac{\beta}{2},\,\ 1\leq a \leq M+N}}
\notag
\end{align}
where
\begin{align}
\vartheta_1''(v_0,\ldots, v_{l-3},v_{l-2},v_{l-1},v_l)=\vartheta_1(v_0, \ldots, v_{l-3}, v_{l-1}+v_{l}\pm M\beta, v_{l-1}, v_l),\notag\\
\vartheta_2'''(v_0,\ldots, v_{l-3},v_{l-2}, v_{l-1}, v_l)=\vartheta_2(v_0, \ldots, v_{l-3}, v_{l-1}+v_l\pm N\beta, v_{l-1}, v_l).\notag
\end{align}

(3) For $l \geq 5$ and $2\leq i \leq l-3$,
\begin{align}
&
\Phi_{M,N}^{(l)}({\bf u}_0, \ldots, {\bf u}_{i}, {\bf u}_{i+1}, {\bf u}_{i+2}, \ldots, {\bf u}_{l}| \vartheta_1, \vartheta_2)
\Bigl|_{u_{i+1,a}=u_{i,a}\pm \frac{\beta}{2}\,,\ 1\leq a \leq 2(M+N)}
\notag\\
=&
\Phi_{M,N}^{(l-1)}({\bf u}_0, \ldots, {\bf u}_{i}, {\bf u}_{i+2}, {\bf u}_{i+3}, \ldots, {\bf u}_l|\vartheta_1'', \vartheta_2''')
\Bigl|_{
{u_{i,a} \to u_{i,a}\pm \frac{\beta}{2}\,, \ 1\leq a \leq 2(M+N)}}
\notag
\end{align}
where
\begin{align}
\vartheta_1''(v_0,\ldots, v_{i},v_{i+2},v_{i+3},\ldots, v_l)=\vartheta_1(v_0, \ldots, v_{i}, v_{i}\pm M\beta, v_{i+2}, \ldots, v_l),\notag\\
\vartheta_2'''(v_0,\ldots, v_{i}, v_{i+2}, v_{i+3},\ldots, v_l)=\vartheta_2(v_0, \ldots, v_{i}, v_{i}\pm N\beta, v_{i+2},\ldots, v_l).\notag
\end{align}
The signs in the equations are taken in the respective order.
\end{lem}
\begin{proof}
Since the calculations are similar, 
we show (3) as an example. 
Setting 
$u_{i+1,a}=u_{i, a}\pm\frac{\beta}{2}$, $1\leq a \leq 2(M+N)$,
terms in ${\rm LHS}_{M,N}^{(l)}$ 
survive only when 
$I_i \cap J_{i+1}=\emptyset$ for the $+$ sign and  
$J_i \cap I_{i+1}=\emptyset$ for the $-$ sign. 
In either case this implies $I_i=I_{i+1}$, $J_i=J_{i+1}$.
Then the factors with $j=i+1$ in the numerator 
cancels the one with $i\to i+1$ in the denominator,
so that ${\rm LHS}_{M,N}^{(l)}$ reduces to 
${\rm LHS}_{M,N}^{(l-1)}$ wherein the variables $\mathbf{u}_i$
and the functions $\vartheta_i$ are changed as shown above. 
For ${\rm RHS}_{M,N}^{(l)}$ the situation is entirely similar.

\end{proof}
{Just as with $\vartheta_i' \in X_l(\nu_1',\ldots, \nu_l')$, 
$\vartheta_i'' \in X_l(\nu_1'',\ldots, \nu_l'')$ and $\vartheta_i''' \in X_l(\nu_1''',\ldots,\nu_l''')$ hold 
for appropriate $\nu_1'',\ldots, \nu_l''$ and $\nu_1''',\ldots, \nu_l'''$.}
\bigskip

\noindent {\it Proof of Proposition \ref{prop:theta1}.}\quad 
Let $l, M, N$ be integers such that $l\geq 3, M>N\geq 0$.
We show $\Phi_{M,N}^{(l)}=0$ by double mathematical induction with respect to $(l,N)$.
Bases of this induction are $\Phi_{M,N}^{(3)}=0$, $N \geq 0$
and $\Phi_{M,0}^{(l)}=0$, $l \geq 3$.

First we make a remark. 
Consider the residue at $u_{i,2}=u_{i,1}+\beta$. 
In
{$\mathrm{LHS}^{(l)}_{M,N}$}, only 
those terms with $1\in I_i$ and $2\in J_i$ survive. 
Similarly in 
{$\mathrm{RHS}^{(l)}_{M,N}$}, only 
those terms with $2\in I_i$ and $1\in J_i$ survive. 
In both cases, the numerator is divisible by the factor
\begin{align*}
\prod_{b=1}^{(M+N)a_j}\theta(u_{i,1}-u_{j,b}+\beta/2)
\end{align*}
for all $j$ such that 
$C_{i,j}=-1$.

For $0 \leq i \leq 2l-2$, we set
\begin{align}
&Z_i=\underset{u_{i+2,4}=u_{i+2,3}+\beta}{\res}Z_{i+1}\,,\quad
0 \leq i \leq l-4,\notag\\
&Z_{l-3}=\underset{u_{l,2}=u_{l,1}+\beta}{\res} Z_{l-2}\,,
\quad
Z_{l-2}=\underset{u_{l-1,2}=u_{l-1,1}+\beta}{\res} Z_{l-1}\,,
\quad
Z_{l-1}=\underset{u_{1,2}=u_{1,1}+\beta}{\res}Z_l\,,
\quad
Z_{l}=\underset{u_{0,2}=u_{0,1}+\beta}{\res} Z_{l+1},\notag
\\
&Z_i=\underset{u_{i-l+1,3}=u_{i-l+1,1}+\beta}{\res} Z_{i+1}\,,
\quad
l+1\leq i \leq 2l-3\,,\quad
Z_{2l-2}=\Phi_{M,N}^{(l)}.\notag
\end{align}
We show $Z_i=0$ by induction with respect to $i$.

By Lemma \ref{lem:recursion1} and 
the induction hypothesis 
$\Phi_{M,N-1}^{(l)}=0$,
we obtain $Z_0=0$.

To show $Z_1=0$, we regard it as a function of $\{u_{2,a}\}_{3\le a\le 2(M+N)}$ and $\tau$.
From the observation made above, we can write
\begin{align*}
Z_1&=f_1(u_{2,3},\ldots,u_{2,2(M+N)})
\frac{\prod\limits_{j=0,1}\prod\limits_{b=3}^{2(M+N)}\theta(u_{j,1}-u_{2,b}+\beta/2)}
{\prod\limits_{b=3}^{2(M+N)}\theta(u_{2,1}-u_{2,b}+\beta)\theta(u_{2,1}-u_{2,b})} \\
&\times
\begin{cases}
\prod\limits_{j=3,4}\prod\limits_{b=3}^{2(M+N)}\theta(u_{j,1}-u_{2,b}+\beta/2)
&\text{if $l=4$},
\\
\prod\limits_{b=3}^{2(M+N)}\theta(u_{3,1}-u_{2,b}+\beta/2)
\prod\limits_{b=3}^{2(M+N)}\theta(u_{3,3}-u_{2,b}+\beta/2)
&\text{if $l\ge5$}\,.
\end{cases}
\end{align*}
The function $f_1$ is symmetric, $1$-periodic 
in $\{u_{2,a}\}_{3\le a\le 2(M+N)}$, and has poles at most
at $u_{2,a}=u_{2,b}+\beta$. 
From $0=Z_0=\mathop{\mathrm{res}}_{u_{2,4}=u_{2,3}+\beta}Z_1$
we see that $f_1$ is entire. On the other hand, from the quasi-periodicity 
of $\Phi^{(l)}_{M,N}$ we have for each $a$
\begin{align*}
f_1(\ldots,u_{2,a}+\tau,\ldots) 
=f_1(\ldots,u_{2,a},\ldots)
\exp\Bigl(
{\scriptstyle \sqrt{-1}
(2\sum\limits_{\substack{\scriptscriptstyle b=3 \\ \scriptscriptstyle
b \neq a}}^{2(M+N)}u_{2,b}+x_1)}
\Bigr)
\end{align*}
where $x_1$ is independent of $\{u_{2,b}\}_{3\le b\le 2(M+N)}$ and $\tau$. 
Such a function must be constant in each $u_{2,a}$.
This is possible only if $f_1=0$, namely $Z_1=0$. 

Next we show $Z_i=0$ for $
{2}\le i\le l-3$ by induction on $i$. 
Suppose we have shown $Z_{i-1}=0$.  
We view $Z_{i}$ as a function of $\{u_{i+1,a}\}_{3\le a\le 2(M+N)}$ and $\tau$. 
It can be written as 
\begin{align*}
Z_i&=f_i(u_{i+1,3},\ldots,u_{i+1,2(M+N)})
\frac{\prod\limits_{b=3}^{2(M+N)}\theta(u_{i,1}-u_{i+1,b}+\beta/2)}
{\prod\limits_{b=3}^{2(M+N)}\theta(u_{i+1,1}-u_{i+1,b}+\beta)
\theta(u_{i+1,1}-u_{i+1,b})} \\
&\times
\begin{cases}
\prod\limits_{j=l-1,l}\prod\limits_{b=3}^{2(M+N)}\theta(u_{j,1}-u_{l-2,b}+\beta/2)
&\text{if $i=l-3$},
\\
\prod\limits_{b=3}^{2(M+N)}\theta(u_{i+2,1}-u_{i+1,b}+\beta/2)
\prod\limits_{b=3}^{2(M+N)}\theta(u_{i+2,3}-u_{i+1,b}+\beta/2)
&\text{if $i\le l-4$}\,.
\end{cases}
\end{align*}
The function $f_i$ is symmetric, $1$-periodic 
in $\{u_{i+1,a}\}_{3\le a\le 2(M+N)}$, and has poles
at most at $u_{i+1,a}=u_{i+1,b}+\beta$. 
From $0=Z_{i-1}=\mathop{\mathrm{res}}_{u_{i+1,4}=u_{i+1,3}+\beta}Z_i$
we see that $f_i$ is entire. 
Further we have 
\begin{align*}
{f_i(\ldots,u_{i+1,a}+\tau,\ldots) 
=f_i(\ldots,u_{i+1,a},\ldots)}
\exp\Bigl(
{\scriptstyle \sqrt{-1}} 
(u_{i+1,a}+\tau/2+2
{\sum\limits_{\substack{\scriptscriptstyle b=3
\\ \scriptscriptstyle b \neq a}}^{2(M+N)}}u_{i+1,b}+x_i)
\Bigr)
{,\quad
2 \leq i \leq l-3,}
\end{align*}
where $x_i$ is independent of $\{u_{i+1,a}\}_{3\le a\le 2(M+N)}$ and $\tau$.
If $f_i$ is not identically zero, then for each $a$
it has exactly one zero 
at $u_{i+1,a}+
2
{\sum\limits_{\substack{\scriptscriptstyle b=3 \\ b \neq a}}^{2(M+N)}}
u_{i+1,b}+x_i\equiv0$ $\bmod \,\Z\oplus\Z\tau$.
Clearly this is impossible. Therefore $f_i=0$, namely $Z_i=0$ 
for $2\leq i \leq l-3$. 

To see $Z_{l-2}=0$ we proceed as follows. 
As a function in $\{u_{l,a}\}_{1\le a\le M+N}$, $Z_{l-2}$
is symmetric, $1$-periodic, and has poles at most at $u_{l,a}=u_{l,b}+\beta$.
From $0=Z_{l-3}=\mathop{\mathrm{res}}_{u_{l,2}=u_{l,1}+\beta}
{Z_{l-2}}$
we see that $Z_{l-2}$ is entire. For each $a$ we have the quasi-periodicity
\begin{align*}
Z_{l-2}\bigl|_{u_{l,a}\to u_{l,a}+\tau}=Z_{l-2}
e^{-2\pi i (\sum_{b=1}^{M+N}u_{l,b}+x_{l-2})}
\end{align*}
where $x_{l-2}$ is independent of $\{u_{l,b}\}_{1\le b\le M+N}$.
By Lemma \ref{lem:theta2}, $Z_{l-2}$ is a function of $\bar{u}_l$. 
On the other hand, it is divisible by 
$\prod_{a=1}^{M+N}\theta(u_{l-2,1}-u_{l,a}+\beta/2)$. This is possible
only when $Z_{l-2}=0$. 

Proof of $Z_i=0$ with $l-1\leq i \leq l$ is entirely similar as those of $Z_{l-2}$.
Proof of $Z_i=0$ with $l+1\le i\le 2l-3$ is entirely similar as those of $Z_i$ with $1\leq i \leq l-3$.

It remains to show $Z_{2l-2}=0$. 
We have $Z_{2l-3}=\underset{u_{l-2,2}=u_{l-2,1}+\beta}{\res}Z_{2l-2}=0$.
Then, $Z_{2l-2}$ has no pole in $u_{l-2,a}$, $1\leq a \leq 2(M+N)$.
Due to Lemma \ref{lem:theta2}, 
$Z_{2l-2}$ is a function of $\overline{u}_{l-2}$. 
On the other hand, by Lemma \ref{lem:recursion2} 
and the hypothesis $\Phi_{M,N}^{(l-1)}=0$,
$Z_{2l-2}$ has the following zeros with respect to $\overline{u}_{l-2}$:
\begin{align}
\overline{u}_{l-2}=\overline{u}_{l-1}+\overline{u}_l\pm\beta(M+N),
\quad
\overline{u}_{l-2}=
\left\{
\begin{array}{cc}
{\overline{u}_{0}+\overline{u}_{1}\pm \beta(M+N)}& 
{\mathrm{if} ~l=4,}\notag
\\
\overline{u}_{l-3}\pm \beta(M+N)& \mathrm{if} ~l\geq 5.\notag
\end{array}\right.
\end{align}
From quasi-periodicity, if $Z_{2l-2}\neq 0$, 
then it can have only two zeros in $\overline{u}_{l-2}$ modulo ${\bf Z} \oplus {\bf Z} \tau$.
It is a contradiction to the above zeros.
Hence we conclude that $Z_{2l-2}=\Phi_{M,N}^{(l)}=0$.
\qed

\section{Proof of Theorem \ref{thm:duality}}\label{section:duality}

In this section we prove the commutativity 
${\bf I}_M \cdot {\bf G}_N(\vartheta)={\bf G}_N(\vartheta) \cdot {\bf I}_M$ 
{for $M,N \geq 1$ and $\vartheta \in X_l(\nu_1,\ldots, \nu_l)$} by a direct calculation.
Since the integrand is real analytic in $\beta$, 
we shall assume $\beta>{\rm Max}(3, 2M-1)$ without loss of generality.

\subsection{Outline of proof}

We start from products ${\bf I}_M \cdot {\bf G}_N(\vartheta)$ and ${\bf G}_N(\vartheta) \cdot {\bf I}_M$, which are given respectively by
\begin{align}
{\bf I}_M \cdot {\bf G}_N(\vartheta)&=\prod_{i=1}^M \int_{|w_i|=R}
\underline{dw_i}
\prod_{j=0}^l \prod_{a=1}^{N a_j}\int_{|z_{j,a}|=\delta_j}
\underline{dz_{j,a}}
\overset{\curvearrowright} 
{\prod_{1\leq i \leq M}}{\bf T}_1(w_i)
\overset{\curvearrowright} 
{\prod_{0\leq j \leq l}}
\overset{\curvearrowright}
{\prod_{~1\leq a \leq N a_j}}
S_j(z_{j,a})\notag\\
&\times
\omega_M(w_1,\ldots,w_M)
h_{N}({\bf u}_0,\ldots, {\bf u}_l;\vartheta),\notag
\end{align}
\begin{align}
{\bf G}_N(\vartheta) \cdot {\bf I}_M&=
\prod_{i=1}^M \int_{|w_i|=r}
\underline{dw_i}
\prod_{j=0}^l \prod_{a=1}^{N a_j}\int_{|z_{j,a}|=\delta_j}
\underline{dz_{j,a}}
\overset{\curvearrowright} 
{\prod_{0\leq j \leq l}}
\overset{\curvearrowright}    
{\prod_{~1\leq a \leq N a_j}}
S_j(z_{j,a})\
\overset{\curvearrowleft} 
{\prod_{1\leq i \leq M}}{\bf T}_1(w_i)
\notag\\
&\times
\omega_M(w_M,\ldots,w_1)
h_{N}({\bf u}_0,\ldots, {\bf u}_l;\vartheta).\notag
\end{align}
The product of operators $\bLambda_k(w_i)S_j(z_{j,a})$ 
and $S_j(z_{j,a})\bLambda_k(w_i)$ converge absolutely in the regions
$|z_{j,a}/w_i| \ll 1$ and $|z_{j,a}/w_i|\gg 1$, respectively.
Accordingly we choose the radii of integration 
circles $R, r$ and 
{$\{\delta_j\}_{j=0}^l$} as
\begin{align}
R\gg 1 \gg r>0 \quad{\rm and}\quad
\delta_0, \delta_l >1,
\quad \delta_1=\delta_2=\cdots=\delta_{l-1}=1\quad {\rm such~that}\quad \delta_0,\delta_l \doteqdot1.\notag
\end{align}
As it turns out, in the course of the proof, 
new singularities arise at $|z_{0,a}|=|z_{1,b}|$ 
and $|z_{l-1,a}|=|z_{l,b}|$, see Sections \ref{section:I1-GN},
\ref{section:I2-GN} and \ref{section:IM-GN} below. 
To avoid these singularities we choose 
$|z_{0,a}|=\delta_0>1$, $|z_{1,b}|=1$ 
and $|z_{l,a}|=\delta_l>1$, $|z_{l-1,b}|=1$.

Consider the product ${\bf I}_M \cdot {\bf G}_N(\vartheta)$.
We move the contour for $\{w_i\}_{i=1}^M$ from the circle $|w_i|=R \gg1$ to the unit circle $|w_i|=1$.
Along the way we take the residues of poles 
{that are contained in the annulus 
$\overline{D}_R=\{w \in \mathbf{C}| 1<|w|<R\}$.}
We shall see that the poles are simple and are of the form 
\begin{align*}
&\overline{d}^{(t)}_{j,s}=\bar{d}^{(t)}_j[z_{j,s}]\,,
{1\leq j \leq l-2,}\quad
\overline{e}^{(t)}_{j,s}=\bar{e}^{(t)}_j[z_{j,s}]\,,
{0 \leq j \leq l-2\quad \mathrm{with}~ 1\leq t \leq M, 1\leq s \leq Na_j},
\\
&
\overline{D}^{(t)}_{j,s}=\overline{D}^{(t)}_j[z_{j,s}]\,,
{1\leq j \leq l-3,}\quad
\overline{E}^{(t)}_{j,s}=\overline{E}^{(t)}_j[z_{j,s}]\,,
1\leq j \leq l-3\quad \mathrm{with}~1\leq t \leq M-1, 1\leq s \leq Na_j,
\\
&
{\overline{d}^{(1)}_{l-1,s}=\bar{d}^{(1)}_{l-1}[z_{l-1,s}]\,,
\overline{d}^{(1)}_{l,s}=\bar{d}^{(1)}_l[z_{l,s}]\,
\quad
\mathrm{with}~1\leq s \leq N,
}
\end{align*}
where
\begin{align}
\overline{d}_{j}^{(t)}[z]=\left\{\begin{array}{cc}
q^{t-\frac{3}{2}-\frac{\beta}{2}(l-j-1)}z,& 1\leq j \leq l-1,\\
q^{t-\frac{3}{2}}z,& j=l,
\end{array}\right.
\quad
\overline{e}_{j}^{(t)}[z]=\left\{\begin{array}{cc}
q^{t-\frac{1}{2}-\frac{\beta}{2}(l-2)}z,& j=0,\\
q^{t-\frac{1}{2}-\frac{\beta}{2}(l-j-1)}z,& 1\leq j \leq l-2,
 \end{array}\right.
\notag
\\
\overline{D}_j^{(t)}[z]=q^{(\beta-1)t-\frac{1}{2}-\frac{\beta}{2}(l-j-1)}z,~1\leq j \leq l-3,
\quad
\overline{E}_j^{(t)}[z]=q^{(\beta-1)t+\frac{1}{2}-\frac{\beta}{2}(l-j-1)}z,~1\leq j \leq l-3.\notag
\end{align}
The position of poles of ${\bf I}_M \cdot {\bf G}_N(\vartheta)$ is depicted 
in Appendix \ref{appendix:IM-GN},
\ref{appendix:I1-GN}, and \ref{appendix:I2-GN}.
{We take the residues of poles at $w_i=\overline{d}_{j,s}^{(t)}, \overline{e}_{j,s}^{(t)}, \overline{D}_{j,s}^{(t)}$, and $\overline{E}_{j,s}^{(t)}$.}
At first glance, from 
TABLE \ref{table:pole1} in Appendix \ref{appendix:normal-ordering},
it may seem that 
only poles $\bar{d}^{(1)}_{j,s}$ and  $\bar{e}^{(1)}_{j,s}$ are present.
However, after taking residues in some $w_i$, poles of the above form 
appear with respect to the remaining variables.
We give the details in Section \ref{section:I2-GN} and Section 
\ref{section:IM-GN}.

We write 
\begin{align}
{\bf I}_M \cdot {\bf G}_N(\vartheta)
=\sum\limits_{\substack{\scriptscriptstyle a,b,c \geq 0
\\ \scriptscriptstyle 0\leq a+b+c \leq M}}
{\rm LHS}_{M,N}(a,b,c).\notag
\end{align}
Here ${\rm LHS}_{M,N}(a,b,c)$ is a 
sum of the terms that are obtained by taking residues 
$a$-times at $w_i=\overline{d}_{1,s}^{(t)}, \overline{e}_{0,s}^{(t)}, \overline{e}_{1,s}^{(t)}$,
{$1\leq t \leq M, 1 \leq s \leq N$},
$\overline{D}_{1,s}^{(t)}$, $\overline{E}_{1,s}^{(t)}$, 
{$1\leq t \leq M-1, 1\leq s \leq N$},
$b$-times at $w_i=\overline{d}_{j,s}^{(t)}, 
\overline{e}_{j,s}^{(t)}$, $2\leq j \leq l-2,
{1\leq t \leq M, 1\leq s \leq 2N}$, 
$\overline{D}_{j,s}^{(t)}$, $\overline{E}_{j,s}^{(t)}$, 
$2\leq j\leq l-3, 
{1\leq t \leq M-1, 1\leq s \leq 2N}$,
and $c$-times at 
{$w_i=\overline{d}_{l-1,s}^{(1)}$, 
$\overline{d}_{l,s}^{(1)}$, $1\leq s \leq N$}.

For the product ${\bf G}_N(\vartheta) \cdot {\bf I}_M$
we proceed similarly, moving 
the contour for $\{w_i\}_{i=1}^M$ from the circle $|w_i|=r\ll 1$ to the unit circle $|w_i|=1$.
Along the way we take the residues of poles that are contained in the 
annulus 
${D}_r=\{w \in \mathbf{C}| r<|w|<1\}$.
All poles are simple.
They are given by
\begin{align*}
&{d}^{(t)}_{j,s}={d}^{(t)}_j[z_{j,s}]\,,
{1\leq j \leq l-2},\quad
{e}^{(t)}_{j,s}={e}^{(t)}_j[z_{j,s}]\,,
{0\leq j \leq l-2\quad \mathrm{with}~1\leq t \leq M, 1\leq s \leq Na_j},
\notag
\\
&{D}^{(t)}_{j,s}={D}^{(t)}_j[z_{j,s}]\,,
{1\leq j \leq l-3},\quad
{E}^{(t)}_{j,s}={E}^{(t)}_j[z_{j,s}]\,, 
{1\leq j \leq l-3 \quad \mathrm{with}~1\leq t \leq M-1, 1\leq s \leq Na_j},\notag
\\
&
{d}^{(1)}_{l-1,s}={d}^{(1)}_{l-1}[z_{l-1,s}]\,,{d}^{(1)}_{l,s}={d}^{(1)}_l[z_{l,s}]\,
\quad
\mathrm{with}~1\leq s \leq N,
\end{align*}
where
\begin{align}
{d}_{j}^{(t)}[z]=\left\{\begin{array}{cc}
q^{-t+\frac{3}{2}+\frac{\beta}{2}(l-j-1)}z,& 1\leq j \leq l-1,\\
q^{-t+\frac{3}{2}}z,& j=l,
\end{array}\right.
\quad
{e}_{j}^{(t)}[z]=\left\{\begin{array}{cc}
q^{-t+\frac{1}{2}+\frac{\beta}{2}(l-2)}z,& j=0,\\
q^{-t+\frac{1}{2}+\frac{\beta}{2}(l-j-1)}z,& 1\leq j \leq l-2,
 \end{array}\right.
\notag
\\
{D}_j^{(t)}[z]=q^{-(\beta-1)t+\frac{1}{2}+\frac{\beta}{2}(l-j-1)}z,~1\leq j \leq l-3,
\quad
{E}_j^{(t)}[z]=q^{-(\beta-1)t-\frac{1}{2}+\frac{\beta}{2}(l-j-1)}z,~1\leq j \leq l-3.\notag
\end{align}
The poles of ${\bf G}_N(\vartheta)\cdot {\bf I}_M$
and those of ${\bf I}_M \cdot {\bf G}_N(\vartheta)$
depicted in Appendix \ref{appendix:IM-GN}, 
\ref{appendix:I1-GN}, and \ref{appendix:I2-GN}
are in positions that are inverted with respect to the circle $|w_i|=1$.
We take the residues of poles at $w_i={d}_{j,s}^{(t)}, {e}_{j,s}^{(t)}, {D}_{j,s}^{(t)}$, and ${E}_{j,s}^{(t)}$.
We denote the sum of these residues as follows.
\begin{align}
{\bf G}_N(\vartheta) \cdot {\bf I}_M
=\sum\limits_{\substack{\scriptscriptstyle a,b,c \geq 0
\\ \scriptscriptstyle 0\leq a+b+c \leq M}}
{\rm RHS}_{M,N}(a,b,c).\notag
\end{align}
Here ${\rm RHS}_{M,N}(a,b,c)$ is sum of the terms that are obtained by taking the residues 
$a$-times at $w_i={d}_{1,s}^{(t)}, {e}_{0,s}^{(t)}, {e}_{1,s}^{(t)}$,
{$1\leq t \leq M, 1\leq s \leq N$}, 
${D}_{1,s}^{(t)}$, ${E}_{1,s}^{(t)}$, 
$1\leq t \leq M-1, 1\leq s \leq N$,
$b$-times at $w_i={d}_{j,s}^{(t)}, {e}_{j,s}^{(t)}$, 
$
1\leq t \leq M-1, 1\leq s \leq 2N, 2\leq j \leq l-2$, 
${D}_{j,s}^{(t)}$, ${E}_{j,s}^{(t)}$, 
$
{1\leq t \leq M-1, 1\leq s \leq 2N}, 2\leq j\leq l-3$,
and $c$-times at $w_i={d}_{l-1,s}^{(1)}$, ${d}_{l,s}^{(1)}$,
$1\leq s \leq N$.
\medskip

Later we give explicit formulae for
${\rm LHS}_{M,N}(a,b,c)$ and 
${\rm RHS}_{M,N}(a,b,c)$ in Sections \ref{section:I1-GN}, \ref{section:I2-GN} and \ref{section:IM-GN}.
\medskip

By definition we have
\begin{align}
{\rm LHS}_{M,N}(0,0,0)
&=\prod_{i=1}^M \int_{|w_i|=1}
\underline{dw_i}
\prod_{j=0}^l \prod_{a=1}^{N a_j}\int_{|z_{j,a}|=\delta_j}
\underline{dz_{j,a}}
\overset{\curvearrowright} 
{\prod_{1\leq i \leq M}}{\bf T}_1(w_i)
\overset{\curvearrowright} 
{\prod_{0\leq j \leq l}}
\overset{\curvearrowright} 
{\prod_{~1\leq a \leq N a_j}}
S_j(z_{j,a})\notag\\
&\times
\omega_M(w_1,\ldots,w_M)
h_{N}({\bf u}_0,\ldots, {\bf u}_l;\vartheta)\,,
\\
{\rm RHS}_{M,N}(0,0,0)
&=\prod_{i=1}^M \int_{|w_i|=1}
\underline{dw_i}
\prod_{j=0}^l \prod_{a=1}^{N a_j}\int_{|z_{j,a}|=\delta_j}
\underline{dz_{j,a}}
\overset{\curvearrowright} 
{\prod_{0\leq j \leq l}}
\overset{\curvearrowright} 
{\prod_{~1\leq a \leq N a_j}}
S_j(z_{j,a})
\overset{\curvearrowleft} 
{\prod_{1\leq i \leq M}}{\bf T}_1(w_i)
\notag\\
&\times
\omega_M(w_M,\ldots, w_1)
h_{N}({\bf u}_0,\ldots, {\bf u}_l;\vartheta).
\notag
\end{align}
In view of the exchange relations which hold as meromorphic functions
\begin{align}
&
{\bf T}_1(w)S_j(z)=S_j(z){\bf T}_1(w),\quad 0 \leq j \leq l\,,
\notag
\end{align}
we have
\begin{align}
{\rm LHS}_{M,N}(0,0,0)={\rm RHS}_{M,N}(0,0,0)\,,\quad M,N\geq 1.
\notag
\end{align}
Therefore, Theorem \ref{thm:duality} will follow if we show
\begin{prop}
\label{prop:duality1}
The followings hold for $M,N \geq 1$.
\begin{align}
&
{\rm LHS}_{M,N}(a,0,c)={\rm LHS}_{M,N}(a,0,c),
\quad a,c \geq 0, 0\leq a+c\leq M,
\notag
\\
&
{\rm LHS}_{M,N}(a,b,c)=0,\quad
{\rm RHS}_{M,N}(a,b,c)=0,\quad
a,c \geq 0, b \geq 1, 1\leq a+b+c\leq M.
\notag
\end{align}
\end{prop}

In the next sections we prove Proposition \ref{prop:duality1}.

\subsection{Proof of $[{\bf I}_1, {\bf G}_N(\vartheta)]=0$}
\label{section:I1-GN}

In this section, we show 
Proposition \ref{prop:duality1} 
in the case $M=1$. 

We start from ${\bf I}_1 \cdot {\bf G}_N(\vartheta)$.
Renaming the integration variable $w_1$ to $w$, 
we move the contour from $|w|=R\gg 1$ to $|w|=1$.
From TABLE \ref{table:pole1} in Appendix \ref{appendix:normal-ordering},
we see that, in the region between $|w|=R$ and $|w|=1$,  
each term $\bLambda_{k}(w)S_j(z_{j,s})$
can have a pole only at $w=\bar{\eta}_{k;j}[z_{j,s}]$, 
{$0\leq j \leq l$}, $1\le s \leq N a_j$ where
\begin{align}
\bar{\eta}_{k;j}[z]=
\begin{cases}
\bar{d}^{(1)}_{j}[z]&(k,j)=
\text{$(\bar{j};j)$ with $1\le j\le l$},\
(l;l-1),\ (\overline{l-1};l),\\
\bar{e}^{(1)}_{j}[z]&(k; j)=(\bar{1};0),\ (\bar{2};0),\ 
\text{$(\overline{j+1};j)$ with $1\le j\le l-2$}.
\end{cases}
\label{def:leta1}
\end{align}
(For $M=1$, poles of the form 
$\bar{D}^{(t)}_{j,s}$, $\bar{E}^{(t)}_{j,s}$ do not appear.)
For $(k; j)$ which appear in the above, we set
\begin{align*}
&\nu_{k;j}\mathcal{A}_{k;j}[z] 
=\mathop{\mathrm{res}}_{w=\overline{\eta}_{k;j}[z]}
\bLambda_{k}(w)S_j(z)\frac{dw}{w}\,,\\
&{\mathcal A}_{k;j}[z]= :S_j(z) \bLambda_{k}(\overline{\eta}_{k;j}[z]):\,.
\end{align*}
The constants $\nu_{k;j}$ will be given in Lemma \ref{lem:nukj} below.

For  each $0\leq i \leq l$, the integrand is 
a symmetric function of $z_{i,1}, z_{i,2}, \ldots, z_{i,Na_i}$.
Taking this symmetry into account and calculating 
residues,  we obtain 
\begin{align}
{\rm LHS}_{1,N}(1,0,0)&
=N
\int_{C_{[\delta]}}\underline{dz}\,
\Bigl(\sum_{k=\overline{1}, \overline{2}}
\nu_{k;0}{\bf S}[(k;0,1)]+\sum_{k=\overline{1}, \overline{2}}\nu_{k;1}
{\bf S}[(k; 1,1)]\Bigr)\bh\,,
\notag
\\
{\rm LHS}_{1,N}(0,1,0)&=
2N
\int_{C_{[\delta]}}\underline{dz}\,
\sum_{j=2}^{l-2}
\sum_{k=\overline{j}, \overline{j+1}}
\nu_{k;j} {\bf S}[(k; j,1)]\bh,\notag
\\
{\rm LHS}_{1,N}(0,0,1)&
=N
\int_{C_{[\delta]}}\underline{dz}\,
\Bigl(\sum_{k=\overline{l-1}, l}\nu_{k;l-1}
{\bf S}[(k;l-1,1)]
+\sum_{k=\overline{l-1}, \overline{l}}\nu_{k;l}
{\bf S}[(k;l,1)]\Bigr)\bh\,.
\notag
\end{align}
Here we set
\begin{align}
\label{Skjb}
{\bf S}[(k;j,b)]= 
\overset{\curvearrowright}{\prod_{
{0 \leq i < j}}}
\overset{\curvearrowright}{\prod_{~~1\leq a \leq Na_i}} S_i(z_{i,a}) 
\cdot
\overset{\curvearrowright}{\prod_{~~1\leq a<b}}S_j(z_{j,a})~
{\mathcal A}_{k;j}[z_{j,b}]
\overset{\curvearrowright}{\prod_{b<a \leq Na_j}}
S_j(z_{j,a})
\cdot \overset{\curvearrowright}{\prod_{
{j<i \leq l}}}
\overset{\curvearrowright}{\prod_{~~1\leq a \leq Na_i}} S_i(z_{i,a})\,,
\end{align}
and $C_{[\delta]}$ stands for the contour
\begin{align*}
&C_{[\delta]}:~|z_{j,a}|=\delta_j\,,\qquad  0\leq j \leq l, 1\leq a \leq Na_j\,.
\end{align*}

Similarly, for ${\bf G}_N(\vartheta) \cdot {\bf I}_1$, 
we move the contour from $|w|=r\ll 1$ to $|w|=1$.
In the region between $|w|=r$ and $|w|=1$,
$S_j(z_{j,s})\bLambda_{k}(w)$
can have a pole only at $w=\eta_{k;j}[z_{j,s}]$, 
{$0\leq j \leq l, 1\le s \leq Na_j$}, where
\begin{align}
\eta_{k;j}[z]=
\begin{cases}
d^{(1)}_{j}[z]&(k; j)=
\text{$(j; j)$ with $1\le j\le l$},\
(\bar{l}; l-1),\ (l-1; l),
\\
e^{(1)}_{j}[z]&(k; j)=(1; 0),\ (2; 0),\ 
\text{$(j+1; j)$ with $1\le j\le l-2$}.
\end{cases}
\label{def:eta1}
\end{align}
For $(k;j)$ which appear in the above, we set
\begin{align*}
&\nu_{k;j}\mathcal{A}_{k;j}[z] 
=\mathop{\mathrm{res}}_{w=\eta_{k;j}[z]}
S_j(z)\bLambda_{k}(w)
\frac{dw}{w}\,,\\
&{\mathcal A}_{k;j}[z]= :S_j(z) \bLambda_{k}(\eta_{k;j}[z]):\,.
\end{align*}
In the notation \eqref{Skjb}, we obtain
\begin{align*}
{\rm RHS}_{1,N}(1,0,0)&
=-N
\int_{C_{[\delta]}}\underline{dz}\,
\bigl(\sum_{k=1,2}\nu_{k;0}
{\bf S}[(k; 0,1)]+\sum_{k=1,2}\nu_{k;1}{\bf S}[(k;1,1)]\bigr)\bh\,,
\\
{\rm RHS}_{1,N}(0,1,0)&=
-2N 
\int_{C_{[\delta]}}\underline{dz}\,
\sum_{j=2}^{l-2}\sum_{k=j, j+1}\nu_{k; j}{\bf S}[(k; j,1)]
\bh\,,
\\
{\rm RHS}_{1,N}(0,0,1)&
=- N
\int_{C_{[\delta]}}\underline{dz}\,
\bigl(
\sum_{k=l-1,\overline{l}}\nu_{k;l-1} {\bf S}[(k;l-1,1)])
+\sum_{k=l-1,l}\nu_{k;l}{\bf S}[(k;l,1)]\bigr)\bh\,.
\end{align*}

The next two lemmas can be shown by a direct calculation.
\begin{lem}\label{lem:pole1}
The following hold.
\begin{align}
&q_1 {\mathcal A}_{\overline{1};0}[z]={\mathcal A}_{2;0}[q z],\quad
q_1 {\mathcal A}_{\overline{2};0}[z]={\mathcal A}_{1;0}[q z],
\notag
\\
&q_1 {\mathcal A}_{\overline{j+1};j}[z]={\mathcal A}_{\overline{j};j}[q z],\quad
q_1^{-1} {\mathcal A}_{j+1;j}[z]={\mathcal A}_{j;j}[q^{-1} z],\quad1\leq j \leq l-1,
\notag
\\
&q_1^{-1} {\mathcal A}_{\overline{l-1};l}[z]={\mathcal A}_{l;l}[q^{-1} z],\quad
q_1^{-1} {\mathcal A}_{\overline{l};l}[z]={\mathcal A}_{l-1;l}[q^{-1} z].
\notag
\end{align}
\end{lem}
\begin{lem}\label{lem:nukj}
The following hold.
\begin{align}
\nu_{k;i}=\left\{\begin{array}{cc}
q_1-1
& (k,i)=(\bar{1}, 0), (\bar{2}, 0), (l, l), (l-1, l),
\\
& \hskip 2cm 
(\overline{j+1}, j), (j, j), 1\leq j \leq l-1,
\\
q_1^{-1}-1
& (k,i)=(1, 0), (2, 0), (\overline{l-1}, l), (\overline{l}, l),
\\
&\hskip 2cm 
(\overline{j}, j), (j+1, j), 1\leq j \leq l-1,
\\
0& {\rm otherwise}.
\end{array}
\right.\notag
\end{align}
\end{lem}

Using Lemma \ref{lem:pole1} one finds 
\begin{align*}
&{\rm LHS}_{1,N}(1,0,0)-{\rm RHS}_{1,N}(1,0,0)
=N\sum_{k=\overline{1},\overline{2}}
\nu_{k;0 }\Delta^{(1,N)}_1[(k;0,1)]
+N\sum_{k=2,\overline{2}}
\nu_{k;1}\Delta^{(1,N)}_1[(k;1,1)]\,,
\\
&{\rm LHS}_{1,N}(0,0,1)-{\rm RHS}_{1,N}(0,0,1)
=N\sum_{k=l-1,\overline{l}}\nu_{k;
{l-1} }
\Delta^{(1,N)}_1[(k; l-1,1)]+
N\sum_{k=l-1, l}\nu_{k;l }
\Delta^{(1,N)}_1[(k;l,1)]\,,
\\
&{\rm LHS}_{1,N}(0,1,0)=
2N\sum_{j=2}^{l-2}
\nu_{\overline{j+1};j}
\Delta^{(1,N)}_1[(\overline{j+1};j,1)]\,,
\\
&{\rm RHS}_{1,N}(0,1,0)=
-2N\sum_{j=2}^{l-2}
\nu_{j+1;j}
\Delta^{(1,N)}_1[(j+1;j,1)]\,.
\end{align*}
Here we set
\begin{align*}
&
\Delta^{(1,N)}_1[(k;0,s)]
=\int_{C[(0,s)]}\underline{dz}
\int_{C_{\delta_0}-C_{q^{-1}\delta_0}}\underline{dz_{0,s}}{\bf S}[(k;0,s)]\bh
\,,\quad 
k=\overline{1},\overline{2}\,,
\\
&
\Delta^{(1,N)}_1[(\overline{j+1}; j,s)]
=\int_{C[(j,s)]}\underline{dz}
\int_{C_{1}-C_{q^{-1}}}\underline{dz_{j,s}}
{\bf S}[(\overline{j+1};j,s)]\bh\,,\quad
1\leq j \leq l-2\,,
\\
&
\Delta^{(1,N)}_1[(j+1; j,s)]
=\int_{C[(j,s)]}\underline{dz}
\int_{C_{1}-C_{q}}\underline{dz_{j,s}}
{\bf S}[(j+1;j,s)]\bh\,,\quad
1\leq j \leq l-2\,,
\\
&
\Delta^{(1,N)}_1[(k;l-1,s)]=
\int_{C[(l-1,s)]}\underline{dz}
\int_{C_1-C_{q^{-1}}}\underline{dz_{l-1,s}}{\bf S}[(k;l-1,s)]\bh\,,
\quad
k=l-1,\overline{l}\,,
\notag
\\
&
\Delta^{(1,N)}_1[(k;l,s)]=
\int_{C[(l,s)]}\underline{dz}\int_{C_{\delta_l}-C_{q^{-1}\delta_l}}
\underline{dz_{l,s}}{\bf S}[(k;l,s)]\bh\,,\quad
k=l-1,l\,,
\end{align*}
and $C[(i,s)]$ stands for the contour
\begin{align*}
&C[(i, s)]:~|z_{j,a}|=\delta_j,
\qquad (j,a)\neq (i, s)\,,\
0\leq j \leq l\,,\ 1\leq a \leq Na_j\,.
\end{align*}
In view of the symmetry, we have $\Delta^{(1,N)}_1[(k; j,s)]
=\Delta^{(1,N)}_1[(k; j,1)]$ for all $s\in[1,Na_j]$. 

\begin{lem}\label{lem:duality1}
The following hold for $N \geq 1$.
\begin{align}
&
\Delta^{(1,N)}_1[(\overline{2};0,1)]=0,\quad
\Delta^{(1,N)}_1[(2;1,1)]=0,\notag
\\
&\Delta^{(1,N)}_1[(\overline{j+1}; j,1)]=0,\quad
\Delta^{(1,N)}_1[(j+1;j,1)]=0,\quad 2 \leq j \leq l-2,
\notag\\
&\Delta^{(1,N)}_1[(\overline{l}; l-1,1)]=0,\quad
\Delta^{(1,N)}_1[(l-1;l ,1)]=0.
\notag
\end{align}
\end{lem}
\begin{proof}
As can be seen from the normal ordering of $S[(\overline{2};0,1)]
{\bf h}_N({\bf u}_0,\ldots, {\bf u}_l;\vartheta)$, 
there is no pole in the region $1\le|z_{0,1}|\le 
{|q^{-1}|}$.
Hence we obtain $\Delta^{{1,N}}_1[(\overline{2};0,1)]=0$.
The other cases are  shown in the same way.
 \end{proof}
 
 \begin{lem}\label{lem:duality2}
The following hold for $N \geq 1$.
\begin{align}
&
\Delta^{(1,N)}_1[(\overline{1};0,1)]+
\Delta^{(1,N)}_1[(\overline{2};1,1)]=0\,,\label{D1Na}
\\
&\Delta^{(1,N)}_1[(l-1; l-1,1)]
+\Delta^{(1,N)}_1[(l;l,1)]=0\,.
\label{D1Nb}
\end{align}
\end{lem}
\begin{proof}
First, we show \eqref{D1Na}. 
Here we find it convenient to use 
$\Delta^{(1,N)}_1[(\overline{1};0,N)]$ in place of 
$\Delta^{(1,N)}_1[(\overline{1};0,1)]$. 

As can be seen from the normal orderings of
 $S[(\overline{1};0,N)]{\bf h}_N({\bf u}_0,\ldots, {\bf u}_l;\vartheta)$ and
$S[(\overline{2};1,1)] {\bf h}_N({\bf u}_0,\ldots, {\bf u}_l;\vartheta)$,
there is a pole of $z_{0,N}$ in the region of integral 
$\int_{C_{\delta_0}-C_{q^{-1}\delta_0}}\underline{dz_{0,N}}$
and there is a pole of $z_{1,1}$ in the region of integral 
$\int_{C_{1}-C_{q^{-1}}}\underline{dz_{1,1}}$.
Taking these residues we obtain 
\begin{align}
\Delta^{(1,N)}_1[(\overline{1};0,N)]
&=N(1-q_1^{-1})
\int_{C[(0,N)]}\underline{dz}
\overset{\curvearrowright}{\displaystyle \prod_{1\leq a \leq N-1}}S_0(z_{0,a})
:\bLambda_{\overline{1}}(q^{-\frac{\beta}{2}(l-2)-\frac{1}{2}}z_{1,1})S_0(q^{-1}z_{1,1})S_1(z_{1,1}):
\notag\\
&\times
\overset{\curvearrowright}{\displaystyle \prod_{2 \leq a \leq N}}S_1(z_{1,a})
{\bf S}_{[2,l]}
\left.{\bf h}_N({\bf u}_0,\ldots, {\bf u}_l;\vartheta)\right|_{z_{0,N}=q^{-1}z_{1,1}}
,\notag
\\
\Delta^{(1,N)}_1[(\overline{2};1,1)]
&=N(1-q_1)\int_{C[(1,1)]}\underline{dz}
\overset{\curvearrowright}{\displaystyle \prod_{1\leq a \leq N-1}}S_0(z_{0,a})
:\bLambda_{\overline{2}}(q^{-\frac{\beta}{2}(l-2)+\frac{1}{2}}z_{0,N})
S_0(z_{0,N})S_1(z_{0,N}):\notag\\
&\times
\overset{\curvearrowright}{\displaystyle \prod_{2 \leq a \leq N}}S_1(z_{1,a}){\bf S}_{[2,l]}
\left.{\bf h}_N({\bf u}_0,\ldots, {\bf u}_l;\vartheta)\right|_{z_{1,1}=z_{0,N}}.\notag
\end{align}
Using the relation
\begin{align*}
:\bLambda_{\overline{1}}(q^{-\frac{\beta}{2}(l-2)-\frac{1}{2}}z)S_0(q^{-1}z)S_1(z): 
=q_1
:\bLambda_{\overline{2}}(q^{-\frac{\beta}{2}(l-2)-\frac{1}{2}}z)S_0(q^{-1}z)S_1(q^{-1}z):
\end{align*}
we arrive at \eqref{D1Na}.

Next, we show \eqref{D1Nb}.

As can be seen from normal orderings of
 $S[(l-1;l-1,N)]{\bf h}_N({\bf u}_0,\ldots, {\bf u}_l;\vartheta)$ and
$S[(l;l,1)] {\bf h}_N({\bf u}_0,\ldots, {\bf u}_l;\vartheta)$,
there is a pole of $z_{l-1,N}$ in the region of the integral 
$\int_{C_{1}-C_{q^{-1}}}\underline{dz_{l-1,N}}$
and there is a pole of $z_{l,1}$ in the region of the integral
$\int_{C_{\delta_l}-C_{q^{-1}\delta_l}}
\underline{dz_{l,1}}$.
Taking these residues gives
\begin{align}
\Delta^{(1,N)}_1[(l-1; l-1,N)]
&=
N(1-q_1)\int_{C[(l-1,N)]}\underline{dz}
{\bf S}_{[0,l-2]} 
\overset{\curvearrowright}{\displaystyle \prod_{1\leq a \leq N-1}}S_{l-1}(z_{l-1,a})
\notag\\
&\times
:\bLambda_{l-1}(q^{\frac{1}{2}}z_{l,1})S_{l-1}(z_{l,1})S_l(z_{l,1}):
\overset{\curvearrowright}{\displaystyle \prod_{2 \leq a \leq N}}S_l(z_{l,a})
\left.{\bf h}_N({\bf u}_0,\ldots, {\bf u}_l;\vartheta)\right|_{z_{l-1,N}=z_{l,1}}
,\notag
\end{align}
\begin{align}
\Delta^{(1,N)}_1[(l;l,1)]&=N(1-q_1^{-1})\int_{C[(l,1)]}\underline{dz}
{\bf S}_{[0,l-2]} 
\overset{\curvearrowright}{\displaystyle \prod_{1\leq a \leq N-1}}S_{l-1}(z_{l-1,a})\notag\\
&\times
:\bLambda_{l}(q^{-\frac{1}{2}}z_{l-1,N})
S_{l-1}(z_{l-1,N})S_l(q^{-1}z_{l-1,N}):
\overset{\curvearrowright}{\displaystyle \prod_{2 \leq a \leq N}}S_l(z_{l,a})
\left.{\bf h}_N({\bf u}_0,\ldots, {\bf u}_l;\vartheta)\right|_{z_{l,1}=q^{-1}z_{l-1,N}}.\notag
\end{align}
Using the relation
\[
:\bLambda_{l-1}(q^{-\frac{1}{2}}z)S_{l-1}(q^{-1}z)S_l(q^{-1}z): 
=q_1^{-1}:\bLambda_{l}(q^{-\frac{1}{2}}z)S_{l-1}(z)S_l(q^{-1}z):
\]
we conclude \eqref{D1Nb}.

\end{proof}

\subsection{Proof of $[{\mathbf I}_2, {\mathbf G}_N(\vartheta)]=0$}
\label{section:I2-GN} 
  
We show Proposition \ref{prop:duality1} in the case of $M=2$ and $N \geq 2$.

We start from ${\bf I}_2 \cdot {\bf G}_N(\vartheta)$.
We set the radii of integration circles $R_1, R_2$ as $|w_1|=R_1$, $|w_2|=R_2$.
We move the contours from $R_1=R_2=R \gg 1$ to $R_1=R_2=1$ 
imposing the conditions:
\begin{align}
1< R_2<R_1< R \quad \textrm{and} \quad 
|q| R_1<R_2,\notag
\end{align}
in the case of $1<R_1, R_2 <R$.
First, in the same way as for $\mathbf{I}_1 \cdot \mathbf{G}_N(\vartheta)$, 
we take the residue in $w_a=\bar{\eta}_{k;j}[z_{j,s}]$ with $(k; j)$ in (\ref{def:leta1}) and $1\le s_1\le Na_{j}$.
After taking the residue in $w_a$,
new poles will appear with respect to the remaining variable $w_b$, $b\neq a$.
We rename the variables $w_a, w_b$ to $w_2, w_1$, respectively.
From Appendix \ref{appendix:normal-ordering},
each term $\bLambda_{k_2}(w_2){\mathcal A}_{k_1; j}[z_{j, s}]$, $1\leq s \leq N a_j$, $k_1, k_2 \in \{\bar{1}, \bar{2}, \ldots,\bar{l}\}$,
can have poles at $w_1=\overline{\eta}_{k_2,k_1; j}[z_{j, s}]$ and
$\overline{\eta}_{k_1; j}^{k_2}[z_{j,s}]$.
Here
\begin{align}
\overline{\eta}_{k_2,k_1; j}[z]=\left\{\begin{array}{cc}
\overline{e}_0^{(2)}[z],&
(k_2,k_1;j)=(\bar{1}, \bar{1};0), (\bar{1}, \bar{2};0), (\bar{2}, \bar{2};0),
\\
\overline{d}_{j}^{(2)}[z],&
(k_2,k_1; j)=(k, \overline{j}; j)~\mathrm{with}
~1\leq j \leq l-2,~ \overline{1} \succeq k \succeq \overline{j+1},
\\
\overline{e}_{j}^{(2)}[z],&
(k_2,k_1;j)=(k, \overline{j+1}; j)~\mathrm{with}~1\leq j \leq l-2,
~ \overline{1} \succeq k \succeq \overline{j+1},
\end{array}\right.
\label{def:leta2}
\end{align}
\begin{align}
\overline{\eta}_{k_1; j}^{k_2}[z]=\left\{\begin{array}{cc}
\overline{D}_{j}^{(1)}[z],&
(k_2,k_1; j)=(k, \overline{j}; j)~\mathrm{with}
~1\leq j \leq l-3,~ \overline{1} \succeq k \succeq \overline{j},
\\
\overline{E}_{j}^{(1)}[z],&
(k_2,k_1; j)=(k, \overline{j+1}; j)~\mathrm{with}~1\leq j \leq l-3,~ \overline{1} \succeq k \succeq \overline{j}.
\end{array}\right.
\label{def:leta3}
\end{align}
For $(k_2, k_1;j)$ in the above, we set
\begin{align*}
&\nu_{k_2,k_1;j}\,\mathcal{A}_{k_2, k_1;j}[z] 
=\mathop{\mathrm{res}}_{w_2=\overline{\eta}_{k_2, k_1;j}[z]}
\left(\mathop{\mathrm{res}}_{w_1=\overline{\eta}_{k_1;j}[z]}
\bLambda_{k_2}(w_2)
\bLambda_{k_1}(w_1)
S_j(z)\omega_2(w_2, w_1)\frac{dw_1}{w_1}\right)
\frac{dw_2}{w_2}\,,\\
&{\mathcal A}_{k_2, k_1; j}[z]= 
:\bLambda_{k_2}(\overline{\eta}_{k_2,k_1;j}[z])
\bLambda_{k_1}(\overline{\eta}_{k_1; j}[z])S_j(z):\,,
\\
&\nu_{k_1;j}^{k_2}\,\mathcal{A}_{k_1;j}^{k_2}[z] 
=\mathop{\mathrm{res}}_{w_2=\overline{\eta}_{k_1;j}^{k_2}[z]}
\left(\mathop{\mathrm{res}}_{w_1=\overline{\eta}_{k_1;j}[z]}
\bLambda_{k_2}(w_2)
\bLambda_{k_1}(w_1)
S_j(z)\omega_2(w_2,w_1)\frac{dw_1}{w_1}\right)
\frac{dw_2}{w_2}\,,\\
&{\mathcal A}_{k_1; j}^{k_2}[z]= 
:\bLambda_{k_2}(\overline{\eta}_{k_1;j}^{k_2}[z])\bLambda_{k_1}(\overline{\eta}_{k_1; j}[z])S_j(z):\,.
\end{align*}
The constants $\nu_{k_2,k_1;j}$ and $\nu_{k_1;j}^{k_2}$ will be given in Lemma \ref{lem:nukkj} below.

For each $0\leq i \leq l$, the integrand is a symmetric function of $z_{i,1}, z_{i,2},\ldots, z_{i, Na_i}$.
Taking this symmetry into account and calculating residues, we obtain ${\rm LHS}_{2, N}(a,b,c)$
as shown in Appendix \ref{appendix:I2-GN}.
\begin{align}
{\rm LHS}_{2,N}(1,0,0)&
=2N
\int_{C_{[\delta]}}\underline{dz}\,\int_{C_1}\underline{dw_1}\,
\Bigl(\sum_{k=\overline{1}, \overline{2}}
\nu_{k;0}{\bf S}[(k;0,1)]\omega_2(w_1, \overline{\eta}_{k;0}[z_{0,1}])\notag\\
&+\sum_{k=\overline{1}, \overline{2}}\nu_{k;1}
{\bf S}[(k; 1,1)]\omega_2(w_1, \overline{\eta}_{k;1}[z_{1,1}])\Bigr)\bh\,,
\notag
\\
{\rm LHS}_{2,N}(0,1,0)&=
4N
\int_{C_{[\delta]}}\underline{dz}\,\int_{C_1}\underline{dw_1}\,
\sum_{i=2}^{l-2}
\sum_{k=\overline{i}, \overline{i+1}}
\nu_{k;i} {\bf S}[(k; i,1)]\omega_2(w_1,\overline{\eta}_{k;i}[z_{i,1}])\bh,
\notag
\\
{\rm LHS}_{2,N}(0,0,1)&
=2N
\int_{C_{[\delta]}}\underline{dz}\,\int_{C_1} \underline{dw_1}
\Bigl(\sum_{k=\overline{l-1}, l}\nu_{k;l-1}
{\bf S}[(k;l-1,1)]\omega_2(w_1, \overline{\eta}_{k;l-1}[z_{l-1;1}])\notag\\
&+\sum_{k=\overline{l-1}, \overline{l}}\nu_{k;l}
{\bf S}[(k;l,1)]\omega_2(w_1, \overline{\eta}_{k;l}[z_{l,1}])\Bigr)\bh\,,
\notag
\end{align}
\begin{align}
{\rm LHS}_{2,N}(2,0,0)&=2N^2 \int_{C_{[\delta]}} \underline{dz}
\sum_{k_1, k_2=\bar{1}, \bar{2}}\nu_{k_1;0} \nu_{k_2;1}
{\bf S}[(k_1;0,1), (k_2;1, 1)]\omega_2(\overline{\eta}_{k_1;0}[z_{0,1}], \overline{\eta}_{k_2;1}[z_{1,1}])
\notag
\\
&\times \bh
\notag
\\
&+N(N-1) \int_{C_{[\delta]}}\underline{dz}
\sum_{i=0,1}
\sum_{k_1, k_2=\bar{1}, \bar{2}}
\nu_{k_1;i}\nu_{k_2;i}
{\bf S}[(k_1;i,1), (k_2;i,2)]\omega_2(\overline{\eta}_{k_1;i}[z_{i,1}], \overline{\eta}_{k_2;i}[z_{i,2}])
\notag\\
&\times
\bh\notag\\
&+2N \int_{C_{[\delta]}}\underline{dz}
\sum_{k_1,k_2=\bar{1}, \bar{2}}\nu_{k_2,k_1;1}
{\bf S}[(k_2,k_1;1,1)]\bh\notag
\\
&+2N \int_{C_{[\delta]}}\underline{dz}
\sum\limits_{\substack{\scriptscriptstyle k_1,k_2=\bar{1}, \bar{2}
\\ \scriptscriptstyle k_2 \succeq k_1}}\nu_{k_2,k_1;0}
{\bf S}[(k_2,k_1;0,1)]\bh\notag
\\
&
+2N \int_{C_{[\delta]}}\underline{dz}
\sum_{k=\bar{1}, \bar{2}}
\nu_{k;1}^{\bar{1}}
{\bf S}[(\boxed{\bar{1}},k;1,1)]\bh,
\notag
\\
{\rm LHS}_{2,N}(0,2,0)&=
2N(2N-1)\int_{C_{[\delta]} }\underline{dz}
\sum_{2\leq i \leq l-2}\sum_{k_1, k_2=\bar{i}, \overline{i+1}}\nu_{k_1;i}\nu_{k_2;i}\mathbf{S}[(k_1;i,1), (k_2;i,2)]\notag\\
&\times \omega_2(\overline{\eta}_{k_1;i}[z_{i,1}], \overline{\eta}_{k_2;i}[z_{i,2}])\bh
\notag
\\
&+8N^2 \int_{C_{[\delta]}}\underline{dz} \sum_{2\leq i<j \leq l-2}\sum_{k_1=\bar{i}, \overline{i+1}}\sum_{k_2=\bar{j}, \overline{j+1}}
\nu_{k_1;i}\nu_{k_2;j}\mathbf{S}[(k_1;i,1), (k_2;j,1)]\notag\\
&\times
\omega_2(\overline{\eta}_{k_1;i}[z_{i,1}], \overline{\eta}_{k_2;j}[z_{j,1}])\bh
\notag\\
&+4N \int_{C_{[\delta]}}\underline{dz}\sum_{2\leq i \leq l-2}\sum_{k_2 \succeq \overline{i+1}} 
\sum_{k_1=\bar{i}, \overline{i+1}} \nu_{k_2,k_1;i}\mathbf{S}[(k_2,k_1;i,1)]\bh\notag\\
&
+4N \int_{C_{[\delta]}}\underline{dz}\sum_{2\leq i \leq l-3}\sum_{k_2 \succeq \overline{i}} 
\sum_{k_1=\bar{i}, \overline{i+1}} \nu^{k_2}_{k_1;i}\mathbf{S}[(\boxed{k_2}, k_1;i,1)]\bh,
\notag
\\
{\rm LHS}_{2,N}(0,0,2)&=
2N(N-1)\int_{C_{[\delta]}}\underline{dz}
\sum_{k_1,k_2=\overline{l-1},l}
\nu_{k_1;l-1}\nu_{k_2;l-1}
\mathbf{S}[(k_1;l-1,1), (k_2;l-1,2)]\notag
\\
&\times
\omega_2(\overline{\eta}_{k_1;l-1}[z_{l-1,1}], \overline{\eta}_{k_2;l-1}[z_{l-1,2}])\bh
\notag\\
&+2N^2\int_{C_{[\delta]}}\underline{dz}
\sum_{k_1=\overline{l-1},l}
\sum_{k_2=\overline{l-1}, \overline{l}}
\nu_{k_1;l-1}\nu_{k_2;l}
\mathbf{S}[(k_1;l-1,1), (k_2;l,1)]\notag\\
&\times \omega_2(\overline{\eta}_{k_1;l-1}[z_{l-1,1}], \overline{\eta}_{k_2;l}[z_{l,1}])\bh
\notag\\
&+2N(N-1)\int_{C_{[\delta]}}\underline{dz}
\sum_{k_1,k_2=\overline{l-1},\overline{l}}
\nu_{k_1;l}\nu_{k_2;l}
\mathbf{S}[(k_1;l,1), (k_2;l,2)]\notag
\\
&\times
\omega_2(\overline{\eta}_{k_1;l}[z_{l,1}], \overline{\eta}_{k_2;l}[z_{l,2}])\bh,
\notag
\end{align}
\begin{align}
\mathrm{LHS}_{2,N}(1,1,0)
&=4N^2 \int_{C_{[\delta]}}\underline{dz}
\sum_{i=0,1}\sum_{2\leq j \leq l-2}\sum_{k_1=\bar{1},\bar{2}}\sum_{k_2=\bar{j}, \overline{j+1}}
\nu_{k_1;i}\nu_{k_2;j}
\mathbf{S}[(k_1;i,1), (k_2;j,1)]\notag\\
&\times \omega_2(\overline{\eta}_{k_1;i}[z_{i,1}], \overline{\eta}_{k_2;j}[z_{j,1}]])\bh,\notag
\\
\mathrm{LHS}_{2,N}(0,1,1)&
=4N^2 \int_{C_{[\delta]}}\underline{dz}
\sum_{2\leq i \leq l-2}\sum_{k_1=\bar{i}, \overline{i+1}}\sum_{k_2=\overline{l-1}, l}\nu_{k_1;i}\nu_{k_2;l-1}
\mathbf{S}[(k_1;i,1), (k_2;l-1,1)]\notag\\
&\times
\omega_2(\overline{\eta}_{k_1; i}[z_{i,1}], \overline{\eta}_{k_2; l-1}[z_{l-1,1}])
\bh
\notag\\
&+
4N^2 \int_{C_{[\delta]}}\underline{dz}
\sum_{2\leq i \leq l-2}\sum_{k_1=\bar{i}, \overline{i+1}}\sum_{k_2=\overline{l-1}, \bar{l}}\nu_{k_1;i}\nu_{k_2;l}
\mathbf{S}[(k_1;i,1), (k_2; l,1)]\notag\\
&\times
\omega_2(\overline{\eta}_{k_1; i}[z_{i,1}], \overline{\eta}_{k_2; l}[z_{l,1}])\bh,
\notag
\\
\mathrm{LHS}_{2,N}(1,0,1)&=
2N^2 \int_{C_{[\delta]}}\underline{dz}
\sum_{i=0,1}\sum_{k_1=\bar{1},\bar{2}}\sum_{k_2=\overline{l-1},l}
\nu_{k_1;i}\nu_{k_2;l-1}
\mathbf{S}[(k_1;i,1), (k_2;l-1,1)]\notag\\
&\times \omega_2(\overline{\eta}_{k_1;i}[z_{i,1}], \overline{\eta}_{k_2;l-1}[z_{l-1,1}])\bh
\notag
\\
&+2N^2 \int_{C_{[\delta]}}\underline{dz}
\sum_{i=0,1}\sum_{k_1=\bar{1},\bar{2}}\sum_{k_2=\overline{l-1},\bar{l}}
\nu_{k_1;i}\nu_{k_2;l}
\mathbf{S}[(k_1;i,1), (k_2;l,1)]\notag\\
&\times \omega_2(\overline{\eta}_{k_1;i}[z_{i,1}], \overline{\eta}_{k_2;l}[z_{l,1}])\bh.
\notag
\end{align}
Here we set
\begin{align}
{\bf S}[(k_2, k_1; j,b)]&= 
\overset{\curvearrowright}{\prod_{0\leq i< j}}
\overset{\curvearrowright}{\prod_{~~1\leq a \leq Na_i}} S_i(z_{i,a}) 
\cdot
\overset{\curvearrowright}{\prod_{~~1\leq a<b}}S_j(z_{j,a})~\notag\\
&\times
{\mathcal A}_{k_2, k_1;j}[z_{j,b}]
\overset{\curvearrowright}{\prod_{b<a \leq Na_j}}
S_j(z_{j,a})
\cdot \overset{\curvearrowright}{\prod_{j<i \leq l}}
\overset{\curvearrowright}{\prod_{~~1\leq a \leq Na_i}} 
S_i(z_{i,a})\,,
\label{Skkjb1}
\\
{\bf S}[(\boxed{k_2}, k_1; j,b)]&= 
\overset{\curvearrowright}{\prod_{0\leq i< j}}
\overset{\curvearrowright}{\prod_{~~1\leq a \leq Na_i}} S_i(z_{i,a}) 
\cdot
\overset{\curvearrowright}{\prod_{~~1\leq a<b}}S_j(z_{j,a})~\notag\\
&\times
{\mathcal A}_{k_1;j}^{k_2}[z_{j,b}]
\overset{\curvearrowright}{\prod_{b<a \leq Na_j}}
S_j(z_{j,a})
\cdot \overset{\curvearrowright}{\prod_{j<i \leq l}}
\overset{\curvearrowright}{\prod_{~~1\leq a \leq Na_i}} 
S_i(z_{i,a})\,,
\label{Skkjb2}
\end{align}
\begin{align}
{\bf S}[(k_1;j_1,b_1), (k_2;j_2,b_2)]&
= 
\overset{\curvearrowright}{\prod_{0\leq i< j_1}}
\overset{\curvearrowright}{\prod_{~~1\leq a \leq Na_i}} S_i(z_{i,a}) 
\cdot
\overset{\curvearrowright}{\prod_{~~1\leq a<b_1}}S_{j_1}(z_{j_1,a})~
{\mathcal A}_{k_1;j_1}[z_{j_1,b_1}]
\overset{\curvearrowright}{\prod_{b_1<a \leq Na_{j_1}}}
S_{j_1}(z_{j_1,a})\notag\\
&\times
\overset{\curvearrowright}{\prod_{j_1<i <j_2}}
\overset{\curvearrowright}{\prod_{~~1\leq a \leq Na_i}} 
S_i(z_{i,a})\,
\cdot
\overset{\curvearrowright}{\prod_{~~1\leq a<b_2}}S_{j_2}(z_{j_2,a})~
{\mathcal A}_{k_2;j_2}[z_{j_2,b_2}]
\overset{\curvearrowright}{\prod_{b_2<a \leq Na_{j_2}}}
S_{j_2}(z_{j_2,a})\notag\\
&\times
\overset{\curvearrowright}{\prod_{j_2<i \leq l}}
\overset{\curvearrowright}{\prod_{~~1\leq a \leq Na_i}} 
S_i(z_{i,a})\,,
\label{Skjbkjb1}
\end{align}
where $(k_1;j_1,b_1)\neq (k_2; j_2,b_2)$, $0 \leq j_1 \leq j_2 \leq l$.
In what follows, we use the following abbreviation.
\begin{align}
\mathbf{S}_{[l_1,l_2]}=
\overset{\curvearrowright}{\prod_{l_1 \leq i \leq l_2}}
\overset{\curvearrowright}{\prod_{~~1\leq a \leq Na_i}}
S_{i}(z_{i,a}), \quad 0\leq l_1 \leq l_2 \leq l.\notag 
\end{align}

\vspace{2mm}

Similarly, we start from ${\bf G}_N(\vartheta) \cdot {\bf I}_2$.
We set the radii of integration circles $r_1, r_2$ as $|w_1|=r_1, |w_2|=r_2$.
We move the contours from $r_1=r_2=r \ll 1$ to $r_1=r_2=1$
imposing the conditions
\begin{align}
r< r_2<r_1< 1 \quad \textrm{and} \quad 
|q| r_1<  r_2,
\notag
\end{align}
in the case of $r<r_1,r_2<1$.
First, in the same way
as for $\mathbf{G}_N(\vartheta)\cdot \mathbf{I}_1$, 
we take the residue in $w_a=\eta_{k; j}[z_{j, s}]$ with $(k; j)$ in (\ref{def:eta1}) and $1\le s\le Na_{j}$.
After taking the residue in $w_a$,
new poles will appear with respect to the remaining variable $w_b$, $b \neq a$.
We rename the variables $w_a, w_b$ to $w_1, w_2$, respectively.
From Appendix \ref{appendix:normal-ordering},
each term ${\mathcal A}_{k_1; j}[z_{j, s}] \bLambda_{k_2}(w_2)$, $1\leq s \leq Na_j$, $k_1, k_2 \in \{1,2,\ldots,l\}$,
can have poles at $w_2=\eta_{k_2,k_1; j}[z_{j, s}]$ and
$\eta_{k_1; j}^{k_2}[z_{j,s}]$ for $0\leq j \leq l-2, 1\leq s \leq Na_{j}$.
Here
\begin{align}
\eta_{k_2,k_1; j}[z]=\left\{\begin{array}{cc}
e_0^{(2)}[z],&
(k_2,k_1;j)=(1, 1;0), (1, 2; 0), (2,2; 0),
\\
d_{j}^{(2)}[z],&
(k_2,k_1; j)=(k, j; j)~\mathrm{with}
~1\leq j \leq l-2,~ 1 \preceq k \preceq j+1,
\\
e_{j}^{(2)}[z],&
(k_2,k_1;j)=(k, j+1; j)~\mathrm{with}~1\leq j \leq l-2,~ 1 \preceq k \preceq j+1,
\end{array}\right.
\label{def:eta2}
\end{align}
\begin{align}
\eta_{k_1; j}^{k_2}[z]=\left\{\begin{array}{cc}
D_{j}^{(1)}[z],&
(k_2,k_1; j)=(k, j; j)~\mathrm{with}
~1\leq j \leq l-3,~ 1 \preceq k \preceq j,
\\
E_{j}^{(1)}[z],&
(k_2,k_1; j)=(k, j+1; j)~\mathrm{with}~1\leq j \leq l-3,~ 1 \preceq k \preceq j.
\end{array}\right.
\label{def:eta3}
\end{align}
For $(k_2, k_1;j)$ which appear in the above, we set
\begin{align*}
&\nu_{k_2,k_1;j}\,
\mathcal{A}_{k_2, k_1;j}[z] 
=\mathop{\mathrm{res}}_{w_2=\eta_{k_2, k_1;j}[z]}
\left(
\mathop{\mathrm{res}}_{w_1=\eta_{k_1;j}[z]}
S_j(z)\bLambda_{k_1}(w_1)
\bLambda_{k_2}(w_2)\omega_2(w_1,w_2)\frac{dw_1}{w_1}\right)\frac{dw_2}{w_2}\,,\\
&{\mathcal A}_{k_2, k_1; j}[z]= 
:S_j(z)\bLambda_{k_1}(\eta_{k_1;j}[z])
\bLambda_{k_2}(\eta_{k_2, k_1; j}[z]):\,,
\qquad
\\
&
\nu_{k_1;j}^{k_2} \mathcal{A}_{k_1;j}^{k_2}[z] 
=\mathop{\mathrm{res}}_{w_2=\eta_{k_1;j}^{k_2}[z]}
\left(\mathop{\mathrm{res}}_{w_1=\eta_{k_1;j}[z]}
S_j(z)
\bLambda_{k_1}(w_1)
\bLambda_{k_2}(w_2)
\omega_2(w_1,w_2)\frac{dw_1}{w_1}\right)
\frac{dw_2}{w_2}\,,
\quad
\\
&{\mathcal A}_{k_1; j}^{k_2}[z]= 
:S_j(z)\bLambda_{k_1}(\eta_{k_1;j}[z])\bLambda_{k_2}(\eta_{k_1; j}^{k_2}[z]):\,.
\end{align*}
The constants $\nu_{k_2, k_1; j}$ and $\nu_{k_1;j}^{k_2}$ will be given in Lemma \ref{lem:nukkj} below.

For each $0\leq i \leq l$, the integrand is a symmetric function of $z_{i,1}, z_{i,2},\ldots, z_{i, Na_i}$.
Taking this symmetry into account and calculating residues, we obtain ${\rm RHS}_{2, N}(a,b,c)$, $1\leq a+b+c \leq 2$
by inverting the pictures in Appendix \ref{appendix:I2-GN}
with respect to the circle $|w_i|=1$.
\begin{align}
{\rm RHS}_{2,N}(1,0,0)&
=-2N
\int_{C_{[\delta]}}\underline{dz}\,\oint_{C_1}\underline{dw_2}\,
\Bigl(\sum_{k=1,2}
\nu_{k;0}{\bf S}[(k;0,1)]\omega_2(\eta_{k;0}[z_{0,1}], w_2)\notag\\
&+\sum_{k=1,2}\nu_{k;1}
{\bf S}[(k; 1,1)]\omega_2(\eta_{k;1}[z_{1,1}], w_2)\Bigr)\bh\,,
\notag
\\
{\rm RHS}_{2,N}(0,1,0)&=
-4N
\int_{C_{[\delta]}}\underline{dz}\,\int_{C_1}\underline{dw_2}\,
\sum_{i=2}^{l-2}
\sum_{k=i, i+1}
\nu_{k;i} {\bf S}[(k; i,1)]\omega_2(\eta_{k;i}[z_{i,1}], w_2)\bh,
\notag
\\
{\rm RHS}_{2,N}(0,0,1)&
=-2N
\int_{C_{[\delta]}}\underline{dz}\,\int_{C_1} \underline{dw_2}
\Bigl(\sum_{k=l-1, \bar{l}}\nu_{k;l-1}
{\bf S}[(k;l-1,1)]\omega_2(\eta_{k;l-1}[z_{l-1;1}], w_2)\notag\\
&+\sum_{k=l-1, l}\nu_{k;l}
{\bf S}[(k;l,1)]\omega_2(\eta_{k;l}[z_{l,1}], w_2)\Bigr)\bh\,,
\notag
\\
{\rm RHS}_{2,N}(2,0,0)&=2N^2 \int_{C_{[\delta]}} \underline{dz}
\sum_{k_1, k_2=1,2}
\nu_{k_1;0} \nu_{k_2;1}
{\bf S}[(k_1;0,1), (k_2;1, 1)]\omega_2(\eta_{k_1;0}[z_{0,1}], \eta_{k_2;1}[z_{1,1}])
\notag
\\
&\times \bh
\notag
\\
&+N(N-1) \int_{C_{[\delta]}}\underline{dz}
\sum_{i=0,1}
\sum_{k_1, k_2=1,2}
\nu_{k_1;i}\nu_{k_2;i}
{\bf S}[(k_1;i,1), (k_2;i,2)]\omega_2(\eta_{k_1;i}[z_{i,1}], \eta_{k_2;i}[z_{i,2}])
\notag\\
&\times
\bh\notag\\
&+2N \int_{C_{[\delta]}}\underline{dz}
\sum_{k_1,k_2=1,2}\nu_{k_2,k_1;1}
{\bf S}[(k_2,k_1;1,1)]\bh\notag
\\
&+2N \int_{C_{[\delta]}}\underline{dz}
\sum\limits_{\substack{\scriptscriptstyle k_1,k_2=1,2 \\
\scriptscriptstyle k_2 \preceq k_1}}\nu_{k_2,k_1;0}
{\bf S}[(k_2,k_1;0,1)]\bh\notag
\\
&
+2N \int_{C_{[\delta]}}\underline{dz}
\sum_{k=1,2}
\nu_{k;1}^{1}
{\bf S}[(\boxed{1},k;1,1)]\bh,
\notag
\\
{\rm RHS}_{2,N}(0,2,0)&=
2N(2N-1)\int_{C_{[\delta]} }\underline{dz}
\sum_{1\leq i \leq l-2}\sum_{k_1, k_2=i, i+1}\nu_{k_1;i}\nu_{k_2;i}\mathbf{S}[(k_1;i,1), (k_2;i,2)]\notag\\
&\times \omega_2(\eta_{k_1;i}[z_{i,1}], \eta_{k_2;i}[z_{i,2}])\bh
\notag
\\
&+8N^2 \int_{C_{[\delta]}}\underline{dz} \sum_{2\leq j<i \leq l-2}\sum_{k_1=i, i+1}\sum_{k_2=j, j+1}
\nu_{k_1;i}\nu_{k_2;j}\mathbf{S}[(k_1;i,1), (k_2;j,1)]\notag\\
&\times
\omega_2(\overline{\eta}_{k_1;i}[z_{i,1}], \overline{\eta}_{k_2;j}[z_{j,1}])\bh
\notag\\
&+4N \int_{C_{[\delta]}}\underline{dz}\sum_{2\leq i \leq l-2}\sum_{k_2 \succeq \overline{i+1}} 
\sum_{k_1=\bar{i}, \overline{i+1}} \nu_{k_2,k_1;i}\mathbf{S}[(k_2,k_1;i,1)]\bh\notag\\
&
+4N \int_{C_{[\delta]}}\underline{dz}\sum_{2\leq i \leq l-3}\sum_{k_2 \succeq \overline{i}} 
\sum_{k_1=\bar{i}, \overline{i+1}} \nu^{k_2}_{k_1;i}\mathbf{S}[(\boxed{k_2}, k_1;i,1)]\bh,
\notag
\end{align}
\begin{align}
{\rm RHS}_{2,N}(0,0,2)&=
2N(N-1)\int_{C_{[\delta]}}\underline{dz}
\sum_{k_1,k_2=l-1,\bar{l}}
\nu_{k_1;l-1}\nu_{k_2;l-1}
\mathbf{S}[(k_1;l-1,1), (k_2;l-1,2)]\notag
\\
&\times
\omega_2(\eta_{k_1;l-1}[z_{l-1,1}], \eta_{k_2;l-1}[z_{l-1,2}])\bh
\notag\\
&+2N^2\int_{C_{[\delta]}}\underline{dz}
\sum_{k_1=l-1,\bar{l}}
\sum_{k_2=l-1, l}
\nu_{k_1;l-1}\nu_{k_2;l}
\mathbf{S}[(k_1;l-1,1), (k_2;l,1)]\notag\\
&\times \omega_2(\eta_{k_1;l-1}[z_{l-1,1}], \eta_{k_2;l}[z_{l,1}])\bh
\notag\\
&+2N(N-1)\int_{C_{[\delta]}}\underline{dz}
\sum_{k_1,k_2=l-1,l}
\nu_{k_1;l}\nu_{k_2;l}
\mathbf{S}[(k_1;l,1), (k_2;l,2)]\notag
\\
&\times
\omega_2(\eta_{k_1;l}[z_{l,1}], \eta_{k_2;l}[z_{l,2}])\bh,
\notag
\\
\mathrm{RHS}_{2,N}(1,1,0)
&=4N^2 \int_{C_{[\delta]}}\underline{dz}
\sum_{i=0,1}\sum_{2\leq j \leq l-2}\sum_{k_1=1,2}\sum_{k_2=j, j+1}
\nu_{k_1;i}\nu_{k_2;j}
\mathbf{S}[(k_1;i,1), (k_2;j,1)]\notag\\
&\times \omega_2(\eta_{k_1;i}[z_{i,1}], \eta_{k_2;j}[z_{j,1}])\bh,\notag
\\
\mathrm{RHS}_{2,N}(0,1,1)&
=4N^2 \int_{C_{[\delta]}}\underline{dz}
\sum_{2\leq i \leq l-2}\sum_{k_1=i, i+1}\sum_{k_2=l-1, \bar{l}}\nu_{k_1;i}\nu_{k_2;l-1}
\mathbf{S}[(k_1;i,1), (k_2;l-1,1)]\notag\\
&\times
\omega_2(\eta_{k_1; i}[z_{i,1}], \eta_{k_2; l-1}[z_{l-1,1}])
\bh
\notag\\
&+
4N^2 \int_{C_{[\delta]}}\underline{dz}
\sum_{2\leq i \leq l-2}\sum_{k_1=i, i+1}
\sum_{k_2=l-1, l}\nu_{k_1;i}\nu_{k_2;l}
\mathbf{S}[(k_1;i,1), (k_2; l,1)]\notag\\\
&\times
\omega_2(\eta_{k_1; i}[z_{i,1}], \eta_{k_2; l}[z_{l,1}])\bh,
\notag
\\
\mathrm{RHS}_{2,N}(1,0,1)&=
2N^2 \int_{C_{[\delta]}}\underline{dz}
\sum_{i=0,1}\sum_{k_1=1,2}\sum_{k_2=l-1,\bar{l}}
\nu_{k_1;i}\nu_{k_2;l-1}
\mathbf{S}[(k_1;i,1), (k_2;l-1,1)]\notag\\
&\times \omega_2(\eta_{k_1;i}[z_{i,1}], \eta_{k_2;l-1}[z_{l-1,1}])\bh
\notag
\\
&+2N^2 \int_{C_{[\delta]}}\underline{dz}
\sum_{i=0,1}\sum_{k_1=1, 2}\sum_{k_2=l-1, l}
\nu_{k_1;i}\nu_{k_2;l}
\mathbf{S}[(k_1;i,1), (k_2;l,1)]\notag\\
&\times \omega_2(\eta_{k_1;i}[z_{i,1}], \eta_{k_2;l}[z_{l,1}])\bh.
\notag
\end{align}

The next lemma can be shown by a direct calculation.
\begin{lem}\label{lem:nukkj}
The following hold.
\begin{align}
&\nu^{k_2}_{k_1; i}=\frac{\Bigl(pq, p q^{-\beta+1}, pq^{\beta-1}, pq^{2\beta-1}; p\Bigr)_\infty}{
\Bigl( p, pq^{-\beta+2}, pq^{\beta}, pq^{2\beta-2}; p\Bigr)_\infty}\times\notag
\end{align}
\begin{align}
&\times
\left\{\begin{array}{cl}
\nu_{k_2;i}\nu_{k_1;i} ~q^{-1}
\frac{\displaystyle (1-q)(1-q^{\beta-1})}{\displaystyle (1-q^2)(1-q^{\beta-2})},
& (k_2, k_1; i)=(\bar{i}, \bar{i}; i), (\bar{i}, \overline{i+1}; i), (i, i; i), (i, i+1; i),\\
&\hskip 2.1cm 1\leq i \leq l-3,
\\
\nu_{k_1;i} ~q^{-1}(1-q)
\frac{\displaystyle 1-q^{\beta-1}}{\displaystyle 1-q^{\beta-2}},
&
(k_2,k_1;i)=(\overline{k}, \bar{i}; i), (\bar{k}, \overline{i+1}, i), (k,i; i), (k,i+1;i),
\\
&\hskip 2.1cm 1\leq i \leq l-3, \bar{1} \succeq \bar{k} \succeq \overline{i-1}, 
1 \preceq k \preceq i-1, 
\\
0,& {\rm otherwise}.
\end{array}
\right.\notag
\end{align}
\begin{align}
&\nu_{k_2, k_1; i}=\frac{\Bigl(p q^{-1}, pq, pq^{\beta-1}, pq^{\beta+1}; p\Bigr)_\infty}{\Bigl( p,pq^2,pq^{\beta-2}, pq^{\beta}; p\Bigr)_\infty}\times\notag
\\
&\times
\left\{\begin{array}{cl}
\nu_{k_2;i}\nu_{k_1;i},
& (k_2,k_1;i)=(\bar{1},\bar{2};0), (1,2;0),\\
&
\\
\nu_{k_2;i}\nu_{k_1;i} ~\frac{\displaystyle (1-q)(1-q^{\beta-1})}{\displaystyle (1-q^2)(1-q^{\beta})},
& (k_2, k_1; i)=(\bar{1}, \bar{1}; 0), (\bar{2}, \overline{2}; 0), 
(1, 1; 0), (2, 2; 0),\\
&\hskip 2.1cm (\overline{i+1}, \overline{i}; i), (\overline{i+1}, \overline{i+1}; i), (i+1,i;i), (i+1, i+1; i),\\
&\hskip 2.1cm 1\leq i \leq l-2,
\\
\nu_{k_2;i}\nu_{k_1;i} ~q^{-1}\frac{\displaystyle (1-q)(1-q^{\beta+1})}{\displaystyle (1-q^2)(1-q^{\beta-2})},
&
(k_2,k_1;i)=(\overline{i}, \bar{i}; i), (\bar{i}, \overline{i+1}, i), (i,i; i), (i,i+1;i),
\\
&
\hskip 2.1cm 1\leq i \leq l-2,\\
\nu_{k_1; i} ~q^{-1}(1-q)\frac{\displaystyle 1-q^{\beta-1}}{\displaystyle 1-q^{\beta-2}},
& (k_2,k_1;i)=(\bar{k},\bar{i};i), (\bar{k}; \overline{i+1};i), (k,i;i), (k, i+1,i)\\
&\hskip 2.1cm 1\leq i \leq l-2, \bar{1} \succeq \bar{k} \succeq \overline{i-1}, 1\preceq k \preceq i-1,
\\
0,& {\rm otherwise}.
\end{array}
\right.\notag
\end{align}
\end{lem}

Using Lemma \ref{lem:pole1} and $\omega(w)=\omega(1/w)=\omega(pw)$ one finds
\begin{align}
&{\rm LHS}_{2,N}(1,0,0)-{\rm RHS}_{2,N}(1,0,0)
=N \sum_{k=\overline{1},\overline{2}}
\nu_{k; 0 }\Delta^{(2,N)}_1[(k; 0,1)]
+
N \sum_{k=2,\bar{2}}\nu_{k; 1}
\Delta^{(2,N)}_1[(k; 1,1)],\notag
\\
&{\rm LHS}_{2,N}(0,0,1)-{\rm RHS}_{2,N}(0,0,1)
=
N \sum_{k=l-1,\overline{l}}
\nu_{k;l-1 }
\Delta^{(2,N)}_1[(k; l-1,1)]+
N \sum_{k=l-1, l}
\nu_{k; l}
\Delta^{(2,N)}_1[(k; l,1)],
\notag\end{align}
\begin{align}
&{\rm LHS}_{2,N}(0,1,0)=
4N
\sum_{i=2}^{l-2}\nu_{\overline{i+1};i}
\Delta^{(2,N)}_1[(\overline{i+1}; i,1)],\notag\\
&{\rm RHS}_{2,N}(0,1,0)=
-4N\sum_{i=2}^{l-2}
\nu_{i+1;i}
\Delta^{(2,N)}_1[(i+1; i,1)].\notag
\end{align}
Here we set
\begin{align}
\Delta^{(2,N)}_1[(k;0,s)]
&
=\int_{C[(0,s)]}\underline{dz}
\int_{C_{\delta_0}-C_{q^{-1}\delta_0}}\underline{dz_{0,s}}
\int_{C_1}\underline{dw}~
\bT_1(w){\bf S}[(k;0,s)]
\omega_2(w_1, \overline{e}_0^{(1)}[z_{0,s}])
\notag\\
&\times
\bh\notag,\quad k=\overline{1},\overline{2}, \quad 1\leq s \leq N,
\notag
\\
\Delta^{(2,N)}_1[(\overline{i+1}; i,s)]
&
=
\int_{C[(i,s)]}\underline{dz_{i,s}}
\int_{C_{1}-C_{q^{-1}}}\underline{dz_{i,s}}
\int_{C_1} \underline{dw}~
\bT_1(w)
{\bf S}[(\overline{i+1};i, s)]
\omega_2(w_1, \overline{e}_i^{(1)}[z_{i,s}])\notag\\
&\times \bh,\quad 1\leq i \leq l-2,\quad 1\leq s \leq Na_i,
\notag
\\
\Delta^{(2,N)}_1[(i+1; i,s)]
&=
\int_{C[(i,s)]}\underline{dz_{i,s}}
\int_{C_{1}-C_{q}}\underline{dz_{i,s}}
\int_{C_1} \underline{dw}~
{\bf S}[(i+1;i,s)]\bT_1(w)\omega_2(e_i^{(1)}[z_{i,s}], w)\notag\\
&\times
\bh,\quad 1\leq i \leq l-2,\quad 1\leq s \leq Na_i,
\notag
\\
\Delta^{(2,N)}_1[(k;l-1,s)]
&=
\int_{C[(l-1,s)]}\underline{dz}\int_{C_1-C_{q^{-1}}}\underline{dz_{l-1,s}}
\int_{C_1}\underline{dw}~{\bf S}[(k;l-1,s)]\bT_1(w)
\omega_2(d_{l-1}^{(1)}[z_{l-1,s}], w)\notag\\
&\times \bh, \quad k=l-1,\overline{l},\quad 1\leq s \leq N,
\notag
\\
\Delta^{(2,N)}_1[(k;l,s)]&
=
\int_{C[(l,s)]}\underline{dz}\int_{C_{\delta_l}-C_{q^{-1}\delta_l}}
\underline{dz_{l,s}}
\int_{C_1}\underline{dw}~
{\bf S}[(k;l,s)]\bT_1(w)\omega_2(d_l^{(1)}[z_{l,s}],w)\notag\\
&\times \bh,\quad k=l-1,l, \quad 1\leq s \leq N.
\notag
\end{align}
 Proposition \ref{prop:duality1} for $M=2$ and $a+b+c=1$ is a direct consequence of the following lemmas,
 which are shown in the same way as Lemmas \ref{lem:duality1} and \ref{lem:duality2}.
\begin{lem}
\label{lem:I2-GN:1}
The following hold for $N \geq 1$.
\begin{align}
&
\Delta^{(2,N)}_1[(\overline{2};0,1)]=0,\quad
\Delta^{(2,N)}_1[(2;1,1)]=0,\notag
\\
&\Delta^{(2,N)}_1[(\overline{i+1}; i,1)]=0,\quad
\Delta^{(2,N)}_1[(i+1;i,1)]=0,\quad 2 \leq i \leq l-2,
\notag\\
&\Delta^{(2,N)}_1[(\overline{l}; l-1,1)]=0,\quad
\Delta^{(2,N)}_1[(l-1;l ,1)]=0.
\notag
\end{align}
\end{lem} 
 \begin{lem}
 \label{lem:I2-GN:2}
The following hold for $N \geq 1$.
\begin{align}
&
\Delta^{(2,N)}_1[(\overline{1};0,1)]+
\Delta^{(2,N)}_1[(\overline{2};1,1)]=0,\notag
\\
&\Delta^{(2,N)}_1[(l-1; l-1,1)]
+\Delta^{(2,N)}_1[(l;l,1)]=0.\notag
\end{align}
\end{lem}

Using Lemma \ref{lem:pole1} and $\omega(w)=\omega(1/w)=\omega(pw)$ one finds
\begin{align}
&{\rm LHS}_{2,N}(2,0,0)-{\rm RHS}_{2,N}(2,0,0)\notag\\
=&
N(N-1)
\sum_{k_{1}, k_{2}=\overline{1},\overline{2}}
\nu_{k_{1};0 }
\nu_{k_{2};0}
\Delta^{(2,N)}_2[(k_{1};0,1), (k_{2};0,2)]
\notag\\
+&
N(N-1)
\nu_{\bar{2}; 1} \nu_{\bar{2}; 1}
\Delta^{(2,N)}_2[(\bar{2};1,1), (\bar{2};1,2)]
-
N(N-1)\nu_{2;1}\nu_{2;1}
\Delta^{(2,N)}_2[(2;1,1), (2;1,1)]
\notag\\
+&
2N
\sum_{k=\bar{1}, \bar{2}}
\nu_{k, \bar{2};1} 
\Delta^{(2,N)}_2[(k, \bar{2};1,1)]
-
2N
\sum_{k=1,2}
\nu_{k, 2; 1}
\Delta^{(2,N)}_2[(k, 2; 1,1)]\notag
\\
+&
2N
\nu_{\bar{2};1}^{\bar{1}} 
\Delta^{(2,N)}_2[(\boxed{\bar{1}}, \bar{2};1,1)]
-
2N
\nu_{2; 1}^{1}
\Delta^{(2,N)}_2[(\boxed{1}, 2; 1,1)]\notag
\\
+ &
2N^2 
\sum_{k=\bar{1}, \bar{2}}
\nu_{k;0}\nu_{\bar{2};1}\Delta_2^{(2,N)}[(k;0,1), (\bar{2};1,1)]
-2N^2\sum_{k=1,2}
\nu_{k;0}\nu_{2;1}\Delta_2^{(2,N)}[(k;0,1), (2;1,1)],\notag
\\
\notag
\\
&{\rm LHS}_{2,N}(0,0,2)-{\rm RHS}_{2,N}(0,0,2)\notag
\\
=&N(N-1)\sum_{k_1, k_2=l-1, \bar{l}}
\nu_{k_1;l-1}\nu_{k_2;l-1}
\Delta^{(2,N)}_2[(k_1; l-1,1), (k_2; l-1,2)]
\notag
\\
+&
2N^2
\sum_{k_1=l-1,\bar{l}}\sum_{k_2=l-1,l}
\nu_{k_{1};l-1}\nu_{k_{2};l}
\Delta^{(2,N)}_2[(k_{1};l-1,1), (k_{2};l,1)]\notag
\\
+&
N(N-1)
 \sum_{k_1, k_2=l-1,l}
 \nu_{k_1;l}\nu_{k_2;l}
 \Delta^{(2,N)}_2[(k_{1};l,1),(k_{2};l,2)],
\notag
\\
\notag\\
&{\rm LHS}_{2,N}(0,2,0)\notag\\
=&2N(2N-1)
\sum_{2\leq i \leq l-2}
\nu_{\overline{i+1}; i}
\nu_{\overline{i+1}; i} 
\Delta^{(2,N)}_2[(\overline{i+1}; i,1), (\overline{i+1}; i,2)]\notag
\\
+&
8N^2
\sum_{2\leq i<j \leq l-2}
\nu_{\overline{i+1};i}
\nu_{\overline{j+1};j}
\Delta^{(2,N)}_2[(\overline{i+1};i,1), (\overline{j+1};j,1)]
\notag
\\
+&4N
 \sum_{2\leq i \leq l-2}
 \sum_{
 k \succeq \overline{i+1}}
 \nu_{k, \overline{i+1}; i}
 \Delta^{(2,N)}_2[(k, \overline{i+1}; i, 1)]
+ 4N
 \sum_{2\leq i \leq l-3}
 \sum_{
 k \succeq \overline{i}}
 \nu_{\overline{i+1}; i}^k
 \Delta^{(2,N)}_2[(\boxed{k}, \overline{i+1}; i, 1)],\notag
 \\
 \notag\\
 &{\rm RHS}_{2,N}(0,2,0)\notag\\
=&2N(2N-1)
\sum_{2\leq i \leq l-2}
\nu_{i+1; i}
\nu_{i+1; i} 
\Delta^{(2,N)}_2[(i+1; i,1), (i+1; i,2)]\notag
\\
+&8N^2
\sum_{2\leq i<j \leq l-2}
\nu_{i+1;i}
\nu_{j+1;j}
\Delta^{(2,N)}_2[(i+1;i,1), (j+1;j,1)]
\notag
\\
+&4N
 \sum_{2\leq i \leq l-2}
 \sum_{
 k \preceq i+1}
 \nu_{k, i+1; i}
 \Delta^{(2,N)}_2[(k, i+1; i, 1)]
+ 4N
 \sum_{2\leq i \leq l-3}
 \sum_{
 k \preceq i}
 \nu_{i+1; i}^k
 \Delta^{(2,N)}_2[(\boxed{k}, i+1; i, 1)],
 \notag 
 \end{align}
\begin{align}
\mathrm{LHS}_{2,N}(1,1,0)
&=
4N^2 \sum_{2 \leq i \leq l-2}
\sum_{k=\bar{1},\bar{2}}
\nu_{k;0}\nu_{\overline{i+1};i}
\Delta_2^{(2,N)}[(k;0,1), (\overline{i+1};i,1)]\notag\\
+&4N^2
\sum_{2 \leq i \leq l-2}
\nu_{\bar{2};1}\nu_{\overline{i+1};i}
\Delta_2^{(2,N)}[(\bar{2};1,1), (\overline{i+1};i,1)],
\notag
\\
\notag
\\
\mathrm{RHS}_{2,N}(1,1,0)
&=
4N^2 \sum_{2 \leq i \leq l-2}
\sum_{k=1, 2}
\nu_{k;0}\nu_{i+1;i}
\Delta_2^{(2,N)}[(k;0,1), (i+1; i,1)]\notag\\
+&4N^2
\sum_{2 \leq i \leq l-2}
\nu_{2;1}\nu_{i+1; i}
\Delta_2^{(2,N)}[(2;1,1), (i+1;i,1)],\notag\\
\notag
\\
\mathrm{LHS}_{2,N}(0,1,1)
&=4N^2 \sum_{2\leq i \leq l-2}\sum_{k=\overline{l-1},l}\nu_{\overline{i+1};i}\nu_{k;l-1}
\Delta_2^{(2,N)}[(\overline{i+1}; i,1), (k;l-1,1)]\notag
\\
+&4N^2 \sum_{2\leq i \leq l-2}\sum_{k=\overline{l-1},\bar{l}}
\nu_{\overline{i+1};i}\nu_{k;l}
\Delta_2^{(2,N)}[(\overline{i+1}; i,1), (k;l,1)],
\notag
\\
\notag
\\
\mathrm{RHS}_{2,N}(0,1,1)
&=4N^2 \sum_{2\leq i \leq l-2}\sum_{k=l-1,\bar{l}}
\nu_{i+1;i}\nu_{k;l-1}
\Delta_2^{(2,N)}[(i+1; i,1), (k;l-1,1)]\notag
\\
+&4N^2 \sum_{2\leq i \leq l-2}\sum_{k=l-1,l}
\nu_{i+1;i}\nu_{k;l}
\Delta_2^{(2,N)}[(i+1; i,1), (k;l,1)],
\notag
\end{align}
\begin{align}
&
\mathrm{LHS}_{2,N}(1,0,1)-\mathrm{RHS}_{2,N}(1,0,1)\notag\\
=&2N^2\sum_{k=\overline{l-1},l}\nu_{\bar{2};1}\nu_{k;l-1}
\Delta_2^{(2,N)}[(\bar{2};1,1), (k;l-1,1)]
+
2N^2\sum_{k=\overline{l-1},\bar{l}}\nu_{\bar{2};1}\nu_{k;l}
\Delta_2^{(2,N)}[(\bar{2};1,1), (k;l,1)]\notag
\\
-&
2N^2\sum_{k=l-1,\bar{l}}\nu_{2;1}\nu_{k;l-1}
\Delta_2^{(2,N)}[(2;1,1), (k;l-1,1)]
-
2N^2\sum_{k=l-1, l}\nu_{2;1}\nu_{k;l}
\Delta_2^{(2,N)}[(2;1,1), (k;l,1)]\notag
\\
+&2N^2 \sum_{k_1=\bar{1},\bar{2}}
\sum_{k_2=\overline{l-1},l} \nu_{k_1;0}\nu_{k_2;l-1}\Delta_2^{(2,N)}[(k_1;0,1), (k_2;l-1,1)]\notag\\
+&2 N^2 \sum_{k_1=\bar{1},\bar{2}}
\sum_{k_2=\overline{l-1}, \bar{l}} \nu_{k_1;0}\nu_{k_2;l}\Delta_2^{(2,N)}[(k_1;0,1), (k_2;l,1)].
\notag
\end{align}

Here we set
\begin{align}
\Delta^{(2,N)}_2[(k_1;0,s), (k_2;0,t)]
&=
\int_{C[(0,s),(0,t)]}\underline{dz}\left(\int_{C_{\delta_0}}\underline{dz_{0,s}}\int_{C_{\delta_0}}\underline{dz_{0,t}}
-\int_{C_{q^{-1}\delta_0}}\underline{dz_{0,s}}\int_{C_{q^{-1}\delta_0}}\underline{dz_{0,t}}
\right)\notag
\\
&\times {\bf S}[(k_1;0,s), (k_2;0,t)]\omega_2(z_{0,s}, z_{0,t})\bh,
\notag\\
&k_1,k_2=\overline{1}, \overline{2},
\quad 1 \leq s \neq t \leq N,
\notag
\end{align}
\begin{align}
 \Delta^{(2,N)}_2[(\overline{i+1}; i,s), (\overline{i+1}; i,t)]&=
\int_{C[(i,s),(i,t)]}\underline{dz}\int_{C_1-C_{q^{-1}}} \underline{dz_{i,s}} 
\int_{C_1-C_{q^{-1}}} \underline{dz_{i,t}}
{\bf S}[(\overline{i+1};i,s), (\overline{i+1};i,t)]\notag\\
&\times
\omega_2(\overline{e}_i^{(1)}[z_{i,s}], \overline{e}_i^{(1)}[z_{i,t}])\bh,\notag
\\
&1\leq i \leq l-2,\quad 1\leq s \neq t \leq Na_i, \notag
 \\
 \Delta^{(2,N)}_2[(i+1; i,s), (i+1; i,t)]&=
\int_{C[(i,s),(i,t)]}\underline{dz}\int_{C_1-C_{q}} \underline{dz_{i,s}} 
\int_{C_1-C_{q}} \underline{dz_{i,t}}
{\bf S}[(i+1;i,s), (i+1;i,t)]\notag\\
&\times
\omega_2({e}_i^{(1)}[z_{i,s}], {e}_i^{(1)}[z_{i,t}])\bh,
\notag\\
& 1\leq i \leq l-2, \quad 1\leq s \neq t\leq Na_i,
 \notag
\\
 \Delta^{(2,N)}_2[(\overline{i+1}; i,s), (\overline{j+1}; j,t)]&=
\int_{C[(i,s),(j,t)]}\underline{dz}\int_{C_1-C_{q^{-1}}} \underline{dz_{i,s}} 
\int_{C_1-C_{q^{-1}}} \underline{dz_{j,t}}
{\bf S}[(\overline{i+1};i,s), (\overline{j+1};j,t)]\notag\\
&\times
\omega_2(\overline{e}_i^{(1)}[z_{i,s}], \overline{e}_j^{(1)}[z_{j,t}])\bh,
\notag\\
& 1\leq i<j \leq l-2,\quad 1\leq s \leq Na_i, 1\leq t \leq Na_j,
 \notag
 \\
 \Delta^{(2,N)}_2[(i+1; i,s), (j+1; j,t)]&=
\int_{C[(i,s),(j,t)]}\underline{dz}\int_{C_1-C_{q}} \underline{dz_{i,s}} 
\int_{C_1-C_{q}} \underline{dz_{j,t}}
{\bf S}[(i+1;i,s), (j+1;j,t)]\notag\\
&\times
\omega_2({e}_i^{(1)}[z_{i,s}], {e}_j^{(1)}[z_{j,t}])\bh,\notag\\
& 1\leq i<j \leq l-2, \quad 1\leq s \leq Na_i, 1\leq t \leq Na_j,
\notag
\\
\Delta^{(2,N)}_2[(k;0,s), (\overline{i+1};i,t)]&=
\int_{C[(i,t)]}\underline{dz} \int_{C_1-C_{q^{-1}}}\underline{dz_{i,t}}
{\bf S}[(k;0,s), (\overline{i+1}; i,t)]\notag\\
&\times \omega_2(\overline{e}_{0}^{(1)}[z_{0,s}], \overline{e}_i^{(1)}[z_{i,t}])\bh,
\notag\\
&1\leq i \leq l-2,\quad k=\overline{1}, \overline{2},\quad 1\leq s \leq N, 1\leq t \leq Na_i,
\notag
\\
\Delta^{(2,N)}_2[(k;0,s), (i+1; i,t)]&=
\int_{C[(i,t)]}\underline{dz} \int_{C_1-C_{q}}\underline{dz_{i,t}}
{\bf S}[(k;0,s), (i+1; i,t)]\notag\\
&\times \omega_2({e}_{0}^{(1)}[z_{0,s}], {e}_i^{(1)}[z_{i,t}])\bh,
\notag\\
&1\leq i \leq l-2,\quad k=1,2,\quad 1\leq s \leq N, 1\leq t \leq Na_i,
\notag
\\
\Delta^{(2,N)}_2[(k_1;l-1,s), (k_2;l-1,t)]
&=\int_{C[(l-1,s),(l-1,t)]} \underline{dz}
\left(\int_{C_{q^{-1}}}\underline{dz_{l-1,s}}\int_{C_{q^{-1}}}\underline{dz_{l-1,t}}-
\int_{C_1} 
\underline{dz_{l-1,s}}\int_{C_1} \underline{dz_{l-1,t}}\right)\notag
\\
&\times
{\bf S}[(k_1;l-1,s),(k_2;l-1,t)]
\omega_2(z_{l-1,s}, z_{l-1,t})\bh,\notag
\\
&k_1,k_2=l-1,\overline{l}, \quad
1\leq s \neq t \leq N,
\notag
\\
\Delta^{(2,N)}_2[(k_1;l-1,s), (k_2;l,t)]&
=\int_{C[(l-1,s),(l,t)]} \underline{dz}
\left(\int_{C_{q^{-1}}}\underline{dz_{l-1,s}}\int_{C_{q^{-1}\delta_l}}\underline{dz_{l,t}}-
\int_{C_1} \underline{dz_{l-1,s}}\int_{C_{\delta_l}} \underline{dz_{l,t}}\right)\notag
\\
&\times
{\bf S}[(k_1;l-1,s),(k_2;l,t)]
\omega_2(z_{l-1,s}, z_{l,t})\bh,\notag\\
&k_1=l-1,\overline{l},\quad k_2=l-1,l,\quad 1\leq s, t \leq N,
\notag
\end{align}
\begin{align}
\Delta^{(2,N)}_2[(k_1;l,s), (k_2;l,t)]&=\int_{C[(l,s),(l,t)]} \underline{dz}
\left(\int_{C_{q^{-1}\delta_l}}\underline{dz_{l,s}}\int_{C_{q^{-1}\delta_l}}\underline{dz_{l,t}}-
\int_{C_{\delta_l}} \underline{dz_{l,s}}\int_{C_{\delta_l}} \underline{dz_{l,t}}\right)\notag
\\
&\times
{\bf S}[(k_1;l,s),(k_2;l,t)]
\omega_2(z_{l,s}, z_{l,t})\bh,\notag
\\
&k_1,k_2=l-1,l, \quad 1\leq s\neq t \leq N,
\notag
\\
\Delta^{(2,N)}_2[(\overline{i+1};i,s), (k;l-1,t)]&
=\int_{C[(i,s)]} \underline{dz}
\int_{C_1-C_{q^{-1}}}\underline{dz_{i,s}}
~{\bf S}[(\overline{i+1}; i,s),(k;l-1,t)]
\notag\\
&\times \omega_2(\overline{e}_i^{(1)}[z_{i,s}], \overline{e}_{l-1}^{(1)}[z_{l-1,t}])\bh,\notag
\\
&1\leq i \leq l-2,\quad k=\overline{l-1},l,\quad 1\leq s \leq Na_i, 1\leq t \leq N,
\notag
\\
\Delta^{(2,N)}_2[(i+1;i,s), (k;l-1,t)]&
=\int_{C[(i,s)]} \underline{dz}
\int_{C_1-C_{q}}\underline{dz_{i,s}}
~{\bf S}[(i+1; i,s),(k; l-1,t)]
\notag\\
&\times \omega_2({e}_i^{(1)}[z_{i,s}], {e}_{l-1}^{(1)}[z_{l-1,t}])\bh,\notag
\\
&1\leq i \leq l-2,\quad k=l-1,\bar{l},\quad 1\leq s \leq Na_i, 1\leq t \leq N,
\notag
\\
\Delta^{(2,N)}_2[(\overline{i+1};i,s), (k;l,t)]&
=\int_{C[(i,s)]} \underline{dz}
\int_{C_1-C_{q^{-1}}}\underline{dz_{i,s}}
~{\bf S}[(\overline{i+1}; i,s),(k;l,t)]
\notag\\
&\times \omega_2(\overline{e}_i^{(1)}[z_{i,s}], \overline{e}_{l}^{(1)}[z_{l,t}])\bh,\notag
\\
&1\leq i \leq l-2,\quad k=\overline{l-1}, \bar{l},\quad 1\leq s \leq Na_i, 1\leq t \leq N,
\notag
\\
\Delta^{(2,N)}_2[(i+1;i,s), (k;l,t)]&
=\int_{C[(i,s)]} \underline{dz}
\int_{C_1-C_{q}}\underline{dz_{i,s}}
~{\bf S}[(i+1; i,s),(k; l,t)]
\notag\\
&\times \omega_2({e}_i^{(1)}[z_{i,s}], {e}_{l}^{(1)}[z_{l,t}])\bh,\notag
\\
&1\leq i \leq l-2,\quad  k=l-1,l,\quad 1\leq s \leq Na_i, 1\leq t \leq N,
\notag
\\
\Delta^{(2,N)}_2[(k_1;0,s), (k_2;l-1,t)]
&=\int_{C[(0,s),(l-1,t)]} \underline{dz}
\left(\int_{C_{\delta_0}}\underline{dz_{0,s}}\int_{C_{1}}\underline{dz_{l-1,t}}-
\int_{C_{q^{-1}\delta_0}} 
\underline{dz_{0,s}}\int_{C_{q}} \underline{dz_{l-1,t}}
\right)\notag
\\
&\times
{\bf S}[(k_1;0,s),(k_2;l-1,t)]
\omega_2(
\overline{e}_0^{(1)}[z_{0,s}], \overline{d}_{l-1}^{(1)}[z_{l-1,t}])\bh,\notag
\\
&k_1=\bar{1}, \bar{2}, \quad k_2=\overline{l-1},l, \quad
1\leq s, t \leq N,
\notag
\\
\Delta^{(2,N)}_2[(k_1;0,s), (k_2;l,t)]
&=\int_{C[(0,s),(l,t)]} \underline{dz}
\left(\int_{C_{\delta_0}}\underline{dz_{0,s}}\int_{C_{\delta_l}}\underline{dz_{l,t}}-
\int_{C_{q^{-1}\delta_0}} 
\underline{dz_{0,s}}\int_{C_{q \delta_l}} \underline{dz_{l,t}}
\right)\notag
\\
&\times
{\bf S}[(k_1;0,s),(k_2;l,t)]
\omega_2(
\overline{e}_0^{(1)}[z_{0,s}], \overline{d}_{l}^{(1)}[z_{l,t}])\bh,\notag
\\
&k_1=\bar{1}, \bar{2},\quad  k_2=\overline{l-1},\bar{l}, \quad
1\leq s, t \leq N,
\notag
\\
\Delta^{(2,N)}_2[(k_2,k_1;0,s)]&=
\int_{C[(0,s)]}\underline{dz} \int_{C_{\delta_0}-C_{q^{-2}\delta_0}}\underline{dz_{0,s}}
{\bf S}[(k_2, k_1; 0,s)]\bh,\notag\\
&
k_1,k_2=\overline{1}, \overline{2}, \quad
k_2 \succeq k_1,\quad 1\leq s \leq N,
\notag
\end{align}
\begin{align}
\Delta^{(2,N)}_2\left[\left(\boxed{k}, \overline{i+1}; i,s\right)\right]
&=
\int_{C[(i,s)]}\underline{dz}\int_{C_1-C_{q^{-1}}}\underline{dz_{i,s}}
{\bf S}\left[\left(\boxed{k}, \overline{i+1}; i,s\right)\right]\bh,
\notag\\
&
1\leq i \leq l-3, \quad
k \succeq \overline{i},\quad 
1\leq s \leq Na_i,
\notag
\\
\Delta^{(2,N)}_2\left[\left(\boxed{k}, i+1; i, s\right)\right]&=
\int_{C[(i,s)]}\underline{dz}\int_{C_1-C_{q}}\underline{dz_{i,s}}
{\bf S}[(\boxed{k}, i+1; i,s)]\bh,
\notag\\
&
1\leq i \leq l-3, \quad
k \preceq i, \quad 1\leq s \leq N a_i,
\notag
\\
\Delta^{(2,N)}_2[(k, \overline{i+1}; i, s)]&=
\int_{C[(i,s)]}\underline{dz}\int_{C_1-C_{q^{-1}}}\underline{dz_{i,s}}
{\bf S}[(k, \overline{i+1}; i,s)]\bh, \notag
\\
&1\leq i \leq l-2,\quad
k \succeq \overline{i+1},\quad 1\leq s \leq Na_i,
\notag
\\
\Delta^{(2,N)}_2[(k, i+1; i,s)]&=
\int_{C[(i,s)]}\underline{dz}\int_{C_1-C_{q}}\underline{dz_{i,s}}
{\bf S}[(k, i+1; i,s)]\bh, \notag\\
&
1\leq i \leq l-2,\quad
k \preceq i+1,\quad 1\leq s \leq Na_i.
\notag
\end{align}
Here $C[(i_1,s_1), (i_2, s_2),\ldots,(i_n, s_n)]$ stands for the contour
\begin{align}
C[(i_1,s_1), (i_2,s_2),\ldots, (i_n,s_n)]: |z_{j,s}|=\delta_j,\quad \quad (j, s)\neq (i_1, s_1), \ldots, (i_n, s_n), 0\leq j \leq l, ~1\leq s \leq Na_j.\notag
\end{align}

Proposition \ref{prop:duality1} for $M=2$ and $a+b+c=2$ is a direct 
consequence of Lemmas 
\ref{lem:I2-GN:3}, 
\ref{lem:I2-GN:4}, 
\ref{lem:I2-GN:5}, 
\ref{lem:I2-GN:6}, and
\ref{lem:I2-GN:7} which will be given below.
The relations $\mathrm{LHS}_{M,N}(a,b,c)=0$,
$\mathrm{RHS}_{M,N}(a,b,c)=0$, $1\leq b \leq 2$, $a+b+c=2$,
in Proposition \ref{prop:duality1}
are direct consequences
of Lemma \ref{lem:I2-GN:3}.
\begin{lem}
\label{lem:I2-GN:3}
The following hold for $N \geq 2$.
\begin{align}
&
\Delta_2^{(2,N)}[(\overline{i+1};i,1), (\overline{i+1};i,2)]=0,\quad
2\leq i \leq l-2,
\notag\\
&
\Delta_2^{(2,N)}[(\overline{i+1};i,1), (\overline{j+1};j,1)]=0,\quad
1\leq i<j  \leq l-2,
\notag\\
&
\Delta_2^{(2,N)}[(k;0,1), (\overline{i+1};i,1)]=0,\quad 2\leq i \leq l-2,\quad k=\bar{1}, \bar{2},
\notag\\
&
\Delta_2^{(2,N)}[(\overline{i+1};i,1), (k; j,1)]=0,\quad 2\leq i \leq l-2,\quad (k; j)=(\overline{l-1}; l-1), (l; l-1), (\overline{l-1};l), (\bar{l};l),
\notag\\
&
\Delta^{(2,N)}_2\left[\left(k, \overline{i+1}; i,1\right)\right]=0,\quad 
2\leq i \leq l-2,~ k \succeq \overline{i+1},
\notag\\
&
\Delta^{(2,N)}_2\left[\left(\boxed{k}, \overline{i+1};i,1\right)\right]=0,
\quad 1\leq i \leq l-3,~ k \succeq \overline{i}.
\notag
\end{align}
\begin{align}
&
\Delta_2^{(2,N)}[(i+1; i,1), (i+1; i,2)]=0,\quad
2\leq i \leq l-2,
\notag\\
&
\Delta_2^{(2,N)}[(i+1;i,1), (j+1; j,1)]=0,\quad
1\leq i<j  \leq l-2,
\notag\\
&
\Delta_2^{(2,N)}[(k;0,1), (i+1; i,1)]=0,\quad 2\leq i \leq l-2,\quad k=1,2,
\notag\\
&
\Delta_2^{(2,N)}[(i+1; i,1), (k; j,1)]=0,\quad 2\leq i \leq l-2,\quad (k; j)=(l-1; l-1), (\bar{l}; l-1), (l-1;l), (l; l),
\notag\\
&
\Delta^{(2,N)}_2\left[\left(k, i+1; i,1\right)\right]=0,\quad 
2\leq i \leq l-2,~ k \preceq i+1,
\notag\\
&
\Delta^{(2,N)}_2\left[\left(\boxed{k}, i+1;i,1\right)\right]=0,
\quad 1\leq i \leq l-3,~ k \preceq i.
\notag
\end{align}
\end{lem}
\begin{proof}
They are shown in the same way as Lemma \ref{lem:duality1}.
\end{proof}

The relation $\mathrm{LHS}_{2,N}(2,0,0)=\mathrm{RHS}_{2,N}(2,0,0)$ in Proposition
\ref{prop:duality1}
is a direct consequence of Lemmas
\ref{lem:I2-GN:4} and \ref{lem:I2-GN:5}
\begin{lem}
\label{lem:I2-GN:4}
The following hold for $N \geq 2$.
\begin{align}
&\sum_{k=1,2}\nu_{k,2;1}\Delta^{(2,N)}_2[(k,2;1,1)]
+N (q_1^{-1}-1)^2
\sum_{k=1,2}\Delta^{(2,N)}_2[(k;0,1), (2;1,1)]=0,
\notag
\end{align}
\begin{align}
&\Delta^{(2,N)}_2[(2;1,1), (2;1,2)]=0.\notag
\end{align}
\end{lem}
\begin{proof}
First, we consider $\Delta^{(2,N)}_2[(2;1,1), (2;1,2)]$.

As can be seen from the normal ordering of
$\mathbf{S}[(2;1,1), (2;1,2)]\omega_2(\overline{e}_1^{(1)}[z_{1,1}], \overline{e}_1^{(1)}[z_{1,2}])\bh$, 
there is no pole of $z_{1,s}$ in the region of the integral
$\int_{C_1-C_q}\underline{dz_{1,s}}$.
Hence, $\Delta^{(2,N)}_2[(2;1,1), (2;1,2)]=0$.

 Next, we consider $\Delta^{(2,N)}_2[(k,2;1,1)]$.

 As can be seen from the normal ordering
 \begin{align}
 S_0(z_{0,s}){\mathcal A}_{k,2;1}(z_{1,1})=
 \frac{(1-q^{\beta-2}\frac{z_{1,1}}{z_{0,s}})(1-q^{\beta-1}\frac{z_{1,1}}{z_{0,s}})}{(1-q^{-1}\frac{z_{1,1}}{z_{0,s}})(1-\frac{z_{1,1}}{z_{0,s}})}
 : S_0(z_{0,s}){\mathcal A}_{k,2;1}(z_{1,1}):,\quad k=1,2,\notag
 \end{align}
 there is a pole at $z_{1,1}=q z_{0,s}$ in the region of the integral
$\int_{C_1-C_q}\underline{dz_{1,1}}$.
 Taking the residue we obtain
 \begin{align}
 &
 \Delta^{(2,N)}_2[(k,2;1,1)]=-N\frac{(1-q^{\beta-1})(1-q^\beta)}{1-q}
 \int_{C[(1,1)]}\underline{dz} \overset{\curvearrowright}{\prod_{1\leq s \leq N-1}}S_0(z_{0,s})
 \notag\\
 &\times
 :\bLambda_{k}(q^{-\frac{1}{2}+\frac{\beta}{2}(l-2)}z_{0,N})
 \bLambda_2(q^{\frac{1}{2}+\frac{\beta}{2}(l-2)}z_{0,N})S_0(z_{0,N})S_1(qz_{0,N}):
 \overset{\curvearrowright}{\prod_{2\leq t \leq N}}S_1(z_{1,t}){\bf S}_{[2,l]}\bh,
 k=1,2.\notag
 \end{align}
 
 Next, we consider $\Delta^{(2,N)}_2[(k;0,1), (2;1,1)]$.
 
 Here we find it convenient to use
 $\Delta^{(2,N)}_2[(k;0,N), (2;1,1)]$ in place of
 $\Delta^{(2,N)}_2[(k;0,1), (2;1,1)]$.

 As can be seen from the normal orderings
 \begin{align}
 &
 {\mathcal A}_{k;0}(z_{0,N}){\mathcal A}_{2;1}(z_{1,1})\omega_2(z_{0,N}, z_{1,1})
 = 
  : {\mathcal A}_{k;0}(z_{0,N}){\mathcal A}_{2;1}(z_{1,1}):\notag\\
  &\times
 \frac{(p \frac{z_{1,1}}{z_{0,N}}, q^\beta \frac{z_{1,1}}{z_{0,N}}, 
 p \frac{z_{0,N}}{z_{1,1}}, q^\beta \frac{z_{0,N}}{z_{1,1}}; p)_\infty}{(pq^{\beta-1}\frac{z_{1,1}}{z_{0,N}},
 q \frac{z_{1,1}}{z_{0,N}}, pq^{\beta-1}\frac{z_{0,N}}{z_{1,1}}, q\frac{z_{0,N}}{z_{1,1}};p )_\infty}
 \times \left\{\begin{array}{cc}
 \frac{1-q^{1-\beta}\frac{z_{0,N}}{z_{1,1}}}{1-q^{\beta-1}\frac{z_{0,N}}{z_{1,1}}},& k=1,\\
 1,&k=2,
 \end{array}
 \right.\notag
 \\
 &{\mathcal A}_{k;0}(z_{0,N})S_1(z_{1,t})=q^{\beta-1}\frac{1-q^{-\beta+1}\frac{z_{1,t}}{z_{0,N}}}{1-\frac{z_{1,t}}{z_{0,N}}}
 : {\mathcal A}_{k;0}(z_{0,N})S_1(z_{1,t}):, k=1,2,\notag
 \end{align}
there is a pole at $z_{1,1}=q z_{0,N}$ in the region of the integral 
$\int_{C_1-C_{q^{-1}}}\underline{dz_{1,1}}$.
Taking the residue we obtain 
\begin{align}
&
\Delta^{(2,N)}_2[(k;0,N), (2;1,1)]=\frac{(qp, q^{\beta+1}, pq^{-1}, q^{\beta-1};p)_\infty}{
(q^\beta p, q^2, pq^{\beta-2}, p; p)_\infty}\int_{C[(1,1)]}\underline{dz}
 \overset{\curvearrowright}{\prod_{1\leq s \leq N-1}}S_0(z_{0,s})\notag\\
&\times 
 :\bLambda_k(q^{-\frac{1}{2}+\frac{\beta}{2}(l-2)}z_{0,N})
 \bLambda_2(q^{\frac{1}{2}+\frac{\beta}{2}(l-2)}z_{0,N})
 S_0(z_{0,N})S_1(qz_{0,N}):
  \overset{\curvearrowright}{\prod_{2\leq t \leq N}}S_1(z_{1,t}){\bf S}_{[2,l]}
  \notag\\
  &\times
  \bh
  \times
  \left\{\begin{array}{cc}
  \frac{1-q^{-\beta}}{1-q^{\beta-2}},& k=1,\\
  1,&k=2.
  \end{array}
  \right.\notag
\end{align}
Adding the above terms we find 
Lemma \ref{lem:I2-GN:4}.
\end{proof}

\begin{lem}
\label{lem:I2-GN:5}
The following hold for $N \geq 2$.
\begin{align}
&(N-1)(q_1-1)^2
\sum_{k_1,k_2=\overline{1},\overline{2}}\Delta^{(2,N)}_2[(k_1;0,1), (k_2;0,2)]\notag
\\
+&(N-1)(q_1-1)^2 \Delta^{(2,N)}_2[(\overline{2};1,1), (\overline{2};1,2)]
+2N(q_1-1)^2
\sum_{k=\overline{1}, \overline{2}}\Delta^{(2,N)}_2[(k;0,1), (\overline{2};1,1)]=0,
\notag
\\
&\sum_{k=\overline{1}, \overline{2}}
\nu_{k, \overline{2};1}\Delta^{(2,N)}_2[(k,\overline{2};1,1)]
+
\sum\limits_{\substack{\scriptscriptstyle 
k_1,k_2=\overline{1}, \overline{2}
\\ \scriptscriptstyle k_1 \succeq k_2}}
\nu_{k_1, k_2 ;0}\Delta^{(2,N)}_2[(k_1, k_2;0,1)]=0.
\notag
\end{align}
\end{lem}

\begin{proof}
First, we consider $\Delta^{(2,N)}_2[(k_1;0,1), (k_2;0,2)]$.

{Here we find it convenient to use $\Delta_2^{(2,N)}[(k_1;0,N-1), (k_2;0,N)]$ 
in place of $\Delta_2^{(2,N)}[(k_1;0,1), (k_2;0,2)]$.}
We rewrite the integral as follows.
\begin{align}
&\Delta^{(2,N)}_2[(k_1;0,N-1), (k_2;0,N)]=
\int_{C[(0,N-1), (0,N)]}\underline{dz}\notag\\
&\times
\left(-\int_{C_{\delta_0}-C_{q^{-1} \delta_0}} \underline{dz_{0,N-1}}
\int_{C_{\delta_0}-C_{q^{-1}\delta_0}}\underline{dz_{0,N}}
+\int_{C_{\delta_0}-C_{q^{-1} \delta_0}} \underline{dz_{0,N-1}}
\int_{C_{\delta_0}}\underline{dz_{0,N}}
+\int_{C_{\delta_0}} \underline{dz_{0,N-1}}
\int_{C_{\delta_0}-C_{q^{-1}\delta_0}}\underline{dz_{0,N}}\right)
\notag
\\
&\times
{\bf S}[(k_1;0,N-1), (k_2;0,N)]
\omega_2(z_{0,N-1}, z_{0,N})
\bh,\quad k_1,k_2=\overline{1}, \overline{2}.
\notag
\end{align}

As can be seen from the normal ordering of $\mathbf{S}[(k_1;0,N-1), (k_2;0,N)]\omega_2(\overline{e}_0^{(1)}[z_{0,N-1}], \overline{e}_0^{(1)}[z_{0,N}])$
$\times \bh$,
there is a pole at $z_{0,s}=q^{-1}z_{1,t}$ in the region of integral 
$\int_{C_{\delta_0}-C_{q^{-1}\delta_0}}\underline{dz_{0,s}}$, $s=N-1,N$.
Taking the residues we obtain 
\begin{align}
&\int_{C[(0,N-1), (0,N)]}\underline{dz}
\int_{C_{\delta_0}-C_{q^{-1} \delta_0}} \underline{dz_{0,N-1}}
\int_{C_{\delta_0}-C_{q^{-1}\delta_0}}\underline{dz_{0,N}}
{\bf S}[(k_1;0,N-1), (k_2;0,N)]\notag\\
&\times
\omega_2(z_{0,N-1}, z_{0,N})\bh
=\delta_{k_1,\overline{1}}\delta_{k_2,\overline{1}}N(N-1) (1-q^{\beta-1})^2
\int_{C[(0,N-1), (0,N)]}\underline{dz}\notag\\
&\times
 \overset{\curvearrowright}{\prod_{1\leq s \leq N-2}}S_0(z_{0,s})
 :\bLambda_{\overline{1}}(q^{-\frac{1}{2}-\frac{\beta}{2}(l-2)}z_{1,1})
 S_0(q^{-1}z_{1,1})S_1(z_{1,1}):
  :\bLambda_{\overline{1}}(q^{-\frac{1}{2}-\frac{\beta}{2}(l-2)}z_{1,2})
 S_0(q^{-1}z_{1,2})S_1(z_{1,2}):\notag
 \\
&\times  \overset{\curvearrowright}{\prod_{3\leq t \leq N}}S_1(z_{1,t})
{\bf S}_{[2,l]}
\omega_2(z_{1,1}, z_{1,2})\bh,\quad k_1, k_2=\overline{1}, \overline{2},
\notag
\end{align}
\begin{align}
&\int_{C[(0,N-1)]}\underline{dz}
\int_{C_{\delta_0}-C_{q^{-1} \delta_0}} \underline{dz_{0,N-1}}
{\bf S}[(k_2; 0,N-1), (k_1; 0,N)]
\omega_2(z_{0,N-1}, z_{0,N})\bh
\notag
\\
=&\int_{C[(0,N)]}\underline{dz}
\int_{C_{\delta_0}-C_{q^{-1} \delta_0}} \underline{dz_{0,N}}
{\bf S}[(k_1;0,N-1), (k_2;0,N)]
\omega_2(z_{0,N-1}, z_{0,N})\bh,
\notag
\end{align}
\begin{align}
&\int_{C[(0,N)]}\underline{dz}
\int_{C_{\delta_0}-C_{q^{-1} \delta_0}} \underline{dz_{0,N}}
{\bf S}[(k_1;0,N-1), (k_2;0,N)]
\omega_2(z_{0,N-1}, z_{0,N})\bh
\notag\\
=&\delta_{k_2,\overline{1}}N(1-q^{\beta-1})
\int_{C[(0,N)]}\underline{dz}
 \overset{\curvearrowright}{\prod_{1\leq s \leq N-2}}S_0(z_{0,s})
 {\mathcal A}_{k_1;0}(z_{0,N-1})
  :\bLambda_{\overline{1}}(q^{-\frac{1}{2}-\frac{\beta}{2}(l-2)}z_{1,1})
 S_0(q^{-1}z_{1,1})S_1(z_{1,1}):\notag
 \\
\times&  \overset{\curvearrowright}{\prod_{2 \leq t \leq N}}S_1(z_{1,t})
{\bf S}_{[2,l]}
\omega_2(z_{0,N-1}, q^{-1} z_{1,1})\bh,\quad k_1, k_2=\overline{1}, \overline{2}.
\notag
\end{align}

Next, we consider
$\Delta^{(2,N)}_2[(\overline{2};1,1), (\overline{2};1,2)]$.

As can be seen from the normal ordering $\mathbf{S}[(\bar{2};1,1), (\bar{2};1,2)]
\omega_2(\overline{e}_1^{(1)}[z_{1,1}], \overline{e}_1^{(1)}[z_{1,2}]) \bh$,
there is a pole at $z_{1,t}=z_{0,s}$, $t=1,2$, $1\leq s \leq N$, in the region of the integral
$\int_{C_1-C_{q^{-1}}}\underline{dz_{1,t}}$.
Taking the residue we obtain
\begin{align}
&
\Delta^{(2,N)}_2[(\overline{2};1,1), (\overline{2};1,2)]=-N(N-1)(1-q^{1-\beta})^2 \int_{C[(1,1), (1,2)]}
\underline{dz}
\overset{\curvearrowright}{\prod_{1\leq s \leq N-2}}S_0(z_{0,s})\notag\\
&\times
 :\bLambda_{\overline{2}}(q^{\frac{1}{2}-\frac{\beta}{2}(l-2)}z_{0,N-1})
 S_0(z_{0,N-1})S_1(z_{0,N-1}):
  :\bLambda_{\overline{2}}(q^{\frac{1}{2}-\frac{\beta}{2}(l-2)}z_{0,N})
 S_0(z_{0,N})S_1(z_{0,N}):\notag\\
&\times
\overset{\curvearrowright}{\prod_{1\leq s \leq N-1}}S_1(z_{1,t})
{\bf S}_{[2,l]}\omega_2(z_{0,N-1}, z_{0,N})\bh.
\notag
\end{align}

Next, we consider
$\Delta^{(2,N)}_2[(k;0,1), (\overline{2};1,2)]$.

As can be seen from the normal ordering of
$\mathbf{S}[(k;0,1), (\bar{2};1,1)]\omega_2(\overline{e}_0^{(1)}[z_{0,1}], \overline{e}_1^{(1)}[z_{1,1}])$
$\bh$, there is a pole at $z_{1,1}=z_{0,s}$ in the region of the integral 
$\int_{C_1-C_{q^{-1}}}\underline{dz_{1,1}}$.
Taking the residue we obtain 
\begin{align}
&
\Delta^{(2,N)}_2[(k; 0,1), (\overline{2};1,1)]=(N-1)(1-q^{1-\beta})\int_{C[(1,1)]}\underline{dz}
\overset{\curvearrowright}{\prod_{1\leq s \leq N-2}}S_0(z_{0,s})\notag\\
&\times
{\mathcal A}_{k;0}(z_{0,N-1})
 :\bLambda_{\overline{2}}(q^{\frac{1}{2}-\frac{\beta}{2}(l-2)}z_{0,N})
 S_0(z_{0,N})S_1(z_{0,N}):
\notag\\
&\times
\overset{\curvearrowright}{\prod_{3 \leq t \leq N}}S_1(z_{1,t})
{\bf S}_{[2,l]}\omega_2(z_{0,N-1}, z_{0,N})\bh,\quad k=\overline{1}, \overline{2}.
\notag
\end{align}

The identity 
$:\bLambda_{\overline{2}}(q^{\frac{1}{2}-\frac{\beta}{2}(l-2)}w)S_0(w)S_1(w):=q_1^{-1}
:\bLambda_{\overline{1}}(q^{\frac{1}{2}-\frac{\beta}{2}(l-2)}w)S_0(w)S_1(q w):$
yields
\begin{align}
&(N-1)\sum_{k_1,k_2=\overline{1},\overline{2}}\nu_{k_1; 0}\nu_{k_2;0}\Delta^{(2,N)}_2[(k_1;0,1), (k_2;0,2)]\notag
\\
+&
(N-1)\nu_{\overline{2};1}^2 \Delta^{(2,N)}_2[(\overline{2};1,1), (\overline{2};1,2)]
+2N
\sum_{k=\overline{1}, \overline{2}}\nu_{k;0}\nu_{\overline{2};1}\Delta^{(2,N)}_2[(k;0,1), (\overline{2};1,1)]=0.
\notag
\end{align}

Next, we consider $\Delta^{(2,N)}_2[(k_1,k_2;0,N)]$ for $k_1,k_2=\overline{1}, \overline{2}$ such that $k_1 \succeq k_2$.

As can be seen from the normal orderings
\begin{align}
&
{\mathcal A}_{\overline{1}, \overline{1};0}(z_{0,N})S_1(z_{1,t})=q^{2(\beta-1)}
\frac{(1-q^{-\beta-1}\frac{z_{1,t}}{z_{0,N}})(1-q^{-\beta}\frac{z_{1,t}}{z_{0,N}})}{
(1-q^{-2}\frac{z_{1,t}}{z_{0,N}})(1-q^{-1}\frac{z_{1,t}}{z_{0,N}})}
:{\mathcal A}_{\overline{1}, \overline{1};0}(z_{0,N})S_1(z_{1,t}):,\notag
\\
&
{\mathcal A}_{\overline{2}, \overline{2};0}(z_{0,N})S_1(z_{1,t})=q^{2(1-\beta)}
\frac{(1-q^{\beta-2}\frac{z_{1,t}}{z_{0,N}})(1-q^{\beta-1}\frac{z_{1,t}}{z_{0,N}})}{
(1-q^{-1}\frac{z_{1,t}}{z_{0,N}})(1-\frac{z_{1,t}}{z_{0,N}})}
:{\mathcal A}_{\overline{2}, \overline{2};0}(z_{0,N})S_1(z_{1,t}):,\notag
\\
&
{\mathcal A}_{\overline{1}, \overline{2};0}(z_{0,N})S_1(z_{1,t})=
\frac{(1-q^{\beta-2}\frac{z_{1,t}}{z_{0,N}})(1-q^{\beta-1}\frac{z_{1,t}}{z_{0,N}})}{
(1-q^{-2}\frac{z_{1,t}}{z_{0,N}})(1-\frac{z_{1,t}}{z_{0,N}})}
:{\mathcal A}_{\overline{1}, \overline{2};0}(z_{0,N})S_1(z_{1,t}):,\notag
\end{align}
there are poles at $z_{0,N}=q^{-1}z_{1,t}, q^{-2}z_{1,t}$ in the region of the integral 
$\int_{C_{\delta_0}-C_{q^{-2}\delta_0}}\underline{dz_{0,N}}$.

Taking the residues we obtain 
\begin{align}
&
\Delta^{(2,N)}_2[(\overline{1}, \overline{1};0,N)]=N \frac{(1-q^\beta)(1-q^{\beta-1})}{1-q}
\int_{C[(0,N)]}\underline{dz}
\overset{\curvearrowright}{\prod_{1\leq s \leq N-1}}S_0(z_{0,s})\notag
\\
&\times
:\bLambda_{\overline{1}}(q^{\frac{1}{2}-\frac{\beta}{2}(l-2)}z_{1,1})
\bLambda_{\overline{1}}(q^{-\frac{1}{2}-\frac{\beta}{2}(l-2)}z_{1,1})S_0(q^{-1}z_{1,1})S_1(z_{1,1}):
\overset{\curvearrowright}{\prod_{2\leq t \leq N}}S_1(z_{1,t}){\bf S}_{[2,l]}\bh
\notag\\
&-qN \frac{(1-q^{\beta-1})(1-q^{\beta-2})}{1-q}
\int_{C[(0,N)]}\underline{dz}\overset{\curvearrowright}{\prod_{1\leq s \leq N-1}}S_0(z_{0,s})\notag
\\
&\times
:\bLambda_{\overline{1}}(q^{-\frac{1}{2}-\frac{\beta}{2}(l-2)}z_{1,1})
\bLambda_{\overline{1}}(q^{-\frac{3}{2}-\frac{\beta}{2}(l-2)}z_{1,1})S_0(q^{-2}z_{1,1})S_1(z_{1,1}):
\overset{\curvearrowright}{\prod_{2\leq t \leq N}}S_1(z_{1,t})
{\bf S}_{[2,l]}\bh,\notag
\end{align}

\begin{align}
&
\Delta^{(2,N)}_2[(\overline{2}, \overline{2};0,N)]=
-N q^{2(1-\beta)}
\frac{(1-q^\beta)(1-q^{\beta-1})}{1-q}
\int_{C[(0,N)]}\underline{dz}
\overset{\curvearrowright}{\prod_{1\leq s \leq N-1}}S_0(z_{0,s})\notag
\\
&\times
:\bLambda_{\overline{2}}(q^{\frac{1}{2}-\frac{\beta}{2}(l-2)}z_{1,1})
\bLambda_{\overline{2}}(q^{-\frac{1}{2}-\frac{\beta}{2}(l-2)}z_{1,1})
S_0(q^{-1}z_{1,1})S_1(z_{1,1}):
\overset{\curvearrowright}{\prod_{2\leq t \leq N}}S_1(z_{1,t}){\bf S}_{[2,l]}\bh,
\notag
\end{align}
\begin{align}
&
\Delta^{(2,N)}_2[(\overline{1}, \overline{2};0,N)]=
N q^{1-\beta}\frac{(1-q^\beta)(1-q^{\beta-1})}{1-q}
\int_{C[(0,N)]}\underline{dz}
\overset{\curvearrowright}{\prod_{1\leq s \leq N-1}}S_0(z_{0,s})\notag
\\
&\times
:\bLambda_{\overline{1}}(q^{-\frac{1}{2}-\frac{\beta}{2}(l-2)}z_{1,1})
\bLambda_{\overline{2}}(q^{-\frac{3}{2}-\frac{\beta}{2}(l-2)}z_{1,1})
S_0(q^{-2}z_{1,1})S_1(z_{1,1}):
\overset{\curvearrowright}{\prod_{1 \leq t \leq N}}S_1(z_{1,t}){\bf S}_{[2,l]}\bh.
\notag
\end{align}

Next, we consider $\Delta^{(2,N)}_2[(k,\overline{2};1,1)]$.

As can be seen from
\begin{align}
S_0(z_{0,s}){\mathcal A}_{k,\overline{2};1}(z_{1,1})=\frac{(1-q^{2-\beta}\frac{z_{1,1}}{z_{0,s}})(1-q^{1-\beta}\frac{z_{1,t}}{z_{0,s}})}{
(1-q\frac{z_{1,t}}{z_{0,s}})(1-\frac{z_{1,1}}{z_{0,s}})}:
S_0(z_{0,s}){\mathcal A}_{k,\overline{2};1}(z_{1,1}):,\quad k=\overline{1}, \overline{2},\notag
\end{align}
there is a pole at $z_{1,1}=z_{0,s}$ in the region of the integral 
$\int_{C_1-C_{q^{-1}}}\underline{dz_{1,1}}$.
Taking the residues we obtain 
\begin{align}
&
\Delta^{(2,N)}_2[(k, \overline{2};1,1)]=
N q^{3-2\beta}\frac{(1-q^{\beta-1})(1-q^{\beta-2})}{1-q}
\int_{C[(1,1)]}\underline{dz}
\overset{\curvearrowright}{\prod_{1\leq s \leq N-1}}S_0(z_{0,s})\notag
\\
&\times
:\bLambda_{k}(q^{\frac{3}{2}-\frac{\beta}{2}(l-2)}z_{0,N})
\bLambda_{\overline{2}}(q^{\frac{1}{2}-\frac{\beta}{2}(l-2)}z_{0,N})
S_0(z_{0,N})S_1(z_{0,N}):
\overset{\curvearrowright}{\prod_{1 \leq t \leq N}}S_1(z_{1,t}){\bf S}_{[2,l]}\bh,
k=\overline{1}, \overline{2}.
\notag
\end{align}
Adding the six terms
 and using
\begin{align}
 &:\bLambda_{\overline{1}}(q^{\frac{1}{2}-\frac{\beta}{2}(l-2)}z)
 \bLambda_{\overline{1}}(q^{-\frac{1}{2}-\frac{\beta}{2}(l-2)}z)
 S_0(q^{-1}z)S_1(z):
 \notag\\
 =&q_1
 :\bLambda_{\overline{1}}(q^{\frac{1}{2}-\frac{\beta}{2}(l-2)}z)
 \bLambda_{\overline{2}}(q^{-\frac{1}{2}-\frac{\beta}{2}(l-2)}z)
 S_0(q^{-1}z)S_1(q^{-1}z):,\notag
\\
&:\bLambda_{\overline{1}}(q^{-\frac{1}{2}-\frac{\beta}{2}(l-2)}z)
 \bLambda_{\overline{1}}(q^{-\frac{3}{2}-\frac{\beta}{2}(l-2)}z)
 S_0(q^{-2}z)S_1(z):
 \notag\\
 =&q_1^2
 :\bLambda_{\overline{2}}(q^{-\frac{1}{2}-\frac{\beta}{2}(l-2)}z)
 \bLambda_{\overline{2}}(q^{-\frac{3}{2}-\frac{\beta}{2}(l-2)}z)
 S_0(q^{-2}z)S_1(q^{-2}z):,
 \notag
\\
&q_1
:\bLambda_{\overline{2}}(q^{-\frac{1}{2}-\frac{\beta}{2}(l-2)}z)
 \bLambda_{\overline{2}}(q^{-\frac{3}{2}-\frac{\beta}{2}(l-2)}z)
 S_0(q^{-2}z)S_1(q^{-1}z):
 \notag\\
 =&
 :\bLambda_{\overline{1}}(q^{-\frac{1}{2}-\frac{\beta}{2}(l-2)}z)
 \bLambda_{\overline{2}}(q^{-\frac{3}{2}-\frac{\beta}{2}(l-2)}z)
 S_0(q^{-2}z)S_1(z):,
\notag 
\end{align} 
we find 
$\sum_{k=\overline{1}, \overline{2}}
\nu_{k, \overline{2};1}\Delta^{(2,N)}_2[(k,\overline{2};1,1)]
+
\sum\limits_{\substack{\scriptscriptstyle 
k_1,k_2=\overline{1}, \overline{2}\\ \scriptscriptstyle 
k_1 \succeq k_2}}
\nu_{k_1, k_2 ;0}\Delta^{(2,N)}_2[(k_1, k_2;0,1)]=0$.

 We obtain Lemma \ref{lem:I2-GN:4}.
\end{proof}

\newpage

The relation $\mathrm{LHS}_{2,N}(0,0,2)=\mathrm{RHS}_{2,N}(0,0,2)$
in Proposition \ref{prop:duality1} is a direct consequence of Lemma \ref{lem:I2-GN:6} below.
\begin{lem}
\label{lem:I2-GN:6}
The following holds for $N \geq 2$.
\begin{align}
&
(N-1)\sum_{k_1,k_2=l-1,\overline{l}}
\Delta^{(2,N)}_2\left[(k_1;l-1,1), (k_2;l-1,2)\right]
\notag\\
+&2N \sum\limits_{\substack{\scriptscriptstyle 
k_1=l-1, \overline{l} \\ \scriptscriptstyle 
k_2=l-1,l}}
\Delta^{(2,N)}_2\left[(k_1;l-1,1), (k_2;l,1)\right]
+
(N-1)\sum_{k_1,k_2=l-1, l}
\Delta^{(2,N)}_2\left[(k_1;l,1), (k_2;l,2)\right]=0.\notag
\end{align}
\end{lem}

\begin{proof}
We rewrite the integral of $\Delta^{(2,N)}_2$ as follows.
\begin{align}
&\int_{C_{q^{-1}\delta_i}}\underline{dz_{i,s}}
\int_{C_{q^{-1}\delta_j}}\underline{dz_{j,t}}-
\int_{C_{\delta_i}}\underline{dz_{i,s}}
\int_{C_{\delta_j}}\underline{dz_{j,t}}\notag\\
&=\int_{C_{q^{-1}\delta_i}-C_{\delta_i}}\underline{dz_{i,s}}
\int_{C_{q^{-1}\delta_j}-C_{\delta_j}}\underline{dz_{j,t}}
+\int_{C_{q^{-1}\delta_i}-C_{\delta_i}}\underline{dz_{i,s}}
\int_{C_{\delta_j}}\underline{dz_{j,t}}
+
\int_{C_{\delta_i}}\underline{dz_{i,s}}
\int_{C_{q^{-1}\delta_j}-C_{\delta_j}}\underline{dz_{j,t}}.
\notag
\end{align}
Taking residues in $\int_{C_{q^{-1}\delta_i}-C_{\delta_i}}\underline{dz_{i,s}}$ we obtain 
\begin{align}
&
\Delta^{(2,N)}_2[(k_1;l-1,N-1),(k_2;l-1,N)]\notag
\\
=&\delta_{k_1,l-1}\delta_{k_2,l-1}N(N-1)
q^{2(1-\beta)}(1-q^{\beta-1})^2
\notag
\\
\times&
\int_{C[(l-1,N-1), (l-1,N)]}\underline{dz}
~{\bf S}_{[0,l-2]}
\overset{\curvearrowright}{\prod_{1\leq s \leq N-2}}
S_{l-1}(z_{l-1, s})
: \bLambda_{l-1}(q^{\frac{1}{2}}z_{l,1})S_{l-1}(z_{l,1})S_l(z_{l,1}):\notag\\
\times&
: \bLambda_{l-1}(q^{\frac{1}{2}}z_{l,2})S_{l-1}(z_{l,2})S_l(z_{l,2}):
\overset{\curvearrowright}{\prod_{2 \leq t \leq N}}
S_{l}(z_{l, t})\omega_2(z_{l,1}, z_{l,2})\bh
\notag
\\
+&\delta_{k_2,l-1}N
q^{1-\beta}(1-q^{\beta-1})
\int_{C[(l-1,N)]}\underline{dz}
~{\bf S}_{[0,l-2]}
\overset{\curvearrowright}{\prod_{1\leq s \leq N-2}}
S_{l-1}(z_{l-1, s})
{\mathcal A}_{k_1;l-1}(z_{l-1,N-1})
\notag\\
\times&
: \bLambda_{l-1}(q^{\frac{1}{2}}z_{l,1})S_{l-1}(z_{l,1})S_l(z_{l,1}):
\overset{\curvearrowright}{\prod_{2 \leq t \leq N}}
S_{l}(z_{l, t})\omega_2(z_{l-1, N-1}, z_{l,1})\bh
\notag
\\
+&\delta_{k_1,l-1}N
q^{1-\beta}(1-q^{\beta-1})
\int_{C[(l-1,N)]}\underline{dz}
~{\bf S}_{[0,l-2]}
\overset{\curvearrowright}{\prod_{1\leq s \leq N-2}}
S_{l-1}(z_{l-1, s})
{\mathcal A}_{k_2;l-1}(z_{l-1,N-1})\notag\\
\times&
: \bLambda_{l-1}(q^{\frac{1}{2}}z_{l,1})S_{l-1}(z_{l,1})S_l(z_{l,1}):
\overset{\curvearrowright}{\prod_{2 \leq t \leq N}}
S_{l}(z_{l, t})\omega_2(z_{l-1, N-1}, z_{l,1})\bh,
\notag
\end{align}
\begin{align}
&
\Delta^{(2,N)}_2[(k_1;l,1),(k_2;l,2)]
\notag\\
=&\delta_{k_1,l}\delta_{k_2,l}N(N-1)
(1-q^{\beta-1})^2
\notag
\\
\times&
\int_{C[(l,1), (l,2)]}\underline{dz}
~{\bf S}_{[0,l-2]}
\overset{\curvearrowright}{\prod_{1\leq s \leq N-2}}
S_{l-1}(z_{l-1, s})
: \bLambda_{l}(q^{-\frac{1}{2}}z_{l-1,N-1})
S_{l-1}(z_{l-1,N-1})S_l(q^{-1}z_{l-1,N-1}):\notag\\
\times&
: \bLambda_{l}(q^{-\frac{1}{2}}z_{l-1,N})
S_{l-1}(z_{l-1,N})S_l(q^{-1}z_{l-1,N}):
\overset{\curvearrowright}{\prod_{3 \leq t \leq N}}
S_{l}(z_{l, t})\omega_2(z_{l-1,N-1}, z_{l-1,N})\bh
\notag
\\
-&\delta_{k_1,l}N
(1-q^{\beta-1})
\int_{C[(l,1)]}\underline{dz}
~{\bf S}_{[0,l-2]}
\overset{\curvearrowright}{\prod_{1\leq s \leq N-1}}
S_{l-1}(z_{l-1, s})
: \bLambda_{l}(q^{-\frac{1}{2}}z_{l-1,N})
S_{l-1}(z_{l-1,N})S_l(q^{-1}z_{l-1,N}):\notag\\
\times&
{\mathcal A}_{k_2;l}(z_{l,2})
\overset{\curvearrowright}{\prod_{3 \leq t \leq N}}
S_{l}(z_{l, t})\omega_2(z_{l-1, N}, q^{-1} z_{l,2})\bh
\notag
\\
-&\delta_{k_2,l}N
(1-q^{\beta-1})
\int_{C[(l,1)]}\underline{dz}
~{\bf S}_{[0,l-2]}
\overset{\curvearrowright}{\prod_{1\leq s \leq N-1}}
S_{l-1}(z_{l-1, s})
: \bLambda_{l}(q^{-\frac{1}{2}}z_{l-1,N})
S_{l-1}(z_{l-1,N})S_l(q^{-1}z_{l-1,N}):\notag\\
\times&
{\mathcal A}_{k_1;l}(z_{l,2})
\overset{\curvearrowright}{\prod_{3 \leq t \leq N}}
S_{l}(z_{l, t})\omega_2(z_{l-1, N}, q^{-1} z_{l,2})\bh,
\notag
\end{align}
\begin{align}
&\Delta^{(2,N)}_2[(k_1;l-1,N),(k_2;l,1)]\notag
\\
=&
-\delta_{k_1,l-1}\delta_{k_2,l}(N-1)^2
q^{1-\beta}(1-q^{\beta-1})^2
\notag
\\
\times&
\int_{C[(l-1,N), (l,1)]}\underline{dz}
~{\bf S}_{[0,l-2]}
\overset{\curvearrowright}{\prod_{1\leq s \leq N-2}}
S_{l-1}(z_{l-1, s})
: \bLambda_{l-1}(q^{\frac{1}{2}}z_{l,2})S_{l-1}(z_{l,2})S_l(z_{l,2}):\notag\\
\times&
: \bLambda_{l}(q^{-\frac{1}{2}}z_{l-1,N-1})
S_{l-1}(z_{l-1,N-1})S_l(q^{-1}z_{l-1,N-1}):
\overset{\curvearrowright}{\prod_{3 \leq t \leq N}}
S_{l}(z_{l, t})\omega_2(z_{l,2}, q^{-1}z_{l-1,N-1})\bh
\notag
\\
+&
\delta_{k_1,l-1}(N-1)
q^{1-\beta}(1-q^{\beta-1})
\int_{C[(l-1,N)]}\underline{dz}
~{\bf S}_{[0,l-2]}
\overset{\curvearrowright}{\prod_{1\leq s \leq N-1}}
S_{l-1}(z_{l-1, s})
\notag\\
\times&
: \bLambda_{l-1}(q^{\frac{1}{2}}z_{l,1})S_{l-1}(z_{l,1})S_l(z_{l,1}):
{\mathcal A}_{k_2;l}(z_{l,2})
\overset{\curvearrowright}{\prod_{3 \leq t \leq N}}
S_{l}(z_{l, t})\omega_2(z_{l,1}, z_{l,2})\bh
\notag
\\
-&\delta_{k_2,l}(N-1)
(1-q^{\beta-1})
\int_{C[(l,1)]}\underline{dz}
~{\bf S}_{[0,l-2]}
\overset{\curvearrowright}{\prod_{1\leq s \leq N-2}}
S_{l-1}(z_{l-1, s})
{\mathcal A}_{k_1; l-1}(z_{l-1,N-1})\notag\\
\times&
: \bLambda_{l}(q^{-\frac{1}{2}}z_{l-1,N})
S_{l-1}(z_{l-1,N})S_l(q^{-1}z_{l-1,N}):
\overset{\curvearrowright}{\prod_{2 \leq t \leq N}}
S_{l}(z_{l, t})\omega_2(z_{l-1, N-1}, q^{-1}z_{l-1,N})\bh.
\notag
\end{align}
Adding the above terms and using
$$\bLambda_{l}(q^{-\frac{1}{2}}z)S_{l-1}(z)S_l(q^{-1}z)=
q_1:\bLambda_{l-1}(q^{-\frac{1}{2}}z)S_{l-1}(q^{-1}z)S_l(q^{-1}z):
$$
we arrive at 
Lemma \ref{lem:I2-GN:6}.
\end{proof}

The relation $\mathrm{LHS}_{2,N}(1,0,1)=\mathrm{RHS}_{2,N}(1,0,1)$ in Proposition \ref{prop:duality1}
is a direct consequence of Lemma \ref{lem:I2-GN:7} below.
\begin{lem}
\label{lem:I2-GN:7}
The following hold for $N \geq 2$.
\begin{align}
&\sum_{k=\overline{l-1},l}\Delta_2^{(2,N)}[(\bar{2};1,1), (k;l-1,1)]
+\sum_{k=\overline{l-1}, \bar{l}}\Delta_2^{(2,N)}[(\bar{2};1,1), (k;l,1)]\notag
\\
&+\sum_{k_1=\bar{1}, \bar{2}}
\sum_{k_2=\overline{l-1},l}\Delta_2^{(2,N)}[(k_1;0,1), (k_2;l-1,1)]
+
\sum_{k_1=\bar{1}, \bar{2}}
\sum_{k_2=\overline{l-1},\bar{l}}\Delta_2^{(2,N)}[(k_1;0,1), (k_2;l,1)]
=0.
\notag
\end{align}
\begin{align}
&\Delta^{(2,N)}_2[(2;1,1), (k;l-1,1)]=0, \quad
k=l-1, \bar{l},\notag
\\
&\Delta^{(2,N)}_2[(2;1,1), (k;l,1)]=0, \quad
k=l-1,l.\notag
\end{align}
\end{lem}
\begin{proof}

First, we consider
$\Delta^{(2,N)}_2[(2;1,1), (k;l-1,1)]$
and
$\Delta^{(2,N)}_2[(2;1,1), (k;l,1)]$.

As can be seen from the normal orderings of
$\mathbf{S}[(2;1,1), (k;l-1,1)] \omega_2(e_1^{(1)}[z_{1,1}], d_{l-1}^{(1)}[z_{l-1,1}])\bh$ and
$\mathbf{S}[(2;1,1), (k;l,1)] \omega_2(e_1^{(1)}[z_{1,1}], d_l^{(1)}[z_{l,1}])$
$\times \bh$,
there is no pole of $z_{1,1}$ in the region of the integral 
$\int_{C_1-C_q}\underline{dz_{1,1}}$.
Hence,
$\Delta^{(2,N)}_2[(2;1,1), (k;l-1,1)]=0$
and
$\Delta^{(2,N)}_2[(2;1,1), (k;l,1)]=0$.

Next, we consider $\Delta_2^{(2,N)}[(\bar{2};1,1), (k;l-1,1)]$.

Here we find it convenient to use $\Delta_2^{(2,N)}[(\bar{2};1,1), (k; l-1,N)]$ in place of 
$\Delta_2^{(2,N)}[(\bar{2};1,1), (k;l-1,1)]$.
As can be seen from the normal  ordering of
$\mathbf{S}[(\bar{2};1,1),(k;l-1,N)]
\omega_2(\overline{e}_1^{(1)}[z_{1,1}], \overline{d}_{l-1}^{(1)}[z_{l-1,N}])\bh$,
there is a pole at $z_{1,1}=z_{0,s}$
in the region of the integral 
$\int_{C_1-C_{q^{-1}}}\underline{dz_{1,1}}$.
Taking the residue we obtain 
\begin{align}
&\Delta_2^{(2,N)}[(\bar{2};1,1), (k;l-1,N)]
\notag\\
=&N(1-q^{1-\beta})\int_{C[(1,1)]}\underline{dz}
\overset{\curvearrowright}{\prod_{1\leq a \leq N-1}}S_0(z_{0,a})
:\bLambda_{\bar{2}}(q^{\frac{1}{2}-\frac{\beta}{2}(l-2)}z_{0,N})S_0(z_{0,N})S_1(z_{0,N}):\notag
\\
&\times
\overset{\curvearrowright}{\prod_{2\leq a \leq N}}S_{1}(z_{1,a})
\mathbf{S}_{[2,l-2]}
\overset{\curvearrowright}{\prod_{1\leq a \leq N-1}}S_{l-1}(z_{l-1,a})
\mathcal{A}_{k;l-1}[z_{l-1,N}]
\mathbf{S}_{[l,l]}\notag
\\
&\times \omega_2(\overline{e}_1^{(1)}[z_{0,N}], \overline{d}_{l-1}^{(1)}[z_{l-1,N}])
\bh,\quad k=\overline{l-1},l.
\notag
\end{align}

We consider $\Delta_2^{(2,N)}[(\bar{2};1,1), (k; l,1)]$.
In the same way as $\Delta_2^{(2,N)}[(\bar{2};1,1), (k;l-1,1)]$, we obtain
\begin{align}
&\Delta_2^{(2,N)}[(\bar{2};1,1), (k;l,1)]\notag\\
=&N(1-q^{1-\beta})\int_{C[(1,1)]}\underline{dz}
\overset{\curvearrowright}{\prod_{1\leq a \leq N-1}}S_0(z_{0,a})
:\bLambda_{\bar{2}}(q^{\frac{1}{2}-\frac{\beta}{2}(l-2)}z_{0,N})S_0(z_{0,N})S_1(z_{0,N}):\notag
\\
&\times
\overset{\curvearrowright}{\prod_{2\leq a \leq N}}S_{1}(z_{1,a})
\mathbf{S}_{[2,l-1]}
\mathcal{A}_{k;l}[z_{l,1}]
\overset{\curvearrowright}{\prod_{2\leq a \leq N}}S_{l}(z_{l,a})
\notag\\
&\times \omega_2(\overline{e}_1^{(1)}[z_{0,N}], \overline{d}_{l}^{(1)}[z_{l,1}])
\bh,\quad k=\overline{l-1}, \bar{l}.\notag
\end{align}

We consider $\Delta_2^{(2,N)}[(k_1;0,N), (k_2;l-1,N)]$ and $\Delta_2^{(2,N)}[(k_1;0,N), (k_2;l,1)]$.

We rewrite the integral as follows.
\begin{align}
&
\int_{C_{\delta_0}}\underline{dz_{0,N}}
\int_{C_{\delta_i}}\underline{dz_{i,1}}
-
\int_{C_{q^{-1} \delta_0}}\underline{dz_{0,N}}
\int_{C_{q \delta_i}}\underline{dz_{i,1}}
\notag
\\
=&
\int_{C_{\delta_0}-C_{q^{-1}\delta_0}}\underline{dz_{0,N}}\int_{C_{\delta_i}}\underline{dz_{i,1}}
+
\int_{C_{\delta_0}}\underline{dz_{0,N}}
\int_{C_{\delta_i}-C_{q \delta_i}}\underline{dz_{i,1}}
-
\int_{C_{\delta_0}-C_{q^{-1}\delta_0}}
\underline{dz_{0,N}}
\int_{C_{\delta_i}-C_{q\delta_i}}\underline{dz_{i,1}}.\notag
\end{align}

First, we consider $\Delta_2^{(2,N)}[(k_1;0,N), (k_2;l-1,N)]$.

As can be seen from the normal ordering of $\mathbf{S}[(k_1;0,N), (k_2;l-1,N)]
\omega_2(\overline{e}_0^{(1)}[z_{0,N}], \overline{d}_{l-1}^{(1)}[z_{l-1,N}])$
$\times\bh$,
there is a pole at $z_{l-1,N}=qz_{l,a}$ in the case of $k_2=l$ and no pole in the case of $k_2=\overline{l-1}$
in the region of the integral 
$\int_{C_1-C_{q}}dz_{l-1,N}$.
Taking the residues we obtain 
\begin{align}
&\Delta_2^{(2,N)}[(k_1;0,N), (k_2;l-1,N)]\notag\\
=&N(q^{1-\beta}-1)\delta_{k_1,\bar{2}}
\int_{C[(0,N)]}\underline{dz}
\overset{\curvearrowright}{\prod_{1\leq a \leq N-1}}S_0(z_{0,a})\notag\\
&\times
:\bLambda_{\bar{2}}(q^{\frac{1}{2}-\frac{\beta}{2}(l-2)}z_{1,1})S_0(z_{1,1})S_1(z_{1,1}):
\overset{\curvearrowright}{\prod_{2\leq a \leq N}}S_{1}(z_{1,a})
\mathbf{S}_{[2,l-2]}
\overset{\curvearrowright}{\prod_{1\leq a \leq N-1}}S_{l-1}(z_{l-1,a})\notag\\
&\times
\mathcal{A}_{k_2;l-1}[z_{l-1,N}]
\mathbf{S}_{[l,l]}\notag
\omega_2(\overline{e}_0^{(1)}[z_{1,1}], \overline{d}_{l-1}^{(1)}[z_{l-1,N}])
\bh
\notag
\\
&+N(1-q^{1-\beta})\delta_{k_2,l}
\int_{C[(l-1,N)]}\underline{dz}
\overset{\curvearrowright}{\prod_{1\leq a \leq N-1}}S_0(z_{0,a})
\mathcal{A}_{k_1;0}[z_{0,N}]\mathbf{S}_{[1,l-2]}\notag
\\
&\times
\overset{\curvearrowright}{\prod_{1\leq a \leq N-1}}S_{l-1}(z_{l-1,a})
:\bLambda_{l}(q^{\frac{1}{2}}z_{l,1})S_{l-1}(q z_{l,1})S_l(z_{l,1}):
\overset{\curvearrowright}{\prod_{2\leq a \leq N}}S_{l}(z_{l,a})\notag
\\
&\times \omega_2(\overline{e}_0^{(1)}[z_{0,N}], \overline{d}_{l-1}^{(1)}[q z_{l,1}])
\bh
\notag
\\
&+N^2 (q^{-\beta+1}-1)(1-q^{1-\beta})\delta_{k_1, \bar{2}}\delta_{k_2,l}
\int_{C[(0,N),(l-1,N)]}\underline{dz}
\overset{\curvearrowright}{\prod_{1\leq a \leq N-1}}S_0(z_{0,a})\notag
\\
&\times
:\bLambda_{\bar{2}}(q^{\frac{1}{2}-\frac{\beta}{2}(l-2)}z_{1,1})S_0(z_{1,1})S_1(z_{1,1}):
\overset{\curvearrowright}{\prod_{2\leq a \leq N}}S_1(z_{1,a})
\mathbf{S}_{[2,l-2]}\notag
\\
&\times
\overset{\curvearrowright}{\prod_{1\leq a \leq N-1}}S_{l-1}(z_{l-1,a})
:\bLambda_{l}(q^{\frac{1}{2}}z_{l,1})S_{l-1}(q z_{l,1})S_l(z_{l,1}):
\overset{\curvearrowright}{\prod_{2\leq a \leq N}}S_{l}(z_{l,a})\notag
\\
&\times \omega_2(\overline{e}_0^{(1)}[z_{1,1}], \overline{d}_{l-1}^{(1)}[q z_{l,1}])
\bh.
\notag
\end{align}

Next, we consider $\Delta_2^{(2,N)}[(k_1;0,N), (k_2;l,1)]$.

As can be seen from the normal ordering of $\mathbf{S}[(k_1;0,N), (k_2;l,1)]
\omega_2(\overline{e}_0^{(1)}[z_{0,N}], \overline{d}_{l}^{(1)}[z_{l,1}])\bh$,
there is a pole at $z_{l,1}=z_{l-1,a}$ in the case of $k_2=\overline{l-1}$ 
and no pole in the case of $k_2=\overline{l}$
in the region of the integral 
$\int_{C_{\delta_l}-C_{q \delta_l}}dz_{l,1}$.
Taking the residues we obtain 
\begin{align}
&
\Delta_2^{(2,N)}[(k_1;0,N), (k_2;l,1)]
\notag\\
&=N(q^{1-\beta}-1)\delta_{k_1,\bar{2}}
\int_{C[(0,N)]}\underline{dz}
\overset{\curvearrowright}{\prod_{1\leq a \leq N-1}}S_0(z_{0,a})
\notag\\
&\times
:\bLambda_{\bar{2}}(q^{\frac{1}{2}-\frac{\beta}{2}(l-2)}z_{1,1})S_0(z_{1,1})S_1(z_{1,1}):
\overset{\curvearrowright}{\prod_{2\leq a \leq N}}S_{1}(z_{1,a})
\mathbf{S}_{[2,l-1]}\notag\\
&\times
\mathcal{A}_{k_2;l}[z_{l,1}]
\overset{\curvearrowright}{\prod_{2\leq a \leq N}}S_{l}(z_{l,a})
\omega_2(\overline{e}_0^{(1)}[z_{1,1}], \overline{d}_{l}^{(1)}[z_{l,1}])
\bh\notag
\\
&+N(1-q^{\beta-1})\delta_{k_2, \overline{l-1}}
\int_{C[(l,1)]}\underline{dz}
\overset{\curvearrowright}{\prod_{1\leq a \leq N-1}}S_0(z_{0,a})
\mathcal{A}_{k_1;0}[z_{0,N}]\mathbf{S}_{[1,l-2]}\notag
\\
&\times
\overset{\curvearrowright}{\prod_{1\leq a \leq N-1}}S_{l-1}(z_{l-1,a})
:\bLambda_{\overline{l-1}}(q^{-\frac{1}{2}}z_{l-1,N})S_{l-1}(z_{l-1,N})S_l(z_{l-1,N}):
\overset{\curvearrowright}{\prod_{2\leq a \leq N}}S_{l}(z_{l,a})\notag
\\
&\times \omega_2(\overline{e}_0^{(1)}[z_{0,N}], \overline{d}_{l}^{(1)}[z_{l-1,N}])
\bh
\notag
\\
&+N^2(q^{-\beta+1}-1)(1-q^{\beta-1})\delta_{k_1, \bar{2}}\delta_{k_2,\overline{l-1}}
\int_{C[(0,N),(l,1)]}\underline{dz}
\overset{\curvearrowright}{\prod_{1\leq a \leq N-1}}S_0(z_{0,a})\notag
\\
&\times
:\bLambda_{\bar{2}}(q^{\frac{1}{2}-\frac{\beta}{2}(l-2)}z_{1,1})S_0(z_{1,1})S_1(z_{1,1}):
\overset{\curvearrowright}{\prod_{2\leq a \leq N}}S_1(z_{1,a})
\mathbf{S}_{[2,l-2]}\notag
\\
&\times
\overset{\curvearrowright}{\prod_{1\leq a \leq N-1}}S_{l-1}(z_{l-1,a})
:\bLambda_{\overline{l-1}}(q^{-\frac{1}{2}}z_{l-1,N})S_{l-1}(z_{l-1,N})S_l(z_{l-1,N}):
\overset{\curvearrowright}{\prod_{2\leq a \leq N}}S_{l}(z_{l,a})\notag
\\
&\times \omega_2(\overline{e}_0^{(1)}[z_{1,1}], \overline{d}_{l}^{(1)}[z_{l-1,N}])
\bh.
\notag
\end{align}
Adding the above terms and using
\begin{align}
:\bLambda_{\overline{l-1}}(q^{-\frac{1}{2}}z)S_{l-1}(z)S_l(z):=q_1
:\bLambda_l(q^{-\frac{1}{2}}z)S_{l-1}(z)S_l(q^{-1}z):\notag
\end{align}
we arrive at 
Lemma \ref{lem:I2-GN:7}.
\end{proof}

\subsection{Proof of $[{\mathbf I}_M, {\mathbf G}_N(\vartheta)]=0$}
\label{section:IM-GN} 

In this section, we will show Proposition \ref{prop:duality1} in the case of $N \geq M \geq 1$.
In what follows, we assume $\beta>\mathrm{Max}(4, 2M-1)$.

\subsubsection{The Left Hand Side}

We start from ${\bf I}_M \cdot {\bf G}_N(\vartheta)$.
We set the radii of integration circles $R_1, R_2, \ldots, R_M$ as $|w_i|=R_i, 1\leq i \leq M$.
We move the contours from $R_1=R_2=\cdots=R_M=R \gg 1$ to $R_1=R_2=\cdots=R_M=1$ imposing the conditions:
\begin{align}
1<R_M<\cdots<R_2<R_1<R
\quad
\mathrm{and}
\quad
|q|R_i<R_j,\quad 1\leq i<j \leq M,\notag
\end{align}
in the case of $1<R_1, R_2, \cdots, R_M<R$.
In the process of this deformation of contours, 
we take the residues at $w_i=\overline{\eta}_{k_{m},\ldots, k_{1};j}[z_{j,s}]$ 
and $w_i=\overline{\eta}_{k_{n},\ldots,k_1;j}^{ k_{m+n},\ldots,k_{n+1}}[z_{j,s}]$
given as follows.

For $m=1$,  $\overline{\eta}_{k_1;j}[z]$ is given in (\ref{def:leta1}).
For $2\leq m \leq M$, we set
\begin{align}
\label{def:leta4}
\overline{\eta}_{k_m,\ldots,k_1;j}[z]=\left\{\begin{array}{cl}
\overline{e}_0^{(m)}[z],& j=0,\quad \bar{1}\succeq k_m \succeq \cdots \succeq k_1 \succeq \bar{2},\\
\overline{d}_j^{(m)}[z],& 1\leq j \leq l-2,\quad \bar{1} \succeq k_m \succeq \cdots \succeq k_2\succeq \overline{j+1}, k_1=\overline{j},\\
\overline{e}_j^{(m)}[z],& 1\leq j \leq l-2,\quad \bar{1} \succeq k_m \succeq \cdots \succeq k_2\succeq \overline{j+1}, k_1=\overline{j+1}.
\end{array}
\right.
\end{align}
For $1\leq m,n \leq M$ such that $2\leq m+n \leq M$, we set
\begin{align}
\label{def:leta5}
\overline{\eta}_{k_{n},\ldots, k_{1};j}^{k_{m+n},\ldots, k_{n+1}}[z]=\left\{\begin{array}{cl}
\overline{D}_j^{(m)}[z],& 1\leq j \leq l-3,
\begin{array}{c}
\bar{1} \succeq k_{m+n} \succeq \cdots \succeq k_{n+1}\succeq \mathrm{Max}(\overline{j}, k_n),\\
\bar{1} \succeq k_n \succeq \cdots \succeq k_2\succeq \overline{j+1}, k_1=\overline{j},
\end{array}
\\
\\
\overline{E}_j^{(m)}[z],& 1\leq j \leq l-3,
\begin{array}{c}
\bar{1} \succeq k_{m+n} \succeq \cdots \succeq k_{n+1}\succeq \mathrm{Max}(\overline{j}, k_n),\\
\bar{1} \succeq k_n \succeq \cdots \succeq k_2\succeq \overline{j+1}, k_1=\overline{j+1}.
\end{array}
\end{array}
\right.
\end{align}

\vspace{2mm}

For $(k_m,\ldots,k_1;j)$ in (\ref{def:leta4}), we set
\begin{align}
&
\nu_{k_m,\ldots,k_1;j}
\mathcal{A}_{k_m,\ldots,k_1;j}[z]
=\overset{\curvearrowleft}{\prod_{1\leq a \leq m}}
\mathop{\mathrm{res}}_{w_a=\overline{\eta}_{k_a,\ldots, k_2, k_1;j}[z]}\frac{dw_a}{w_a}
\overset{\curvearrowleft}{\prod_{1\leq a \leq m}}\bLambda_{k_a}(w_a)S_j(z)
\omega_{m}(w_m,\ldots,w_2,w_1),
\notag
\\
&
\mathcal{A}_{k_m,\ldots,k_1;j}[z]=:
\prod_{1\leq a \leq m}
\bLambda_{k_a}(\overline{\eta}_{k_a,\ldots,k_1;j}[z])S_j(z):.\notag
\end{align}
For $\left(\begin{array}{c}
k_{m+n},\ldots, k_{n+1}\\
k_n, \ldots, k_1; j\end{array}
\right)$ in (\ref{def:leta5}), 
we set
\begin{align}
&\nu_{k_n,\ldots,k_1;j}^{k_{m+n},\ldots,k_{n+1}}
\mathcal{A}_{k_n,\ldots,k_1;j}^{k_{m+n},\ldots,k_{n+1}}[z]\notag\\
=&\overset{\curvearrowleft}{\prod_{n+1\leq a \leq m+n}}
\mathop{\mathrm{res}}_{w_a=\overline{\eta}_{k_n,\ldots, k_1;j}^{k_a,
\ldots,k_{n+1}}[z]}
\frac{dw_a}{w_a}
\overset{\curvearrowleft}{\prod_{1\leq a \leq n}}
\mathop{\mathrm{res}}_{w_a=\overline{\eta}_{k_a,\ldots, k_1;j}[z]}\frac{dw_a}{w_a}
\overset{\curvearrowleft}{\prod_{1\leq a \leq m+n}}\bLambda_{k_a}(w_a)S_j(z)
\omega_{m+n}(w_{m+n},\ldots, w_2, w_1),
\notag
\\
&\mathcal{A}_{k_n,\ldots,k_1;j}^{k_{m+n},\ldots,k_{n+1}}[z]=
:\prod_{n+1\leq a \leq m+n}
\bLambda_{k_a}(\overline{\eta}_{k_n,\ldots,k_1;j}^{k_a,\ldots,k_{n+1}}[z])
\prod_{1\leq a \leq n}
\bLambda_{k_a}(\overline{\eta}_{k_a,\ldots,k_1;j}[z])S_j(z):.\notag
\end{align}

To make formulae of $\mathrm{LHS}_{M,N}(a,b,c)$ and $\mathrm{RHS}_{M,N}(a,b,c)$ easier to read, we prepare notation.

We extend the index set $I=\{1,\ldots,l,\overline{l},\ldots,\overline{1}\}$ 
by adjoining its copy 
$\boxed{I}=\{\boxed{1},\ldots, \boxed{l}, \boxed{\overline{l}},\ldots, \boxed{\overline{1}}\}$.
The set $I \cup \boxed{I}$ is ordered as
$\boxed{1} \prec \boxed{2} \prec \ldots \prec \boxed{l} \prec 1 \prec \ldots  \prec l \prec \bar{l} \prec \ldots \prec \bar{1} \prec \boxed{\bar{l}} \prec \ldots \prec \boxed{\bar{2}} \prec \boxed{\bar{1}}.$
We define $\boxed{\boxed{k}}=k$ for 
$k \in I \cup \boxed{I}$.

A partition $\lambda=(\lambda_1,\lambda_2,\ldots,\lambda_m)$
is a sequence of non-negative integers in 
non-increasing order: $\lambda_1 \geq \lambda_2\geq \ldots \geq \lambda_m$.
The non-zero $\lambda_i$ are called parts of $\lambda$.
The length of $\lambda$, denoted by $l(\lambda)$, 
is the number of its parts.
The weight of $\lambda$, denoted by $|\lambda|$, 
is the sum of its parts.
We use the notation which indicates the number of times each integer occurs as a part:
$\lambda=(1^{m_1} 2^{m_2} \ldots i^{m_i}\ldots)$
means that exactly $m_i$ of the parts of $\lambda$ are equal to $i$.
The number $m_i=m_i(\lambda)=|\{\left.j \right|\lambda_j=i \}|$ is called the multiplicity of $i$ in $\lambda$.
The Young diagram of a partition $\lambda$ is the set of points 
$(i,j) \in {\mathbf Z}^2$ such that $1\leq j \leq \lambda_i$.
The first coordinate $i$ (the row index) increases as one goes downwards, and 
the second coordinate $j$ (the column index) increases as one goes from left to right.
We often refer to a node by a `box'.
We do not make distinction between 
a partition $\lambda$  and the Young diagram it represents.
For $\lambda=(1^{m_1} 2^{m_2} \cdots r^{m_r})$ and $n \geq l(\lambda)$, we define
\begin{align}
\com{n}{\lambda}=
\left(
\begin{array}{c}
~~~n\\
m_1, m_2, \cdots, m_r
\end{array}\right)
.\notag
\end{align}

A Young tableau is a pair $(\lambda,k)$ consisting of a 
Young diagram  $\lambda$ 
and a filling $k=(k_{i,j})$ in the alphabet $I_l \cup \boxed{I_l}$ 
assigned to each box $(i,j)$ of $\lambda$. 
Let $T$ be the set of all Young tableaux.
For brevity of notation we shall write $(\lambda,k)\in T$ simply as
 $\lambda\in T$, suppressing the filling $k$.

We define subsets $\mathcal{\overline{D}}_j, \mathcal{\overline{E}}_j, \boxed{\mathcal{\overline{D}}_j}, \boxed{\mathcal{\overline{E}}_j}
\subset T$ as follows.
\begin{align}
\overline{\mathcal D}_j&=
\left\{\lambda \in T\left|
\bar{1} \succeq k_{i,\lambda_i}\succeq \cdots \succeq k_{i,2} \succeq \overline{j+1}, k_{i,1}=\overline{j}
\textrm{~for~} 1\leq i \leq l(\lambda)
\right.
\right\},~1\leq j \leq l-2,
\notag\\
\overline{\mathcal D}_{l-1}&=\left\{\lambda \in T\left|
k_{i,1}=\overline{l-1}, l \textrm{~and~} \lambda_i=1 \textrm{~for~} 1\leq i \leq l(\lambda)
\right.\right\},
\notag\\
\overline{\mathcal D}_{l}&=\left\{\lambda \in T \left|
k_{i,1}=\overline{l-1}, \overline{l} \textrm{~and~} \lambda_i=1 \textrm{~for~} 1\leq i \leq l(\lambda)
\right.\right\},
\notag\\
\overline{\mathcal E}_0&=\left\{\lambda \in T\left|
\bar{1}\succeq k_{i,\lambda_i}\succeq \cdots \succeq k_{i,1}=\overline{1}, \overline{2} \textrm{~for~} 1\leq i \leq l(\lambda)
\right.\right\},\notag
\\
\overline{\mathcal E}_j&=
\left\{\lambda \in T\left|
\bar{1} \succeq k_{i,\lambda_i} \succeq \cdots \succeq k_{i,2} \succeq \overline{j+1}, k_{i,1}=\overline{j+1} \textrm{~for~} 1\leq i \leq l(\lambda)
\right.\right\},\quad 1\leq j \leq l-2,
\notag\\
\overline{\mathcal D}_0&=\overline{\mathcal E}_{l-1}=\overline{\mathcal E}_l=\emptyset.
\notag
\\
\boxed{\overline{\mathcal D}_j}&
=\left\{\lambda \in T\left|
\begin{array}{l}
\text{$\exists\underline{\lambda_i}\in\Z$, 
$\lambda_i-\left[\frac{l-j-1}{2}\right] \leq \underline{\lambda_i}\leq \lambda_i$, $1\leq i \leq l(\lambda)$, such that}
\\
\boxed{\bar{1}}\succeq k_{i,\lambda_i}\succeq \cdots \succeq k_{i,\underline{\lambda_i}+1}
\succeq \textrm{Max}\left(\boxed{\overline{j}}, \boxed{k_{i,\underline{\lambda_i}}}\right),
\\
\bar{1} \succeq k_{i,\underline{\lambda_i}} \succeq \cdots \succeq k_{i,2} \succeq \overline{j+1},
k_{i,1}=\bar{j} ~\mathrm{for}~ 1\leq i \leq l(\lambda) 
\\
\end{array}
\right.\right\}
\setminus \overline{\mathcal D}_j,\quad 1\leq j \leq l-3,
\notag
\\
\boxed{\overline{\mathcal E}_j}&
=
\left\{\lambda \in T\left|
\begin{array}{l}
\text{$\exists\underline{\lambda_i}\in\Z$, 
$\lambda_i-\left[\frac{l-j-1}{2}\right] \leq \underline{\lambda_i}\leq \lambda_i$, $1\leq i \leq l(\lambda)$, such that
}
\\
\boxed{\bar{1}}\succeq k_{i,\lambda_i}\succeq \cdots \succeq k_{i,\underline{\lambda_i}+1}
\succeq \textrm{Max}\left(\boxed{\overline{j}}, \boxed{k_{i,\underline{\lambda_i}}}\right),
\\
\bar{1} \succeq k_{i,\underline{\lambda_i}} \succeq \cdots \succeq k_{i,2} \succeq \overline{j+1},
k_{i,1}=\overline{j+1}
~\mathrm{for}~ 1\leq i \leq l(\lambda)
\\
\end{array}
\right.\right\}
\setminus \overline{\mathcal E}_j,\quad 1\leq j \leq l-3,
\notag
\\
\boxed{\overline{\mathcal D}_0}&=
\boxed{\overline{\mathcal D}_{l-2}}=
\boxed{\overline{\mathcal D}_{l-1}}
=\boxed{\overline{\mathcal D}_l}=
\boxed{\overline{\mathcal E}_0}=
\boxed{\overline{\mathcal E}_{l-2}}=
\boxed{\overline{\mathcal E}_{l-1}}
=\boxed{\overline{\mathcal E}_l}=\emptyset.
\notag
\end{align}
Here $[x]$ is the floor function. 

We set
\begin{align}
&
\overline{\mathcal F}_0=\overline{\mathcal E}_0,\quad
\overline{\mathcal F}_j=\overline{\mathcal D}_j \cup 
\overline{\mathcal E}_j \cup \boxed{\overline{\mathcal D}_j} \cup \boxed{\overline{\mathcal E}_j}, 1\leq j \leq l-3,
\notag\\
&
\overline{\mathcal F}_{l-2}=\overline{\mathcal D}_{l-2}\cup \overline{\mathcal E}_{l-2},\quad
\overline{\mathcal F}_{l-1}=\overline{\mathcal D}_{l-1},
\quad
\overline{\mathcal F}_l=\overline{\mathcal D}_l.\notag
\end{align}

For $\lambda \in \overline{\mathcal D}_j \cup \overline{\mathcal E}_j$, $0\leq j \leq l, 1\leq i \leq Na_j$, we set
\begin{align}
&\overline{\eta}_{\lambda; (j,i)}^a[z]
=\overline{\eta}_{k_{i,a},\ldots, k_{i,1};j}[z],\quad 1\leq a \leq \lambda_i,
\notag\\
&\notag \\
&\nu_{\lambda;j}=\prod_{i=1}^{l(\lambda)}\nu_{k_{i,\lambda_i},\ldots, k_{i,1};j},\quad
{\mathcal A}_{\lambda;j}[z_{j,1},\ldots, z_{j,l(\lambda)}]=\overset{\curvearrowright}{\prod_{1\leq i \leq l(\lambda)}}
{\mathcal A}_{k_{i, \lambda_i},\ldots, k_{i,1}; j}[z_{j,i}].
\notag
\end{align}
For $\lambda \in \boxed{\overline{\mathcal D}_j} \cup \boxed{\overline{\mathcal E}_j}$, $1\leq j \leq l-3, 1\leq i \leq Na_j$, we set
\begin{align}
&\overline{\eta}_{\lambda;(j,i)}^a[z]=
\left\{\begin{array}{cc}
\overline{\eta}_{k_{i,a},\ldots, k_{i,1};j}[z],&1\leq a \leq \underline{\lambda_i},\\
\\
\overline{\eta}_{k_{i,\underline{\lambda_i}},\ldots, k_{i,1};j}^{b_{i,a},\ldots, b_{i,\underline{\lambda_i}+1}}[z],
&\underline{\lambda_i}<a \leq \lambda_i,
\end{array}\right.\notag
\\
&\nu_{\lambda;j}=\prod_{i=1}^{l(\lambda)}
\nu_{k_{i,\underline{\lambda_i}},\ldots, k_{i,1};j}^{
{b_{i,\lambda_i}},\ldots, {b_{i,\underline{\lambda_i}+1}}},\quad
{\mathcal A}_{\lambda; j}
[z_{j,1},\ldots, z_{j,l(\lambda)}]=
\overset{\curvearrowright}{\prod_{1\leq i \leq l(\lambda)}}
{\mathcal A}_{k_{i,\underline{\lambda_i}},\ldots, k_{i,1};j}^{
{b_{i,\lambda_i}},\ldots, {b_{i,\underline{\lambda_i}+1}}}[z_{j,i}]
,\notag
\end{align}
where we set
$b_{i,\lambda_i}=
\boxed{k_{i,\lambda_i}},\ldots, b_{i,\underline{\lambda_i}+1}=
\boxed{k_{i, \underline{\lambda_i}+1}} \in I$.
In the case of $\lambda_i=\underline{\lambda_i}$,
$\eta_{k_{i,\underline{\lambda_i}},\ldots, k_{i, 1};j}^{b_{i,\lambda_i},\ldots, b_{i, \underline{\lambda_i}+1}}[z]$,
$\nu_{k_{i,\underline{\lambda_i}},\ldots, k_{i,1};j}^{
{b_{i,\lambda_i}},\ldots, {b_{i,\underline{\lambda_i}+1}}}$, and
${\mathcal A}_{k_{i,\underline{\lambda_i}},\ldots, k_{i,1};j}^{
{b_{i,\lambda_i}},\ldots, {b_{i,\underline{\lambda_i}+1}}}[z]$ 
represents
$\eta_{k_{i,\lambda_i},\ldots, k_{i, 1};j}[z]$,
$\nu_{k_{i,\lambda_i},\ldots, k_{i,1};j}$, and
${\mathcal A}_{k_{i,\lambda_i},\ldots, k_{i,1};j}[z]$, respectively.

For $0\leq a, b, c \leq M$ such that $1\leq a+b+c \leq M$, we obtain $\mathrm{LHS}_{M,N}(a,b,c)$ as follows.
\begin{align}
&\mathrm{LHS}_{M,N}(a,b,c)
=
\left(\begin{array}{c}
M\\
a+b+c
\end{array}
\right)(a+b+c)! \int_{C_{[\delta]}}\underline{dz}\prod_{i=1}^{M-(a+b+c)}\int_{C_1}\underline{dw_i}
\notag\\
&\times
\left(\sum\limits_{\substack{\scriptscriptstyle 
\overline{\lambda}[0] \in \overline{\mathcal F}_0 \\ \scriptscriptstyle 
|\overline{\lambda}[0]|=a}}
\com{N}{\overline{\lambda}[0]}
\nu_{\overline{\lambda}[0];0}
+
\sum_{\overline{\lambda}[0] \in \overline{\mathcal F}_0,
\overline{\lambda}[1]\in \overline{\mathcal F}_1
\atop{|\overline{\lambda}[0]|+|\overline{\lambda}[1]|=a}}
\com{N}{\overline{\lambda}[0]}\com{N}{\overline{\lambda}[1]}
\nu_{\overline{\lambda}[0];0}
\nu_{\overline{\lambda}[1];1}
+
\sum_{\overline{\lambda}[1] \in \overline{\mathcal F}_1
\atop{|\overline{\lambda}[1]|=a}}\com{N}{\overline{\lambda}[1]}
\nu_{\overline{\lambda}[1];1}\right)
\notag\\
&\times
\left(\sum_{n=0}^{Min(b,l-3)}
\sum_{2\leq j_1<\cdots<j_n \leq l-2}
\sum_{\overline{\lambda}[j_1]\in \overline{\mathcal F}_{j_1},\ldots, \overline{\lambda}[j_n]\in \overline{\mathcal F}_{j_n}
\atop{|\overline{\lambda}[j_1]|+\cdots+|\overline{\lambda}[j_n]|=b}}
\prod_{s=1}^n
\com{2N}{\overline{\lambda}[j_s]}\prod_{s=1}^n \nu_{\overline{\lambda}[j_s];j_s}
\right)
\times\left(
\sum_{\overline{\lambda}[l-1]\in \overline{\mathcal F}_{l-1}
\atop{|\overline{\lambda}[l-1]|=c}}
\com{N}{\overline{\lambda}[l-1]}\nu_{\overline{\lambda}[l-1];l-1}\right.
\notag
\\
&+
\left.
\sum_{\overline{\lambda}[l-1]\in \overline{\mathcal F}_{l-1},
\overline{\lambda}[l] \in \overline{\mathcal F}_{l}
\atop{|\overline{\lambda}[l-1]|+|\overline{\lambda}[l]|=c}}
\com{N}{\overline{\lambda}[l-1]}\com{N}{\overline{\lambda}[l]}
\nu_{\overline{\lambda}[l-1];l-1}\nu_{\overline{\lambda}[l];l}
+
\sum_{\overline{\lambda}[l] \in \overline{\mathcal F}_{l}
\atop{|\overline{\lambda}[l]|=c}}
\com{N}{\overline{\lambda}[l]}\nu_{\overline{\lambda}[l];l}
\right)~\overset{\curvearrowright}{\prod_{1\leq i \leq M-(a+b+c)}}{\bf T}_1(w_i)
\notag\\
&\times
\mathbf{S}[(\overline{\lambda}[0];0), (\overline{\lambda}[1];1),
\ldots,
(\overline{\lambda}[j_s];j_s),\ldots, (\overline{\lambda}[l-1], l-1), (\overline{\lambda}[l];l)]
\notag\\
&\times
\overline{\omega}_M
\left[w_1,\ldots, w_{M-(a+b+c)}; 
(\overline{\lambda}[0];0), (\overline{\lambda}[1];1), \ldots, (\overline{\lambda}[j_s];j_s),\ldots, (\overline{\lambda}[l-1];l-1), (\overline{\lambda}[l],l)\right]
\bh.
\notag
\end{align}

For $0\leq j_1<j_2<\cdots<j_n \leq l$ and $\overline{\lambda}[j_s] \in \overline{\mathcal F}_{j_s}$, $1\leq s \leq n$,
we set
\begin{align}
&\mathbf{S}\left[(\overline{\lambda}[j_1];j_1), (\overline{\lambda}[j_2];j_2),\ldots,(\overline{\lambda}[j_n];j_n)\right]\notag\\
=&\mathbf{S}_{[0,j_1-1]} {\mathcal A}_{\overline{\lambda}[j_1]; j_1}[z_{j_1,1},\ldots, z_{j_1,l(\overline{\lambda}[j_1])}]
\overset{\curvearrowright}{\prod_{l(\overline{\lambda}[j_1])<s \leq N a_{j_1}}} 
S_{j_1}(z_{j_1,s})\notag
\\
&\times
\mathbf{S}_{[j_1+1,j_2-1]} {\mathcal A}_{\overline{\lambda}[j_2]; j_2}[z_{j_2,1},\ldots, z_{j_2,l(\overline{\lambda}[j_2])}]
\overset{\curvearrowright}{\prod_{l(\overline{\lambda}[j_2])<s \leq N a_{j_2}}} 
S_{j_2}(z_{j_2,s})\notag
\\
&
\cdots \cdots \cdots
\notag
\end{align}
\begin{align}
&\times
\mathbf{S}_{[j_{n-1}+1,j_n-1]} {\mathcal A}_{\overline{\lambda}[j_n]; j_n}[z_{j_n,1},\ldots, z_{j_n,l(\overline{\lambda}[j_n])}]
\overset{\curvearrowright}{\prod_{l(\overline{\lambda}[j_n])<s \leq N a_{j_n}}} 
S_{j_n}(z_{j_n,s})\mathbf{S}_{[j_{n}+1,l]}.
\notag
\end{align}
For $0\leq j_1<j_2<\cdots<j_n \leq l$, 
and $|\overline{\lambda}[j_1]|+\cdots+|\overline{\lambda}[j_n]|=a+b+c$,
we set
\begin{align}
&\overline{\omega}_{M}[w_1, w_2,\ldots,w_{M-(a+b+c)}; (\overline{\lambda}[j_1];j_1), (\overline{\lambda}[j_2];j_2),\ldots, (\overline{\lambda}[j_n];j_n))]\notag\\
&=\omega_{M-(a+b+c)}(w_1,\ldots,w_{M-(a+b+c)})\prod_{i=1}^{M-(a+b+c)}
\prod_{s=1}^n
\prod_{i_s=1}^{l(\overline{\lambda}[j_s])}
\prod_{1\leq d \leq \overline{\lambda}[j_s]_{i_s}}
\omega_2\left(w_i, \overline{\eta}_{\overline{\lambda}[j_s];(j_s,i_s)}^d [z_{j_s, i_s}]\right)\notag\\
&\times
\prod_{1\leq s \leq n} 
\prod_{1\leq i_1<i_2 \leq l(\overline{\lambda}[j_s])}
\prod_{1\leq d_1 \leq \overline{\lambda}[j_s]_{i_1}
\atop{1\leq d_2 \leq \overline{\lambda}[j_s]_{i_2}}}
\omega_2\left(
\overline{\eta}_{\overline{\lambda}[j_s]; (j_s,i_1)}^{d_1}[z_{j_s,i_1}], 
\overline{\eta}_{\overline{\lambda}[j_s]; (j_s,i_2)}^{d_2}[z_{j_s,i_2}]\right)\notag
\\
&\times
\prod_{1\leq s< t\leq n}
\prod_{1\leq i_1 \leq l(\overline{\lambda}[j_s])
\atop{1\leq i_2 \leq l(\overline{\lambda}[j_t])}}
\prod_{1\leq d_1 \leq \overline{\lambda}[j_s]_{i_1}
\atop{1\leq d_2 \leq \overline{\lambda}[j_t]_{i_2}}}
\omega_2\left(
\overline{\eta}_{\overline{\lambda}[j_s]; (j_s,i_1)}^{d_1}[z_{j_s,i_1}], 
\overline{\eta}_{\overline{\lambda}[j_t]; (j_t,i_2)}^{d_2}[z_{j_t,i_2}]\right).
\notag
\end{align}
In the case of $|\lambda[j_s]|=0$, we read
\begin{align}
&\mathbf{S}[(\lambda[j_1]; j_1),\ldots, (\lambda[j_s];j_s),\ldots, (\lambda[j_n];j_n)]~\mathrm{and}\notag\\
&\overline{\omega}_{M}[w_1,\ldots,w_{M-(a+b+c)}; (\lambda[j_1];j_1), \ldots, (\lambda[j_s];j_s),\ldots, (\lambda[j_n];j_n))]
~\mathrm{as}
\notag
\\
&\mathbf{S}[(\lambda[j_1]; j_1),\ldots, (\lambda[j_{s-1}];j_{s-1}),
(\lambda[j_{s+1}];j_{s+1}),
\ldots, (\lambda[j_n];j_n)]~{\rm and}
\notag\\
&\overline{\omega}_{M}[w_1,\ldots,w_{M-(a+b+c)}; (\lambda[j_1];j_1), \ldots, (\lambda[j_{s-1}];j_{s-1}),
(\lambda[j_{s+1}];j_{s+1}),
\ldots, (\lambda[j_n];j_n))], \mathrm{respectively}.\notag
\end{align}

\subsubsection{The Right Hand Side}

We start from ${\bf G}_N(\vartheta)\cdot {\bf I}_M$.
We set the radii of integration circles $r_1, r_2, \ldots, r_M$ as $|w_i|=r_i, 1\leq i \leq M$.
We move the contours from $r_1=r_2=\cdots=r_M=r \ll 1$ to $r_1=r_2=\cdots=r_M=1$ imposing the conditions:
\begin{align}
r<r_M<\cdots<r_2<r_1<1
\quad
\mathrm{and}
\quad
|q|r_i<r_j,\quad 1\leq i<j \leq M,\notag
\end{align}
in the case of $1<r_1, r_2, \cdots, r_M<r$.
In the process of this deformation of contours,
we take the residues at $w_i={\eta}_{k_{m},\ldots, k_{1};j}[z_{j,s}]$ 
and $w_i={\eta}_{k_{n},\ldots,k_1;j}^{ k_{m+n},\ldots,k_{n+1}}[z_{j,s}]$
given as follows.

For $m=1$, ${\eta}_{k_1;j}[z]$ is given in (\ref{def:eta1}).
For $2\leq m \leq M$, we set
\begin{align}
\label{def:eta4}
{\eta}_{k_m,\ldots,k_1;j}[z]=\left\{\begin{array}{cl}
{e}_0^{(m)}[z],& j=0,\quad 1\preceq k_m \preceq \cdots \preceq k_1 \preceq 2,\\
{d}_j^{(m)}[z],& 1\leq j \leq l-2,\quad 1 \preceq k_m \preceq \cdots \preceq k_2\preceq j+1, k_1=j,\\
{e}_j^{(m)}[z],& 1\leq j \leq l-2,\quad 1 \preceq k_m \preceq \cdots \preceq k_2\preceq j+1, k_1=j+1.
\end{array}
\right.
\end{align}
For $1\leq m,n \leq M$ such that $2\leq m+n \leq M$, we set
\begin{align}
\label{def:eta5}
\eta_{k_{n},\ldots, k_{1};j}^{k_{m+n},\ldots, k_{n+1}}[z]=\left\{\begin{array}{cl}
{D}_j^{(m)}[z],& 1\leq j \leq l-3,
\begin{array}{c}
1 \preceq k_{m+n} \preceq \cdots \preceq k_{n+1}\preceq \mathrm{Min}(j, k_n),\\
1 \preceq k_n \preceq \cdots \preceq k_2\preceq j+1, k_1=j,
\end{array}
\\
\\
{E}_j^{(m)}[z],& 1\leq j \leq l-3,
\begin{array}{c}
1 \preceq k_{m+n} \preceq \cdots \preceq k_{n+1}\preceq \mathrm{Min}(j, k_n),\\
1 \preceq k_n \preceq \cdots \preceq k_2\preceq j+1, k_1=j+1.
\end{array}
\end{array}
\right.
\end{align}

\vspace{2mm}

For $(k_m,\ldots,k_1;j)$ in (\ref{def:eta4}), we set
\begin{align}
&
\nu_{k_m,\ldots,k_1;j}
\mathcal{A}_{k_m,\ldots,k_1;j}[z]
=\overset{\curvearrowleft}{\prod_{1\leq a \leq m}}
\mathop{\mathrm{res}}_{w_a={\eta}_{k_a,\ldots, k_2, k_1;j}[z]}
\frac{dw_a}{w_a}S_j(z)
\overset{\curvearrowright}{\prod_{1\leq a \leq m}}\bLambda_{k_a}(w_a)
\omega_{m}(w_1,w_2, \ldots, w_m),
\notag
\\
&
\mathcal{A}_{k_m,\ldots,k_1;j}[z]=:S_j(z)\prod_{1\leq a \leq m}
\bLambda_{k_a}({\eta}_{k_a,\ldots,k_1;j}[z]):.\notag
\end{align}
For $\left(\begin{array}{c}
k_{m+n},\ldots, k_{n+1}\\
k_n, \ldots, k_1; j\end{array}
\right)$ in (\ref{def:eta5}), 
we set
\begin{align}
&\nu_{k_n,\ldots,k_1;j}^{k_{m+n},\ldots,k_{n+1}}
\mathcal{A}_{k_n,\ldots,k_1;j}^{k_{m+n},\ldots,k_{n+1}}[z]\notag\\
=&\overset{\curvearrowleft}{\prod_{n+1\leq a \leq m+n}}
\mathop{\mathrm{res}}_{w_a={\eta}_{k_n,\ldots, k_1;j}^{k_a,
\ldots,k_{n+1}}[z]}
\frac{dw_a}{w_a}
\overset{\curvearrowleft}{\prod_{1\leq a \leq n}}
\mathop{\mathrm{res}}_{w_a={\eta}_{k_a,\ldots, k_1;j}[z]}\frac{dw_a}{w_a}
S_j(z)\overset{\curvearrowright}{\prod_{1\leq a \leq m+n}}\bLambda_{k_a}(w_a)
\omega_{m+n}(w_{1}, w_2, \ldots, w_2, w_{m+n}),
\notag
\\
&\mathcal{A}_{k_n,\ldots,k_1;j}^{k_{m+n},\ldots,k_{n+1}}[z]=
:S_j(z)\prod_{n+1\leq a \leq m+n}
\bLambda_{k_a}({\eta}_{k_n,\ldots,k_1;j}^{k_a,\ldots,k_{n+1}}[z])
\prod_{1\leq a \leq n}
\bLambda_{k_a}({\eta}_{k_a,\ldots,k_1;j}[z]):.\notag
\end{align}

We define subsets $\mathcal{{D}}_j, \mathcal{{E}}_j, \boxed{\mathcal{{D}}_j}, \boxed{\mathcal{{E}}_j}
\subset T$ as follows.
\begin{align}
{\mathcal D}_j&=
\left\{\lambda \in T\left|
1 \preceq k_{i,\lambda_i}\preceq \cdots \preceq k_{i,2} \preceq j+1, k_{i,1}=j
\textrm{~for~} 1\leq i \leq l(\lambda)
\right.
\right\},~1\leq j \leq l-2,
\notag\\
{\mathcal D}_{l-1}&=\left\{\lambda \in T\left|
k_{i,1}=l-1, \bar{l} \textrm{~and~} \lambda_i=1 \textrm{~for~} 1\leq i \leq l(\lambda)
\right.\right\},
\notag\\
{\mathcal D}_{l}&=\left\{\lambda \in T \left|
k_{i,1}=l-1, l \textrm{~and~} \lambda_i=1 \textrm{~for~} 1\leq i \leq l(\lambda)
\right.\right\},
\notag\\
{\mathcal E}_0&=\left\{\lambda \in T\left|
1\preceq k_{i,\lambda_i}\preceq \cdots \preceq k_{i,1}=1, 2 \textrm{~for~} 1\leq i \leq l(\lambda)
\right.\right\},\notag
\\
{\mathcal E}_j&=
\left\{\lambda \in T\left|
1 \preceq k_{i,\lambda_i} \preceq \cdots \preceq k_{i,2} \preceq  j+1, k_{i,1}=j+1 \textrm{~for~} 1\leq i \leq l(\lambda)
\right.\right\},\quad 1\leq j \leq l-2,
\notag\\
{\mathcal D}_0&={\mathcal E}_{l-1}={\mathcal E}_l=\emptyset.
\notag
\\
\boxed{{\mathcal D}_j}&
=\left\{\lambda \in T\left|
\begin{array}{l}
\text{
$\exists\underline{\lambda_i}\in\Z$, 
$\lambda_i-\left[\frac{l-j-1}{2}\right] \leq \underline{\lambda_i}\leq \lambda_i$, $1\leq i \leq l(\lambda)$, such that
}\\
\boxed{1}\preceq k_{i,\lambda_i}\preceq \cdots \preceq k_{i,\underline{\lambda_i}+1}
\preceq \textrm{Min}\left(\boxed{j}, \boxed{k_{i,\underline{\lambda_i}}}\right),
\\
1 \preceq k_{i,\underline{\lambda_i}} \preceq \cdots \preceq k_{i,2} \preceq j+1,
k_{i,1}=j ~\mathrm{for}~ 1\leq i \leq l(\lambda)
\\
\end{array}
\right.\right\}
\setminus {\mathcal D}_j,\quad 1\leq j \leq l-3,
\notag
\\
\boxed{{\mathcal E}_j}&
=
\left\{\lambda \in T\left|
\begin{array}{l}
\text{
$\exists\underline{\lambda_i}\in\Z$,
$\lambda_i-\left[\frac{l-j-1}{2}\right] \leq \underline{\lambda_i}\leq \lambda_i$, $1\leq i \leq l(\lambda)$, such that
}\\
\boxed{1}\preceq k_{i,\lambda_i}\preceq \cdots \preceq k_{i,\underline{\lambda_i}+1}
\preceq \textrm{Min}\left(\boxed{j}, \boxed{k_{i,\underline{\lambda_i}}}\right),
\\
1 \preceq k_{i,\underline{\lambda_i}} \preceq \cdots \preceq k_{i,2} \preceq j+1,
k_{i,1}=j+1 ~\mathrm{for}~ 1\leq i \leq l(\lambda)
\\
\end{array}
\right.\right\}
\setminus {\mathcal E}_j,\quad 1\leq j \leq l-3,
\notag
\\
\boxed{{\mathcal D}_0}&=
\boxed{{\mathcal D}_{l-2}}=
\boxed{{\mathcal D}_{l-1}}
=\boxed{{\mathcal D}_l}=
\boxed{{\mathcal E}_0}=
\boxed{{\mathcal E}_{l-2}}=
\boxed{{\mathcal E}_{l-1}}
=\boxed{{\mathcal E}_l}=\emptyset.
\notag
\end{align}

We set
\begin{align}
&
{\mathcal F}_0={\mathcal E}_0,\quad
{\mathcal F}_j={\mathcal D}_j \cup 
{\mathcal E}_j \cup \boxed{{\mathcal D}_j} \cup \boxed{{\mathcal E}_j}, 1\leq j \leq l-3,
\notag\\
&
{\mathcal F}_{l-2}={\mathcal D}_{l-2}\cup \overline{\mathcal E}_{l-2},\quad
{\mathcal F}_{l-1}={\mathcal D}_{l-1},
\quad
{\mathcal F}_l={\mathcal D}_l.\notag
\end{align}

For $\lambda \in {\mathcal D}_j \cup {\mathcal E}_j$, $0\leq j \leq l, 1\leq i \leq Na_j$, we set
\begin{align}
&{\eta}_{\lambda; (j,i)}^a[z]={\eta}_{k_{i,a},\ldots, k_{i,1};j}[z],\quad 1\leq a \leq \lambda_i,
\notag\\
&\nu_{\lambda;j}=\prod_{i=1}^{l(\lambda)}\nu_{k_{i,\lambda_i},\ldots, k_{i,1};j},\quad
{\mathcal A}_{\lambda;j}[z_{j,1},\ldots, z_{j,l(\lambda)}]=\overset{\curvearrowright}{\prod_{1\leq i \leq l(\lambda)}}
{\mathcal A}_{k_{i, \lambda_i},\ldots, k_{i,1}; j}[z_{j,i}].
\notag
\end{align}
For $\lambda \in \boxed{{\mathcal D}_j} \cup \boxed{{\mathcal E}_j}$, $1\leq j \leq l-3, 1\leq i \leq Na_j$, we set
\begin{align}
&{\eta}_{\lambda;(j,i)}^a[z]=
\left\{\begin{array}{cc}
{\eta}_{k_{i,a},\ldots, k_{i,1};j}[z],&1\leq a \leq \underline{\lambda_i},\\
{\eta}_{k_{i,\underline{\lambda_i}},\ldots, k_{i,1};j}^{b_{i,a},\ldots, b_{i,\underline{\lambda_i}+1}}[z],
&\underline{\lambda_i}<a \leq \lambda_i,
\end{array}\right.\notag
\\
&\nu_{\lambda;j}=\prod_{i=1}^{l(\lambda)}
\nu_{k_{i,\underline{\lambda_i}},\ldots, k_{i,1};j}^{
{b_{i,\lambda_i}},\ldots, {b_{i,\underline{\lambda_i}+1}}},\quad
{\mathcal A}_{\lambda; j}
[z_{j,1},\ldots, z_{j,l(\lambda)}]=
\overset{\curvearrowright}{\prod_{1\leq i \leq l(\lambda)}}
{\mathcal A}_{k_{i,\underline{\lambda_i}},\ldots, k_{i,1};j}^{
{b_{i,\lambda_i}},\ldots, {b_{i,\underline{\lambda_i}+1}}}[z_{j,i}]
,\notag
\end{align}
where we set
$b_{i,\lambda_i}=
\boxed{k_{i,\lambda_i}},\ldots, b_{i,\underline{\lambda_i}+1}=
\boxed{k_{i, \underline{\lambda_i}+1}} \in I$.

For $0\leq a, b, c \leq M$ such that $1\leq a+b+c \leq M$, we obtain $\mathrm{RHS}_{M,N}(a,b,c)$ as follows.
\begin{align}
&\mathrm{RHS}_{M,N}(a,b,c)
=(-1)^{a+b+c}
\left(\begin{array}{c}
M\\
a+b+c
\end{array}
\right)
(a+b+c)! \int_{C_{[\delta]}}\underline{dz}\prod_{i=a+b+c+1}^{M}\int_{C_1}\underline{dw_i}
\notag\\
&\times
\left(\sum_{\lambda[0] \in {\mathcal F}_0
\atop{|\lambda[0]|=a}}\com{N}{\lambda[0]}
\nu_{\lambda[0];0}
+
\sum_{\lambda[0] \in {\mathcal F}_0,
\lambda[1]\in {\mathcal F}_1
\atop{|\lambda[0]|+|\lambda[1]|=a}}
\com{N}{\lambda[0]}\com{N}{\lambda[1]}
\nu_{\lambda[0];0}\nu_{\lambda[1];1}
+
\sum_{\lambda[1] \in {\mathcal F}_1
\atop{|\lambda[1]|=a}}\com{N}{\lambda[1]}
\nu_{\lambda[1];1}\right)
\notag\\
&\times
\left(\sum_{n=0}^{Min(b,l-3)}
\sum_{2\leq j_1<\cdots<j_n \leq l-2}
\sum_{\lambda[j_1]\in {\mathcal F}_{j_1},\ldots, \lambda[j_n]\in {\mathcal F}_{j_n}
\atop{|\lambda[j_1]|+\cdots+|\lambda[j_n]|=b}}
\prod_{s=1}^n
\com{2N}{\lambda[j_s]}\prod_{s=1}^n \nu_{\lambda[j_s];j_s}
\right)
\times\left(
\sum_{\lambda[l-1]\in {\mathcal F}_{l-1}
\atop{|\lambda[l-1]|=c}}
\com{N}{\lambda[l-1]}\nu_{\lambda[l-1];l-1}\right.\notag\\
&+
\left.\sum_{\lambda[l-1]\in {\mathcal F}_{l-1},
\lambda[l]\in {\mathcal F}_{l}
\atop{|\lambda[l-1]|+|\lambda[l]|=c}}
\com{N}{\lambda[l-1]}\com{N}{\lambda[l]}
\nu_{\lambda[l-1];l-1}\nu_{\lambda[l];l}
+
\sum_{\lambda[l]\in {\mathcal F}_{l}
\atop{|\lambda[l]|=c}}
\com{N}{\lambda[l]}\nu_{\lambda[l];l}
\right)
\notag\\
&\times
\mathbf{S}[(\lambda[0];0), (\lambda[1];1),
\ldots,
(\lambda[j_s];j_s),\ldots, (\lambda[l-1], l-1), (\lambda[l];l)]
~\overset{\curvearrowright}{\prod_{a+b+c+1\leq i \leq M}}{\bf T}_1(w_i)
\notag\\
&\times
{\omega}_M
\left[
(\lambda[0];0), (\lambda[1];1), \ldots, (\lambda[j_s];j_s),\ldots, (\lambda[l-1];l-1), (\lambda[l],l); w_{a+b+c+1},\ldots, w_{M}\right]
\bh.
\notag
\end{align}
For $0\leq j_1<j_2<\cdots<j_n \leq l$ and $|\lambda[j_1]|+\cdots+|\lambda[j_n]|=a+b+c$, we set
\begin{align}
&{\omega}_{M}[(\lambda[j_1];j_1), (\lambda[j_2];j_2),\ldots, (\lambda[j_n];j_n)); 
w_{a+b+c+1},\ldots, w_M ]\notag\\
&=\omega_{M-(a+b+c)}(
w_{a+b+c-1},\ldots,w_{M})
\prod_{
i=a+b+c+1}^{M}
\prod_{s=1}^n
\prod_{i_s=1}^{l(\lambda[j_s])}
\prod_{1\leq d \leq \lambda[j_s]_{i_s}}
\omega_2\left({\eta}_{\lambda[j_s];(j_s,i_s)}^d [z_{j_s, i_s}], w_i\right)\notag\\
&\times
\prod_{1\leq s \leq n} 
\prod_{1\leq i_1<i_2 \leq l(\lambda[j_s])}
\prod_{1\leq d_1 \leq \lambda[j_s]_{i_1}
\atop{1\leq d_2 \leq \lambda[j_s]_{i_2}}}
\omega_2\left(
{\eta}_{\lambda[j_s]; (j_s,i_1)}^{d_1}[z_{j_s,i_1}], 
{\eta}_{\lambda[j_s]; (j_s,i_2)}^{d_2}[z_{j_s,i_2}]\right)\notag
\\
&\times
\prod_{1\leq s< t\leq n}
\prod_{1\leq i_1 \leq l(\lambda[j_s])
\atop{1\leq i_2 \leq l(\lambda[j_t])}}
\prod_{1\leq d_1 \leq \lambda[j_s]_{i_1}
\atop{1\leq d_2 \leq \lambda[j_t]_{i_2}}}
\omega_2\left(
{\eta}_{\lambda[j_s]; (j_s,i_1)}^{d_1}[z_{j_s,i_1}], 
{\eta}_{\lambda[j_t]; (j_t,i_2)}^{d_2}[z_{j_t,i_2}]\right).
\notag
\end{align}

\vspace{2mm}

\subsubsection{The Cancellation.}

Lemma \ref{lem:pole1} yields the following.

For $1\leq b \leq M$, $0\leq a, c \leq M$, and $1\leq a+b+c \leq M$,
we obtain $\mathrm{LHS}_{M,N}(a,b,c)$ as follows.
\begin{align}
&\mathrm{LHS}_{M,N}(a,b,c)
=
\left(\begin{array}{c}
M\\
a+b+c
\end{array}
\right)(a+b+c)! 
\notag\\
&\times
\left(\sum_{\overline{\lambda}[0] \in \overline{\mathcal F}_0
\atop{|\overline{\lambda}[0]|=a}}\com{N}{\overline{\lambda}[0]}
\nu_{\overline{\lambda}[0];0}
+
\sum_{\overline{\lambda}[0] \in \overline{\mathcal F}_0,
\lambda[1]\in \overline{\mathcal F}_1
\atop{|\overline{\lambda}[0]|+|\overline{\lambda}[1]|=a}}
\com{N}{\overline{\lambda}[0]}\com{N}{\overline{\lambda}[1]}
\nu_{\overline{\lambda}[0];0}\nu_{\overline{\lambda}[1];1}
+
\sum_{\overline{\lambda}[1] \in \overline{\mathcal F}_1
\atop{|\overline{\lambda}[1]|=a}}\com{N}{\overline{\lambda}[1]}
\nu_{\overline{\lambda}[1];1}\right)
\notag\\
&\times
\left(\sum_{n=1}^{Min(b,l-3)}
\sum_{2\leq j_1<\cdots<j_n \leq l-2}
\sum_{\overline{\lambda}[j_1]\in \overline{\mathcal E}_{j_1} \cup 
{\tiny
\boxed{\scriptstyle  \overline{\mathcal E}_{j_1}}}
,\ldots, 
\overline{\lambda}[j_n]\in \overline{\mathcal E}_{j_n}
\cup{\tiny \boxed{\scriptstyle  
\overline{\mathcal E}_{j_n}}
}
\atop{|\overline{\lambda}[j_1]|+\cdots+|\overline{\lambda}[j_n]|=b}}
\prod_{s=1}^n
\com{2N}{\overline{\lambda}[j_s]}\prod_{s=1}^n \nu_{\overline{\lambda}[j_s];j_s}
\right)\notag\\
&
\times\left(
\sum_{\overline{\lambda}[l-1]\in \overline{\mathcal F}_{l-1}
\atop{|\overline{\lambda}[l-1]|=c}}
\com{N}{\overline{\lambda}[l-1]}\nu_{\overline{\lambda}[l-1];l-1}\right.
+
\left.\sum_{\overline{\lambda}[l-1]\in \overline{\mathcal F}_{l-1},
\lambda[l]\in \overline{\mathcal F}_{l}
\atop{|\overline{\lambda}[l-1]|+|\overline{\lambda}[l]|=c}}
\com{N}{\overline{\lambda}[l-1]}\com{N}{\overline{\lambda}[l]}
\nu_{\overline{\lambda}[l-1];l-1}\nu_{\overline{\lambda}[l];l}\right.\notag\\
&+\left.
\sum_{\overline{\lambda}[l]\in \overline{\mathcal F}_{l}
\atop{|\overline{\lambda}[l]|=c}}
\com{N}{\overline{\lambda}[l]}\nu_{\overline{\lambda}[l];l}
\right)
\times
\Delta_{a+b+c}^{(M,N)}[(\overline{\lambda}[0];0), (\overline{\lambda}[1];1),
\ldots,
(\overline{\lambda}[j_s];j_s),\ldots, (\overline{\lambda}[l-1], l-1), (\overline{\lambda}[l];l)].
\notag
\end{align}
For $1\leq b \leq M$, $0\leq a,c \leq M$, and $1\leq a+b+c \leq M$,
we obtain $\mathrm{RHS}_{M,N}(a,b,c)$ as follows.
\begin{align}
&\mathrm{RHS}_{M,N}(a,b,c)
=(-1)^{a+b+c}
\left(\begin{array}{c}
M\\
a+b+c
\end{array}
\right)(a+b+c)!
\notag
\end{align}
\begin{align}
&\times
\left(\sum_{\lambda[0] \in {\mathcal F}_0
\atop{|\lambda[0]|=a}}\com{N}{\lambda[0]}
\nu_{\lambda[0];0}
+
\sum_{\lambda[0] \in {\mathcal F}_0,
\lambda[1]\in {\mathcal F}_1
\atop{|\lambda[0]|+|\lambda[1]|=a}}
\com{N}{\lambda[0]}\com{N}{\lambda[1]}
\nu_{\lambda[0];0}\nu_{\lambda[1];1}
+
\sum_{\lambda[1] \in {\mathcal F}_1
\atop{|\lambda[1]|=a}}\com{N}{\lambda[1]}
\nu_{\lambda[1];1}\right)
\notag
\\
&\times
\left(\sum_{n=1}^{Min(b,l-3)}
\sum_{2\leq j_1<\cdots<j_n \leq l-2}
\sum_{\lambda[j_1]\in {\mathcal E}_{j_1}\cup
{\tiny \boxed{\scriptstyle  {\mathcal E}_{j_1}}}
,\ldots, \lambda[j_n]\in {\mathcal E}_{j_n}\cup 
{\tiny \boxed{\scriptstyle  
{\mathcal E}_{j_n}}}
\atop{|\lambda[j_1]|+\cdots+|\lambda[j_n]|=b}}
\prod_{s=1}^n
\com{2N}{\lambda[j_s]}\prod_{s=1}^n \nu_{\lambda[j_s];j_s}
\right)\notag
\\
&\times
\left(
\sum_{\lambda[l-1]\in {\mathcal F}_{l-1}
\atop{|\lambda[l-1]|=c}}
\com{N}{\lambda[l-1]}\nu_{\lambda[l-1];l-1}\right.+
\left.\sum_{\lambda[l-1]\in {\mathcal F}_{l-1},
\lambda[l]\in {\mathcal F}_{l}
\atop{|\lambda[l-1]|+|\lambda[l]|=c}}
\com{N}{\lambda[l-1]}\com{N}{\lambda[l]}
\nu_{\lambda[l-1];l-1}\nu_{\lambda[l];l}\right.\notag\\
&+\left.
\sum_{\lambda[l]\in {\mathcal F}_{l}
\atop{|\lambda[l]|=c}}
\com{N}{\lambda[l]}\nu_{\lambda[l];l}
\right)
\times
\Delta_{a+b+c}^{(M,N)}[(\lambda[0];0), (\lambda[1];1),
\ldots,
(\lambda[j_s];j_s),\ldots, (\lambda[l-1], l-1), (\lambda[l];l)].
\notag
\end{align}

For $1 \leq a \leq M$, we obtain
\begin{align}
&\mathrm{LHS}_{M,N}(a,0,0)-{\rm RHS}_{M,N}(a,0,0)=
\comb{M}{a}a!
\times
\left(\sum_{\overline{\lambda}[0] \in \overline{\mathcal E}_0
\atop{|\overline{\lambda}[0]|=a}}
\com{N }{\overline{\lambda}[0]} \nu_{\overline{\lambda}[0];0}\Delta^{(M,N)}_a[(\overline{\lambda}[0];0)]\right.
\notag
\\
&+\sum_{\overline{\lambda}[1] \in \overline{\mathcal E}_1 \cup 
{\tiny \boxed{\scriptstyle  
\overline{\mathcal E}_1}}
\atop{|\overline{\lambda}[1]|=a}}
\com{N}{\overline{\lambda}[1]}  
\nu_{\overline{\lambda}[1];1}\Delta^{(M,N)}_a[(\overline{\lambda}[1];1)]
-(-1)^a
\sum_{\lambda[1] \in {\mathcal E}_1 \cup 
{\tiny \boxed{\scriptstyle {\mathcal E}_1}}
\atop{|{\lambda}[1]|=a}}
\com{N}{{\lambda}[1]} \nu_{{\lambda}[1];1}\Delta^{(M,N)}_a[({\lambda[1]};1)]
\notag
\\
&+\sum_{\overline{\lambda}[0] \in \overline{\mathcal E}_0, 
\overline{\lambda}[1] \in \overline{\mathcal E}_1 \cup 
{\tiny \boxed{\scriptstyle  \overline{\mathcal E}_1}}
\atop{|\overline{\lambda}[0]|+|\overline{\lambda}[1]|=a}}
\com{N}{\overline{\lambda}[0]}
\com{N}{\overline{\lambda}[1]}
\nu_{\overline{\lambda}[0];0}
\nu_{\overline{\lambda}[1];1}
\Delta^{(M,N)}_a
[(\overline{\lambda}[0];0), (\overline{\lambda}[1];1)]
\notag\\
&\left.-(-1)^a
\sum_{{\lambda}[0] \in {\mathcal E}_0, {\lambda}[1] \in {\mathcal  E}_1 
\cup {\tiny \boxed{\scriptstyle  
{\mathcal E}_1}}
\atop{|{\lambda}[0]|+|{\lambda}[1]|=a}}
\com{N}{{\lambda}[0]}
\com{N}{{\lambda}[1]} 
\nu_{{\lambda}[0];0}
\nu_{{\lambda}[1];1}
\Delta^{(M,N)}_a
[({\lambda}[0];0), ({\lambda}[1];1)]\right).
\notag
\end{align}
For $1 \leq c \leq M$, we obtain
\begin{align}
&{\rm LHS}_{M,N}(0,0,c)-{\rm RHS}_{M,N}(0,0,c)=
\comb{M}{c} c!\notag\\
&\times
\left(
\sum_{\lambda[l-1] \in {\mathcal D}_{l-1}
\atop{|\lambda[l-1]|=c}}
\com{N}{\lambda[l-1]} 
\nu_{\lambda[l-1]; l-1}\Delta^{(M,N)}_c[(\lambda[l-1];l-1)]+
\sum_{\lambda[l] \in {\mathcal D}_{l}
\atop{|\lambda[l] |=c}}
\com{N}{\lambda[l]} 
\nu_{\lambda[l];l}\Delta^{(M,N)}_c[(\lambda[l]; l)]\right.
\notag
\end{align}
\begin{align}
&+\left.
\sum_{\lambda[l-1] \in {\mathcal D}_{l-1}, \lambda[l] \in {\mathcal D}_l
\atop{|\lambda[l-1]|+|\lambda[l]|=c}}
\com{N}{\lambda[l-1]} \com{N}{\lambda[l]}
\nu_{\lambda[l-1];l-1}\nu_{\lambda[l];l}
\Delta^{(M,N)}_c[(\lambda[l-1];l-1), (\lambda[l];l)]\right).
\notag
\end{align}
For $1 \leq a,c \leq M$ such that $2 \leq a+c \leq M$, we obtain
\begin{align}
&
\mathrm{LHS}_{M,N}(a,0,c)-{\rm RHS}_{M,N}(a,0,c)=\comb{M}{a+c}(a+c)!
\notag
\\
&\times
\left(
\sum_{\overline{\lambda}[0] \in \overline{\mathcal E}_0
\atop{|\overline{\lambda}[0]|=a}}
\sum_{\overline{\lambda}[l-1] \in \overline{\mathcal D}_{l-1}
\atop{|\overline{\lambda}[l-1]|=c}}
\com{N}{\overline{\lambda}[0]}
\com{N}{\overline{\lambda}[l-1]}
\nu_{\overline{\lambda}[0];0} 
\nu_{\overline{\lambda}[l-1];l-1} 
\Delta_{a+c}^{(M,N)}[(\overline{\lambda}[0];0), (\overline{\lambda}[l-1];l-1)]
\right.
\notag\\
&+
\sum_{\overline{\lambda}[0] \in \overline{\mathcal E}_0
\atop{|\overline{\lambda}[0]|=a}}
\sum_{\overline{\lambda}[l] \in \overline{\mathcal D}_{l}
\atop{|\overline{\lambda}[l]|=c}}
\com{N}{\overline{\lambda}[0]}
\com{N}{\overline{\lambda}[l]}
\nu_{\overline{\lambda}[0];0} 
\nu_{\overline{\lambda}[l];l} 
\Delta_{a+c}^{(M,N)}[(\overline{\lambda}[0];0), (\overline{\lambda}[l];l)]
\notag\\
&+
\sum_{\overline{\lambda}[0] \in \overline{\mathcal E}_0
\atop{|\overline{\lambda}[0]|=a}}
\sum_{\overline{\lambda}[l-1]\in \overline{\mathcal D}_{l-1},
\overline{\lambda}[l] \in \overline{\mathcal D}_{l}
\atop{|\overline{\lambda}[l-1]|+|\overline{\lambda}[l]|=c}}
\com{N}{\overline{\lambda}[0]}
\com{N}{\overline{\lambda}[l-1]}
\com{N}{\overline{\lambda}[l]}
\nu_{\overline{\lambda}[0];0} 
\nu_{\overline{\lambda}[l-1];l-1}
\nu_{\overline{\lambda}[l];l} \notag\\
&\times
\Delta_{a+c}^{(M,N)}[(\overline{\lambda}[0];0), (\overline{\lambda}[l-1];l-1), (\overline{\lambda}[l];l)]
\notag\\
&+
\sum_{\overline{\lambda}[1] \in \overline{\mathcal E}_1\cup 
{\tiny \boxed{\scriptstyle  \overline{\mathcal E}_1}}
\atop{|\overline{\lambda}[1]|=a}}
\sum_{\overline{\lambda}[l-1] \in \overline{\mathcal D}_{l-1}
\atop{|\overline{\lambda}[l-1]|=c}}
\com{N}{\overline{\lambda}[1]}
\com{N}{\overline{\lambda}[l-1]}
\nu_{\overline{\lambda}[1];1} 
\nu_{\overline{\lambda}[l-1];l-1} 
\Delta_{a+c}^{(M,N)}[(\overline{\lambda}[1];1), (\overline{\lambda}[l-1];l-1)]\notag
\\
&+
\sum_{\overline{\lambda}[1] \in \overline{\mathcal E}_1\cup 
{\tiny \boxed{\scriptstyle  \overline{\mathcal E}_1}}
\atop{|\overline{\lambda}[1]|=a}}
\sum_{\overline{\lambda}[l] \in \overline{\mathcal D}_{l}
\atop{|\overline{\lambda}[l]|=c}}
\com{N}{\overline{\lambda}[1]}
\com{N}{\overline{\lambda}[l]}
\nu_{\overline{\lambda}[1];1} 
\nu_{\overline{\lambda}[l];l} 
\Delta_{a+c}^{(M,N)}[(\overline{\lambda}[1];1), (\overline{\lambda}[l];l)]
\notag\\
&+
\sum_{\overline{\lambda}[1] \in \overline{\mathcal E}_1 \cup 
{\tiny \boxed{\scriptstyle  \overline{\mathcal E}_1}}
\atop{|\overline{\lambda}[1]|=a}}
\sum_{\overline{\lambda}[l-1]\in \overline{\mathcal D}_{l-1},
\overline{\lambda}[l] \in \overline{\mathcal D}_{l}
\atop{|\overline{\lambda}[l-1]|+|\overline{\lambda}[l]|=c}}
\com{N}{\overline{\lambda}[1]}
\com{N}{\overline{\lambda}[l-1]}
\com{N}{\overline{\lambda}[l]}
\nu_{\overline{\lambda}[1];1} 
\nu_{\overline{\lambda}[l-1];l-1}
\nu_{\overline{\lambda}[l];l} \notag\\
&\times
\Delta_{a+c}^{(M,N)}[(\overline{\lambda}[1];1), (\overline{\lambda}[l-1];l-1), (\overline{\lambda}[l];l)]
\notag\\
&+
\sum_{\overline{\lambda}[0] \in \overline{\mathcal E}_0, \overline{\lambda}[1] \in \overline{\mathcal E}_1\cup {
\tiny \boxed{\scriptstyle  
{\overline{\mathcal E}_1}}}
\atop{|\overline{\lambda}[0]|+|\overline{\lambda}[1]|=a}}
\sum_{\overline{\lambda}[l-1] \in \overline{\mathcal D}_{l-1}
\atop{|\overline{\lambda}[l-1]|=c}}
\com{N}{\overline{\lambda}[0]}
\com{N}{\overline{\lambda}[1]}
\com{N}{\overline{\lambda}[l-1]}
\nu_{\overline{\lambda}[0];0}
\nu_{\overline{\lambda}[1];1} 
\nu_{\overline{\lambda}[l-1];l-1}\notag\\
&\times
\Delta_{a+c}^{(M,N)}[(\overline{\lambda}[0];0), (\overline{\lambda}[1];1), (\overline{\lambda}[l-1];l-1)]\notag
\\
&+
\sum_{\overline{\lambda}[0]\in \overline{\mathcal E}_0,
\overline{\lambda}[1] \in \overline{\mathcal E}_1\cup {
\tiny \boxed{\scriptstyle  
\overline{\mathcal E}_1}}
\atop{|\overline{\lambda}[0]|+|\overline{\lambda}[1]|=a}}
\sum_{\overline{\lambda}[l] \in \overline{\mathcal D}_{l}
\atop{|\overline{\lambda}[l]|=c}}
\com{N}{\overline{\lambda}[0]}
\com{N}{\overline{\lambda}[1]}
\com{N}{\overline{\lambda}[l]}
\nu_{\overline{\lambda}[0];0} 
\nu_{\overline{\lambda}[1];1} 
\nu_{\overline{\lambda}[l];l} 
\Delta_{a+c}^{(M,N)}[
(\overline{\lambda}[0];0),
(\overline{\lambda}[1];1), (\overline{\lambda}[l];l)]
\notag
\end{align}
\begin{align}
&+
\sum_{\overline{\lambda}[0]\in \overline{\mathcal E}_0,
\overline{\lambda}[1] \in \overline{\mathcal E}_1 \cup {
\tiny
\boxed{\scriptstyle  
\overline{\mathcal E}_1}}
\atop{|\overline{\lambda}[0]|+|\overline{\lambda}[1]|=a}}
\sum_{\overline{\lambda}[l-1]\in \overline{\mathcal D}_{l-1},
\overline{\lambda}[l] \in \overline{\mathcal D}_{l}
\atop{|\overline{\lambda}[l-1]|+|\overline{\lambda}[l]|=c}}
\com{N}{\overline{\lambda}[0]}
\com{N}{\overline{\lambda}[1]}
\com{N}{\overline{\lambda}[l-1]}
\com{N}{\overline{\lambda}[l]}
\nu_{\overline{\lambda}[0];0} 
\nu_{\overline{\lambda}[1];1} 
\nu_{\overline{\lambda}[l-1];l-1}
\nu_{\overline{\lambda}[l];l} \notag\\
&\times
\Delta_{a+c}^{(M,N)}[(\overline{\lambda}[0];0), (\overline{\lambda}[1];1), (\overline{\lambda}[l-1];l-1), (\overline{\lambda}[l];l)]
\notag
\\
&-(-1)^{a+c}
\sum_{{\lambda}[1] \in {\mathcal E}_1\cup {
\tiny
\boxed{\scriptstyle  
{\mathcal E}_1}}
\atop{|{\lambda}[1]|=a}}
\sum_{{\lambda}[l-1] \in {\mathcal D}_{l-1}
\atop{|{\lambda}[l-1]|=c}}
\com{N}{{\lambda}[1]}
\com{N}{{\lambda}[l-1]}
\nu_{{\lambda}[1];1} 
\nu_{{\lambda}[l-1];l-1} 
\Delta_{a+c}^{(M,N)}[({\lambda}[1];1), ({\lambda}[l-1];l-1)]\notag
\\
&-(-1)^{a+c}
\sum_{{\lambda}[1] \in {\mathcal E}_1\cup {\tiny
\boxed{\scriptstyle  
{\mathcal E}_1}}
\atop{|{\lambda}[1]|=a}}
\sum_{\lambda[l] \in {\mathcal D}_{l}
\atop{|{\lambda}[l]|=c}}
\com{N}{{\lambda}[1]}
\com{N}{{\lambda}[l]}
\nu_{{\lambda}[1];1} 
\nu_{{\lambda}[l];l} 
\Delta_{a+c}^{(M,N)}[({\lambda}[1];1), ({\lambda}[l];l)]
\notag\\
&-(-1)^{a+c}
\sum_{{\lambda}[1] \in {\mathcal E}_1 \cup {\tiny
\boxed{\scriptstyle  
{\mathcal E}_1}}
\atop{|{\lambda}[1]|=a}}
\sum_{{\lambda}[l-1]\in {\mathcal D}_{l-1},
{\lambda}[l] \in {\mathcal D}_{l}
\atop{|{\lambda}[l-1]|+|{\lambda}[l]|=c}}
\com{N}{{\lambda}[1]}
\com{N}{{\lambda}[l-1]}
\com{N}{{\lambda}[l]}
\nu_{{\lambda}[1];1} 
\nu_{{\lambda}[l-1];l-1}
\nu_{{\lambda}[l];l} \notag\\
&\times
\Delta_{a+c}^{(M,N)}[({\lambda}[1];1), ({\lambda}[l-1];l-1), ({\lambda}[l];l)]
\notag
\\
&-(-1)^{a+c}
\sum_{{\lambda}[0] \in {\mathcal E}_0, {\lambda}[1] \in {\mathcal E}_1\cup 
{\tiny \boxed{\scriptstyle  
{\mathcal E}_1}}
\atop{|{\lambda}[0]|+|{\lambda}[1]|=a}}
\sum_{{\lambda}[l-1] \in {\mathcal D}_{l-1}
\atop{|{\lambda}[l-1]|=c}}
\com{N}{{\lambda}[0]}
\com{N}{{\lambda}[1]}
\com{N}{{\lambda}[l-1]}
\nu_{{\lambda}[0];0}
\nu_{{\lambda}[1];1} 
\nu_{{\lambda}[l-1];l-1}\notag\\
&\times
\Delta_{a+c}^{(M,N)}[({\lambda}[0];0), ({\lambda}[1];1), ({\lambda}[l-1];l-1)]\notag
\\
&-(-1)^{a+c}
\sum_{{\lambda}[0]\in {\mathcal E}_0,
{\lambda}[1] \in {\mathcal E}_1\cup {\tiny
\boxed{\scriptstyle  
{\mathcal E}_1}}
\atop{|{\lambda}[0]|+|{\lambda}[1]|=a}}
\sum_{\lambda[l] \in {\mathcal D}_{l}
\atop{|{\lambda}[l]|=c}}
\com{N}{{\lambda}[0]}
\com{N}{{\lambda}[1]}
\com{N}{{\lambda}[l]}
\nu_{{\lambda}[0];0} 
\nu_{{\lambda}[1];1} 
\nu_{{\lambda}[l];l} 
\Delta_{a+c}^{(M,N)}[
({\lambda}[0];0),
({\lambda}[1];1), ({\lambda}[l];l)]
\notag\\
&-(-1)^{a+c}
\sum_{{\lambda}[0]\in {\mathcal E}_0,
{\lambda}[1] \in {\mathcal E}_1 \cup {\tiny 
\boxed{\scriptstyle  
{\mathcal E}_1}}
\atop{|{\lambda}[0]|+|{\lambda}[1]|=a}}
\sum_{{\lambda}[l-1]\in {\mathcal D}_{l-1},
{\lambda}[l] \in {\mathcal D}_{l}
\atop{|{\lambda}[l-1]|+|{\lambda}[l]|=c}}
\com{N}{{\lambda}[0]}
\com{N}{{\lambda}[1]}
\com{N}{{\lambda}[l-1]}
\com{N}{{\lambda}[l]}
\nu_{{\lambda}[0];0} 
\nu_{{\lambda}[1];1} 
\nu_{{\lambda}[l-1];l-1}
\nu_{{\lambda}[l];l} \notag\\
&\times
\left.
\Delta_{a+c}^{(M,N)}[({\lambda}[0];0), ({\lambda}[1];1), ({\lambda}[l-1];l-1), ({\lambda}[l];l)]\right).
\notag
\end{align}

For $\overline{\lambda}[0] \in \overline{\mathcal F}_0$,
$\overline{\lambda}[1] \in \overline{\mathcal F}_1$,
$\overline{\lambda}[j_s] \in \overline{\mathcal E}_{j_s}\cup
\boxed{\overline{\mathcal E}_{j_s}}$, $1\leq s \leq n$,
$\overline{\lambda}[l-1] \in \overline{\mathcal F}_{l-1}$,
$\overline{\lambda}[l] \in \overline{\mathcal F}_l$
such that $|\overline{\lambda}[0]|+|\overline{\lambda}[1]|=a\geq 0$, 
$|\overline{\lambda}[j_1]|+\cdots+|\overline{\lambda}[j_n]|=b\geq 1$, 
$|\overline{\lambda}[l-1]|+|\overline{\lambda}[l]|=c \geq 0$,
$M \geq a+b+c \geq 1$,
we set
\begin{align}
&
\Delta^{(M,N)}_{a+b+c}
[(\overline{\lambda}[0];0),
(\overline{\lambda}[1];1),
(\overline{\lambda}[j_1];j_1),
\ldots, 
(\overline{\lambda}[j_n];j_n), (\overline{\lambda}[l-1];l-1), (\overline{\lambda}[l];l)]
\notag\\
=&\prod_{i=1}^{M-(a+b+c)}\int_{C_1}\underline{dw_i}
\int_{C[(\overline{\lambda}[0];0),
(\overline{\lambda}[1];1),
(\overline{\lambda}[j_1];j_1),
\ldots, (\overline{\lambda}[j_n];j_n), (\overline{\lambda}[l-1];l-1), (\overline{\lambda}[l];l)]}
\underline{dz}
~
\prod_{t=1}^n \prod_{s=1}^{l(\overline{\lambda}[j_t])}
\int_{C_1-C_{q^{-1}}}\underline{dz_{j_t; s}}
\notag
\end{align}
\begin{align}
&\times
\overset{\curvearrowright}{\prod_{1\leq i \leq M-(a+b+c)}}\bT_1(w_i)
{\bf S}[(\overline{\lambda}[0];0), (\overline{\lambda}[1];1), 
(\overline{\lambda}[j_1]; j_1),\ldots,
(\overline{\lambda}[j_n]; j_n),
(\overline{\lambda}[l-1];l-1),
(\overline{\lambda}[l];l)]\notag
\\
&\times
\overline{\omega}_M
[w_1,\ldots, w_{M-(a+b+c)};
(\overline{\lambda}[0];0), (\overline{\lambda}[1];1),
(\overline{\lambda}[j_1]; j_1),\ldots,
(\overline{\lambda}[j_n]; j_n),
(\overline{\lambda}[l-1];l-1),
(\overline{\lambda}[l];l)]\notag\\
&\times
\bh.
\notag
\end{align}
For ${\lambda}[0] \in {\mathcal F}_0$,
${\lambda}[1] \in {\mathcal F}_1$,
${\lambda}[j_s] \in {\mathcal E}_{j_s}\cup
\boxed{{\mathcal E}_{j_s}}$, $1\leq s \leq n$,
${\lambda}[l-1] \in {\mathcal F}_{l-1}$,
${\lambda}[l] \in {\mathcal F}_l$
such that $|{\lambda}[0]|+|{\lambda}[1]|=a\geq 0$, 
$|{\lambda}[j_1]|+\cdots+|{\lambda}[j_n]|=b\geq 1$, 
$|{\lambda}[l-1]|+|{\lambda}[l]|=c \geq 0$,
$M \geq a+b+c \geq 1$,
we set
\begin{align}
&
\Delta^{(M,N)}_{a+b+c}
[({\lambda}[0];0),
({\lambda}[1];1),
({\lambda}[j_1];j_1),
\ldots, 
({\lambda}[j_n];j_n), ({\lambda}[l-1];l-1), ({\lambda}[l];l)]
\notag\\
=&\prod_{i=a+b+c+1}^{M}
\int_{C_1}\underline{dw_i}
\int_{C[({\lambda}[0];0),
({\lambda}[1];1),
({\lambda}[j_1];j_1),
\ldots, ({\lambda}[j_n];j_n), ({\lambda}[l-1];l-1), ({\lambda}[l];l)]}
\underline{dz}
~
\prod_{t=1}^n \prod_{s=1}^{l({\lambda}[j_t])}
\int_{C_1-C_{q}}\underline{dz_{j_t; s}}
\notag\\
&\times
{\bf S}[({\lambda}[0];0), ({\lambda}[1];1), 
({\lambda}[j_1]; j_1),\ldots,
({\lambda}[j_n]; j_n),
({\lambda}[l-1];l-1),
({\lambda}[l];l)]
\overset{\curvearrowright}{\prod_{a+b+c+1\leq i \leq M}}\bT_1(w_i)
\notag
\\
&\times
{\omega}_M
[({\lambda}[0];0), ({\lambda}[1];1),
({\lambda}[j_1]; j_1),\ldots,
({\lambda}[j_n]; j_n),
({\lambda}[l-1];l-1),
({\lambda}[l];l); w_{a+b+c+1},\ldots, w_M]\notag\\
&\times
\bh.
\notag
\end{align}
For $\ol{\lambda}[0] \in \ol{\mathcal E}_0$ such that $
M \geq |\ol{\lambda}[0]|=a \geq 1$, we set
\begin{align}
\Delta^{(M,N)}_a
[(\overline{\lambda}[0]; 0)]&
=\prod_{i=1}^{M-a}\int_{C_1}\underline{dw_i}
\int_{C[(\overline{\lambda}[0];0)]}\underline{dz}
\left(\prod_{s=1}^{l(\overline{\lambda}[0])}\int_{C_{\delta_0}}\underline{dz_{0,s}}-
\prod_{s=1}^{l(\overline{\lambda}[0])}
\int_{C_{\delta_0 
q^{-\overline{\lambda}[0]_s}}}
\underline{dz_{0,s}}
\right)\notag\\
&\times\overset{\curvearrowright}{\prod_{1\leq i \leq M-a}}\bT_1(w_i)
{\bf S}[(\overline{\lambda}[0]; 0)]
\ol{\omega}_M[w_1,\ldots, w_{M-a}; (\overline{\lambda}[0]; 0)]
\bh.
\notag
\end{align}
For $\ol{\lambda}[1] \in \ol{\mathcal E}_1 \cup \boxed{\ol{\mathcal E}_1}$
such that $
M \geq |\ol{\lambda}[1]|=a \geq 1$, we set
\begin{align}
\Delta^{(M,N)}_a
[(\overline{\lambda}[1]; 1)]&
=\prod_{i=1}^{M-a}\int_{C_1}\underline{dw_i}
\int_{C[(\overline{\lambda}[1];1)]}\underline{dz}
\prod_{s=1}^{l(\overline{\lambda}[1])}\int_{C_1-C_q^{-1}}\underline{dz_{1,s}}
\notag\\
&\times\overset{\curvearrowright}{\prod_{1\leq i \leq M-a}}
\bT_1(w_i)
{\bf S}[(\overline{\lambda}[1]; 1)]
\ol{\omega}_M[(w_1,\ldots, w_{M-a}; \overline{\lambda}[1]; 1)]
\bh.\notag
\end{align}
For $\lambda[1]\in {\mathcal E}_1 \cup \boxed{{\mathcal E}_1}$ such that 
$
M \geq |\lambda[1]|=a\geq 1$, 
we set
\begin{align}
\Delta^{(M,N)}_a
[(\lambda[1]; 1)]&
=\prod_{i=a+1}^{M}\int_{C_1}\underline{dw_i}
\int_{C[(\lambda[1]; 1)]}\underline{dz}
\prod_{s=1}^{l({\lambda[1]})}\int_{C_1-C_q}\underline{dz_{1,s}}
\notag\\
&\times
{\bf S}[(\lambda[1]; 1)]
\overset{\curvearrowright}{\prod_{a+1\leq i \leq M}}\bT_1(w_i)
{\omega}_M[(\lambda[1]; 1), w_{a+1},\ldots, w_M ]
\bh.
\notag
\end{align}
For $\ol{\lambda}[0]\in \ol{\mathcal E}_0$, $\ol{\lambda}[1] \in \ol{\mathcal E}_1 \cup \boxed{\ol{\mathcal E}_1}$
such that 
$M \geq |\ol{\lambda}[0]|+|\ol{\lambda}[1]|=a\geq 2$, we set
\begin{align}
\Delta^{(M,N)}_a
[(\overline{\lambda}[0]; 0), (\overline{\lambda}[1]; 1)]&
=\prod_{i=1}^{M-a}\int_{C_1}\underline{dw_i}
\int_{C[(\overline{\lambda}[1]; 1)]}
\underline{dz}
\prod_{s=1}^{l(\overline{\lambda}[1])}
\int_{C_1-C_{q^{-1}}}\underline{dz_{1,s}}
\notag\\
&\times\overset{\curvearrowright}{\prod_{1\leq i \leq M-a}}\bT_1(w_i)
{\bf S}[(\overline{\lambda}[0]; 0),(\overline{\lambda}[1];1)]
\ol{\omega}_M[(w_1,\ldots, w_{M-a}; (\overline{\lambda}[0];0), (\overline{\lambda}[1];1)]
\notag\\
&\times \bh.\notag
\end{align}
For ${\lambda}[0]\in {\mathcal E}_0$, ${\lambda}[1] \in {\mathcal E}_1 \cup \boxed{{\mathcal E}_1}$
such that $
M \geq |{\lambda}[0]|+|{\lambda}[1]|=a\geq 2$, we set
\begin{align}
\Delta^{(M,N)}_a
[({\lambda}[0]; 0), ({\lambda}[1]; 1)]&
=\prod_{i=a+1}^{M}\int_{C_1}\underline{dw_i}
\int_{C[({\lambda}[0];0), ({\lambda}[1];1)]}
\underline{dz}
\prod_{s=1}^{l({\lambda}[0])}
\int_{C_{\delta_0}}\underline{dz_{0,s}}
\prod_{s=1}^{l({\lambda}[1])}
\int_{C_1-C_{q}}\underline{dz_{1,s}}
\notag\\
&\times
{\bf S}[({\lambda}[0]; 0),({\lambda}[1];1)]
\overset{\curvearrowright}{\prod_{a+1\leq i \leq M}}\bT_1(w_i)
{\omega}_M[({\lambda}[0];0), ({\lambda}[1];1); w_{a+1},\ldots, w_M]\notag\\
&\times \bh.\notag
\end{align}

For ${\lambda}[l-1] \in {\mathcal D}_{l-1}$ such that 
$M \geq |{\lambda}[l-1]|=c \geq 1$, we set
\begin{align}
 \Delta^{(M,N)}_c[(\lambda[l-1]; l-1)]&=
 \prod_{i=c+1}^{M}
 \int_{C_1}\underline{dw_i}
\int_{C[({\lambda}[l-1];l-1)]}
\underline{dz}
\left(\prod_{s=1}^{l({\lambda}[l-1])}
\int_{C_{q^{-1}}}\underline{dz_{l-1,s}}-
\prod_{s=1}^{l({\lambda}[l-1])}
\int_{C_{1}}
\underline{dz_{l-1,s}}
\right)
\notag\\
&\times
{\bf S}[({\lambda}[l-1]; l-1)]
\overset{\curvearrowright}{\prod_{c+1\leq i \leq M}}\bT_1(w_i)
{\omega}_M[(\lambda[l-1];l-1);w_{c+1},\ldots, w_M]\notag\\
&\times
\bh.\notag
\end{align}
For ${\lambda}[l] \in {\mathcal D}_{l}$ such that 
$M \geq |{\lambda}[l]|=c \geq 1$, we set
\begin{align}
 \Delta^{(M,N)}_c[(\lambda[l]; l)]
 &=
 \prod_{i=c+1}^{M}\int_{C_1}\underline{dw_i}
\int_{C[({\lambda}[l]; l)]}
\underline{dz}
\left(\prod_{s=1}^{l({\lambda}[l])}
\int_{C_{\delta_l q^{-1}}}
\underline{dz_{l,s}}-
\prod_{s=1}^{l({\lambda}[l])}
\int_{C_{\delta_l}}
\underline{dz_{l,s}}
\right)
\notag\\
&\times
{\bf S}[({\lambda}[l]; l)]
\overset{\curvearrowright}{\prod_{c+1 \leq i \leq M}}\bT_1(w_i)
{\omega}_M[(\lambda[l]; l); w_{c+1},\ldots, w_M]
\bh.
 \notag
 \end{align}
 For $\lambda[l-1] \in {\mathcal D}_{l-1}, \lambda[l] \in {\mathcal D}_l$ such that 
 $M \geq |\lambda[l-1]|+|\lambda[l]|=c\geq 2$, we set
 \begin{align}
 &\Delta^{(M,N)}_c[(\lambda[l-1]; l-1), (\lambda[l]; l)]
 =\prod_{i=c+1}^{M}\int_{C_1}\underline{dw_i}
\int_{C[({\lambda}[l-1];l-1), (\lambda[l];l)]}
\underline{dz}\notag\\
&\times
\left(\prod_{s=1}^{l({\lambda[l-1]})}
\int_{C_{q^{-1}}}\underline{dz_{l-1,s}}
\prod_{s=1}^{l(\lambda[l])}\int_{C_{\delta_l q^{-1}}}
\underline{dz_{l,s}}\right.-\left.
\prod_{s=1}^{l({\lambda[l-1]})}
\int_{C_{1}}\underline{dz_{l-1,s}}
\prod_{s=1}^{l(\lambda[l])}\int_{C_{\delta_l}}
\underline{dz_{l,s}}
\right)\notag\\
&\times
{\bf S}[({\lambda}[l-1];l-1), (\lambda[l];l)]
\overset{\curvearrowright}{\prod_{c+1 \leq i \leq M}}\bT_1(w_i)
{\omega}_M[(\lambda[l-1];l-1), (\lambda[l];l); w_{c+1},\ldots, w_M]
\notag\\
&\times \bh.\notag
 \end{align}
For $\ol{\lambda}[0] \in \ol{\mathcal E}_0$, $\lambda[l-1] \in \ol{\mathcal D}_{l-1}$ such that $|\ol{\lambda}[0]|=a \geq 1$, 
$|\ol{\lambda}[l-1]|=c \geq 1$, 
$M \geq a+c \geq 2$, we set
\begin{align}
&\Delta_{a+c}^{(M,N)}[
(\ol{\lambda}[0]; 0), (\ol{\lambda}[l-1];l-1)]
=
\prod_{i=1}^{M-(a+c)} \int_{C_1} \underline{dw_i} 
\int_{C[(\ol{\lambda}[0];0), (\ol{\lambda}[l-1];l-1)]}
\underline{dz}
\notag
\\
&\times
\left(
\prod_{s=1}^{l(\ol{\lambda}[0])}\int_{C_{\delta_0}}\underline{dz_{0,s}}
\prod_{s=1}^{l(\ol{\lambda}[l-1])}\int_{C_1}\ul{dz_{l-1,s}}
-
\prod_{s=1}^{l(\ol{\lambda}[0])}\int_{C_{\delta_0q^{-1}}}\underline{dz_{0,s}}
\prod_{s=1}^{l(\ol{\lambda}[l-1])}\int_{C_q}\ul{dz_{l-1,s}}
\right)\notag\\
&\times \overset{\curvearrowright}{\prod_{1\leq i \leq M-(a+c)}}
{\bT}_1(w_i) {\bf S}[(\ol{\lambda}[0];0), (\ol{\lambda}[l-1];l-1)]
\notag\\
&\times
\ol{\omega}_M
[w_1,\ldots, w_{M-(a+c)}; 
(\ol{\lambda}[0];0), (\ol{\lambda}[l-1];l-1)]
\bh.\notag
\end{align}
For $\ol{\lambda}[0] \in \ol{\mathcal E}_0$, $\lambda[l] \in \ol{\mathcal D}_{l}$ such that $|\ol{\lambda}[0]|=a \geq 1$, 
$|\ol{\lambda}[l]|=c \geq 1$, $M \geq a+c \geq 2$, 
we set
\begin{align}
&\Delta_{a+c}^{(M,N)}[
(\ol{\lambda}[0]; 0), (\ol{\lambda}[l];l)]
=
\prod_{i=1}^{M-(a+c)} \int_{C_1} \underline{dw_i} 
\int_{C[(\ol{\lambda}[0];0), (\ol{\lambda}[l];l)]}
\underline{dz}
\notag
\\
&\times
\left(
\prod_{s=1}^{l(\ol{\lambda}[0])}\int_{C_{\delta_0}}\underline{dz_{0,s}}
\prod_{s=1}^{l(\ol{\lambda}[l])}\int_{C_{\delta_l}}\ul{dz_{l,s}}
-
\prod_{s=1}^{l(\ol{\lambda}[0])}\int_{C_{\delta_0 q^{-1}}}\underline{dz_{0,s}}
\prod_{s=1}^{l(\ol{\lambda}[l])}\int_{C_{\delta_l q}}
\ul{dz_{l,s}}
\right)\notag\\
&\times \overset{\curvearrowright}{\prod_{1\leq i \leq M-(a+c)}}
{\bT}_1(w_i) {\bf S}[(\ol{\lambda}[0];0), (\ol{\lambda}[l];l)]
\notag\\
&\times
\ol{\omega}_M
[w_1,\ldots, w_{M-(a+c)}; 
(\ol{\lambda}[0];0), (\ol{\lambda}[l];l)]
\bh.\notag
\end{align}
For $\ol{\lambda}[0] \in \ol{\mathcal E}_0$, 
$\ol{\lambda}[l-1] \in \ol{\mathcal D}_{l-1}$,
$\ol{\lambda}[l] \in \ol{\mathcal D}_{l}$ such that $|\ol{\lambda}[0]|=a \geq 1$, 
$|\ol{\lambda}[l-1]|+|\ol{\lambda}[l]|=c \geq 2$, $M \geq a+c \geq 3$,
we set
\begin{align}
&\Delta_{a+c}^{(M,N)}[
(\ol{\lambda}[0]; 0), (\ol{\lambda}[l-1];l-1)), (\ol{\lambda}[l];l)]
=
\prod_{i=1}^{M-(a+c)} \int_{C_1} \underline{dw_i} 
\int_{C[(\ol{\lambda}[0];0), (\ol{\lambda}[l-1];l-1), (\ol{\lambda}[l];l)]}
\underline{dz}
\notag
\\
&\times
\left(
\prod_{s=1}^{l(\ol{\lambda}[0])}\int_{C_{\delta_0}}\underline{dz_{0,s}}
\prod_{s=1}^{l(\ol{\lambda}[l-1])}\int_{C_1}\ul{dz_{l-1,s}}
\prod_{s=1}^{l(\ol{\lambda}[l])}\int_{C_{\delta_l}}\ul{dz_{l,s}}
-
\prod_{s=1}^{l(\ol{\lambda}[0])}\int_{C_{\delta_0 q^{-1}}}\underline{dz_{0,s}}
\prod_{s=1}^{l(\ol{\lambda}[l-1])}\int_{C_q}\ul{dz_{l-1;s}}
\prod_{s=1}^{l(\ol{\lambda}[l])}\int_{C_{\delta_l q}}
\ul{dz_{l,s}}
\right)\notag\\
&\times \overset{\curvearrowright}{\prod_{1\leq i \leq M-(a+c)}}
{\bT}_1(w_i) {\bf S}[(\ol{\lambda}[0];0), (\ol{\lambda}[l-1];l-1), (\ol{\lambda}[l];l)]
\notag\\
&\times
\ol{\omega}_M
[w_1,\ldots, w_{M-(a+c)}; 
(\ol{\lambda}[0];0), (\ol{\lambda}[l-1];l-1), (\ol{\lambda}[l];l)
]
\bh.\notag
\end{align}
For $\ol{\lambda}[1] \in \ol{\mc E}_1 \cup \boxed{\ol{\mc E}_1}$,
$\ol{\lambda}[l-1] \in \ol{\mc D}_{l-1}$ such that 
$|\ol{\lambda}[1]|=a\geq 1$,
$|\ol{\lambda}[l-1]|=c \geq 1$, 
$M \geq a+c \geq 2$, we set
\begin{align}
&
\Delta_{a+c}^{(M,N)}[(\ol{\lambda}[1];1), (\ol{\lambda}[l-1];l-1)]\notag
\\
=&\prod_{i=1}^{M-(a+c)}\int_{C_1}\ul{dw_i} \int_{C[(\ol{\lambda}[1];1)]}\ul{dz}
\prod_{s=1}^{l(\ol{\lambda}[1])}
\int_{C_1-C_{q^{-1}}}\ul{dz_{1,s}}  \overset{\curvearrowright}{\prod_{1\leq i \leq M-(a+c)}}
{\bT}_1(w_i) {\bf S}[(\ol{\lambda}[1];1), (\ol{\lambda}[l-1];l-1)]
\notag
\\
&\times
\ol{\omega}_M
[w_1,\ldots, w_{M-(a+c)}; 
(\ol{\lambda}[1];1), (\ol{\lambda}[l-1];l-1)]
\bh.\notag
\end{align}
For $\ol{\lambda}[1] \in \ol{\mc E}_1 \cup \boxed{\ol{\mc E}_1}$,
$\ol{\lambda}[l] \in \ol{\mc D}_{l}$ such that $|\ol{\lambda}[1]|=a\geq 1$,
$|\ol{\lambda}[l]|=c \geq 1$, $M \geq a+c \geq 2$,
we set
\begin{align}
&
\Delta_{a+c}^{(M,N)}[(\ol{\lambda}[1];1), (\ol{\lambda}[l];l)]\notag
\\
=&\prod_{i=1}^{M-(a+c)}\int_{C_1}\ul{dw_i} \int_{C[(\ol{\lambda}[1];1)]}\ul{dz}
\prod_{s=1}^{l(\ol{\lambda}[1])}
\int_{C_1-C_{q^{-1}}}\ul{dz_{1,s}}  \overset{\curvearrowright}{\prod_{1\leq i \leq M-(a+c)}}
{\bT}_1(w_i) {\bf S}[(\ol{\lambda}[1];1), (\ol{\lambda}[l];l)]
\notag
\\
&\times
\ol{\omega}_M
[w_1,\ldots, w_{M-(a+c)}; 
(\ol{\lambda}[1];1), (\ol{\lambda}[l];l)]
\bh.\notag
\end{align}
For $\ol{\lambda}[1] \in \ol{\mc E}_1 \cup \boxed{\ol{\mc E}_1}$,
$\ol{\lambda}[l-1] \in \ol{\mc D}_{l-1}$ $\ol{\lambda}[l] \in \ol{\mc D}_l$ such that $|\ol{\lambda}[1]|=a\geq 1$,
$|\ol{\lambda}[l-1]|+|\ol{\lambda}[l]|=c \geq 2$, 
$M \geq a+c \geq 3$, we set
\begin{align}
&
\Delta_{a+c}^{(M,N)}[(\ol{\lambda}[1];1), (\ol{\lambda}[l-1];l-1), (\ol{\lambda}[l];l)]\notag
\\
=&\prod_{i=1}^{M-(a+c)}\int_{C_1}\ul{dw_i} \int_{C[(\ol{\lambda}[1];1)]}\ul{dz}
\prod_{s=1}^{l(\ol{\lambda}[1])}
\int_{C_1-C_{q^{-1}}}\ul{dz_{1,s}}  \overset{\curvearrowright}{\prod_{1\leq i \leq M-(a+c)}}
{\bT}_1(w_i) {\bf S}[(\ol{\lambda}[1];1), (\ol{\lambda}[l-1];l-1), (\ol{\lambda}[l];l)]
\notag
\\
&\times
\ol{\omega}_M
[w_1,\ldots, w_{M-(a+c)}; 
(\ol{\lambda}[1];1), (\ol{\lambda}[l-1];l-1), (\ol{\lambda}[l];l)]
\bh.\notag
\end{align}
For $\ol{\lambda}[0] \in \ol{\mc E}_0$, $\ol{\lambda}[1] \in \ol{\mc E}_1 \cup \boxed{\ol{\mc E}_1}$,
$\ol{\lambda}[l-1] \in \ol{\mc D}_{l-1}$ such that $|\ol{\lambda}[0]|+|\ol{\lambda}[1]|=a\geq 2$,
$|\ol{\lambda}[l-1]|=c \geq 1$, 
$M \geq a+c \geq 3$, we set
\begin{align}
&
\Delta_{a+c}^{(M,N)}[(\ol{\lambda}[0];0), (\ol{\lambda}[1];1), (\ol{\lambda}[l-1];l-1)]\notag
\\
=&\prod_{i=1}^{M-(a+c)}\int_{C_1}\ul{dw_i} \int_{C[(\ol{\lambda}[1];1)]}\ul{dz}
\prod_{s=1}^{l(\ol{\lambda}[1])}
\int_{C_1-C_{q^{-1}}}\ul{dz_{1,s}}  \overset{\curvearrowright}{\prod_{1\leq i \leq M-(a+c)}}
{\bT}_1(w_i) {\bf S}[(\ol{\lambda}[0];0), (\ol{\lambda}[1];1), (\ol{\lambda}[l-1];l-1)]
\notag
\\
&\times
\ol{\omega}_M
[w_1,\ldots, w_{M-(a+c)}; 
(\ol{\lambda}[0];0), (\ol{\lambda}[1];1), (\ol{\lambda}[l-1];l-1)]
\bh.\notag
\end{align}
For $\ol{\lambda}[0] \in \ol{\mc E}_0$,
$\ol{\lambda}[1] \in \ol{\mc E}_1 \cup \boxed{\ol{\mc E}_1}$,
$\ol{\lambda}[l] \in \ol{\mc D}_{l}$ such that $|\ol{\lambda}[0]|+|\ol{\lambda}[1]|=a\geq 2$,
$|\ol{\lambda}[l]|=c \geq 1$, 
$M\geq a+c \geq 3$, we set
\begin{align}
&
\Delta_{a+c}^{(M,N)}[(\ol{\lambda}[0];0), (\ol{\lambda}[1];1), (\ol{\lambda}[l];l)]\notag
\\
=&\prod_{i=1}^{M-(a+c)}\int_{C_1}\ul{dw_i} \int_{C[(\ol{\lambda}[1];1)]}\ul{dz}
\prod_{s=1}^{l(\ol{\lambda}[1])}
\int_{C_1-C_{q^{-1}}}\ul{dz_{1,s}}  \overset{\curvearrowright}{\prod_{1\leq i \leq M-(a+c)}}
{\bT}_1(w_i) {\bf S}[(\ol{\lambda}[0];0), (\ol{\lambda}[1];1), (\ol{\lambda}[l];l)]
\notag
\\
&\times
\ol{\omega}_M
[w_1,\ldots, w_{M-(a+c)}; 
(\ol{\lambda}[0];0), (\ol{\lambda}[1];1), (\ol{\lambda}[l];l)]
\bh.\notag
\end{align}
For $\ol{\lambda}[0]\in \ol{\mc E}_0$, $\ol{\lambda}[1] \in \ol{\mc E}_1 \cup \boxed{\ol{\mc E}_1}$,
$\ol{\lambda}[l-1] \in \ol{\mc D}_{l-1}$\jm{,}
$\ol{\lambda}[l] \in \ol{\mc D}_l$ such that $|\ol{\lambda}[0]|+|\ol{\lambda}[1]|=a\geq 2$,
$|\ol{\lambda}[l-1]|+|\ol{\lambda}[l]|=c \geq 2$, 
$M \geq a+c \geq 4$, we set
\begin{align}
&
\Delta_{a+c}^{(M,N)}[(\ol{\lambda}[0]; 0), (\ol{\lambda}[1];1), (\ol{\lambda}[l-1];l-1), (\ol{\lambda}[l];l)]\notag
\\
=&\prod_{i=1}^{M-(a+c)}\int_{C_1}\ul{dw_i} \int_{C[(\ol{\lambda}[1];1)]}\ul{dz}
\prod_{s=1}^{l(\ol{\lambda}[1])}
\int_{C_1-C_{q^{-1}}}\ul{dz_{1,s}}  \overset{\curvearrowright}{\prod_{1\leq i \leq M-(a+c)}}
{\bT}_1(w_i) \notag\\
&\times
{\bf S}[(\ol{\lambda}[0];0), (\ol{\lambda}[1];1), (\ol{\lambda}[l-1];l-1), (\ol{\lambda}[l];l)]
\notag
\\
&\times
\ol{\omega}_M
[w_1,\ldots, w_{M-(a+c)}; 
(\ol{\lambda}[0];0), (\ol{\lambda}[1];1), (\ol{\lambda}[l-1];l-1), (\ol{\lambda}[l];l)]
\bh.\notag
\end{align}
For ${\lambda}[1] \in {\mc E}_1 \cup \boxed{{\mc E}_1}$,
${\lambda}[l-1] \in {\mc D}_{l-1}$ such that $|{\lambda}[1]|=a\geq 1$,
$|{\lambda}[l-1]|=c \geq 1$, 
$M \geq a+c \geq 2$, we set
\begin{align}
&
\Delta_{a+c}^{(M,N)}[({\lambda}[1];1), ({\lambda}[l-1];l-1)]\notag
\\
=&\prod_{i=a+c+1}^{M}\int_{C_1}\ul{dw_i} \int_{C[({\lambda}[1];1)]}\ul{dz}
\prod_{s=1}^{l({\lambda}[1])}
\int_{C_1-C_{q}}\ul{dz_{1,s}} 
{\bf S}[({\lambda}[1];1), ({\lambda}[l-1];l-1)]
 \overset{\curvearrowright}{\prod_{a+c+1\leq i \leq M}}{\bT}_1(w_i) 
\notag
\\
&\times
{\omega}_M
[(\lambda[1];1), (\lambda [l-1]; l-1); w_{a+c+1},\ldots, w_M]
\bh.\notag
\end{align}
For ${\lambda}[1] \in {\mc E}_1 \cup \boxed{{\mc E}_1}$,
${\lambda}[l] \in {\mc D}_{l}$ such that $|{\lambda}[1]|=a\geq 1$,
$|{\lambda}[l]|=c \geq 1$,
$M \geq a+c \geq 2$, we set
\begin{align}
&
\Delta_{a+c}^{(M,N)}[({\lambda}[1];1), ({\lambda}[l];l)]\notag
\\
=&\prod_{i=a+c+1}^{M}
\int_{C_1}\ul{dw_i} \int_{C[({\lambda}[1];1)]}\ul{dz}
\prod_{s=1}^{l({\lambda}[1])}
\int_{C_1-C_{q}}\ul{dz_{1,s}}
{\bf S}[({\lambda}[1];1), ({\lambda}[l];l)]
\overset{\curvearrowright}{\prod_{a+c+1\leq i \leq M}}
{\bT}_1(w_i) 
\notag
\\
&\times
{\omega}_M
[({\lambda}[1];1), ({\lambda}[l];l); w_{a+c+1},\ldots, w_M]
\bh.\notag
\end{align}
For ${\lambda}[1] \in {\mc E}_1 \cup \boxed{{\mc E}_1}$,
$\lambda [l-1] \in {\mc D}_{l-1}$,
${\lambda}[l] \in {\mc D}_l$ such that $|{\lambda}[1]|=a\geq 1$,
$|{\lambda}[l-1]|+|{\lambda}[l]|=c \geq 2$, 
$M \geq a+c \geq 3$, we set
\begin{align}
&
\Delta_{a+c}^{(M,N)}[({\lambda}[1];1), ({\lambda}[l-1];l-1), ({\lambda}[l];l)]\notag
\\
=&\prod_{i=a+c+1}^{M}\int_{C_1}\ul{dw_i} 
\int_{C[({\lambda}[1];1)]}\ul{dz}
\prod_{s=1}^{l({\lambda}[1])}
\int_{C_1-C_{q}}\ul{dz_{1,s}} 
{\bf S}[({\lambda}[1];1), ({\lambda}[l-1];l-1), ({\lambda}[l];l)]
\overset{\curvearrowright}{\prod_{a+c+1\leq i \leq M}}
{\bT}_1(w_i) 
\notag
\\
&\times
{\omega}_M
[ ({\lambda}[1];1), ({\lambda}[l-1];l-1), ({\lambda}[l];l);
w_{a+c+1},\ldots, w_{M}]
\bh.\notag
\end{align}
For ${\lambda}[0] \in {\mc E}_0$, ${\lambda}[1] \in {\mc E}_1 \cup \boxed{{\mc E}_1}$,
${\lambda}[l-1] \in {\mc D}_{l-1}$ such that $|{\lambda}[0]|+|{\lambda}[1]|=a\geq 2$,
$|{\lambda}[l-1]|=c \geq 1$, 
$M \geq a+c \geq 3$, we set
\begin{align}
&
\Delta_{a+c}^{(M,N)}[({\lambda}[0];0), ({\lambda}[1];1), ({\lambda}[l-1];l-1)]\notag
\\
=&\prod_{i=a+c+1}^{M}\int_{C_1}\ul{dw_i} 
\int_{C[({\lambda}[1];1)]}\ul{dz}
\prod_{s=1}^{l({\lambda}[1])}
\int_{C_1-C_{q}}\ul{dz_{1,s}} 
{\bf S}[({\lambda}[0];0), ({\lambda}[1];1), ({\lambda}[l-1];l-1)]
\overset{\curvearrowright}{\prod_{a+c+1\leq i \leq M}}
{\bT}_1(w_i) 
\notag
\\
&\times
{\omega}_M
[({\lambda}[0];0), ({\lambda}[1];1), ({\lambda}[l-1];l-1);
w_{a+c+1},\ldots, w_M]
\bh.\notag
\end{align}
For ${\lambda}[0] \in {\mc E}_0$,
${\lambda}[1] \in {\mc E}_1 \cup \boxed{{\mc E}_1}$,
${\lambda}[l] \in {\mc D}_{l}$ such that $|{\lambda}[0]|+|{\lambda}[1]|=a\geq 2$,
$|{\lambda}[l]|=c \geq 1$, 
$M \geq a+c \geq 3$, we set
\begin{align}
&
\Delta_{a+c}^{(M,N)}[({\lambda}[0];0), ({\lambda}[1];1), ({\lambda}[l];l)]\notag
\\
=&\prod_{i=a+c+1}^{M}\int_{C_1}\ul{dw_i} \int_{C[({\lambda}[1];1)]}\ul{dz}
\prod_{s=1}^{l({\lambda}[1])}
\int_{C_1-C_{q}}\ul{dz_{1,s}}
{\bf S}[({\lambda}[0];0), ({\lambda}[1];1), ({\lambda}[l];l)]
\overset{\curvearrowright}{\prod_{a+c+1\leq i \leq M}}
{\bT}_1(w_i) 
\notag
\\
&\times
{\omega}_M
[ ({\lambda}[0];0), ({\lambda}[1];1), ({\lambda}[l];l);
w_{a+c+1},\ldots, w_M]
\bh.\notag
\end{align}
For ${\lambda}[0]\in {\mc E}_0$, ${\lambda}[1] \in {\mc E}_1 \cup \boxed{{\mc E}_1}$,
${\lambda}[l-1] \in {\mc D}_{l-1}$, ${\lambda}[l] \in {\mc D}_l$ such that $|{\lambda}[0]|+|{\lambda}[1]|=a\geq 2$,
$|{\lambda}[l-1]|+|{\lambda}[l]|=c \geq 2$, 
$M \geq a+c \geq 4$, we set
\begin{align}
&
\Delta_{a+c}^{(M,N)}[({\lambda}[0]), ({\lambda}[1];1), ({\lambda}[l-1];l-1), ({\lambda}[l];l)]\notag
\\
=&\prod_{i=a+c+1}^{M}\int_{C_1}\ul{dw_i} 
\int_{C[({\lambda}[1];1)]}\ul{dz}
\prod_{s=1}^{l({\lambda}[1])}
\int_{C_1-C_{q}}\ul{dz_{1,s}}  \notag
\\
&\times
{\bf S}[({\lambda}[0];0), ({\lambda}[1];1), ({\lambda}[l-1];l-1), ({\lambda}[l];l)]
\overset{\curvearrowright}{\prod_{a+c+1\leq i \leq M}}
{\bT}_1(w_i) \notag
\\
&\times
{\omega}_M
[ ({\lambda}[0];0), ({\lambda}[1];1), ({\lambda}[l-1];l-1), ({\lambda}[l];l);
w_{a+c+1},\ldots, w_M]
\bh.\notag
\end{align}

Here $C[(\lambda; j)]$, $\lambda \in \ol{\mc F}_j \cup {\mc F}_j, 0\leq j \leq l$, stands for the contour
\begin{align}
C[(\lambda; j)]
=C[(j,1), (j,2),\ldots, (j, l(\lambda))]\,,
\notag
\end{align}
and $C[(\lambda_1; j_1), \ldots, (\lambda_n; j_n)]$, 
$\lambda_s \in \ol{\mc F}_{j_s} \cup {\mc F}_{j_s}, 1 \leq s \leq n, 0\leq j_1<j_2<\cdots <j_n \leq l$, 
stands for the contour
\begin{align}
C[(\lambda_1; j_1),\ldots, (\lambda_n;j_n)]
=C[(j_1,1), (j_1,2),\ldots, (j_1, l(\lambda_1)),\ldots,
(j_n,1), (j_n, 2), \ldots, (j_n,l(\lambda_n))].
\notag
\end{align}

\vspace{5mm}

Proposition \ref{prop:duality1} for $0\leq a,c \leq M$, 
$0\leq b \leq M$,
$0\leq a+b+c \leq M$ is a direct consequence of Lemmas
\ref{lem:IM-GN:1},
\ref{lem:IM-GN:2},
\ref{lem:IM-GN:3},
\ref{lem:IM-GN:4},
\ref{lem:IM-GN:5},
\ref{lem:IM-GN:6}
and
\ref{lem:IM-GN:7} given below.
\begin{lem}
\label{lem:IM-GN:1}
The following hold for $1\leq b \leq M$, $0 \leq a,c \leq M$, $1\leq a+b+c \leq M$.
\begin{align}
{\rm LHS}_{M,N}(a,b,c)=0,\quad
{\rm RHS}_{M,N}(a,b,c)=0.\notag
\end{align}
\end{lem}
\begin{proof}
First, we consider $\mathrm{LHS}_{M,N}(a,b,c)$.
As can be seen from the formula for $b \geq 1$,
all terms have the integral 
$\int_{C_1-C_{q^{-1}}}\ul{dz_{j_1,1}}$, $2 \leq j_1 \leq l-2$.
As can be seen from the normal ordering, there is no pole of $z_{j_1,1}$ in the region of
the integral 
$\int_{C_1-C_{q^{-1}}}\ul{dz_{j_1,1}}$, $2 \leq j_1 \leq l-2$.
Hence, we obtain $\mathrm{LHS}_{M,N}(a,b,c)=0$.
Next, we consider $\mathrm{RHS}_{M,N}(a,b,c)$.
As can be seen from the formula for $b \geq 1$,
all terms have the integral 
$\int_{C_1-C_{q}}\ul{dz_{j_1,1}}$, $2 \leq j_1 \leq l-2$.
As can be seen from the normal ordering, there is no pole of $z_{j_1,1}$ in the region of
the integral 
$\int_{C_1-C_{q}}\ul{dz_{j_1,1}}$, $2 \leq j_1 \leq l-2$.
Hence, we obtain $\mathrm{RHS}_{M,N}(a,b,c)=0$.
\end{proof}

\hspace{2mm}

The relation $\mathrm{LHS}_{M,N}(0,0,c)=\mathrm{RHS}_{M,N}(0,0,c)$, $1\leq c \leq M$, in Proposition \ref{prop:duality1}
is a direct consequence of Lemma \ref{lem:IM-GN:2}.
\begin{lem}
\label{lem:IM-GN:2}
The following holds for $M \geq m \geq 3$.
\begin{align}
&
\sum_{\lambda \in {\mc D}_{l-1}
\atop{|\lambda|=m}}
\com{N}{\lambda}
\nu_{\lambda;l-1}
\Delta_m^{(M,N)}[(\lambda;l-1)]
+
\sum_{\mu \in {\mc D}_{l}
\atop{|\mu|=m}}
\com{N}{\mu}
\nu_{\mu;l}
\Delta_m^{(M,N)}[(\mu;l)]\notag
\\
&
+\sum_{\lambda \in {\mc D}_{l-1}, \mu \in {\mc D}_l
\atop{|\lambda|+|\mu|=m}}
\com{N}{\lambda}\com{N}{\mu}
\nu_{\lambda;l-1}\nu_{\mu;l}
\Delta_m^{(M,N)}[(\lambda;l-1), (\mu;l)]=0.
\notag
\end{align}
\end{lem}

\begin{proof}

First, we prepare some identities.
\begin{align}
&\int_{C_{q^{-1}}-C_1}\ul{dz_{l-1}}~
{\mc A}_{k;l-1}(z_{l-1})S_l(z_{l})
\notag\\
=&\delta_{k,l-1}(q_1-1)
:\bLambda_{l-1}(q^{\frac{1}{2}}z_{l})S_{l-1}(z_{l})S_l(z_{l}):,\quad |z_l|=\delta_l>1, \quad k=l-1, \bar{l},\notag
\end{align}
\begin{align}
&\int_{C_{q^{-1}}-C_1}\ul{dz_{l-1}}~
{\mc A}_{k_1;l-1}(z_{l-1}){\mc A}_{k_2;l}(z_{l})=0,\quad
|z_l|=\delta_l>1,\quad k_1=l-1, \bar{l}, \quad k_2=l-1,l,
\notag
\end{align}
\begin{align}
&\int_{C_{\delta_l q^{-1}}-C_{\delta_l}}
\ul{dz_{l}}~
S_{l-1}(z_{l-1}){\mc A}_{k;l}(z_{l})
\notag\\
=&\delta_{k,l}(q_1^{-1}-1) 
:\bLambda_{l}(q^{-\frac{1}{2}}z_{l-1})S_{l-1}(z_{l-1})S_l(q^{-1}z_{l-1}):,
\quad |z_{l-1}|=1, \quad k=l-1, l,\notag
\end{align}
\begin{align}
&\int_{C_{\delta_l q^{-1}}-C_{\delta_l}}\ul{dz_{l}}~
{\mc A}_{k_1;l-1}(z_{l-1}){\mc A}_{k_2;l}(z_{l})=0,\quad
|z_{l-1}|=1,\quad k_1=l-1, \bar{l}, \quad k_2=l-1,l.
\notag
\end{align}

We consider 
$\Delta_m^{(M,N)}[(k_1; l-1, 1), (k_{2}; l-1,2), \ldots, (k_m; l-1, m)]$.
We rewrite the integral as follows.
\begin{align}
&
\prod_{s=1}^{l(\lambda[l-1])}
\int_{C_{q^{-1}}}\ul{dz_{l-1,s}}
-
\prod_{s=1}^{l(\lambda[l-1])}\int_{C_1}\ul{dz_{l-1,s}}
=
\sum_{n=1}^{l(\lambda[l-1])}
\sum_{S \subset \{1,2,\ldots,l(\lambda[l-1])\}
\atop{|S|=n}}
\prod_{s=1
\atop{s \in S}}^{l(\lambda[l-1])} \int_{C_{q^{-1}}-C_1} \ul{dz_{l-1, s}}
\prod_{s=1 \atop{s \notin S}}^{l(\lambda[l-1])} \int_{C_1}\ul{dz_{l-1,s}}.
\notag 
\end{align}
Let $(k_i)_{1\leq i \leq l(\lambda[l-1])}$ be the fillings of $\lambda[l-1] \in {\mc D}_{l-1}$.
Set $l(\lambda[l-1])=m$.
Taking the residues we obtain 
\begin{align}
&\Delta_m^{(M,N)}[(k_1; l-1,1), (k_{2}; l-1,2),\dots, (k_m; l-1,m)]\notag\\
=&
\sum_{n=1}^m 
\comb{N}{n} n!
\sum_{S \subset\{1,2,\ldots,m\}
\atop{|S|=n}}
\prod_{s=1
\atop{s \in S}}^m \delta_{k_s, l-1} (q_1-1)^n
\prod_{i=m+1}^M \int_{C_1} \ul{dw_i}
\int_{C[(l-1,s_1), (l-1,s_2), \ldots, (l-1,s_n)]}
\underline{dz} \notag
\\
&
\times
\mathbf{S}_{[0,l-2]}
\overset{\curvearrowright}{\prod_{1\leq s \leq n}}
:\bLambda_{l-1}(q^{\frac{1}{2}}z_{l,s})S_{l-1}(z_{l,s})S_l(z_{l,s}):
\overset{\curvearrowright}{\prod_{1\leq s \leq m
\atop{s \notin S}}}
{\mc A}_{k_{s}; l-1}[z_{l-1,s}]
\overset{\curvearrowright}{\prod_{m+1\leq s \leq N}}
S_{l-1}(z_{l-1,s})
\notag\\
&\times
\overset{\curvearrowright}{\prod_{n+1\leq s \leq N}}S_l(z_{l,s})
\overset{\curvearrowright}{\prod_{m+1\leq i \leq M}}{\bT}_1(w_i)
\left.\omega_M[(\lambda[l-1]; l-1); w_{m+1},\ldots, w_M]\bh
\right|_{z_{l-1,s_t}=z_{l,t} (1\leq t \leq n)}.\notag
\end{align}
Here $s_1,s_2,\ldots, s_n$ are given by
$S=\{s_1,s_2,\ldots, s_n
\left|1\leq s_1<s_2<\cdots<s_n \leq m\right.\}$.

\vspace{2mm}

We consider 
$\Delta_m^{(M,N)}[(l_1; l, 1),(l_2; l, 2), \ldots, (l_m; l, m)]$.
We rewrite the integral as follows.
\begin{align}
&
\prod_{t=1}^{l(\lambda[l])}
\int_{C_{\delta_l q^{-1}}}\ul{dz_{l,t}}
-
\prod_{t=1}^{l(\lambda[l])}\int_{C_{\delta_l}}\ul{dz_{l,t}}
=
\sum_{n=1}^{l(\lambda[l])}
\sum_{T \subset \{1,2,\ldots,l(\lambda[l])\}
\atop{|T|=n}}
\prod_{t=1
\atop{t \in T}}^{l(\lambda[l])} \int_{C_{\delta_l q^{-1}}-C_{\delta_l}} \ul{dz_{l, t}}
\prod_{t=1 \atop{t \notin T}}^{l(\lambda[l])} \int_{C_{\delta_l}}\ul{dz_{l, t}}.
\notag 
\end{align}

Let $(l_i)_{1\leq i \leq l(\lambda[l])}$ be the fillings of $\lambda[l] \in {\mc D}_{l}$.
Set $l(\lambda[l])=m$.
Taking the residues we obtain 
\begin{align}
&
\Delta_m^{(M,N)}[(l_1; l, 1),(l_2; l, 2), \ldots, (l_m; l, m)]
\notag\\
=&
\sum_{n=1}^m
\comb{N}{n} n!
\sum_{T \subset \{1,2,\ldots, m\}
\atop{|T|=n}}  
\prod_{t=1
\atop{t \in T}}^m \delta_{l_t, l} (q_1^{-1}-1)^n
\prod_{i=m+1}^M
\int_{C_1}
\ul{dw_i}
\int_{C[(l,t_1),(l,t_2), \ldots, (l,t_n)]}
\underline{dz} \notag\\
&\times \mathbf{S}_{[0,l-2]}
\overset{\curvearrowright}{\prod_{n+1\leq s \leq N}}
S_{l-1}(z_{l-1,s})
\overset{\curvearrowright}{\prod_{1\leq s \leq n}}
:\bLambda_{l}(q^{-\frac{1}{2}}z_{l-1,s})
S_{l-1}(z_{l-1,s})S_l(q^{-1}z_{l-1,s}):
\notag\\
&\times
\overset{\curvearrowright}
{\prod_{1\leq s \leq m
\atop{s \notin T}}}
{\mc A}_{l_{s}; l}[z_{l,s}]
\overset{\curvearrowright}{\prod_{m+1\leq s \leq N}}S_l(z_{l,s})
\overset{\curvearrowright}{\prod_{m+1\leq i \leq M}}{\bT}_1(w_i)\notag
\\
&\times
\left.\omega_M[(\lambda[l]; l); w_{m+1},\ldots, w_M]\bh
\right|_{z_{l, t_s}=q^{-1} z_{l-1,s} (1\leq s \leq n)}.\notag
\end{align}
Here $t_1,t_2,\ldots, t_n$ are given by
$T=\{t_1,t_2,\ldots, t_n
\left|1\leq t_1<t_2<\cdots<t_n \leq m\right.\}$.

We consider
$$\Delta_{m_1+m_2}^{(M,N)}[(k_1;l-1,1), (k_{2};l-1,2), \ldots, (k_{m_1};l-1,m_1),
(l_1;l,1),(l_2;l,2),\ldots,(l_{m_2};l,m_2)].$$
We rewrite the integral as follows.
\begin{align}
&
\prod_{s=1}^{l(\lambda[l-1])}
\int_{C_{q^{-1}}}\ul{dz_{l-1,s}}
\prod_{t=1}^{l(\lambda[l])}
\int_{C_{\delta_l q^{-1}}}\ul{dz_{l,t}}
-
\prod_{s=1}^{l(\lambda[l-1])}\int_{C_{\delta_l}}\ul{dz_{l-1,s}}
\prod_{t=1}^{l(\lambda[l])}\int_{C_{\delta_l}}\ul{dz_{l,t}}
=
\sum_{n_1=1}^{l(\lambda[l-1])}
\sum_{n_2=1}^{l(\lambda[l])}
\notag\\
&\times
\sum_{S \subset \{1,2,\ldots,l(\lambda[l-1])\}
\atop{|S|=n_1}}
\sum_{T \subset \{1,2,\ldots,l(\lambda[l])\}
\atop{|T|=n_2}}
\prod_{s=1
\atop{s \in S}}^{l(\lambda[l-1])} \int_{C_{q^{-1}}-C_{1}} \ul{dz_{l-1, s}}
\prod_{s=1 \atop{s \notin S}}^{l(\lambda[l-1])} 
\int_{C_{1}}\ul{dz_{l-1, s}}
\prod_{t=1
\atop{t \in T}}^{l(\lambda[l])} \int_{C_{\delta_l q^{-1}}-C_{\delta_l}} \ul{dz_{l, t}}
\prod_{t=1 \atop{t \notin T}}^{l(\lambda[l])} \int_{C_{\delta_l}}\ul{dz_{l, t}}.
\notag 
\end{align}
 Taking the residues we obtain 
 \begin{align}
&
 \Delta_{m_1+m_2}^{(M,N)}
 [(k_1;l-1,1), (k_{2};l-1,2), \ldots, (k_{m_1};l-1,m_1), (l_1;l,1),(l_2;l,2),\ldots,(l_{m_2};l,m_2)]\notag\\
=&
\sum_{n_1=1}^{m_1}
\sum_{n_2=1}^{m_2}
\comb{N-m_2}{n_1} n_1!
\comb{N-m_1}{n_2} n_2!
\sum_{S \subset \{1,2,\ldots,l(\lambda[l-1])\}
\atop{|S|=n_1}}
\sum_{T \subset \{1,2,\ldots,l(\lambda[l])\}
\atop{|T|=n_2}}
\prod_{s=1\atop{s \in S}}^{m_1}\delta_{k_s, l-1}
\prod_{t=1\atop{t \in T}}^{m_2} \delta_{l_t,l}
\notag\\
&\times
(q_1-1)^{n_1}(q_1^{-1}-1)^{n_2}
\prod_{i=m_1+m_2+1}^M
\int_{C_1} \ul{dw_i}
\int_{C[(l-1,s_1),(l-1,s_2),\ldots,(l-1,s_{n_1}), 
(l,t_1), (l,t_2),\ldots, (l,t_{n_2})]}\ul{dz}
~\mathbf{S}_{[0,l-2]} 
\notag\\
&\times
\overset{\curvearrowright}{\prod_{1\leq s \leq n_1}}
:\bLambda_{l-1}(q^{\frac{1}{2}}z_{l,m_2+s})S_{l-1}(z_{l,m_2+s})S_l(z_{l,m_2+s}):
\overset{\curvearrowright}{\prod_{1\leq s \leq m_1
\atop{s \notin S}}}
{\mc A}_{k_{s}; l-1}[z_{l-1,s}]
\overset{\curvearrowright}{\prod_{m_1+n_2+1\leq s \leq N}}
S_{l-1}(z_{l-1,s})
\notag\\
&\times
\overset{\curvearrowright}{\prod_{m_1+1\leq s \leq m_1+n_2}}
:\bLambda_{l}(q^{-\frac{1}{2}}z_{l-1,s})
S_{l-1}(z_{l-1,s})S_l(q^{-1}z_{l-1,s}):
\notag\\
&\times
\overset{\curvearrowright}
{\prod_{1\leq s \leq m_2
\atop{s \notin T}}}
{\mc A}_{l_{s}; l}[z_{l,s}]
\overset{\curvearrowright}{\prod_{m_2+n_1+1\leq s \leq N}}S_l(z_{l,s})
\overset{\curvearrowright}{\prod_{m_1+m_2+1\leq i \leq M}}{\bT}_1(w_i)
\notag
\\
&\times
\omega_M[(\lambda[l-1]; l-1), (\lambda[l]; l); w_{m_1+m_2+1},\ldots, w_M]
\bh
\left|_{
\begin{array}{c}
{\scriptstyle  z_{l-1,s_t}=z_{l, m_2+t}~(1\leq t \leq n_1)}\\
{\scriptstyle  z_{l,t_s}=q^{-1}z_{l-1,m_1+s}~(1\leq s \leq n_2)}
\end{array}}
\right..
\notag
 \end{align}
 Here $s_1, s_2, \ldots, s_{n_1}$ and $t_1,t_2,\ldots, t_{n_2}$ are given by
$S=\{s_1, s_2,\ldots, s_{n_1}
\left|1\leq s_1<s_2<\cdots< s_{n_1} \leq m_1\right.\}$
and
$T=\{t_1,t_2,\ldots, t_{n_2}
\left|1\leq t_1<t_2<\cdots<t_{n_2} \leq m_2\right.\}$,
respectively.

Adding up the above and using
\begin{align}
q_1:\bLambda_{l-1}(q^{\frac{1}{2}}z)S_{l-1}(z)S_l(z):=:\bLambda_l(q^{\frac{1}{2}}z)S_{l-1}(qz)S_l(z):\notag
\end{align}
we find
\begin{align}
&
\sum_{\lambda \in {\mc D}_{l-1}
\atop{l(\lambda)=m}}
\comb{N}{m}
\nu_{\lambda;l-1}
\Delta_m^{(M,N)}[(\lambda;l-1)]
+
\sum_{\mu \in {\mc D}_{l}
\atop{l(\mu)=m}}
\comb{N}{m}
\nu_{\mu;l}
\Delta_m^{(M,N)}[(\mu;l)]\notag
\\
&
+
\sum_{1\leq m_1, m_2 \leq m-1\atop{m_1+m_2=m}}
\sum_{\lambda \in {\mc D}_{l-1}, \mu \in {\mc D}_l
\atop{l(\lambda)=m_1, l(\mu)=m_2}}
\comb{N}{m_1}\comb{N}{m_2}
\nu_{\lambda; l-1} \nu_{\mu; l}
\Delta_{m_1+m_2}^{(M,N)}[(\lambda;l-1), (\mu;l)]=0,\quad
2 \leq m \leq M.
\notag
\end{align}
Here we used
$\com{N}{\lambda}=\comb{N}{m}$, $\lambda=(1^{m})$.

\end{proof}

The relation $\mathrm{LHS}_{M,N}(a,0,0)=\mathrm{RHS}_{M,N}(a,0,0)$, $1\leq a \leq M$, in Proposition \ref{prop:duality1}
is a direct consequence of Lemmas \ref{lem:IM-GN:3}, \ref{lem:IM-GN:4},  and \ref{lem:IM-GN:5} given below.
\begin{lem}
\label{lem:IM-GN:3}
The following hold for $M \geq m \geq 3$.
\begin{align}
&
\Delta^{(M,N)}_m[(\lambda;1)]=0,\quad
\lambda \in \boxed{{\mathcal E}_1}, ~|\lambda |=m,
\notag
\\
&
\Delta^{(M,N)}_m[(\overline{\lambda};1)]=0,\quad
\overline{\lambda} \in \boxed{\overline{\mathcal E}_1}, ~|\overline{\lambda}|=m,
\notag
\\
&
\Delta^{(M,N)}_m[(\mu;0),(\lambda;1)]=0,\quad
\mu \in {\mathcal E}_0, \lambda \in \boxed{{\mathcal E}_1}, 
~|\mu|+|\lambda|=m,
\notag
\\
&
\Delta^{(M,N)}_m[(\overline{\mu};0),(\overline{\lambda};1)]=0,
\quad
\overline{\mu} \in \overline{\mathcal E}_0, \overline{\lambda} \in \boxed{\overline{\mathcal E}_1},
~
|\overline{\mu}|+|\overline{\lambda} |=m.
\notag
\end{align}
\end{lem}
\begin{proof}
We can show them in the same way as Lemma \ref{lem:duality1}.
\end{proof}

\vspace{2mm}

\begin{lem}
\label{lem:IM-GN:4}
The following holds for $M\geq m \geq 3$.
\begin{align}
\sum_{\lambda \in {\mathcal E}_1
\atop{|\lambda|=m
}}
\com{N}{\lambda} \nu_{\lambda;1}
\Delta^{(M,N)}_m [(\lambda;1)]
+\sum_{
\mu \in {\mathcal E}_0,
\lambda \in {\mathcal E}_1
\atop{|\mu|+|\lambda|=m}}
\com{N}{\mu} \com{N}{\lambda} \nu_{\mu;0}\nu_{\lambda;1}
\Delta^{(M,N)}_m[(\mu;0),(\lambda;1)]=0.\notag
\end{align}
\end{lem}

\begin{proof}

We prepare the calculation rules
$(\alpha 1), (\alpha 2), (\alpha 3)$, and $(\alpha 4)$.
Let
\begin{align}
&
I^{k_a,\ldots, k_2, k_1;0}_{l_b,\ldots,l_2,2;1}=\int_{C_1-C_{q}}\ul{dz_1}{\mc A}_{k_a,\ldots, k_1;0}[z_0]
{\mc A}_{l_b,\ldots,l_2,2;1}[z_1]\prod_{i=1}^a \prod_{j=1}^b \omega(q^{j-i}z_1/z_0)f(z_1),\notag
\\
&a \geq 0, b \geq 1, \quad 1\preceq k_a \preceq \cdots \preceq k_1 \preceq 2,\quad
1\preceq l_b \preceq \cdots \preceq l_2 \preceq l_1=2.\notag
\end{align}
Here $f(z)$ is a holomorphic function in the region $|q|<|z|<1$.
If $a=0$, it is understood
${\mc A}_{k_a,\ldots,k_1;0}[z]=S_0(z)$ and $I_{l_b,\ldots,l_2,2;1}^{k_0;0}=\int_{C_1-C_q}\ul{dz_1}S_0(z_0){\mc A}_{l_b,\ldots,l_2,2;1}(z_1)f(z_1)$.
\vspace{2mm}~\\
$(\alpha 1)$
If $a=0, b=1$,
as can be seen from the normal ordering of 
$S_0(z_0){\mc A}_{2;1}[z_1]$,
there is no pole in the region of the integral
$\int_{C_1-C_q}\ul{dz_{1}}$.
Hence, we obtain
\begin{align}
I_{2;1}^{k_0;0}=0.\notag
\end{align}
$(\alpha 2)$
If $a \geq 1$, $b=1$,
as can be seen from the normal ordering of 
${\mc A}_{k_a,\ldots, k_1;0}[z_0]{\mc A}_{2;1}[z_1]$
$\times \prod_{i=1}^a \omega(q^{1-i}z_1/z_0)$,
there is a pole at $z_1=qz_0$ in the region of the integral
$\int_{C_1-C_q}\ul{dz_{1}}$.
Taking the residue we obtain
\begin{align}
I_{2;1}^{k_a,\ldots, k_2, k_1;0}
=Const. :{\mc A}_{k_a,\ldots, k_1;0}[z_0]{\mc A}_{2;1}[qz_0]:f(qz_0),
\quad
a \geq 1, 1\preceq k_a \preceq \cdots \preceq k_1 \preceq 2.\notag
\end{align}
$(\alpha 3)$
If $a=0$, $b \geq 2$, as can be seen from the normal ordering of
$S_0(z_0){\mc A}_{l_b,\ldots, l_2, 2;1}[z_1]$,
there is a pole at $z_1=qz_0$ in the region of the integral
$\int_{C_1-C_q}\ul{dz_{1}}$.
Taking the residue we obtain
\begin{align}
I_{l_b,\ldots,l_2,2;1}^{k_0;0}=
Const. :S_0(z_0){\mc A}_{l_b,\ldots, l_2, 2;1}[qz_0]:f(q z_0),
\quad
b \geq 2,\quad 1\preceq l_b \preceq \cdots \preceq l_2 \preceq l_1=2
.\notag
\end{align}
$(\alpha 4)$
If $a \geq 1, b\geq 2$, as can be seen from the normal ordering of 
${\mc A}_{k_a,\ldots, k_1;0}[z_0]
{\mc A}_{l_b,\ldots, l_2, 2;1}[z_1]$\\
$\times \prod_{i=1}^a \prod_{j=1}^b \omega(q^{j-i}z_1/z_0)$,
there is no pole in the region of the integral
$\int_{C_1-C_q}\ul{dz_1}$.
We obtain
\begin{align}
I^{k_a,\ldots, k_1;0}_
{l_b,\ldots, l_2, 2;1}=0,\quad
a \geq 1, b \geq 2, \quad
1\preceq k_a \preceq \cdots \preceq k_1 \preceq 2,\quad
1\preceq l_b \preceq \cdots \preceq l_2 \preceq l_1=2.
\notag
\notag
\end{align}

\vspace{2mm}

Let $(\mu, (k_{i,j})) \in {\mc E}_0$ and $(\lambda, (l_{i,j})) \in {\mc E}_1$. 
Set $n(\lambda)=n\left(\{i|\lambda_i \geq 2 \}\right)$ for a Young diagram $\lambda$.
We consider $\Delta_m^{(M,N)}[(\lambda;1)]$ with $(\lambda, (l_{ij})) \in {\mc E}_1$ such that $n(\lambda)<l(\lambda)$.
As can be seen from $(\alpha 1)$, there is no pole in the region of the integral 
$\int_{C_1-C_q}\ul{dz_{1,l(\lambda)}}$.
We obtain
\begin{align}
\Delta_m^{(M,N)}[(\lambda;1)]=0,\quad n(\lambda)<l(\lambda).
\label{lem:IM-GN:4:eq0}
\end{align} 

We consider Lemma \ref{lem:IM-GN:4} for the case $m=3$.

We consider $\Delta_3^{(M,N)}[(\lambda;1)]$ with $|\lambda|=3$.
Using $(\ref{lem:IM-GN:4:eq0})$ we obtain
\begin{align}
&
\Delta_3^{(M,N)}[(2;1,1), (2;1,2), (2;1,3)]=0,\notag
\\
&
\Delta_3^{(M,N)}[(l_{1 2},2;1,1), (2;1,2)]=0,\quad l_{1 2}=1,2.\notag
\end{align}
Using $(\alpha 3)$ 
to $\Delta_3^{(M,N)}[(l_{1 3}, l_{1 2},2;1,1)]$
we obtain
\begin{align}
\label{lem:IM-GN:4:eq1}
:{\mc A}_{l_{13}, l_{12}, 2;1}[qz_{0,s}]S_0(z_{0,s}):,\quad
1 \preceq l_{1 3} \preceq l_{1 2} \preceq 2,
\end{align}
as a core part of the residue.
In what follows we sometimes omit repeating coefficients and write only core parts of the residues.

We consider $\Delta_3^{(M,N)}[(\mu; 0), (\lambda;1)]$ with $|\mu|=1, |\lambda |=2$.
Using $(\alpha 3)$ and $(\alpha 4)$ to
$\Delta_3^{(M,N)}[(k_{1 1};0,1),$ $(l_{1 2}, 2;1,1)]$
we obtain
\begin{align}
:{\mc A}_{k_{1 1};0}[z_{0,s}]: \times :{\mc A}_{l_{1 2},2; 1}[qz_{0,t}]S_0(z_{0,t}):, \quad 1 \preceq k_{11}, l_{12} \preceq 2.
\label{lem:IM-GN:4:eq2}
\end{align}
Using $(\alpha 1)$ and $(\alpha 2)$ to $\Delta_3^{(M,N)}[(k_{11};0,1), (2;1,1), (2;1,2)]$ 
we obtain
\begin{align}
\Delta_3^{(M,N)}[(k_{11};0,1), (2;1,1), (2;1,2)]=0,\quad 1 \preceq k_{11} \preceq 2.
\notag
\end{align}

We consider $\Delta_3^{(M,N)}[(\mu; 0), (\lambda;1)]$ with $|\mu|=2, |\lambda |=1$.
Using $(\alpha 1)$ and $(\alpha 2)$ to
$\Delta_3^{(M,N)}[(k_{12}, k_{11};0,1),$ $(2;1,1)]$
we obtain
\begin{align}
:{\mc A}_{k_{12}, k_{11}; 0}[z_{0,s}] {\mc A}_{2;1}[qz_{0,s}]:,\quad
1 \preceq k_{1 2} \preceq k_{1 1} \preceq 2.
\label{lem:IM-GN:4:eq3}
\end{align}
Using $(\alpha 1)$ and $(\alpha 2)$ to
$\Delta_3^{(M,N)}[(k_{11};0,1), (k_{21};0,2), (2;1,1)]$
we obtain
\begin{align}
&:{\mc A}_{k_{11};0}[z_{0,s}]: \times :{\mc A}_{k_{21}; 0}[z_{0,t}] {\mc A}_{2;1}[q z_{0,t}]:,\quad 1 \preceq k_{11}, k_{21} \preceq 2.
\label{lem:IM-GN:4:eq4}
\end{align}

The operators (\ref{lem:IM-GN:4:eq1}) and (\ref{lem:IM-GN:4:eq3}) (resp. (\ref{lem:IM-GN:4:eq2}) and (\ref{lem:IM-GN:4:eq4})) 
are the same. They appear exactly twice.
To summarize, in $\Delta_3^{(M,N)}[(\lambda;1)]$ and $\Delta_3^{(M,N)}[(\mu;0), (\lambda;1)]$,
core parts of the residues are given as follows.
\begin{align}
&:\bLambda_{l_2}(q^{\frac{\beta}{2}(l-2)-\frac{3}{2}}z_{0,s})
\bLambda_{l_1}(q^{\frac{\beta}{2}(l-2)-\frac{1}{2}}z_{0,s})
\bLambda_2(q^{\frac{\beta}{2}(l-2)+\frac{1}{2}}z_{0,s})
S_0(z_{0,s})S_1(q z_{0,s}):,\quad
1 \preceq l_{2} \preceq l_{1} \preceq 2, \notag
\\
&:\bLambda_{l_2}(q^{\frac{\beta}{2}(l-2)-\frac{1}{2}}z_{0,s})S_0(z_{0,s}):\times
:\bLambda_{l_1}(q^{\frac{\beta}{2}(l-2)-\frac{1}{2}}z_{0,t})
\bLambda_{2}(q^{\frac{\beta}{2}(l-2)+\frac{1}{2}})S_0(z_{0,t})S_1(qz_{0,t}):,\quad
1\preceq l_1, l_2 \preceq 2.\notag
\end{align}
The residues corresponding  to (\ref{lem:IM-GN:4:eq1}) and (\ref{lem:IM-GN:4:eq3}) 
(resp. (\ref{lem:IM-GN:4:eq2}) and (\ref{lem:IM-GN:4:eq4})) 
have opposite signs. Their sum becomes zero.
Hence, we obtain Lemma \ref{lem:IM-GN:4} in the case $m=3$.

\vspace{2mm}

We consider Lemma \ref{lem:IM-GN:4} for the case $m=4$.
We consider $\Delta_4^{(M,N)}[(\lambda;1)]$ with $|\lambda|=4$.

Using $(\alpha 3)$ and $(\ref{lem:IM-GN:4:eq0})$ 
we obtain
\begin{align}
&\Delta_4^{(M,N)}[(l_{12}, l_{11},2;1,1), (2;1,2)]=0,\quad 1\preceq l_{12} \preceq l_{11} \preceq 2,\notag
\\
&\Delta_4^{(M,N)}[(l_{12},2;1,1), (2;1,2), (2;1,3)]=0,\quad 1\preceq l_{12}\preceq 2,\notag\\
&\Delta_4^{(M,N)}[(2;1,1), (2;1,2), (2;1,3), (2;1,4)]=0.\notag
\end{align}
Using $(\alpha 3)$ to $\Delta_4^{(M,N)}[(l_{13}, l_{12}, l_{11},2; 1,1)]$
we obtain
\begin{align}
:S_0(z_{0,s}){\mc A}_{l_{13}, l_{12}, l_{11},2; 1}[qz_{0,s}]:,\quad 1\preceq l_{13} \preceq l_{12} \preceq l_{11} \preceq 2.
\notag
\end{align}
Using $(\alpha 3)$ to $\Delta_4^{(M,N)}[(l_{11}, 2;1,1), (l_{21},2;1,2)]$ we obtain
\begin{align}
:S_0(z_{0,s}){\mc A}_{l_{11},2;1}[qz_{0,s}]:\times
:S_0(z_{0,t}){\mc A}_{l_{21},2;1}[qz_{0,t}]:,\quad
1\preceq l_{11}, l_{21} \preceq 2.
\notag
\end{align}

We consider $\Delta_4^{(M,N)}[(\mu;0), (\lambda;1)]$ with $|\mu|=1$, $|\lambda|=3$.

Using $(\alpha 3)$ and $(\alpha 4)$ to $\Delta_4^{(M,N)}[(k_{11};0,1), (l_{12},l_{11},2;1,1)]$ we obtain
\begin{align}
:{\mc A}_{k_{11}:0}[z_{0,s}]:\times
:S_0(z_{0,t}){\mc A}_{l_{12},l_{11},2;1}[q z_{0,t}]:,\quad 1 \preceq k_{11} \preceq 2, 1\preceq l_{12} \preceq l_{11} \preceq 2.
\notag
\end{align}
Using $(\alpha 3)$ and $(\alpha 4)$ to $\Delta_4^{(M,N)}[(k_{11};0,1), (l_{12},2;1,1), (2;1,2)]$ we obtain
\begin{align}
:{\mc A}_{k_{11}; 0}[z_{0,s}]{\mc A}_{2;1}[q z_{0,s}]: \times :S_0(z_{0,t}){\mc A}_{l_{12},2;1}[qz_{0,t}]:,\quad 1 \preceq k_{11}, l_{12} \preceq 2.
\notag
\end{align} 

We consider $\Delta_4^{(M,N)}[(\mu;0), (\lambda;1)]$ with $|\mu|=|\lambda|=2$.

Using $(\alpha 3)$ and $(\alpha 4)$ to $\Delta_4^{(M,N)}[(k_{12}, k_{11},;0,1), (l_{12},2; 1,1))]$
we obtain 
\begin{align}
:{\mc A}_{k_{12}, k_{11};0}[z_{0,s}]:\times
:S_0(z_{0,t}){\mc A}_{l_{12},2;1}[qz_{0,t}]:,\quad
1\preceq k_{12} \preceq k_{11} \preceq 2, 1\preceq l_{12} \preceq 2.\notag
\end{align}
Using $(\alpha 3)$ to $\Delta_4^{(M,N)}[(k_{11};0,1), (k_{21};0,2), (l_{12},2;1,1)]$
we obtain
\begin{align}
:{\mc A}_{k_{11};0}[z_{0,s}]:\times
:{\mc A}_{k_{21};0}[z_{0,t}]:\times
:{\mc A}_{l_{12},2;1}[qz_{0,u}]S_0(z_{0,u}):,\quad 1\preceq k_{11}, k_{21}, l_{12} \preceq 2.\notag
\end{align}

We consider $\Delta_4[(\mu;0), (\lambda;1)]$ with $|\mu|=3$, $|\lambda|=1$.

Using $(\alpha 1)$ and $(\alpha 2)$ to $\Delta_4^{(M,N)}[(k_{13}, k_{12}, k_{11}; 0,1), (2;1,1)]$ we obtain
\begin{align}
:{\mc A}_{k_{13}, k_{12}, k_{11} ;0}[z_{0,s}] {\mc A}_{2;1}[qz_{0,s}]:,\quad 
1 \preceq k_{13} \preceq k_{12} \preceq k_{11} \preceq 2.\notag
\end{align}
Using $(\alpha 1)$ and $(\alpha 2)$ to 
$\Delta_4^{(M,N)}[(k_{12} ,k_{11}; 0,1), (k_{21}; 0,2), (2;1,1)]$ we obtain
\begin{align}
&:{\mc A}_{k_{12}, k_{11}; 0}[z_{0,s}] 
{\mc A}_{2;1}[qz_{0,s}]:\times
:{\mc A}_{k_{21}; 0}[z_{0,t}]:,\quad 1 \preceq k_{12} \preceq k_{11} \preceq 2,  1\preceq k_{21} \preceq 2,
\notag
\\
&:{\mc A}_{k_{12},k_{11}; 0}[z_{0,s}]:
\times
:{\mc A}_{k_{21}; 0}[z_{0,t}]{\mc A}_{2;1}[qz_{0,t}]:,\quad 1 \preceq k_{12} \preceq k_{11} \preceq 2,  1\preceq k_{21} \preceq 2.\notag
\end{align}
Using $(\alpha 1)$ and $(\alpha 2)$ to $\Delta_4^{(M,N)}[(k_{11};0,1), (k_{21};0,2), (k_{31};0,3), (2;1,1)]$ we obtain
\begin{align}
:{\mc A}_{k_{11};0}[z_{0,s}] {\mc A}_{2;1}[qz_{0,s}]:
\times
:{\mc A}_{k_{21};0}[z_{0,t}]:
\times
:{\mc A}_{k_{31};0}[z_{0,u}]:,
\quad 1 \preceq k_{11}, k_{21}, k_{31} \preceq 2.\notag
\end{align}

In $\Delta_4^{(M,N)}[(\lambda;1)]$ with $|\lambda|=4$,
$\Delta_4^{(M,N)}[(\mu;0), (\lambda;1)]$ with $|\mu|=1$, $|\lambda|=3$ such that $l(\lambda)=1$, 
and $\Delta_4^{(M,N)}[(\mu;0), (\lambda;1)]$ with
$|\mu|=|\lambda|=2$,
the core parts of the residues are given by
(\ref{lem:IM-GN:4:eq5}) and (\ref{lem:IM-GN:4:eq6}) for $m=4$.
In $\Delta_4^{(M,N)}[(\mu;0), (\lambda;1)]$ with $|\mu|=3$, $|\lambda|=1$, 
and $\Delta_4^{(M,N)}[(\mu;0), (\lambda;1)]$ with
$|\mu|=1$, $|\lambda|=3$ such that $l(\lambda)=2$,
the core parts of the residues are given by
(\ref{lem:IM-GN:4:eq7}) and (\ref{lem:IM-GN:4:eq8}) for $m=4$.
The core parts (\ref{lem:IM-GN:4:eq5}) and (\ref{lem:IM-GN:4:eq7}) (resp.  (\ref{lem:IM-GN:4:eq6}) and (\ref{lem:IM-GN:4:eq8})) are the same.
The residues corresponding to (\ref{lem:IM-GN:4:eq5}) and (\ref{lem:IM-GN:4:eq7}) 
(resp.  (\ref{lem:IM-GN:4:eq6}) and (\ref{lem:IM-GN:4:eq8})) 
have opposite signs.
Hence, the sum of the residues becomes zero.
We obtain Lemma \ref{lem:IM-GN:4} for $m=4$.

We consider Lemma \ref{lem:IM-GN:4} for general $m \geq 3$.

We calculate the residues in $\Delta_m^{(M,N)}[(\lambda;1)]$ and $\Delta_m^{(M,N)}[(\mu;0), (\lambda;1)]$
using $(\alpha 1)$, $(\alpha 2)$, $(\alpha 3)$, and $(\alpha 4)$, in the same way as the case $m=3,4$.
The core parts are summarized as follows.
\begin{align}
&
:\overset{}{\prod_{1\leq i \leq l(\lambda)}}
{\mc A}_{l_{i,\lambda_i},\ldots, l_{i,2},2;1}[qz_{0,i}]S_0(z_{0,i}):,
\label{lem:IM-GN:4:eq5}
\\
&
\quad
(\lambda, (l_{ij})) \in {\mc E}_1,\quad |\lambda|=m,\quad \lambda_i \geq 2, 1\leq i \leq l(\lambda),\notag
\\
&
:\overset{}{\prod_{1\leq i \leq l(\mu)}}
{\mc A}_{k_{i,\mu_i},\ldots, k_{i,1};0}[z_{0,i}] \times
\overset{}{\prod_{1\leq i \leq l(\lambda)}}
{\mc A}_{l_{i,\lambda_i},\ldots, l_{i,2},2;1}[qz_{0, i+l(\mu)}]S_0(z_{0, i+l(\mu)}):,
\label{lem:IM-GN:4:eq6}
\\
&
(\mu, (k_{i,j})) \in {\mc E}_0,
\lambda=(l_{ij}) \in {\mc E}_1, 
\quad
|\mu|+|\lambda|=m, 
\quad
\lambda_i \geq 2, 1\leq i \leq l(\lambda).\notag
\end{align}

\begin{align}
&
:\overset{}{\prod_{1\leq i \leq l(\lambda)}}
{\mc A}_{l_{i,\lambda_i},\ldots, l_{i,2};0}[z_{0,i}]{\mc A}_{2;1}[q z_{0,i}]:,
\label{lem:IM-GN:4:eq7}
\\
&
\quad
(\lambda, (l_{ij})) \in {\mc E}_1,\quad |\lambda|=m,\quad \lambda_i \geq 2, 1\leq i \leq l(\lambda),\notag
\\
&
:\overset{}{\prod_{1\leq i \leq l(\mu)}}
{\mc A}_{k_{i,\mu_i},\ldots, k_{i,1};0}[z_{0,i}]\times
\overset{}{\prod_{1\leq i \leq l(\lambda)}}
{\mc A}_{l_{i,\lambda_i},\ldots, l_{i,2};0}[z_{0, i+l(\mu)}]{\mc A}_{2;1}[qz_{0, i+l(\mu)}]:,
\label{lem:IM-GN:4:eq8}
\\
&
(\mu, (k_{i,j})) \in {\mc E}_0,
(\lambda, (l_{ij})) \in {\mc E}_1, 
\quad
|\mu|+|\lambda|=m, 
\quad
\lambda_i \geq 2, 1\leq i \leq l(\lambda).\notag
\end{align}
The core parts (\ref{lem:IM-GN:4:eq5}) and (\ref{lem:IM-GN:4:eq7}) (resp.  (\ref{lem:IM-GN:4:eq6}) and (\ref{lem:IM-GN:4:eq8})) are the same.
The residues corresponding to (\ref{lem:IM-GN:4:eq5}) and (\ref{lem:IM-GN:4:eq7}) 
(resp.  (\ref{lem:IM-GN:4:eq6}) and (\ref{lem:IM-GN:4:eq8})) 
have opposite signs.
Hence, the sum of the residues becomes zero.
We obtain Lemma \ref{lem:IM-GN:4} for general $m$.
\end{proof}

\vspace{2mm}

\begin{lem}
\label{lem:IM-GN:5}
The following holds for 
$M \geq m \geq 3$.
\begin{align}
&
\sum_{\overline{\mu} \in \overline{\mathcal E}_0
\atop{|\overline{\mu}|=m}}
\com{N}{\overline{\mu}} \nu_{\overline{\mu};0}
\Delta^{(M,N)}_m[(\overline{\mu}; 0)]
+\sum_{\overline{\lambda} \in \overline{\mathcal E}_1
\atop{|\overline{\lambda}|=m}}
\com{N}{\overline{\lambda}} \nu_{\overline{\lambda};1}
\Delta^{(M,N)}_m[(\overline{\lambda}; 1)]
\notag\\
&+
\sum_{\overline{\mu} \in \overline{\mathcal E}_0,
\overline{\lambda} \in \overline{\mathcal E}_1
\atop{|\overline{\mu}|+|\overline{\lambda}|=m}}
\com{N}{\overline{\mu}} \com{N}{\overline{\lambda}} \nu_{\overline{\mu};0}
\nu_{\overline{\lambda};1}\Delta^{(M,N)}_m[(\overline{\mu}; 0), (\overline{\lambda};1)]=0.
\notag
\end{align}
\end{lem}

\begin{proof}

We prepare the rules $(\beta 1), (\beta 2), (\beta 3)$, and $(\beta 4)$ on the integral
\begin{align}
\int_{C_{\delta_0}-C_{\delta_0 q^{-n}}}\ul{dz_0}~{\mc A}_{k_n,\ldots, k_1;0}[z_0]S_1(z_1)f(z_0),
\quad n \geq 1,\quad \bar{1} \succeq k_n \succeq \cdots \succeq k_1 \succeq \bar{2}.
\notag
\end{align}
Here $f(z)$ is a holomorphic function in the region $|\delta_0|<|z|<|\delta_0 q^{-n}|$.
\vspace{3mm}~\\
$(\beta 1)$ If $n \geq 1$, $k_1=\cdots=k_n=\bar{1}$,
as can be seen from the normal ordering of
\begin{align}
{{\mc A}_{\underbrace{\scriptstyle \bar{1},\ldots,\bar{1}}_{n} ;0}
[z_0]S_1(z_1),\notag}
\end{align}
there are poles at $z_0=q^{-s}z_1$, $1\leq s \leq n$, in the region of the integral $\int_{C_{\delta_0}-C_{\delta_0}q^{-n}}\ul{dz_0}$.

Taking the residues we obtain
\begin{align}
{{\mc A}_{\underbrace{\scriptstyle \bar{1},\ldots,\bar{1}}_{n} ;0}}
[q^{-s}z_1]S_1(z_1):
=q_1^s
:
{{\mc A}_{\underbrace{\scriptstyle \ol{1},\ldots,\ol{1}}_{n-s},
\underbrace{\scriptstyle \ol{2},\ldots,\ol{2}}_s ;0}}
[q^{-s}z_1]S_1(q^{-s}z_1):,\quad 1\leq s \leq n.\notag
\end{align}
$(\beta 2)$ If $n \geq 2$, $k_1=\bar{2}, k_2=\cdots=k_n=\bar{1}$, as can be seen from
the normal ordering of
\begin{align}
{\mc A}_{
{\underbrace{\scriptstyle \ol{1},\ldots,\ol{1}}_{n-1},\ol{2};0}}
[z_0]S_1(z_1),\notag
\end{align}
there are poles at $z_0=q^{-s}z_1, 2\leq s \leq n$,
in the region of the integral 
$\int_{C_{\delta_0}-C_{\delta_0}q^{-n}}\ul{dz_0}$.

Taking the residues we obtain
\begin{align}
:{{\mc A}_{\underbrace{\scriptstyle 
\bar{1},\ldots,\bar{1}}_{n-1},\ol{2};0}}
[q^{-s-1}z_1]S_1(z_1):
=q_1^{s-1}
:{{\mc A}_{
\underbrace{\scriptstyle \ol{1},\ldots,\ol{1}}_{n-s},
\underbrace{\scriptstyle \ol{2},\ldots,\ol{2}}_s ;0}
}[q^{-s}z_1]S_1(q^{-s+1}z_1):,\quad 2\leq s \leq n.\notag
\end{align}
\vspace{3mm}~\\
$(\beta 3)$ If $n \geq 3$, $2\leq r \leq n-1$,
$k_1=\cdots=k_r=\bar{2}$, $k_{r+1}=\cdots=k_n=\bar{1}$, as can be seen from
the normal ordering of
\begin{align}
{\mc A}_{
{\underbrace{\scriptstyle \ol{1},\ldots,\ol{1}}_{n-r}, 
\underbrace{\scriptstyle \ol{2},\ldots, \ol{2}}_r;0}}
[z_0]S_1(z_1),\notag
\end{align}
there are
\begin{center}
(i)~ poles at $z_0=q^{-r-s}z_1$, $1\leq s \leq n-r$,
(ii)~ poles at $z_0=q^{-s}z_1$, $1 \leq s \leq r-1$,
\end{center}
in the region of the integral 
$\int_{C_{\delta_0}-C_{\delta_0}q^{-n}}\ul{dz_0}$.

Taking the residues we obtain
\begin{align}
\mathrm{(i)}~&:{\mc A}_{
{\underbrace{\scriptstyle \bar{1},\ldots,\bar{1}}_{n-r},
\underbrace{\scriptstyle \ol{2},\ldots, \ol{2}}_r;0}}
[q^{-r-s}z_1]S_1(z_1):
\notag\\
&=q_1^s
:{\mc A}_{
{\underbrace{\scriptstyle \ol{1},\ldots,\ol{1}}_{n-r-s},
\underbrace{\scriptstyle \ol{2},\ldots,\ol{2}}_{r+s};0}}
[q^{-r-s}z_1]S_1(q^{-s}z_1):,\quad 2\leq r \leq n-1, 1\leq s \leq n-r,
\notag
\\
\mathrm{(ii)}~&:{\mc A}_{
{\underbrace{\scriptstyle \bar{1},\ldots,\bar{1}}_{n-r},
\underbrace{\scriptstyle \ol{2},\ldots, \ol{2}}_r;0}}
[q^{-s}z_1]S_1(z_1):,
\quad 2\leq r \leq n-1, 1\leq s \leq r-1.
\notag
\end{align}
$(\beta 4)$
If $n \geq 1$, $k_1=\cdots=k_n=\ol{2}$,
as can be seen from the normal ordering of
\begin{align}
{\mc A}_{
{\underbrace{\scriptstyle \ol{2},\ldots,\ol{2}}_n ;0}}
[z_0]S_1(z_1),\notag
\end{align}
there are poles at $z_0=q^{-s}z_1$, $1\leq s \leq n-1$,
in the region of the integral 
$\int_{C_{\delta_0}-C_{\delta_0}q^{-n}}\ul{dz_0}$.

Taking the residues we obtain
\begin{align}
{\mc A}_{
{\underbrace{\scriptstyle \ol{2},\ldots,\ol{2}}_n;0}}
[q^{-s}z_1]S_1(z_1),\quad 1\leq s \leq n-1.\notag
\end{align}

\vspace{2mm}

Next, we prepare the rules $(\gamma 1), (\gamma 2), (\gamma 3)$, $(\gamma 4)$ and $(\gamma 5)$ on the integral
\begin{align}
&\int_{C_1-C_{q^{-1}}}\ul{dz_1}
~{\mc A}_{k_a,\ldots,k_1;0}[z_0] {\mc A}_{l_b,\ldots,l_2,\ol{2};1}[z_1]
\prod_{i=1}^a \prod_{j=1}^b \omega(q^{-i+j}z_1/z_0)f(z_1),\notag
\\& a \geq 0, b\geq 2,\quad 
\ol{1} \succeq k_a \succeq \cdots \succeq k_1 \succeq \ol{2},\quad
\ol{1} \succeq l_b \succeq \cdots \succeq l_2 \succeq \ol{2}=l_1.\notag
\end{align}
Here $f(z)$ is a holomorphic function in the region $1<|z|<|q^{-1}|$.
If $a=0$, it is understood that ${\mc A}_{k_a,\ldots,k_1;0}[z_0]=S_0(z_0)$.
\vspace{3mm}~\\
$(\gamma 1)$ If $a \geq 0$, $b=1$, $\ol{1}\succeq k_a \succeq \cdots \succeq k_1 \succeq \ol{2}$, and $l_1=\ol{2}$,
as can be seen from the normal ordering of
\begin{align}
{\mc A}_{k_a,\ldots,k_1;0}[z_0]{\mc A}_{\ol{2};1}[z_1]\prod_{i=1}^a\omega(q^{-i+1}z_1/z_0),\notag
\end{align}
there is
\begin{center}
(i)~ a pole at $z_1=z_0$ for $a=0$,
(ii)~ no pole for $a=1$,
(iii)~ a pole at $z_1=z_0$ for $a \geq 2$,
\end{center}
in the region of the integral $\int_{C_1-C_{q^{-1}}}\ul{dz_1}$.

Taking the residues we obtain
\begin{align}
\mathrm{(i)}~&:S_0(z_0){\mc A}_{\ol{2};1}[z_0]:~\mathrm{in~the~case~of}~a=0,
\notag\\
\mathrm{(ii)}~&~0\quad \mathrm{in~the~case~of}~a=1,
\notag
\\
\mathrm{(iii)}~&:{\mc A}_{k_a,\ldots,k_1;0}[z_0]{\mc A}_{\ol{2};1}[z_0]:~
\mathrm{in~the~case~of}~a\geq 2, \ol{1}\succeq k_a \succeq \cdots \succeq k_1 \succeq  \ol{2}.
\notag
\end{align}
\vspace{2mm}~\\
$(\gamma 2)$ If $a=0$, $b \geq 2$, and $\ol{1} \succeq l_b \succeq \cdots \succeq l_2 \succeq \ol{2}=l_1$,
as can be seen from the normal ordering of
\begin{align}
S_0(z_0){\mc A}_{l_b,\ldots,l_2,\ol{2};1}[z_1],\notag
\end{align}
there is a pole at $z_1=z_0$ in the region of the integral $\int_{C_1-C_{q^{-1}}}\ul{dz_1}$.

Taking the residues we obtain
\begin{align}
:S_0(z_0){\mc A}_{l_b,\ldots,l_2,\ol{2};1}[z_0]:,\quad
\ol{1}\succeq l_b \succeq \cdots \succeq l_2 \succeq l_1=\ol{2}.\notag
\end{align}
\vspace{2mm}~\\
$(\gamma 3)$
In the case of $a=b\geq 2$,
$\ol{1} \succeq k_a \succeq \cdots \succeq k_1 \succeq \ol{2}$,
and $\ol{1} \succeq l_a \succeq \cdots \succeq l_2 \succeq \ol{2}=l_1$,
as can be seen from the normal ordering
\begin{align}
{\mc A}_{k_a,\ldots,k_1;0}[z_0]
{\mc A}_{l_a,\ldots,l_2,\ol{2};1}[z_1]
\prod_{i=1}^a \prod_{j=1}^a \omega(q^{-i+j}z_1/z_0),\notag 
\end{align}
there is no pole in the region of the integral $\int_{C_1-C_{q^{-1}}}\ul{dz_1}$.
Taking the residue we obtain zero.
\vspace{2mm}~\\
$(\gamma 4)$
If $a>b\geq 1$,
$\ol{1} \succeq k_a \succeq \cdots \succeq k_b \succeq \cdots \succeq k_1 \succeq \ol{2}$,
and 
$\ol{1} \succeq l_b \succeq \cdots \succeq l_2 \succeq \ol{2}=l_1$,
as can be seen from the normal ordering
\begin{align}
{\mc A}_{k_a,\ldots, k_{b+1}, k_b,\ldots, k_1;0}[z_0]
{\mc A}_{l_b,\ldots,l_2,\ol{2};1}[z_1]
\prod_{i=1}^a \prod_{j=1}^b \omega(q^{-i+j}z_1/z_0),\notag 
\end{align}
there is 
\begin{center}
(i)~ a pole at $z_1=z_0$ in the case of $k_{b+1} \succeq l_b$,
(ii)~ no pole in the case of $k_{b+1} \prec l_b$,
\end{center}
in the region of the integral $\int_{C_1-C_{q^{-1}}}\ul{dz_1}$.
Taking the residues we obtain
\begin{align}
\mathrm{(i)}~&:{\mc A}_{k_a,\ldots, k_{b+1}, k_b,\ldots, k_1;0}[z_0]
{\mc A}_{l_b,\ldots,l_2,\ol{2};1}[z_0]:~\mathrm{in~the~case~of}~a>b \geq 2,~k_{b+1} \succeq l_b,
\notag\\
\mathrm{(ii)}~&~0\quad \mathrm{in~the~case~of}~a>b \geq 2,~ k_{b+1} \prec l_b.
\notag
\end{align}
\vspace{2mm}~\\
$(\gamma 5)$
If $1 \leq a<b$,
$\ol{1} \succeq k_a \succeq \cdots \succeq k_1 \succeq \ol{2}$,
and 
$\ol{1} \succeq l_b \succeq \cdots \succeq l_a \succeq \cdots \succeq l_2 \succeq \ol{2}=l_1$,
as can be seen from the normal ordering
\begin{align}
{\mc A}_{k_a, \ldots, k_1;0}[z_0]
{\mc A}_{l_b,\ldots, l_{a+1}, l_a, \ldots, l_2,\ol{2};1}[z_1]
\prod_{i=1}^a \prod_{j=1}^b \omega(q^{-i+j}z_1/z_0),\notag 
\end{align}
there is 
\begin{center}
(i)~ a pole at $z_1=z_0$ in the case of $k_{a} \preceq l_{a+1}$,
(ii)~ no pole in the case of $k_{a} \succ l_{a+1}$,
\end{center}
in the region of the integral $\int_{C_1-C_{q^{-1}}}\ul{dz_1}$.
Taking the residues we obtain
\begin{align}
\mathrm{(i)}~&:{\mc A}_{k_a, \ldots, k_1;0}[z_0]
{\mc A}_{l_b,\ldots,l_{a+1}, l_a,\ldots, l_2,\ol{2};1}[z_0]:~\mathrm{in~the~case~of}~1\leq a<b,~k_{a} \preceq l_{a+1},
\notag\\
\mathrm{(ii)}~&~0\quad \mathrm{in~the~case~of}~1\leq a<b,~ k_{a} \succ l_{a+1}.
\notag
\end{align}
\vspace{2mm}

We show Lemma \ref{lem:IM-GN:5} for the case $m=3$.

\vspace{2mm}

We consider $\Delta_3^{(M,N)}[(\ol{\mu}; 0)]$ with $|\ol{\mu}|=3$.
We consider
$\Delta_3^{(M,N)}[(k_3, k_2, k_1;0,1)]$ with $\ol{1}\succeq k_3 \succeq k_2 \succeq k_1 \succeq \ol{2}$.
Using $(\beta 1)$ to $\Delta_3^{(M,N)}[(\ol{1},\ol{1},\ol{1};0,1)]$ and shifting the variable $z_{1,s}$,
we obtain
\begin{align}
:{\mc A}_{l_2, l_1, \ol{2};0}[z_{1,s}]S_1(z_{1,s}):,\quad
\ol{1}\succeq l_2 \succeq l_1 \succeq \ol{2}.
\label{lem:IM-GN:5:eq1}
\end{align}
Using $(\beta 2)$ to 
$\Delta_3^{(M,N)}[(\ol{1},\ol{1},\ol{2};0,1)]$
we obtain
\begin{align}
:{\mc A}_{\ol{1}, \ol{2}, \ol{2};0}[q^{-1}z_{1,s}]S_1(z_{1,s}):,\quad
:{\mc A}_{\ol{2}, \ol{2}, \ol{2};0}[q^{-1}z_{1,s}]S_1(z_{1,s}):.
\label{lem:IM-GN:5:eq2}
\end{align}
Using $(\beta 3)$ to 
$\Delta_3^{(M,N)}[(\ol{1},\ol{2},\ol{2};0,1)]$
we obtain
\begin{align}
:{\mc A}_{\ol{1}, \ol{2}, \ol{2};0}[q^{-1}z_{1,s}]S_1(z_{1,s}):,\quad
:{\mc A}_{\ol{2}, \ol{2}, \ol{2};0}[q^{-2}z_{1,s}]S_1(z_{1,s}):.
\label{lem:IM-GN:5:eq3}
\end{align}
Using $(\beta 4)$ to 
$\Delta_3^{(M,N)}[(\ol{2},\ol{2},\ol{2};0,1)]$
we obtain
\begin{align}
:{\mc A}_{\ol{2}, \ol{2}, \ol{2};0}[q^{-1}z_{1,s}]S_1(z_{1,s}):,\quad
:{\mc A}_{\ol{2}, \ol{2}, \ol{2};0}[q^{-2}z_{1,s}]S_1(z_{1,s}):.
\label{lem:IM-GN:5:eq4}
\end{align}
The $:{\mc A}_{\ol{1}, \ol{2}, \ol{2};0}[q^{-1}z_{1,s}]S_1(z_{1,s}):$,
$:{\mc A}_{\ol{2}, \ol{2}, \ol{2};0}[q^{-1}z_{1,s}]S_1(z_{1,s}):$,
and $:{\mc A}_{\ol{2}, \ol{2}, \ol{2};0}[q^{-2}z_{1,s}]S_1(z_{1,s}):$ appear exactly twice
in (\ref{lem:IM-GN:5:eq2}),  (\ref{lem:IM-GN:5:eq3}), and  (\ref{lem:IM-GN:5:eq4}).
Taking the residues reveals that
the two residues have opposite signs.
The sum of the residues becomes zero.
Hence, combining all the terms arising from $\Delta_3^{(M,N)}[(k_3, k_2, k_1;0,1)]$
with $\ol{1}\succeq k_3 \succeq k_2 \succeq k_1 \succeq \ol{2}$
results in (\ref{lem:IM-GN:5:eq1}) only.

Next, we consider $\Delta_{1}^{(M,N)}[(k;0,1)]$ and $\Delta_{2}^{(M,N)}[(k_2,k_1;0,1)]$.
Using $(\beta 1)$ and $(\beta 4)$ to $\Delta_1^{(M,N)}[(k;0,1)]$, $k=\ol{1}, \ol{2}$, we obtain
\begin{align}
\delta_{k,\ol{2}} :{\mc A}_{\ol{2};0}[z_{1,s}]S_1(z_{1,s}):,\quad \ol{1} \succeq k \succeq \ol{2}.
\label{lem:IM-GN:5:eq5}
\end{align}

Using $(\beta 1)$ to $\Delta_2^{(M,N)}[(\ol{1}, \ol{1};0,1)]$ we obtain
\begin{align}
:{\mc A}_{\ol{1},\ol{2};0}[z_{1,s}]S_1(z_{1,s}):,\quad
:{\mc A}_{\ol{2},\ol{2};0}[z_{1,s}]S_1(z_{1,s}):.
\label{lem:IM-GN:5:eq6}
\end{align}
Using $(\beta 2)$ to $\Delta_2^{(M,N)}[(\ol{1}, \ol{2};0,1)]$ we obtain
\begin{align}
:{\mc A}_{\ol{2},\ol{2};0}[q^{-1}z_{1,s}]S_1(z_{1,s}):.
\label{lem:IM-GN:5:eq7}
\end{align}
Using $(\beta 4)$ to $\Delta_2^{(M,N)}[(\ol{2}, \ol{2};0,1)]$ we obtain
\begin{align}
:{\mc A}_{\ol{2},\ol{2};0}[q^{-1}z_{1,s}]S_1(z_{1,s}):.
\label{lem:IM-GN:5:eq8}
\end{align}
In (\ref{lem:IM-GN:5:eq7}) and  (\ref{lem:IM-GN:5:eq8})
$:{\mc A}_{\ol{2}, \ol{2}, ;0}[q^{-1}z_{1,s}]S_1(z_{1,s}):$ appears exactly twice.
Taking the residues and determining the constants reveal that
the two residues have opposite signs.
The sum of them becomes zero.
Hence, combining the residues arising from $\Delta_2^{(M,N)}[(k_2, k_1;0,1)]$, 
$\ol{1}\succeq k_2 \succeq k_1 \succeq \ol{2}$
results in 
\begin{align}
:{\mc A}_{l, \ol{2};0}[z_{1,s}]S_1(z_{1,s}):,\quad l=\ol{1}, \ol{2}.
\label{lem:IM-GN:5:eq9}
\end{align}

We consider $\Delta_3^{(M,N)}[(k_{12},k_{11};0,1), (k_{21};0,2)]$ and $\Delta_3^{(M,N)}[(k_{11};0,1), (k_{21};0,2), (k_{31};0,3)]$.
We rewrite the integrals of $\Delta_3^{(M,N)}[(k_{12},k_{11};0,1), (k_{21};0,2)]$ as follows.
\begin{align}
&\int_{C_{\delta_0}}\ul{dz_{0,1}}
\int_{C_{\delta_0}}\ul{dz_{0,2}}-
\int_{C_{\delta_0}q^{-2}}\ul{dz_{0,1}}
\int_{C_{\delta_0}q^{-1}}\ul{dz_{0,2}}\notag
\\&=
-\int_{C_{\delta_0}-C_{\delta_0 q^{-2}}}\ul{dz_{0,1}}
\int_{C_{\delta_0}-C_{\delta_0 q^{-1}}}\ul{dz_{0,2}}+
\int_{C_{\delta_0}-C_{\delta_0}q^{-2}}\ul{dz_{0,1}}
\int_{C_{\delta_0}}\ul{dz_{0,2}}
+
\int_{C_{\delta_0}}\ul{dz_{0,1}}
\int_{C_{\delta_0}-C_{\delta_0}q^{-1}}\ul{dz_{0,2}}.
\notag
\end{align}
Using (\ref{lem:IM-GN:5:eq5}) and (\ref{lem:IM-GN:5:eq9}) to $\Delta_3^{(M,N)}[(k_{12},k_{11};0,1), (k_{21};0,2)]$
we obtain
\begin{align}
&:{\mc A}_{l_{12}, \ol{2};0}[z_{1,s}]S_1(z_{1,s}):\times :{\mc A}_{\ol{2};0}[z_{1,t}]S_1(z_{1,t}):,
\quad l_{12}=\ol{1}, \ol{2},
\label{lem:IM-GN:5:eq10}
\\
&:{\mc A}_{l_{12}, \ol{2};0}[z_{1,s}]S_1(z_{1,s}):\times :{\mc A}_{k_{21};0}[z_{1,t}]:,
\quad l_{12}=\ol{1}, \ol{2}, ~k_{21}=\ol{1}, \ol{2},
\label{lem:IM-GN:5:eq11}
\\
&~:{\mc A}_{k_{12}, k_{11};0}[z_{1,s}]:\times :{\mc A}_{\ol{2}; 0}[z_{1,t}]S_1(z_{1,t}):,\quad
\ol{1} \succeq k_{12} \succeq k_{11} \succeq \ol{2}.
\label{lem:IM-GN:5:eq12}
\end{align}

We rewrite the integrals of $\Delta_3^{(M,N)}[(k_{11};0,1), (k_{21};0,2), (k_{31};0,3)]$ as follows.
\begin{align}
&\int_{C_{\delta_0}}\ul{dz_{0,1}}
\int_{C_{\delta_0}}\ul{dz_{0,2}}
\int_{C_{\delta_0}}\ul{dz_{0,3}}
-
\int_{C_{\delta_0}q^{-1}}\ul{dz_{0,1}}
\int_{C_{\delta_0}q^{-1}}\ul{dz_{0,2}}
\int_{C_{\delta_0q^{-1}}}\ul{dz_{0,3}}\notag
\\&=
\int_{C_{\delta_0}-C_{\delta_0 q^{-1}}}\ul{dz_{0,1}}
\int_{C_{\delta_0}-C_{\delta_0 q^{-1}}}\ul{dz_{0,2}}
\int_{C_{\delta_0}-C_{\delta_0}q^{-1}}\ul{dz_{0,3}}-
\int_{C_{\delta_0}-C_{\delta_0 q^{-1}}}\ul{dz_{0,1}}
\int_{C_{\delta_0}-C_{\delta_0 q^{-1}}}\ul{dz_{0,2}}
\int_{C_{\delta_0}}\ul{dz_{0,3}}
\notag
\\
&-
\int_{C_{\delta_0}-C_{\delta_0 q^{-1}}}\ul{dz_{0,1}}
\int_{C_{\delta_0}}\ul{dz_{0,2}}
\int_{C_{\delta_0}-C_{\delta_0 q^{-1}}}\ul{dz_{0,3}}
-
\int_{C_{\delta_0}}\ul{dz_{0,1}}
\int_{C_{\delta_0}-C_{\delta_0 q^{-1}}}\ul{dz_{0,2}}
\int_{C_{\delta_0}-C_{\delta_0 q^{-1}}}\ul{dz_{0,3}}
\notag
\\
&
+
\int_{C_{\delta_0}-C_{\delta_0 q^{-1}}}\ul{dz_{0,1}}
\int_{C_{\delta_0}}\ul{dz_{0,2}}
\int_{C_{\delta_0}}\ul{dz_{0,3}}
+
\int_{C_{\delta_0}}\ul{dz_{0,1}}
\int_{C_{\delta_0}-C_{|\delta_0 q^{-1}}}\ul{dz_{0,2}}
\int_{C_{\delta_0}}\ul{dz_{0,3}}
\notag\\
&
+
\int_{C_{\delta_0}}\ul{dz_{0,1}}
\int_{C_{\delta_0}}\ul{dz_{0,2}}
\int_{C_{\delta_0}-C_{\delta_0 q^{-1}}}\ul{dz_{0,3}}.
\notag
\end{align}
Using (\ref{lem:IM-GN:5:eq6}) to $\Delta_3^{(M,N)}[(k_{11};0,1), (k_{21};0,2), (k_{31};0,3)]$
we obtain
\begin{align}
&:{\mc A}_{\ol{2};0}[z_{1,s}]S_1(z_{1,s}):\times
:{\mc A}_{\ol{2};0}[z_{1,t}]S_1(z_{1,t}):\times
:{\mc A}_{\ol{2};0}[z_{1,u}]S_1(z_{1,u}):,
\label{lem:IM-GN:5:eq13}
\\
&:{\mc A}_{\ol{2};0}[z_{1,s}]S_1(z_{1,s}):
\times :{\mc A}_{\ol{2};0}[z_{1,t}]S_1(z_{1,t}):
\times :{\mc A}_{k_{31};0}[z_{1,u}]:, \quad k_{31}=\ol{1}, \ol{2},
\label{lem:IM-GN:5:eq14}
\\
&
:{\mc A}_{\ol{2};0}[z_{1,s}]S_1(z_{1,s}):\times
:{\mc A}_{k_{21};0}[z_{1,t}]:\times
:{\mc A}_{\ol{2};0}[z_{1,u}]S_1(z_{1,u}):, \quad k_{21}=\ol{1}, \ol{2},
\notag
\\
&
:{\mc A}_{k_{11};0}[z_{1,s}]:\times
:{\mc A}_{\ol{2};0}[z_{1,t}]S_1(z_{1,t}):\times
:{\mc A}_{\ol{2};0}[z_{1,u}]S_1(z_{1,u}):, \quad k_{11}=\ol{1}, \ol{2},
\notag
\\
&:{\mc A}_{\ol{2};0}[z_{1,s}]S_1(z_{1,s}):
\times :{\mc A}_{k_{21}; 0}[z_{1,t}]:
\times:{\mc A}_{k_{31};0}[z_{1,u}]:,
\quad
k_{21}, k_{31}=\ol{1}, \ol{2},
\label{lem:IM-GN:5:eq15}
\\
&
:{\mc A}_{k_{11};0}[z_{1,s}]:\times
:{\mc A}_{\ol{2};0}[z_{1,t}]S_1(z_{1,t}):\times
:{\mc A}_{k_{31};0}[z_{1,u}]:,\quad k_{11}, k_{31}=\ol{1}, \ol{2},
\notag
\\
&
:{\mc A}_{k_{11};0}[z_{1,s}]:\times
:{\mc A}_{k_{21};0}[z_{1,t}]:\times
:{\mc A}_{\ol{2};0}[z_{1,u}]S_1(z_{1,u}):, \quad k_{11}, k_{21}=\ol{1}, \ol{2}.
\notag
\end{align}

\vspace{2mm}

To summarize, in $\Delta_3^{(M,N)}[(\ol{\mu};0)]$ with $|\ol{\mu}|=3$,
the core parts of the residues are given by
(\ref{lem:IM-GN:5:eq1}), (\ref{lem:IM-GN:5:eq10}), 
(\ref{lem:IM-GN:5:eq11}), (\ref{lem:IM-GN:5:eq12}), 
(\ref{lem:IM-GN:5:eq13}), (\ref{lem:IM-GN:5:eq14}), 
and (\ref{lem:IM-GN:5:eq15}).
To summarize, the core parts of the residues are summarized by 
(\ref{lem:IM-GN:5:eq25}) and (\ref{lem:IM-GN:5:eq26}) for the case $m=3$.

\vspace{2mm}

Next, we consider $\Delta_3^{(M,N)}[(\ol{\lambda};1)]$ with $|\ol{\lambda}|=3$.

Using $(\gamma 2)$ to $\Delta_3^{(M,N)}[(l_{13}, l_{12}, \ol{2};1,1)]$
we obtain
\begin{align}
\label{lem:IM-GN:5:eq16}
:S_0(z_{0,s}){\mc A}_{l_{13}, l_{12}, \ol{2};1}[z_{0,s}]:,\quad
\ol{1} \succeq l_{13} \succeq l_{12} \succeq \ol{2}.
\end{align}
Using $(\gamma 2)$ to $\Delta_3^{(M,N)}[(l_{12}, \ol{2};1,1), (\ol{2};1,2)]$
we obtain
\begin{align}
:S_0(z_{0,s}){\mc A}_{l_{12}, \ol{2};1}[z_{0,s}]:\times
:S_0(z_{0,t}){\mc A}_{\ol{2};1}(z_{0,t}):,\quad \ol{1} \succeq l_{12} \succeq \ol{2}.
\label{lem:IM-GN:5:eq17}
\end{align}
Using $(\gamma 2)$ to $\Delta_3^{(M,N)}[(\ol{2};1,1), (\ol{2};1,2), (\ol{2};1,3)]$
we obtain
\begin{align}
:S_0(z_{0,s}){\mc A}_{\ol{2};1}[z_{0,s}]:\times
:S_0(z_{0,t}){\mc A}_{\ol{2};1}[z_{0,t}]:\times
:S_0(z_{0,u}){\mc A}_{\ol{2};1}[z_{0,u}]:.
\label{lem:IM-GN:5:eq18}
\end{align}

Next, we consider $\Delta_3^{(M,N)}[(\ol{\mu}; 0), (\ol{\lambda};1)]$ with $|\ol{\mu}|=1, |\ol{\lambda}|=2$.

Using $(\gamma 2)$ to $\Delta_3^{(M,N)}[(k_{11};0,1), (l_{12}, \ol{2};1,1)]$
we obtain
\begin{align}
\label{lem:IM-GN:5:eq19}
:{\mc A}_{k_{11};0}[z_{0,s}]:\times
:S_0(z_{0,t}){\mc A}_{l_{12}, \ol{2};1}[z_{0,t}]:,\quad
\ol{1} \succeq k_{11} \succeq \ol{2}.
\end{align}
Using $(\gamma 5)$ to $\Delta_3^{(M,N)}[(k_{11};0,1), (l_{12}, \ol{2};1,1)]$
we obtain
\begin{align}
\label{lem:IM-GN:5:eq20}
:{\mc A}_{k_{11};0}[z_{0,s}]
{\mc A}_{l_{12}, \ol{2};1}[z_{0,s}]:,\quad
\ol{1} \succeq k_{11} \succeq l_{12} \succeq \ol{2}.
\end{align}

Using $(\gamma 1)$ and $(\gamma 3)$ to $\Delta_3^{(M,N)}[(k_{11};0,1), (\ol{2};1,1), (\ol{2};1,2)]$
we obtain 
\begin{align}
:{\mc A}_{k_{11};0}[z_{0,s}]:\times
:S_0(z_{0,t}){\mc A}_{\ol{2};1}[z_{0,t}]:\times
:S_0(z_{0,u}){\mc A}_{\ol{2};1}[z_{0,u}]:,\quad
\ol{1} \succeq k_{11} \succeq \ol{2}.
\label{lem:IM-GN:5:eq21}
\end{align}

Next, we consider $\Delta_3^{(M,N)}[(\ol{\mu};0), (\ol{\lambda};1)]$ with $|\ol{\mu}|=2$, $|\ol{\lambda}|=1$.

Using $(\gamma 1)$ to $\Delta_m^{(M,N)}[(k_{11};0,1), (k_{21};0,2), (\ol{2};1,1)]$ we obtain
\begin{align}
:{\mc A}_{\ol{2};1}[z_{0,s}]S_0(z_{0,s}):\times
:{\mc A}_{k_{11};0}[z_{0,t}]:\times
:{\mc A}_{k_{21};0}[z_{0,u}]:,\quad \ol{1} \succeq k_{11}, k_{21} \succeq \ol{2}.
\label{lem:IM-GN:5:eq22}
\end{align}
Using $(\gamma 2)$ to $\Delta_3^{(M,N)}[(k_{12}, k_{11};0,1), (\ol{2};1,1)]$
we obtain
\begin{align}
:S_0(z_{0,s}){\mc A}_{\ol{2};1}[z_{0,s}]:\times
:{\mc A}_{k_{12}, k_{11};0}[z_{0,t}]:,\quad 
\ol{1} \succeq k_{12} \succeq k_{11} \succeq \ol{2}.
\label{lem:IM-GN:5:eq23}
\end{align}
Using $(\gamma 4)$ to $\Delta_3^{(M,N)}[(k_{12}, k_{11};0,1), (\ol{2};1,1)]$
we obtain
\begin{align}
:{\mc A}_{k_{12}, k_{11};0}[z_{0,s}]
{\mc A}_{\ol{2};1}[z_{0,s}]:,\quad 
\ol{1} \succeq k_{12} \succeq k_{11} \succeq \ol{2}.
\label{lem:IM-GN:5:eq24}
\end{align}

Using $:{\mc A}_{l_1;0}[z]{\mc A}_{l_2,\ol{2};1}[z]:=:{\mc A}_{l_2,l_1;0}[z]{\mc A}_{\ol{2};1}[z]:$
we can cancel the residues emerging from $(\ref{lem:IM-GN:5:eq20})$
and $(\ref{lem:IM-GN:5:eq24})$.
Hence, in $\Delta_3^{(M,N)}[(\ol{\lambda};1)]$ and $\Delta_3^{(M,N)}[(\ol{\mu};0), (\ol{\lambda};1)]$,
the core parts of the residues are given by 
(\ref{lem:IM-GN:5:eq16}),
(\ref{lem:IM-GN:5:eq17}),
(\ref{lem:IM-GN:5:eq18}),
(\ref{lem:IM-GN:5:eq19}),
(\ref{lem:IM-GN:5:eq21}),
(\ref{lem:IM-GN:5:eq22}),
and (\ref{lem:IM-GN:5:eq23}).
To summarize, the core parts of the residues are summarized by 
(\ref{lem:IM-GN:5:eq25}), (\ref{lem:IM-GN:5:eq26}), (\ref{lem:IM-GN:5:eq27}), and (\ref{lem:IM-GN:5:eq28}) for the case $m=3$.
The core parts (\ref{lem:IM-GN:5:eq25}) and (\ref{lem:IM-GN:5:eq27}) 
(resp.  (\ref{lem:IM-GN:5:eq26}) and (\ref{lem:IM-GN:5:eq28})) are the same.
The residues corresponding to (\ref{lem:IM-GN:5:eq25}) and (\ref{lem:IM-GN:5:eq27}) 
(resp.  (\ref{lem:IM-GN:5:eq26}) and (\ref{lem:IM-GN:5:eq28}))
have opposite signs.
Hence, the sum of the residues becomes zero.
We obtain Lemma \ref{lem:IM-GN:5} for the case $m=3$.

We consider Lemma \ref{lem:IM-GN:5} for general $m \geq 3$.
We calculate the residues of $\Delta_m^{(M,N)}[(\ol{\mu};0)]$
using $(\beta 1)$, $(\beta 2)$, $(\beta 3)$, and $(\beta 4)$ in the same way as the case $m=3$.
The core parts of the residues are summarized as follows.
\begin{align}
&
:\overset{}{\prod_{1\leq i \leq l(\ol{\lambda})}}
{\mc A}_{l_{i,\ol{\lambda}_i},\ldots, l_{i,2}, \ol{2};0}[z_{0,i}]S_1(z_{0,i}):,
\quad
(\ol{\lambda}, (l_{ij})) \in \ol{\mc E}_1, \quad
|\ol{\lambda}|=m,
\label{lem:IM-GN:5:eq25}
\\
&
:\overset{}{\prod_{1\leq i \leq l(\ol{\mu})}}
{\mc A}_{k_{i,\ol{\mu}_i},\ldots, k_{i,1};0}[z_{0,i}]\times
\overset{}{\prod_{1\leq i \leq l(\ol{\lambda})}}
{\mc A}_{l_{i,\ol{\lambda}_i},\ldots, l_{i,2},\ol{2};0}[z_{0, i+l(\ol{\mu})}]S_1(z_{0, i+l(\ol{\mu})}):,
\label{lem:IM-GN:5:eq26}
\\
&
(\ol{\mu}, (k_{i,j})) \in \ol{\mc E}_0,~
(\ol{\lambda}, (l_{ij})) \in \ol{\mc E}_1, \quad
|\ol{\mu}|+|\ol{\lambda}|=m.
\notag
\end{align}

We calculate the residues of $\Delta_m^{(M,N)}[(\ol{\lambda};1)]$ and $\Delta_m^{(M,N)}[(\ol{\mu};0), (\ol{\lambda};1)]$
using $(\gamma 1)$, $(\gamma 2)$, $(\gamma 3)$, $(\gamma 4)$ and $(\gamma 5)$ in the same way as the case $m=3$.
The core parts of the residues are summarized as follows.

\begin{align}
&
:\overset{}{\prod_{1\leq i \leq l(\ol{\lambda})}}
{\mc A}_{l_{i, \ol{\lambda}_i},\ldots, l_{i,2}, \ol{2};1}[z_{0,i}]S_0(z_{0,i}):,
\quad
(\ol{\lambda}, (l_{ij})) \in \ol{\mc E}_1, \quad
|\ol{\lambda}|=m,
\label{lem:IM-GN:5:eq27}
\\
&
:\overset{}{\prod_{1\leq i \leq l(\ol{\mu})}}
{\mc A}_{k_{i,\ol{\mu}_i},\ldots, k_{i,1};0}[z_{0,i}]\times
\overset{}{\prod_{1\leq i \leq l(\ol{\lambda})}}
{\mc A}_{l_{i,\ol{\lambda}_i},\ldots, l_{i,2},\ol{2};1}[z_{0, i+l(\ol{\mu})}]S_0(z_{0, i+l(\ol{\mu})}):,
\label{lem:IM-GN:5:eq28}
\\
&
(\ol{\mu}, (k_{i,j})) \in \ol{\mc E}_0,~
(\ol{\lambda}, (l_{ij})) \in \ol{\mc E}_1, \quad
|\ol{\mu}|+|\ol{\lambda}|=m.
\notag
\end{align}

The operators (\ref{lem:IM-GN:5:eq25}) and (\ref{lem:IM-GN:5:eq27}) 
(resp.  (\ref{lem:IM-GN:5:eq26}) and (\ref{lem:IM-GN:5:eq28}))
are the same.
The residues corresponding to (\ref{lem:IM-GN:5:eq25}) and (\ref{lem:IM-GN:5:eq27}) 
(resp.  (\ref{lem:IM-GN:5:eq26}) and (\ref{lem:IM-GN:5:eq28}))
have opposite signs.
Hence, the sum of the residues becomes zero.
We obtain Lemma \ref{lem:IM-GN:5} for general $m$.

\end{proof}

The relation $\mathrm{LHS}_{M,N}(a,0,c)=\mathrm{RHS}_{M,N}(a,0,c)$, $1\leq a,c \leq M, 2\leq a+c \leq M$ in Proposition \ref{prop:duality1}
is a direct consequence of Lemma \ref{lem:IM-GN:6} and \ref{lem:IM-GN:7},
which can be shown in the same way as Lemma
\ref{lem:IM-GN:3},
\ref{lem:IM-GN:4}, and \ref{lem:IM-GN:5}.

\begin{lem}
\label{lem:IM-GN:6}
The following hold.
\begin{align}
&
\Delta_{a+c}^{(M,N)}[(\ol{\lambda}[1];1), (\ol{\lambda}[l-1];l-1)]=0,
\notag\\
&
\ol{\lambda}[1] \in \boxed{\ol{\mc E}_1}, 
\ol{\lambda}[l-1] \in \ol{\mc D}_{l-1}, \quad
|\ol{\lambda}[1]|=a\geq 2, |\ol{\lambda}[l-1]|=c\geq 1, M \geq a+c \geq 3,
\notag
\\
\notag
\\
&\Delta_{a+c}^{(M,N)}[(\ol{\lambda}[1];1), (\ol{\lambda}[l];l)]=0,
\notag
\\
&
\ol{\lambda}[1] \in \boxed{\ol{\mc E}_1}, 
\ol{\lambda}[l] \in \ol{\mc D}_{l}, \quad
|\ol{\lambda}[1]|=a \geq 2, |\ol{\lambda}[l]|=c \geq 1, M \geq a+c \geq 3,
\notag
\\
\notag
\\
&\Delta_{a+c}^{(M,N)}[((\ol{\lambda}[1];1), (\ol{\lambda}[l-1];l-1), (\ol{\lambda}[l];l)]=0,\notag
\\
&\ol{\lambda}[1] \in \boxed{\ol{\mc E}_1}, 
\ol{\lambda}[l-1] \in \ol{\mc D}_{l-1},
\ol{\lambda}[l] \in \ol{\mc D}_l,
\quad |\ol{\lambda}[1]|=a \geq 2, |\ol{\lambda}[l-1]|+|\ol{\lambda}[l]|=c\geq 2, M \geq a+c \geq 4,
\notag
\\
\notag
\\
&\Delta_{a+c}^{(M,N)}[(\ol{\lambda}[0];0), (\ol{\lambda}[1];1), (\ol{\lambda}[l-1];l-1)]=0,\notag
\\
&
\ol{\lambda}[0] \in \ol{\mc E}_0,
\ol{\lambda}[1] \in \boxed{\ol{\mc E}_1}, 
\ol{\lambda}[l-1] \in \ol{\mc D}_{l-1},
\quad |\ol{\lambda}[0]|+|\ol{\lambda}[1]|=a \geq 3, |\ol{\lambda}[l-1]|=c \geq 1, M \geq a+c \geq 4,
\notag
\\
\notag
\\
&\Delta_{a+c}^{(M,N)}[(\ol{\lambda}[0];0), (\ol{\lambda}[1];1), (\ol{\lambda}[l];l)]=0,\notag
\\
&
\ol{\lambda}[0] \in \ol{\mc E}_0,
\ol{\lambda}[1] \in \boxed{\ol{\mc E}_1}, 
\ol{\lambda}[l] \in \ol{\mc D}_l,
\quad |\ol{\lambda}[0]|+|\ol{\lambda}[1]|=a \geq 3, |\ol{\lambda}[l]|=c \geq 1, M \geq a+c \geq 4,
\notag
\\
\notag
\\
&\Delta_{a+c}^{(M,N)}[(\ol{\lambda}[0];0), (\ol{\lambda}[1];1), (\ol{\lambda}[l-1];l-1), (\ol{\lambda}[l];l)]=0,\notag
\\
&
\ol{\lambda}[0] \in \ol{\mc E}_0,
\ol{\lambda}[1] \in \boxed{\ol{\mc E}_1}, 
\ol{\lambda}[l-1] \in \ol{\mc D}_{l-1},
\ol{\lambda}[l] \in \ol{\mc D}_l,
\quad |\ol{\lambda}[0]|+|\ol{\lambda}[1]|=a \geq 3, |\ol{\lambda}[l-1]|+|\ol{\lambda}[l]|=c \geq 2,\notag
\\
& M \geq a+c\geq 5,
\notag
\\
\notag
\\
&\Delta_{a+c}^{(M,N)}[({\lambda}[1];1), ({\lambda}[l-1];l-1)]=0,\notag
\\
&{\lambda}[1] \in \boxed{{\mc E}_1}, 
{\lambda}[l-1] \in {\mc D}_{l-1}, 
\quad |{\lambda}[1]|=a \geq 2, |{\lambda}[l-1]|=c \geq 1,
M \geq a+c \geq 3, \notag
\\
\notag
\\
&\Delta_{a+c}^{(M,N)}[({\lambda}[1];1), ({\lambda}[l];l)]=0,\quad
{\lambda}[1] \in \boxed{{\mc E}_1}, 
{\lambda}[l] \in \ol{\mc D}_{l}, |{\lambda}[1]|=a \geq 2, |{\lambda}[l]|=c \geq 1, M \geq a+c \geq 3,
\notag
\\
\notag
\\
&\Delta_{a+c}^{(M,N)}[(({\lambda}[1];1), ({\lambda}[l-1];l-1), ({\lambda}[l];l)]=0,\notag
\\
&{\lambda}[1] \in \boxed{{\mc E}_1}, 
{\lambda}[l-1] \in {\mc D}_{l-1},
{\lambda}[l] \in {\mc D}_l,
\quad |{\lambda}[1]|=a \geq 2,
 |{\lambda}[l-1]|+|{\lambda}[l]|=c \geq 2, M \geq a+c \geq 4,
 \notag
\\
\notag
\\
&\Delta_{a+c}^{(M,N)}[({\lambda}[0];0), ({\lambda}[1];1), ({\lambda}[l-1];l-1)]=0,\notag
\\
&
{\lambda}[0] \in {\mc E}_0,
{\lambda}[1] \in \boxed{{\mc E}_1}, 
{\lambda}[l-1] \in {\mc D}_{l-1},
\quad |{\lambda}[0]|+|{\lambda}[1]|=a \geq 3, 
|{\lambda}[l-1]|=c \geq 1, M \geq a+c \geq 4,
\notag
\end{align}

\begin{align}
&\Delta_{a+c}^{(M,N)}[({\lambda}[0];0), ({\lambda}[1];1), ({\lambda}[l];l)]=0,\notag
\\
&
{\lambda}[0] \in {\mc E}_0,
{\lambda}[1] \in \boxed{{\mc E}_1}, 
{\lambda}[l] \in {\mc D}_l,
\quad |{\lambda}[0]|+|{\lambda}[1]|=a \geq 3, |{\lambda}[l]|=c \geq 1, M \geq a+c \geq 4, \notag
\\
\notag
\\
&\Delta_{a+c}^{(M,N)}[({\lambda}[0];0), ({\lambda}[1];1), ({\lambda}[l-1];l-1), ({\lambda}[l];l)]=0,\notag
\\
&
{\lambda}[0] \in {\mc E}_0,
{\lambda}[1] \in \boxed{{\mc E}_1}, 
{\lambda}[l-1] \in {\mc D}_{l-1},
{\lambda}[l] \in {\mc D}_l,
\quad |{\lambda}[0]|+|{\lambda}[1]|=a, \geq 3, 
|{\lambda}[l-1]|+|{\lambda}[l]|=c \geq 2,
\notag\\
&M \geq a+c \geq 5
.\notag
\end{align}
\end{lem}

\begin{lem}\label{lem:IM-GN:7}
The following hold for 
$1\leq a, c \leq M, 2\leq a+c \leq M$.
\begin{align}
&
\sum_{
{\lambda}[1] \in {\mathcal E}_1
\atop{|{\lambda}[1]|=a}}
\sum_{{\lambda}[l-1] \in {\mathcal D}_{l-1}
\atop{|{\lambda}[l-1]|=c}}
\com{N}{{\lambda}[1]}
\com{N}{{\lambda}[l-1]}
\nu_{{\lambda}[1];1} 
\nu_{{\lambda}[l-1];l-1} 
\Delta_{a+c}^{(M,N)}[({\lambda}[1];1), ({\lambda}[l-1];l-1)]\notag
\\
&+
\sum_{{\lambda}[1] \in {\mathcal E}_1
\atop{|{\lambda}[1]|=a}}
\sum_{\lambda[l] \in {\mathcal D}_{l}
\atop{|{\lambda}[l]|=c}}
\com{N}{{\lambda}[1]}
\com{N}{{\lambda}[l]}
\nu_{{\lambda}[1];1} 
\nu_{{\lambda}[l];l} 
\Delta_{a+c}^{(M,N)}[({\lambda}[1];1), ({\lambda}[l];l)]
\notag\\
&+
\sum_{{\lambda}[1] \in {\mathcal E}_1
\atop{|{\lambda}[1]|=a}}
\sum_{{\lambda}[l-1]\in {\mathcal D}_{l-1},
{\lambda}[l] \in {\mathcal D}_{l}
\atop{|{\lambda}[l-1]|+|{\lambda}[l]|=c}}
\com{N}{{\lambda}[1]}
\com{N}{{\lambda}[l-1]}
\com{N}{{\lambda}[l]}
\nu_{{\lambda}[1];1} 
\nu_{{\lambda}[l-1];l-1}
\nu_{{\lambda}[l];l} \notag\\
&\times
\Delta_{a+c}^{(M,N)}[({\lambda}[1];1), ({\lambda}[l-1];l-1), ({\lambda}[l];l)]
\notag
\\
&+
\sum_{{\lambda}[0] \in {\mathcal E}_0, {\lambda}[1] \in {\mathcal E}_1
\atop{|{\lambda}[0]|+|{\lambda}[1]|=a}}
\sum_{{\lambda}[l-1] \in {\mathcal D}_{l-1}
\atop{|{\lambda}[l-1]|=c}}
\com{N}{{\lambda}[0]}
\com{N}{{\lambda}[1]}
\com{N}{{\lambda}[l-1]}
\nu_{{\lambda}[0];0}
\nu_{{\lambda}[1];1} 
\nu_{{\lambda}[l-1];l-1}\notag\\
&\times
\Delta_{a+c}^{(M,N)}[({\lambda}[0];0), ({\lambda}[1];1), ({\lambda}[l-1];l-1)]\notag
\\
&+
\sum_{{\lambda}[0]\in {\mathcal E}_0,
{\lambda}[1] \in {\mathcal E}_1
\atop{|{\lambda}[0]|+|{\lambda}[1]|=a}}
\sum_{\lambda[l] \in {\mathcal D}_{l}
\atop{|{\lambda}[l]|=c}}
\com{N}{{\lambda}[0]}
\com{N}{{\lambda}[1]}
\com{N}{{\lambda}[l]}
\nu_{{\lambda}[0];0} 
\nu_{{\lambda}[1];1} 
\nu_{{\lambda}[l];l} 
\Delta_{a+c}^{(M,N)}[
({\lambda}[0];0),
({\lambda}[1];1), ({\lambda}[l];l)]
\notag\\
&+
\sum_{{\lambda}[0]\in {\mathcal E}_0,
{\lambda}[1] \in {\mathcal E}_1 
\atop{|{\lambda}[0]|+|{\lambda}[1]|=a}}
\sum_{{\lambda}[l-1]\in {\mathcal D}_{l-1},
{\lambda}[l] \in {\mathcal D}_{l}
\atop{|{\lambda}[l-1]|+|{\lambda}[l]|=c}}
\com{N}{{\lambda}[0]}
\com{N}{{\lambda}[1]}
\com{N}{{\lambda}[l-1]}
\com{N}{{\lambda}[l]}
\nu_{{\lambda}[0];0} 
\nu_{{\lambda}[1];1} 
\nu_{{\lambda}[l-1];l-1}
\nu_{{\lambda}[l];l} \notag\\
&\times
\Delta_{a+c}^{(M,N)}[({\lambda}[0];0), ({\lambda}[1];1), ({\lambda}[l-1];l-1), ({\lambda}[l];l)]=0.
\notag
\end{align}

\begin{align}
&
\sum_{\overline{\lambda}[0] \in \overline{\mathcal E}_0
\atop{|\overline{\lambda}[0]|=a}}
\sum_{\overline{\lambda}[l-1] \in \overline{\mathcal D}_{l-1}
\atop{|\overline{\lambda}[l-1]|=c}}
\com{N}{\overline{\lambda}[0]}
\com{N}{\overline{\lambda}[l-1]}
\nu_{\overline{\lambda}[0];0} 
\nu_{\overline{\lambda}[l-1];l-1} 
\Delta_{a+c}^{(M,N)}[(\overline{\lambda}[0];0), (\overline{\lambda}[l-1];l-1)]
\notag
\end{align}

\begin{align}
&+
\sum_{\overline{\lambda}[0] \in \overline{\mathcal E}_0
\atop{|\overline{\lambda}[0]|=a}}
\sum_{\overline{\lambda}[l] \in \overline{\mathcal D}_{l}
\atop{|\overline{\lambda}[l]|=c}}
\com{N}{\overline{\lambda}[0]}
\com{N}{\overline{\lambda}[l]}
\nu_{\overline{\lambda}[0];0} 
\nu_{\overline{\lambda}[l];l} 
\Delta_{a+c}^{(M,N)}[(\overline{\lambda}[0];0), (\overline{\lambda}[l];l)]
\notag\\
&+
\sum_{\overline{\lambda}[0] \in \overline{\mathcal E}_0
\atop{|\overline{\lambda}[0]|=a}}
\sum_{\overline{\lambda}[l-1]\in \overline{\mathcal D}_{l-1},
\overline{\lambda}[l] \in \overline{\mathcal D}_{l}
\atop{|\overline{\lambda}[l-1]|+|\overline{\lambda}[l]|=c}}
\com{N}{\overline{\lambda}[0]}
\com{N}{\overline{\lambda}[l-1]}
\com{N}{\overline{\lambda}[l]}
\nu_{\overline{\lambda}[0];0} 
\nu_{\overline{\lambda}[l-1];l-1}
\nu_{\overline{\lambda}[l];l} \notag\\
&\times
\Delta_{a+c}^{(M,N)}[(\overline{\lambda}[0];0), (\overline{\lambda}[l-1];l-1), (\overline{\lambda}[l];l)]
\notag\\
&+
\sum_{\overline{\lambda}[1] \in \overline{\mathcal E}_1
\atop{|\overline{\lambda}[1]|=a}}
\sum_{\overline{\lambda}[l-1] \in \overline{\mathcal D}_{l-1}
\atop{|\overline{\lambda}[l-1]|=c}}
\com{N}{\overline{\lambda}[1]}
\com{N}{\overline{\lambda}[l-1]}
\nu_{\overline{\lambda}[1];1} 
\nu_{\overline{\lambda}[l-1];l-1} 
\Delta_{a+c}^{(M,N)}[(\overline{\lambda}[1];1), (\overline{\lambda}[l-1];l-1)]\notag
\\
&+
\sum_{\overline{\lambda}[1] \in \overline{\mathcal E}_1
\atop{|\overline{\lambda}[1]|=a}}
\sum_{\overline{\lambda}[l] \in \overline{\mathcal D}_{l}
\atop{|\overline{\lambda}[l]|=c}}
\com{N}{\overline{\lambda}[1]}
\com{N}{\overline{\lambda}[l]}
\nu_{\overline{\lambda}[1];1} 
\nu_{\overline{\lambda}[l];l} 
\Delta_{a+c}^{(M,N)}[(\overline{\lambda}[1];1), (\overline{\lambda}[l];l)]
\notag
\\
&+
\sum_{\overline{\lambda}[1] \in \overline{\mathcal E}_1
\atop{|\overline{\lambda}[1]|=a}}
\sum_{\overline{\lambda}[l-1]\in \overline{\mathcal D}_{l-1},
\overline{\lambda}[l] \in \overline{\mathcal D}_{l}
\atop{|\overline{\lambda}[l-1]|+|\overline{\lambda}[l]|=c}}
\com{N}{\overline{\lambda}[1]}
\com{N}{\overline{\lambda}[l-1]}
\com{N}{\overline{\lambda}[l]}
\nu_{\overline{\lambda}[1];1} 
\nu_{\overline{\lambda}[l-1];l-1}
\nu_{\overline{\lambda}[l];l} \notag\\
&\times
\Delta_{a+c}^{(M,N)}[(\overline{\lambda}[1];1), (\overline{\lambda}[l-1];l-1), (\overline{\lambda}[l];l)]
\notag
\\
&+
\sum_{\overline{\lambda}[0] \in \overline{\mathcal E}_0, \overline{\lambda}[1] \in \overline{\mathcal E}_1
\atop{|\overline{\lambda}[0]|+|\overline{\lambda}[1]|=a}}
\sum_{\overline{\lambda}[l-1] \in \overline{\mathcal D}_{l-1}
\atop{|\overline{\lambda}[l-1]|=c}}
\com{N}{\overline{\lambda}[0]}
\com{N}{\overline{\lambda}[1]}
\com{N}{\overline{\lambda}[l-1]}
\nu_{\overline{\lambda}[0];0}
\nu_{\overline{\lambda}[1];1} 
\nu_{\overline{\lambda}[l-1];l-1}\notag\\
&\times
\Delta_{a+c}^{(M,N)}[(\overline{\lambda}[0];0), (\overline{\lambda}[1];1), (\overline{\lambda}[l-1];l-1)]\notag
\\
&+
\sum_{\overline{\lambda}[0]\in \overline{\mathcal E}_0,
\overline{\lambda}[1] \in \overline{\mathcal E}_1
\atop{|\overline{\lambda}[0]|+|\overline{\lambda}[1]|=a}}
\sum_{\overline{\lambda}[l] \in \overline{\mathcal D}_{l}
\atop{|\overline{\lambda}[l]|=c}}
\com{N}{\overline{\lambda}[0]}
\com{N}{\overline{\lambda}[1]}
\com{N}{\overline{\lambda}[l]}
\nu_{\overline{\lambda}[0];0} 
\nu_{\overline{\lambda}[1];1} 
\nu_{\overline{\lambda}[l];l} 
\Delta_{a+c}^{(M,N)}[
(\overline{\lambda}[0];0),
(\overline{\lambda}[1];1), (\overline{\lambda}[l];l)]
\notag\\
&+
\sum_{\overline{\lambda}[0]\in \overline{\mathcal E}_0,
\overline{\lambda}[1] \in \overline{\mathcal E}_1
\atop{|\overline{\lambda}[0]|+|\overline{\lambda}[1]|=a}}
\sum_{\overline{\lambda}[l-1]\in \overline{\mathcal D}_{l-1},
\overline{\lambda}[l] \in \overline{\mathcal D}_{l}
\atop{|\overline{\lambda}[l-1]|+|\overline{\lambda}[l]|=c}}
\com{N}{\overline{\lambda}[0]}
\com{N}{\overline{\lambda}[1]}
\com{N}{\overline{\lambda}[l-1]}
\com{N}{\overline{\lambda}[l]}
\nu_{\overline{\lambda}[0];0} 
\nu_{\overline{\lambda}[1];1} 
\nu_{\overline{\lambda}[l-1];l-1}
\nu_{\overline{\lambda}[l];l} \notag\\
&\times
\Delta_{a+c}^{(M,N)}[(\overline{\lambda}[0];0), (\overline{\lambda}[1];1), (\overline{\lambda}[l-1];l-1), (\overline{\lambda}[l];l)]=0.
\notag
\end{align}
\end{lem}

\newpage

\appendix

\section{Quantum toroidal algebra $U_{\ssq}(\ddot{D}_l)$}
\label{appendix:bosonization-quantum toroidal D}

\subsection{Presentation}
Assume $l\ge4$. 
Let $\ssq\in\C^{\times}$ be such that
$\ssq^n\neq 1$ if $n\in\Z\backslash\{0\}$.
We set $[r]=(\ssq^r-\ssq^{-r})/(\ssq-\ssq^{-1})$.
{We retain the notation in Section \ref{section:deformed-W}.}

The quantum toroidal algebra $U_{\ssq}(\ddot{D}_l)$ 
of type $D_l$ is a unital associative algebra generated by elements
$x^\pm_{i,j}$, $h_{i,r}$, $k_{i,0}^{\pm1}$,
with $0\le i\le l$, $j\in\Z$, $r\in\Z\backslash\{0\}$,
and an invertible central element $C$.
We use the generating series
\begin{align*}
&x^\pm_i(z) =\sum_{j\in \Z} x^\pm_{i,j}z^{-j}, \quad 
k_i^{\pm}(z) = k_{i,0}^{\pm1}
\exp\bigl(\pm(\ssq-\ssq^{-1})\sum_{r=1}^\infty h_{i,\pm r}z^{\mp r}\bigr)\,.
\end{align*}
The defining relations are as follows. 
\begin{align*}
&[h_{i,r},h_{j,s}]=
\delta_{r+s,0}\,\frac{1}{r} [(\alb_i,\alb_j)r]\frac{C^{r}-C^{-r}}{\ssq-\ssq^{-1}}\,,
\\ 
&(z-\ssq^{\pm(\alb_i,\alb_j)}w)k_i^\mp(z)x_j^\pm(w)
+(w-\ssq^{\pm(\alb_i,\alb_j)}z)x^\pm_j(w)k_i^\mp(z)=0
\,,\\ 
&(z-\ssq^{\pm(\alb_i,\alb_j)}w)k_i^\pm\bigl(C^{-1}z\bigr)
x^\pm_j(w)
+(w-\ssq^{\pm(\alb_i,\alb_j)}z)
x^\pm_j(w)k_i^\pm\bigl(C^{-1}z\bigr)=0\,,
\\
&(z-\ssq^{\pm(\alb_i,\alb_j)}w)x^\pm_i(z)x^\pm_j(w)
+(w-\ssq^{\pm(\alb_i,\alb_j)}z)x^\pm_j(w)x^\pm_i(z)=0\,,
\\
&[x^+_i(z),x^-_j(w)]=\frac{\delta_{i,j}}{\ssq-\ssq^{-1}}
\Bigl\{\delta\Bigl(\frac{Cw}{z}\Bigr)k_i^+(w)
-\delta\Bigl(\frac{Cz}{w}\Bigr)k_i^-(z)\Bigr\}\,,
\\
&\mathop{\mathrm{Sym}}_{z_1,z_2}
[x^\pm_i(z_1),[x^\pm_i(z_2),x^\pm_j(w)]_\ssq]_{\ssq^{-1}}=0
\qquad \text{if $(\alb_i,\alb_j)=-1$}\,,
\\
&[x^\pm_i(z),x^\pm_j(w)]=0\qquad
\qquad \text{if $(\alb_i,\alb_j)=0$}\,.
\end{align*}
Here $[X,Y]_{\ssq^{\pm1}}=XY-\ssq^{\pm1}YX$.

Algebra  $\uq$ has a topological Hopf structure 
{but we do not use it here.}

Algebra  $\uq$ 
has a $\Z\times \bar{Q}$-grading, which we call 
degree and weight, and denote by $\deg$ and $\wt$ respectively, 
given by the assignment
\begin{align*}
&\deg x^\pm_{i,j}=j\,,\quad \deg h_{i,r}=r\,,\quad 
\deg k_{i,0}=\deg C=0\,,\\
&\wt x^\pm_{i,j}=\pm \alb_i\,,\quad \wt h_{i,r}=0\,,\quad 
\wt k_{i,0}=\wt C=0\,.
\end{align*}

\subsection{{``Level one''} representations.}

Consider a Heisenberg algebra with generators
$\{h_{i,r}\}_{{0\le i\le l},   
r\in\Z\backslash\{0\}}$ satisfying
\begin{align*}
[h_{i,r},h_{j,s}]=\delta_{r+s,0}\frac{[r]}{r}[(\alb_i,\alb_j)r]\,.
\end{align*}
Let $\mathcal{H}=\C[h_{i,-r}]_{{0\le i\le l},  
r>0}$ be its Fock representation.

Denote by $\C[\bar{Q}]$ {(resp. $\C[\bar P]$)}
the group ring of $\bar{Q}$ (resp. $\bar P$). 
For $\gamma\in \bar{P}$ 
we write the corresponding element in $\C[\bar{P}]$
as $e^\gamma$.
The coset $\bar{P}/\bar{Q}$ has representatives 
$\{\omb_0,\omb_1,\omb_{l-1},\omb_l\}$.
For $\omb \in \{\omb_0,\omb_1,\omb_{l-1},\omb_l\}$
define
\begin{align*}
\FF_{\omb}=
\C[ \bar{Q}]e^{\omb}\otimes \mathcal{H}\,.
\end{align*}
Define a $\Z\times \bar{Q}$-grading on $\FF_{\omb}$ by
setting, for $\gamma\in\bar{Q}$, $f\in\mathcal{H}$,
\begin{align*}
&\deg (e^{\gamma+\omb}\otimes f)
=-\frac{1}{2}(\gamma,\gamma)-(\gamma,\omb)
+\deg f\,,\\
&\wt (e^{\gamma+\omb}\otimes f)=\gamma\,.
\end{align*}

Let $\varepsilon: \bar{Q}\times \bar{Q}\to \{1,-1\}$ 
be the asymmetry function 
{\cite{Kac}} 
satisfying
\begin{align*}
&\varepsilon(\alb+\alb',\gamma) 
=\varepsilon(\alb,\gamma) \varepsilon(\alb',\gamma) \,,
\quad
\varepsilon(\alb,\gamma+\gamma')
=\varepsilon(\alb,\gamma)\varepsilon(\alb,\gamma')\,,
\\
&\varepsilon(\alb,\alb)=(-1)^{(\alb,\alb)/2}\,.
\end{align*}
We have $\varepsilon(\alb,\gamma)\varepsilon(\gamma,\alb)=(-1)^{(\alb,\gamma)}$. 
For $\alb\in \bar{Q}$, define operators $e^{\alb}$, 
$\partial_{\alb}$ and $c_{\alb}$ on $\FF_\omb$ by
\begin{align*}
&e^{\alb}\bigl(e^{\gamma+\omb}\otimes f\bigr) 
=e^{\alb+\gamma+\omb}\otimes f\,,
\\
&\partial_{\alb}\bigl(e^{\gamma+\omb}\otimes f\bigr) 
=(\alb,\gamma+\omb)
e^{\gamma+\omb}\otimes f \,,\\
&c_{\alb} \bigl(e^{\gamma+\omb}
\otimes f\bigr) 
=\varepsilon(\alb,\gamma)
 e^{\gamma+\omb}\otimes f\,,
\end{align*}
where $\gamma\in \bar{Q}$ and $f\in \mathcal{H}$.

\begin{thm}\label{thm:levone}
Let $\omb\in \{\omb_0,\omb_1,\omb_{l-1},\omb_l\}$. 
Then the following assignment makes $\FF_{\omb}$ a $\Z\times \bar{Q}$
graded representation of $\uq$:
\begin{align*}
&x^+_i(z)\mapsto e^{{\bar{\al}_i}}z^{\partial_{\bar{\al}_i}+1}c_{\bar{\al}_i}
\exp\Bigl(\sum_{r>0}\frac{\ssq^{-r}}{[r]}h_{i,{-r}}z^r\Bigr) 
\exp\Bigl(-\sum_{r>0}\frac{1}{[r]}h_{i,{r}}z^{-r}\Bigr) \,,
\\
&x^-_i(z)\mapsto -e^{-{\bar{\al}_i}}z^{-\partial_{\bar{\al}_i}+1}
c_{\bar{\al}_i}
\exp\Bigl(-\sum_{r>0}\frac{1}{[r]}h_{i,{-r}}z^r\Bigr) 
\exp\Bigl(\sum_{r>0}\frac{\ssq^r}{[r]}h_{i,{r}}z^{-r}\Bigr) \,,
\\
&k_{i,0}\mapsto \ssq^{\partial_{\bar{\al}_i}}\,,
\quad h_{i,r}\mapsto h_{i,r}\,,
\quad C\mapsto \ssq\,.
\end{align*}
Here $0\le i\le l$ and $r\in\Z\backslash\{0\}$.
\end{thm}
Theorem \ref{thm:levone} 
is proved in \cite{Frenkel-Jing-Wang} for $\omb=\omb_0$. The same proof works
for other $\omb_\nu$. 
We only note that the action of
$x^\pm_i(z)$ on $\FF_\omb$ has integer powers in $z$
due to the condition $\omb\in \bar{P}$.

\section{Normal ordering rules} \label{appendix:normal-ordering}
The following normal ordering rule holds for $\mathbf{\Lambda}_i(z)$.
\begin{align}
\omega(w_2/w_1)\bLambda_i(w_1)\bLambda_i(w_2)&
=
:\bLambda_i(w_1)\bLambda_i(w_2):
\frac{\big(\frac{w_2}{w_1}, q^\beta \frac{w_2}{w_1}, 
\frac{w_1}{w_2}, q^\beta \frac{w_1}{w_2}; p\big)_\infty}
{
\big(q^{\beta-1}\frac{w_2}{w_1}, q \frac{w_2}{w_1}, 
q^{\beta-1} \frac{w_1}{w_2}, q \frac{w_1}{w_2}; p\big)_\infty
},\notag\\
&~ i=1,\ldots,l,\bar{l},\ldots,\bar{1},
\notag
\\
\omega(w_2/w_1)\bLambda_i(w_1)\bLambda_j(w_2)
&=
:\bLambda_i(w_1)\bLambda_j(w_2):
\frac{\big(
p \frac{w_2}{w_1}, 
p q^{\beta} \frac{w_2}{w_1}, 
\frac{w_1}{w_2}, 
q^\beta \frac{w_1}{w_2}; p\big)_\infty}
{
\big(p q^{\beta-1}\frac{w_2}{w_1}, 
p q \frac{w_2}{w_1}, 
q^{\beta-1} \frac{w_1}{w_2}, 
q \frac{w_1}{w_2}; p\big)_\infty
},
\notag
\\
&~~ i \prec j, i\neq \bar{j}, i,j=1,\ldots,l,\bar{l},\ldots,\bar{1},
\notag
\\
\omega(w_2/w_1)\bLambda_j(w_1)\bLambda_i(w_2)
&=
:\bLambda_j(w_1)\bLambda_i(w_2):
\frac{\big(
\frac{w_2}{w_1}, q^\beta \frac{w_2}{w_1}, 
p
\frac{w_1}{w_2}, 
p 
q^{\beta} \frac{w_1}{w_2}; p\big)_\infty}
{
\big(q^{\beta-1}\frac{w_2}{w_1}, 
q \frac{w_2}{w_1}, 
p q^{\beta-1} \frac{w_1}{w_2}, 
p q \frac{w_1}{w_2}; p\big)_\infty
},\notag\\
&~~ i \prec j, i\neq \bar{j}, i,j=1,\ldots,l,\bar{l},\ldots,\bar{1},
\notag
\\
\omega(w_2/w_1)\bLambda_i(w_1)\bLambda_{\bar{i}}(w_2)
&=
:\bLambda_i(w_1)\bLambda_{\bar{i}}(w_2):
\frac{\big(
p \frac{w_2}{w_1}, 
p q^{\beta} \frac{w_2}{w_1}, 
\frac{w_1}{w_2}, 
q^\beta \frac{w_1}{w_2}; p\big)_\infty}
{
\big(p q^{\beta-1}\frac{w_2}{w_1}, 
p q \frac{w_2}{w_1}, 
q^{\beta-1} \frac{w_1}{w_2}, 
q \frac{w_1}{w_2}; p\big)_\infty
}\notag\\
&\times
\frac{\big(1-q^{(l-1-i)\beta+1}\frac{w_2}{w_1}\big)
\big(1-q^{(l-i)\beta-1}\frac{w_2}{w_1}\big)}{
\big(1-q^{(l-i)\beta}\frac{w_2}{w_1}\big)\big(
1-q^{(l-i-1)\beta}\frac{w_2}{w_1}\big)},~ i=1,2,\ldots,l.
\notag
\\
\omega(w_2/w_1)\bLambda_{\bar{i}}(w_1)\bLambda_i(w_2)
&=
:\bLambda_{\bar{i}}(w_1)\bLambda_i(w_2):
\frac{\big(
\frac{w_2}{w_1}, q^\beta \frac{w_2}{w_1}, 
p \frac{w_1}{w_2}, 
p q^{\beta} \frac{w_1}{w_2}; p\big)_\infty}
{
\big(q^{\beta-1}\frac{w_2}{w_1}, 
q \frac{w_2}{w_1}, 
p q^{\beta-1} \frac{w_1}{w_2}, 
p q \frac{w_1}{w_2}; p\big)_\infty}\notag\\
&\times
\frac{\big(1-q^{\beta(i-l)+1}\frac{w_2}{w_1}\big)
\big(1-q^{(i+1-l)\beta-1}\frac{w_2}{w_1}\big)}
{\big(1-q^{(i-l)\beta}\frac{w_2}{w_1}\big)\big(1-q^{(i+1-l)\beta}\frac{w_2}{w_1}\big)},~~ i=1,2,\ldots,l.
\notag
\end{align}
Here we used $\omega(w)$ defined in (\ref{def:omega}).

In Table 1 and Table 2 below, 
we give the list of contractions
$\mathcal{C}(\mathbf{\Lambda}_i(w),S_j(z))$ and
$\mathcal{C}(S_j(z),\mathbf{\Lambda}_i(w))$ respectively, 
using the shorthand notation in Section \ref{sec:basicW}.
For example, the top-left corners of Table 1 and Table 2 mean respectively
\begin{align*}
&\bLambda_1(w)S_0(z) =:\bLambda_1(w)S_0(z) :\, 
q^{\beta-1}
\frac{w-q^{\frac{\beta}{2}(l-4)+\frac{1}{2}}z}
{w-q^{\frac{\beta}{2}(l-2)-\frac{1}{2}}z}\,,\\
&S_0(z)\bLambda_1(w)
 =:S_0(z) \bLambda_1(w):\, 
\frac{z-q^{-\frac{\beta}{2}(l-4)-\frac{1}{2}}w}
{z-q^{-\frac{\beta}{2}(l-2)+\frac{1}{2}}w}\,.
\end{align*}
Contribution from the zero mode (the power $q^{\beta-1}$ in the first line)
is not shown in the tables.

In the following tables, $3\le i\le l-2$. 

\begin{table}[h]
\begin{tabular}{|c|c|c|c|c|}
\hline
& $S_0$ & $S_1$ & $S_2$ & ~~~$S_3$~~~
\\
\hline
$\bLambda_1$ 
& $q^{\frac{\beta}{2}(l-2)-\frac{1}{2}}-q^{\frac{\beta}{2}(l-4)+\frac{1}{2}}$ 
& $q^{\frac{\beta}{2}(l-2)+\frac{1}{2}}-q^{\frac{\beta}{2}l-\frac{1}{2}}$ & $0$ & $0$\\
\hline
$\bLambda_2$ 
& $q^{\frac{\beta}{2}(l-2)-\frac{1}{2}}-q^{\frac{\beta}{2}(l-4)+\frac{1}{2}}$ 
& $q^{\frac{\beta}{2}(l-2)-\frac{1}{2}}-q^{\frac{\beta}{2}(l-4)+\frac{1}{2}}$ 
& $q^{\frac{\beta}{2}(l-3)+\frac{1}{2}}-q^{\frac{\beta}{2}(l-1)-\frac{1}{2}}$ 
& $0$\\
\hline
$\bLambda_{\bar{2}}$
& $q^{-\frac{\beta}{2}(l-2)+\frac{1}{2}}-q^{-\frac{\beta}{2}(l-4)-\frac{1}{2}}$ 
& $q^{-\frac{\beta}{2}(l-2)+\frac{1}{2}}-q^{-\frac{\beta}{2}(l-4)-\frac{1}{2}}$ 
& $q^{-\frac{\beta}{2}(l-3)-\frac{1}{2}}-q^{-\frac{\beta}{2}(l-1)+\frac{1}{2}}$ 
& $0$\\
\hline
$\bLambda_{\bar{1}}$ 
& $q^{-\frac{\beta}{2}(l-2)+\frac{1}{2}}-q^{-\frac{\beta}{2}(l-4)-\frac{1}{2}}$ 
& $q^{-\frac{\beta}{2}(l-2)-\frac{1}{2}}-q^{-\frac{\beta}{2}l+\frac{1}{2}}$ & $0$ & $0$
\\
\hline
\hline
& $S_{i-2}$ & $S_{i-1}$ & $S_i$ & $S_{i+1}$
\\
\hline
$\bLambda_i$ 
& $0$ 
& $q^{\frac{\beta}{2}(l-i)-\frac{1}{2}}-q^{\frac{\beta}{2}(l-i-2)+\frac{1}{2}}$ 
& $q^{\frac{\beta}{2}(l-i-1)+\frac{1}{2}}-q^{\frac{\beta}{2}(l-i+1)-\frac{1}{2}}$ & $0$\\
\hline
$\bLambda_{\bar{i}}$ 
& $0$ 
& $q^{-\frac{\beta}{2}(l-i)+\frac{1}{2}}-q^{-\frac{\beta}{2}(l-i-2)-\frac{1}{2}}$ 
& $q^{-\frac{\beta}{2}(l-i-1)-\frac{1}{2}}-q^{-\frac{\beta}{2}(l-i+1)+\frac{1}{2}}$
& $0$\\
\hline
\end{tabular}
\vspace{2mm}
\caption{${\mathcal C}(\bLambda_i(w), S_j(z))$}
\label{table:pole1}
\end{table}

\begin{table}
\begin{tabular}{|c|c|c|c|c|}
\hline
& $S_0$ & $S_1$ & $S_2$ & ~~~$S_3$~~~
\\
\hline
$\bLambda_1$ 
& $q^{-\frac{\beta}{2}(l-2)+\frac{1}{2}}-q^{-\frac{\beta}{2}(l-4)-\frac{1}{2}}$ 
& $q^{-\frac{\beta}{2}(l-2)-\frac{1}{2}}-q^{-\frac{\beta}{2}l+\frac{1}{2}}$ & $0$ & $0$\\
\hline
$\bLambda_2$ 
& $q^{-\frac{\beta}{2}(l-2)+\frac{1}{2}}-q^{-\frac{\beta}{2}(l-4)-\frac{1}{2}}$ 
& $q^{-\frac{\beta}{2}(l-2)+\frac{1}{2}}-q^{-\frac{\beta}{2}(l-4)-\frac{1}{2}}$ 
& $q^{-\frac{\beta}{2}(l-3)-\frac{1}{2}}-q^{-\frac{\beta}{2}(l-1)+\frac{1}{2}}$ 
& $0$\\
\hline
$\bLambda_{\bar{2}}$
& $q^{\frac{\beta}{2}(l-2)-\frac{1}{2}}-q^{\frac{\beta}{2}(l-4)+\frac{1}{2}}$ 
& $q^{\frac{\beta}{2}(l-2)-\frac{1}{2}}-q^{\frac{\beta}{2}(l-4)+\frac{1}{2}}$ 
& $q^{\frac{\beta}{2}(l-3)+\frac{1}{2}}-q^{\frac{\beta}{2}(l-1)-\frac{1}{2}}$ 
& $0$\\
\hline
$\bLambda_{\bar{1}}$ 
& $q^{\frac{\beta}{2}(l-2)-\frac{1}{2}}-q^{\frac{\beta}{2}(l-4)+\frac{1}{2}}$ 
& $q^{\frac{\beta}{2}(l-2)+\frac{1}{2}}-q^{\frac{\beta}{2}l-\frac{1}{2}}$ & $0$ & $0$
\\
\hline
\hline
& $S_{i-2}$ & $S_{i-1}$ & $S_i$ & $S_{i+1}$
\\
\hline
$\bLambda_i$ 
& $0$ 
& $q^{-\frac{\beta}{2}(l-i)+\frac{1}{2}}-q^{-\frac{\beta}{2}(l-i-2)-\frac{1}{2}}$ 
& $q^{-\frac{\beta}{2}(l-i-1)-\frac{1}{2}}-q^{-\frac{\beta}{2}(l-i+1)+\frac{1}{2}}$ & $0$\\
\hline
$\bLambda_{\bar{i}}$ 
& $0$ 
& $q^{\frac{\beta}{2}(l-i)-\frac{1}{2}}-q^{\frac{\beta}{2}(l-i-2)+\frac{1}{2}}$ 
& $q^{\frac{\beta}{2}(l-i-1)+\frac{1}{2}}-q^{\frac{\beta}{2}(l-i+1)-\frac{1}{2}}$
& $0$\\
\hline
\end{tabular}
\vspace{2mm}
\caption{${\mathcal C}(S_j(z), \bLambda_i(w))$}
\label{table:pole2}
\end{table}

\section{Locations of poles for ${\bf I}_M \cdot {\bf G}_N(\vartheta)$}
 \label{appendix:IM-GN}
 
%

\newpage

\section{Transformation of integrals in ${\bf I}_1 \cdot {\bf G}_N(\vartheta)$}
\label{appendix:I1-GN}

%

\newpage
\section{Transformation of integrals in ${\bf I}_2 \cdot {\bf G}_N(\vartheta)$}
\label{appendix:I2-GN}

%

\vskip 1cm

{\bf Acknowledgments.\ }

MJ and TK are partially supported by JSPS KAKENHI Grant Numbers 25K07041,
19K03509, 23K03137 
and 25K06942. 
TK thanks Rikkyo University, where a part of this work was done, 
for hospitality.

\end{document}